\documentclass[10pt, oneside, a4paper]{article}
\usepackage{amsmath, amsthm, amssymb, amsfonts}
\usepackage{geometry}
\usepackage{graphicx}
\usepackage{natbib}
\usepackage{setspace}
\usepackage{algorithm}
\usepackage{algorithmicx}
\usepackage{algpseudocode}
\usepackage{longtable}
\usepackage{caption}

\newtheorem{theorem}{Theorem}[section]
\newtheorem{lemma}[theorem]{Lemma}
\newtheorem{proposition}[theorem]{Proposition}
\newtheorem{corollary}[theorem]{Corollary}
\theoremstyle{definition}
\newtheorem{definition}{Definition}[section]
\newtheorem{example}{Example}[section]
\newtheorem{assumption}{Assumption}
\newtheorem{remark}[theorem]{Remark}

\input{tcilatex}

\begin{document}

\title{Controlling the Size of Autocorrelation Robust Tests\thanks{%
Financial support of the second author by Austrian Science Fund (FWF)
Project P27398-G11 and by the Danish National Research Foundation (Grant
DNRF 78, CREATES) is gratefully acknowledged. We thank the associate editor
and the referees as well as Majid Al-Sadoon, Frank Kleibergen, the
participants in the Workshop \textquotedblleft New Approaches to the
Identification of Macroeconomics Models\textquotedblright\ at Oxford
University, seminar participants at the Tinbergen Institute, at Aarhus
University, at the University of Vienna, and at TU Vienna for helpful
comments. Address correspondence to Benedikt P\"{o}tscher, Department of
Statistics, University of Vienna, A-1090 Oskar-Morgenstern Platz 1. E-Mail:
benedikt.poetscher@univie.ac.at.}}
\author{Benedikt M. P\"{o}tscher and David Preinerstorfer \\
Department of Statistics, University of Vienna\\
ECARES, Universit\'{e} libre de Bruxelles}
\date{First version: November 2016\\
This version: August 2018}
\maketitle

\begin{abstract}
Autocorrelation robust tests are notorious for suffering from size
distortions and power problems. We investigate under which conditions the
size of autocorrelation robust tests can be controlled by an appropriate
choice of critical value.
\end{abstract}

\section{Introduction}

Autocorrelation robust tests have gained prominence in econometrics in the
late 1980s and early 1990s, mainly through the work by \cite{NW87, NW94}, 
\cite{A91}, and \cite{A92}. These tests are Wald-type tests that use a
nonparametric variance estimator which tries to take the autocorrelation in
the data into account. For more on the history of such tests, that can
actually be traced back at least to \cite{GR57}, see Section 1 of \cite%
{PP2016}. Critical values are usually obtained from the asymptotic null
distribution, which is a chi-square distribution under the framework used in
the before mentioned papers. Considerable size distortions of the resulting
tests have been reported in several Monte Carlo studies. In an attempt to
ameliorate this problem, \cite{KVB2000}, \cite{KiefVogl2002, KV2002, KV2005}
have suggested an alternative asymptotic framework (\textquotedblleft fixed
bandwidth asymptotics\textquotedblright ), which resulted in another
proposal for a critical value. While this proposal often leads to some
reduction in the size distortion (typically at the expense of some loss in
power), it does not eliminate the problem. A finite-sample theory explaining
the aforementioned size distortions and power deficiencies of
autocorrelation robust tests (in either of their two variants) has been
developed in \cite{PP2016} in the framework of testing affine restrictions
on the regression parameters of a linear regression model, the errors of
which follow a Gaussian stationary process: Under the very mild assumption
that the model for the autocorrelations encompasses the stationary
autoregressive model of order one, Section 3 of \cite{PP2016} shows, in
particular, that the following cases typically arise for a given
autocorrelation robust test statistic (which of the cases arises depending
on observable quantities only):

\begin{enumerate}
\item The test has size one for every choice of the critical value $C$. This
case arises, for example, when one wants to test a hypothesis concerning the
intercept of a regression.

\item There exists an observable threshold $c$ such that the size of the
test is $1$, if the critical value $C$ is less than $c$. If the critical
value $C$ is larger than the threshold $c$, the infimal power of the test is
zero.
\end{enumerate}

We note that these results are to a large extent driven by a
\textquotedblleft concentration phenomenon\textquotedblright\ that arises
for strong correlations in the error process. For a more detailed discussion
and some intuition see Sections 3.1 and 5.2 as well as pp.~275-276 of \cite%
{PP2016}.

While the two cases just mentioned do not exhaust all possibilities, their
union is shown to be \emph{generic} in \cite{PP2016}, entailing that general
autocorrelation robust tests are typically fraught with considerable size
and/or power problems. \cite{PP2016} also show that there is a case not
covered by Cases 1 and 2 (and which hence is nongeneric) for which the
following holds, \emph{provided} the model for the autocorrelation is
precisely the stationary autoregressive model of order one (which, of
course, is a quite restrictive assumption):

\begin{enumerate}
\item[3.] For every prescription of a significance value $\alpha $, the size
of the test can be made less than or equal to $\alpha $ by a proper choice
of the critical value $C$. For this $C$, infimal power is positive (and more
can be said about power properties, but we do not discuss this here).
\end{enumerate}

While this case is nongeneric, \cite{PP2016} also show how one can -- under
the restrictive assumption on the correlation structure just mentioned --
often force the original testing problem into the situation of Case 3 by
augmenting the model with auxiliary regressors and then using the
autocorrelation robust test on the augmented model. We note that the above
mentioned results in \cite{PP2016} are given for commonly used
autocorrelation robust test statistics where the bandwidth is nonrandom. 
\cite{Prein2014b} extends these results to the case of random bandwidths and
when prewithening is used. As shown in Section 3.3 of \cite{PP2016}, similar
results also hold if a parametrically-based long-run variance estimator or
feasible generalized least squares are used.

A fundamental question that is not satisfactorily answered by the results in 
\cite{PP2016} is under which conditions one can control the size of the test
(i.e., under which conditions one can find a critical value $C$ in such a
way that the resulting test has size less than or equal to the prescribed
significance value $\alpha $, $\alpha <1$): While in Case 1 we have a
definite (negative) answer, the result in Case 2 as it stands tells us only
that one can not choose $C$ smaller than the threshold $c$ if one wants to
control the size by $\alpha $ (thus allowing one to rule out proposals for $%
C $ if they fail this criterion), but it does not tell us anything about the
behavior of the size if $C$ is chosen larger than the threshold $c$. And the
result in Case 3, which guarantees size control, is limited by a quite
restrictive assumption on the correlation structure.

In the present paper we answer the question under which conditions the size
of an autocorrelation robust test can be controlled. In particular, we
provide sufficient conditions that guarantee that size control is possible
even if we allow for \emph{arbitrary} stationary autocorrelation in the
data. We then show that these conditions are broadly satisfied, namely that
they are satisfied \emph{generically} (i.e., are -- given the restriction to
be tested -- satisfied except for a Lebesgue null set in an appropriate
universe of relevant design matrices). We also discuss how the critical
value that leads to size control can be determined numerically and provide
the \textsf{R}-package \textbf{acrt} (\cite{csart}) for its computation. The
usefulness of the proposed algorithms and their implementation in the 
\textsf{R}-package are illustrated on testing problems involving
macroeconomic data taken from the FRED-MD database of \cite{mccracken}. In
particular, we show that -- even in situations where size control is
possible -- standard critical values in the literature based on asymptotic
considerations are often substantially smaller than the size-controlling
critical values devised in the present paper, and hence fail to deliver size
control. Furthermore, in a subsequent paper \cite{PP4} we show that the
sufficient conditions for size control we provide here are in a sense also
necessary (for a large class of test statistics). While the main emphasis in
the current paper is on stationary autocorrelation in the data, the general
theory for size control developed in Section \ref{general} applies to
arbitrary models for the covariance structure of the data (including e.g.,
the case of heteroskedasticity or the case of spatial correlation). The
present paper is only concerned with the size of autocorrelation robust
tests and its control. Power properties of such tests and ways of improving
power without loosing the size control property is the topic of a companion
paper. While a trivial remark, we would like to note that the size control
results given in this paper can obviously be translated into results stating
that the minimal coverage probability of the associated confidence set
obtained by \textquotedblleft inverting\textquotedblright\ the test is not
less than the nominal level.

The development in the present paper concentrates on controlling the size of
autocorrelation robust tests by an appropriate choice of a `fixed' (i.e.,
nonrandom) critical value. Extensions to size control by use of random
critical values are the subject of ongoing research.

We also note that, on a broader scale, our results contribute to an
important recent literature concerning size properties of tests and minimal
coverage properties of confidence sets, see, e.g., \cite{andrews_guggen_2009}%
, \cite{guggenberger2012}, \cite{mikusheva2007}, \cite{moreira2009}, or \cite%
{romano2012}).

The paper is organized as follows: Section \ref{frame} lays out the
framework and introduces the basic concepts. Section \ref{autocorr} contains
the size control results for commonly used autocorrelation robust tests: In
Section \ref{main} we present a theorem (Theorem \ref{HAR_F_all}) that lists
easy to verify sufficient conditions for size control under a very broad
model for the autocorrelation structure in the data. Several examples
illustrate this theorem. A similar result for a variant of the test
statistic considered in Theorem \ref{HAR_F_all} is provided in Theorem \ref%
{HAR_F_all_Eicker}. In Section \ref{gen_size_contol} we make precise in
which sense the conditions for size control in Theorems \ref{HAR_F_all} and %
\ref{HAR_F_all_Eicker} are generically satisfied. After a short section
(Section \ref{power}) commenting on power properties, we discuss extensions
of the results in Sections \ref{main} and \ref{gen_size_contol} to other
classes of test statistics. Computational issues regarding the determination
of the critical values effecting size control are discussed in Section \ref%
{numerical}. This section also contains a numerical illustration of the
algorithms suggested. The results in Section \ref{autocorr} rest on a
general theory for size control of tests under general forms of
nonsphericity, which is developed in Section \ref{general}. In Section \ref%
{structure under Topelitz} we apply this general theory to the case where
the nonsphericity is a consequence of stationary autocorrelation in the data
and we obtain size control results more general (but also more complex) than
the ones considered in Section \ref{autocorr}. In Appendix \ref{App E} we
show that the Gaussianity assumption underlying the paper can be
substantially relaxed, whereas Appendices \ref{App A}-\ref{App D} contain
the proofs and some auxiliary results for the main body of the paper.
Appendices \ref{App F}-\ref{App H} contain material relevant for the
numerical results in Section \ref{numerical} including a detailed
description of the algorithms used. Appendix \ref{App I} investigates a
proposal for choosing critical values suggested by a referee, while Appendix %
\ref{App J} discusses what happens if, instead of stationary solutions,
starting value solutions of autoregressive models are used as the error
process in the linear model; in particular, in this appendix we show that
under this alternative covariance model size control is always possible, but
that standard critical values proposed in the literature fail to deliver
size control. Appendix \ref{App K} provides some comments on the case of
stochastic regressors. Finally, Appendix \ref{app:tables} contains some
tables pertaining to Section~\ref{numerical}.

\section{Framework\label{frame}}

Consider the linear regression model 
\begin{equation}
\mathbf{Y}=X\beta +\mathbf{U},  \label{lm}
\end{equation}%
where $X$ is a (real) nonstochastic regressor (design) matrix of dimension $%
n\times k$ and where $\beta \in \mathbb{R}^{k}$ denotes the unknown
regression parameter vector. We always assume $\limfunc{rank}(X)=k$ and $%
1\leq k<n$. We furthermore assume that the $n\times 1$ disturbance vector $%
\mathbf{U}=(\mathbf{u}_{1},\ldots ,\mathbf{u}_{n})^{\prime }$ is normally
distributed with mean zero and unknown covariance matrix $\sigma ^{2}\Sigma $%
, where $\Sigma $ varies in a prescribed (nonempty) set $\mathfrak{C}$ of
symmetric and positive definite $n\times n$ matrices and where $0<\sigma
^{2}<\infty $ holds ($\sigma $ always denoting the positive square root).%
\footnote{%
Since we are concerned with finite-sample results only, the elements of $%
\mathbf{Y}$, $X$, and $\mathbf{U}$ (and even the probability space
supporting $\mathbf{Y}$ and $\mathbf{U}$) may depend on sample size $n$, but
this will not be expressed in the notation. Furthermore, the obvious
dependence of$\ \mathfrak{C}$ on $n$ will also not be shown in the notation.}
The set $\mathfrak{C}$ will be referred to as the covariance model. We shall
always assume that $\mathfrak{C}$ allows $\sigma ^{2}$ and $\Sigma $ to be
uniquely determined from $\sigma ^{2}\Sigma $.\footnote{%
That is, $\mathfrak{C}$ has the property that $\Sigma \in \mathfrak{C}$
implies $\delta \Sigma \notin \mathfrak{C}$ for every $\delta \neq 1$.}
[This entails virtually no loss of generality and can always be achieved,
e.g., by imposing some normalization assumption on the elements of $%
\mathfrak{C}$ such as normalizing the first diagonal element of $\Sigma $ or
the norm of $\Sigma $ to one, etc.] The leading case discussed in Section %
\ref{autocorr} will concern the situation where $\mathfrak{C}$ results from
the assumption that the elements $\mathbf{u}_{1},\ldots ,\mathbf{u}_{n}$ of
the $n\times 1$ disturbance vector $\mathbf{U}$ are distributed like
consecutive elements of a zero mean weakly stationary Gaussian process with
an unknown spectral density.

The linear model described in (\ref{lm}) together with the Gaussianity
assumption on $\mathbf{U}$ induces a collection of distributions on the
Borel-sets of $\mathbb{R}^{n}$, the sample space of $\mathbf{Y}$. Denoting a
Gaussian probability measure with mean $\mu \in \mathbb{R}^{n}$ and
(possibly singular) covariance matrix $A$ by $P_{\mu ,A}$, the induced
collection of distributions is then given by 
\begin{equation}
\left\{ P_{\mu ,\sigma ^{2}\Sigma }:\mu \in \mathrm{\limfunc{span}}%
(X),0<\sigma ^{2}<\infty ,\Sigma \in \mathfrak{C}\right\} .  \label{lm2}
\end{equation}%
Since every $\Sigma \in \mathfrak{C}$ is positive definite by assumption,
each element of the set in the previous display is absolutely continuous
with respect to (w.r.t.) Lebesgue measure on $\mathbb{R}^{n}$.

We shall consider the problem of testing a linear (better: affine)
hypothesis on the parameter vector $\beta \in \mathbb{R}^{k}$, i.e., the
problem of testing the null $R\beta =r$ against the alternative $R\beta \neq
r$, where $R$ is a $q\times k$ matrix always of rank $q\geq 1$ and $r\in 
\mathbb{R}^{q}$. Set $\mathfrak{M}=\limfunc{span}(X)$. Define the affine
space 
\begin{equation*}
\mathfrak{M}_{0}=\left\{ \mu \in \mathfrak{M}:\mu =X\beta \text{ and }R\beta
=r\right\}
\end{equation*}%
and let 
\begin{equation*}
\mathfrak{M}_{1}=\left\{ \mu \in \mathfrak{M}:\mu =X\beta \text{ and }R\beta
\neq r\right\} .
\end{equation*}%
Adopting these definitions, the above testing problem can then be written
more precisely as 
\begin{equation}
H_{0}:\mu \in \mathfrak{M}_{0},\ 0<\sigma ^{2}<\infty ,\ \Sigma \in 
\mathfrak{C}\quad \text{ vs. }\quad H_{1}:\mu \in \mathfrak{M}_{1},\
0<\sigma ^{2}<\infty ,\ \Sigma \in \mathfrak{C}.  \label{testing problem}
\end{equation}%
The previously introduced concepts and notation will be used throughout the
paper.

The assumption of Gaussianity is made mainly in order not to obscure the
structure of the problem by technicalities. Substantial generalizations away
from Gaussianity are possible and will be discussed in Appendix \ref{App E}.
The assumption of nonstochastic regressors can be relaxed somewhat: If $X$
is random and, e.g., independent of $\mathbf{U}$, the results of the paper
apply after one conditions on $X$. For arguments supporting conditional
inference see, e.g., \cite{RO1979}. In case $X$ is random with $X$ and $%
\mathbf{U}$ being dependent, our \emph{theory} (at least in its current
state) is not applicable. However, we show in Appendix \ref{App K} by means
of numerical examples that in such a case the critical values generated by
the \emph{algorithm} proposed in Section \ref{numerical} can lead to tests
that perform at least as well as (and often better than) tests that are
based on standard critical values suggested in the literature.\footnote{%
That is, the algorithm is used acting as if $X$ and $\mathbf{U}$ were
independent.} And this is so despite the fact that the asymptotic theory
used in the literature to justify the latter critical values \emph{is}
applicable to cases where $X$ and $\mathbf{U}$ are dependent (in particular
to the examples considered in Appendix \ref{App K}).

We next collect some further terminology and notation used throughout the
paper. A (nonrandomized) \textit{test} is the indicator function of a
Borel-set $W$ in $\mathbb{R}^{n}$, with $W$ called the corresponding \textit{%
rejection region}. The \textit{size} of such a test (rejection region) is
the supremum over all rejection probabilities under the null hypothesis $%
H_{0}$, i.e., $\sup_{\mu \in \mathfrak{M}_{0}}\sup_{0<\sigma ^{2}<\infty
}\sup_{\Sigma \in \mathfrak{C}}P_{\mu ,\sigma ^{2}\Sigma }(W)$. Throughout
the paper we let $\hat{\beta}(y)=\left( X^{\prime }X\right) ^{-1}X^{\prime
}y $, where $X$ is the design matrix appearing in (\ref{lm}) and $y\in 
\mathbb{R}^{n}$. The corresponding ordinary least squares (OLS) residual
vector is denoted by $\hat{u}(y)=y-X\hat{\beta}(y)$. We use $\Pr $ as a
generic symbol for a probability measure. Lebesgue measure on the Borel-sets
of $\mathbb{R}^{n}$ will be denoted by $\lambda _{\mathbb{R}^{n}}$, whereas
Lebesgue measure on an affine subspace $\mathcal{A}$ of $\mathbb{R}^{n}$
(but viewed as a measure on the Borel-sets of $\mathbb{R}^{n}$) will be
denoted by $\lambda _{\mathcal{A}}$, with zero-dimensional Lebesgue measure
being interpreted as point mass. The set of real matrices of dimension $%
l\times m$ is denoted by $\mathbb{R}^{l\times m}$ (all matrices in the paper
will be real matrices) and Lebesgue measure on this set, equipped with its
Borel $\sigma $-field, is denoted by $\lambda _{\mathbb{R}^{l\times m}}$.
The Euclidean norm of a vector is denoted by $\left\Vert \cdot \right\Vert $%
, but the same symbol is also used to denote a norm of a matrix. Let $%
B^{\prime }$ denote the transpose of a matrix $B\in \mathbb{R}^{l\times m}$
and let $\mathrm{\limfunc{span}}(B)$ denote the subspace in $\mathbb{R}^{l}$
spanned by its columns. For a symmetric and nonnegative definite matrix $B$
we denote the unique symmetric and nonnegative definite square root by $%
B^{1/2}$. For a linear subspace $\mathcal{L}$ of $\mathbb{R}^{n}$ we let $%
\mathcal{L}^{\bot }$ denote its orthogonal complement and we let $\Pi _{%
\mathcal{L}}$ denote the orthogonal projection onto $\mathcal{L}$. We use
the convention that the adjoint of a $1\times 1$ dimensional matrix $D$,
i.e., $\mathrm{\limfunc{adj}}(D)$, equals one. Given an $m$-dimensional
vector $v$ we write $\mathrm{\limfunc{diag}}(v)$ for the $m\times m$
diagonal matrix with main diagonal given by $v$. The $j$-th standard basis
vector in $\mathbb{R}^{n}$ is written as $e_{j}(n)$. With $e_{+}$ we denote
the $n\times 1$ vector of ones, i.e., $e_{+}=(1,\ldots ,1)^{\prime }$ and we
define the $n\times 1$ vector $e_{-}=(-1,1,\ldots ,(-1)^{n})^{\prime }$.
Furthermore, we let $\mathbb{N}$ denote the set of all positive integers. A
sum (product, respectively) over an empty index set is to be interpreted as $%
0$ ($1$, respectively). Finally, for a subset $A$ of a topological space we
denote by $\limfunc{cl}(A)$, $\limfunc{int}(A)$, and $\limfunc{bd}(A)$ the
closure, interior, and boundary of $A$ (w.r.t. the ambient space),
respectively.

\section{Size control of tests of affine restrictions in regression models
with stationary autocorrelated errors: First results \label{autocorr}}

In this section we are concerned with size control of autocorrelation robust
tests that have been designed for use in case of stationary disturbances. We
thus assume throughout this section that the elements $\mathbf{u}_{1},\ldots
,\mathbf{u}_{n}$ of the $n\times 1$ disturbance vector $\mathbf{U}$ are
distributed like consecutive elements of a zero mean weakly stationary
Gaussian process that has an unknown spectral density, which is not almost
everywhere equal to zero. Consequently, the covariance matrix of the
disturbance vector is positive definite and can be written as $\sigma
^{2}\Sigma (f)$ where%
\begin{equation*}
\Sigma (f)=\left[ \int_{-\pi }^{\pi }e^{-\iota (j-l)\omega }f(\omega
)d\omega \right] _{j,l=1}^{n},
\end{equation*}%
with $f$ varying in $\mathfrak{F}$, a prescribed (nonempty) family of \emph{%
normalized} (i.e., $\int_{-\pi }^{\pi }f(\omega )d\omega =1$) spectral
densities, and where $0<\sigma ^{2}<\infty $ holds. Here $\iota $ denotes
the imaginary unit. The set $\mathfrak{F}$ may, for example, be $\mathfrak{F}%
_{\mathrm{all}}$, the set of \emph{all} normalized spectral densities, or a
subset thereof (e.g., the set of normalized spectral densities corresponding
to stationary autoregressive or autoregressive moving average models of a
certain order, or to fractional autoregressive moving average models, etc.).
We define the associated covariance model $\mathfrak{C}(\mathfrak{F}%
)=\left\{ \Sigma (f):f\in \mathfrak{F}\right\} $ and note that the test
problem (\ref{testing problem})\ now becomes

\begin{equation}
H_{0}:\mu \in \mathfrak{M}_{0},\ 0<\sigma ^{2}<\infty ,\ f\in \mathfrak{F}%
\quad \text{ vs. }\quad H_{1}:\mu \in \mathfrak{M}_{1},\ 0<\sigma
^{2}<\infty ,\ f\in \mathfrak{F}.  \label{testing problem 0}
\end{equation}

\begin{remark}
\label{Fall=FAR(n-1)}(i) As is well-known, the covariance model $\mathfrak{C}%
(\mathfrak{F}_{\mathrm{all}})$ is precisely the set of \emph{all} $n\times n$
symmetric and positive definite Toeplitz matrices with ones on the main
diagonal, cf. Lemma \ref{Toeplitz} and Remark \ref{Toeplitz_rem} in Appendix %
\ref{App C}. It is thus maximal in the sense that it coincides with the set
of \emph{all} positive definite $n\times n$ correlation matrices that can be
generated from sections of length $n$ of stationary processes (possessing a
spectral density or not).

(ii) Furthermore, as is well-known, $\mathfrak{C}(\mathfrak{F}_{\mathrm{all}%
})$ coincides with $\mathfrak{C}(\mathfrak{F}_{\mathrm{AR(}p\mathrm{)}})$ if 
$p\geq n-1$, where $\mathfrak{F}_{\mathrm{AR(}p\mathrm{)}}$ is the set of
all normalized spectral densities corresponding to stationary autoregressive
processes of order not larger than $p$, cf. Remark \ref{Toeplitz_rem} in
Appendix \ref{App C}. As the testing problem depends on $\mathfrak{F}$ only
via $\mathfrak{C}(\mathfrak{F})$, the testing problems with $\mathfrak{F}=%
\mathfrak{F}_{\mathrm{all}}$ and $\mathfrak{F}=\mathfrak{F}_{\mathrm{AR(}p%
\mathrm{)}}$ (with $p\geq n-1$), respectively, coincide. Thus we can use
that \textquotedblleft parameterization\textquotedblright\ of the covariance
model (represented by $\mathfrak{F}_{\mathrm{all}}$ and $\mathfrak{F}_{%
\mathrm{AR(}p\mathrm{)}}$, respectively), which is more convenient for our
purpose. More generally, if $\mathfrak{C}(\mathfrak{F}_{1})=\mathfrak{C}(%
\mathfrak{F}_{2})$ holds for two sets of normalized spectral densities $%
\mathfrak{F}_{1}$ and $\mathfrak{F}_{2}$, the same argument can be made.
\end{remark}

Commonly used autocorrelation robust tests for the null hypothesis $H_{0}$
given by (\ref{testing problem 0}) are based on test statistics $T_{w}:%
\mathbb{R}^{n}\rightarrow \mathbb{R}$ of the form 
\begin{equation}
T_{w}(y)=\left\{ 
\begin{array}{cc}
(R\hat{\beta}(y)-r)^{\prime }\hat{\Omega}_{w}^{-1}\left( y\right) (R\hat{%
\beta}(y)-r) & \text{if }\det \hat{\Omega}_{w}\left( y\right) \neq 0 \\ 
0 & \text{if }\det \hat{\Omega}_{w}\left( y\right) =0%
\end{array}%
\right.  \label{tlr}
\end{equation}%
for $y\in \mathbb{R}^{n}$, where $R$ and $r$ are as in Section \ref{frame},
and where%
\begin{equation}
\hat{\Omega}_{w}\left( y\right) =nR(X^{\prime }X)^{-1}\hat{\Psi}%
_{w}(y)(X^{\prime }X)^{-1}R^{\prime },  \label{omega_0}
\end{equation}%
\begin{equation}
\hat{\Psi}_{w}(y)=\sum\limits_{j=-(n-1)}^{n-1}w(j,n)\hat{\Gamma}_{j}(y).
\label{lrve}
\end{equation}%
Here $\hat{\Gamma}_{j}(y)=n^{-1}\sum_{t=j+1}^{n}\hat{v}_{t}(y)\hat{v}%
_{t-j}(y)^{\prime }$ if $j\geq 0$ and $\hat{\Gamma}_{j}\left( y\right) =\hat{%
\Gamma}_{-j}(y)^{\prime }$ else, where $\hat{v}_{t}(y)=\hat{u}_{t}(y)x_{t%
\mathbf{\cdot }}^{\prime }$, with $\hat{u}_{t}(y)$ denoting the $t$-th
coordinate of the least squares residual vector $\hat{u}(y)=y-X\hat{\beta}%
(y) $ and with $x_{t\cdot }$ denoting the $t$-th row vector of $X$.
Rejection is for large values of $T_{w}$. We make the following standard
assumption on the weights.

\begin{assumption}
\label{AW} The weights $w(j,n)$ for $j=-(n-1),\ldots ,n-1$ are
data-independent and satisfy $w(0,n)=1$ as well as $w\left( -j,n\right)
=w\left( j,n\right) $. Furthermore, the symmetric $n\times n$ Toeplitz
matrix $\mathcal{W}_{n}$ with elements $w\left( i-j,n\right) $ is positive
definite.
\end{assumption}

This assumption implies that $\hat{\Psi}_{w}(y)$, and hence $\hat{\Omega}%
_{w}\left( y\right) $, is always nonnegative definite, see Lemma 3.1 in \cite%
{PP2016}. In many applications the weights take the form $w(j,n)=w_{0}\left(
|j|/M_{n}\right) $, where the lag-window $w_{0}$ is an even function with $%
w_{0}(0)=1$ and where $M_{n}>0$ is a truncation lag (bandwidth) parameter.
In this case the first part of the above assumption means that we are
considering deterministic bandwidths only. For extensions of the results in
this section to more general classes of tests statistics, including the case
of data-dependent bandwidth choices and prewhitening, see Subsection \ref%
{extension}. Assumption \ref{AW} is known to be satisfied, e.g., for the
(modified) Bartlett, Parzen, or the Quadratic Spectral lag-window, but is
not satisfied, e.g., for the rectangular lag-window (with $M_{n}>1$). It is
also satisfied for many exponentiated lag-windows as used in \cite%
{SPJ06,SPJ07} and \cite{SPJ11}. For more discussion of Assumption \ref{AW}
see \cite{PP2016}.

Assigning the test statistic $T_{w}$ the value zero on the set where $\hat{%
\Omega}_{w}\left( y\right) $ is singular is, of course, arbitrary. Given
Assumption \ref{AW}, the set where $\det \hat{\Omega}_{w}\left( y\right) =0$
holds can be shown to be%
\begin{equation*}
\mathsf{B}=\left\{ y\in \mathbb{R}^{n}:\func{rank}(B(y))<q\right\} ,
\end{equation*}%
where 
\begin{eqnarray}
B(y) &=&R(X^{\prime }X)^{-1}X^{\prime }\limfunc{diag}\left( \hat{u}%
_{1}(y),\ldots ,\hat{u}_{n}(y)\right)  \notag \\
&=&R(X^{\prime }X)^{-1}X^{\prime }\limfunc{diag}\left( e_{1}^{\prime }(n)\Pi
_{\limfunc{span}(X)^{\bot }}y,\ldots ,e_{n}^{\prime }(n)\Pi _{\limfunc{span}%
(X)^{\bot }}y\right) ,  \label{Def_B}
\end{eqnarray}%
see Lemma 3.1 of \cite{PP2016}. Clearly, $\func{span}(X)\subseteq \mathsf{B}$
and $\mathsf{B}+\func{span}(X)=\mathsf{B}$ always hold. Furthermore, Lemma
3.1 in \cite{PP2016} shows that the set $\mathsf{B}$ is a $\lambda _{\mathbb{%
R}^{n}}$-null set if and only if Assumption \ref{R_and_X} given below holds,
and is the entire space $\mathbb{R}^{n}$ otherwise. It thus transpires that
under Assumptions \ref{AW} and \ref{R_and_X} the chosen assignment is
irrelevant for size (and power) properties of the test (since all relevant
distributions $P_{\mu ,\sigma ^{2}\Sigma (f)}$ are absolutely continuous
w.r.t. $\lambda _{\mathbb{R}^{n}}$ due to the fact that every $\Sigma (f)$
is positive definite); the case where Assumption \ref{R_and_X} is violated
is hopeless for autocorrelation robust tests based on $T_{w}$ in that these
then break down in the sense that the quadratic form appearing in (\ref{tlr}%
) is not defined for any $y\in \mathbb{R}^{n}$ (as then $\mathsf{B}=\mathbb{R%
}^{n}$). [If one insists on using $T_{w}$ as defined by (\ref{tlr}) in this
case, $T_{w}$ then reduces to the trivial test statistic that is identically
zero, an uninteresting situation.\footnote{%
Of course, size control is then trivially possible, but leads to a test that
never rejects.}] We stress that the subsequently given Assumption \ref%
{R_and_X} can be readily checked in any given application and is not very
restrictive since -- for given restriction matrix $R$ -- it holds
generically in $X$.

\begin{assumption}
\label{R_and_X}Let $1\leq i_{1}<\ldots <i_{s}\leq n$ denote all the indices
for which $e_{i_{j}}(n)\in \limfunc{span}(X)$ holds where $e_{j}(n)$ denotes
the $j$-th standard basis vector in $\mathbb{R}^{n}$. If no such index
exists, set $s=0$. Let $X^{\prime }\left( \lnot (i_{1},\ldots i_{s})\right) $
denote the matrix which is obtained from $X^{\prime }$ by deleting all
columns with indices $i_{j}$, $1\leq i_{1}<\ldots <i_{s}\leq n$ (if $s=0$ no
column is deleted). Then $\limfunc{rank}\left( R(X^{\prime }X)^{-1}X^{\prime
}\left( \lnot (i_{1},\ldots i_{s})\right) \right) =q$ holds.
\end{assumption}

For later use we note that under Assumption \ref{R_and_X} the set $\mathsf{B}
$ not only is a $\lambda _{\mathbb{R}^{n}}$-null set, but is a finite union
of proper linear subspaces of $\mathbb{R}^{n}$, see Lemma \ref{Fact2} in
Section \ref{suff_con}. Also note that if $\mathsf{B}=\func{span}(X)$ holds,
then Assumption \ref{R_and_X} must hold (since $\func{span}(X)$ is a $%
\lambda _{\mathbb{R}^{n}}$-null set due to $k<n$).

The test statistic $T_{w}$ defined in (\ref{tlr}) is based on a long-run
covariance estimator for $\mathbf{u}_{t}x_{t\mathbf{\cdot }}^{\prime }$.%
\footnote{%
This test statistic makes sense for stochastic as well as nonstochastic
regressors under appropriate assumptions. Note that the asymptotic theory
for $T_{w}$ in a bandwidth-to zero as well as in a fixed-bandwidth scenario
relies only on limit theorems for partial sums of $\mathbf{u}_{t}x_{t\mathbf{%
\cdot }}^{\prime }$ and of $x_{t\mathbf{\cdot }}^{\prime }x_{t\mathbf{\cdot }%
}$ which are satisfied for various classes of stochastic as well as
nonstochastic regressors.} For nonstochastic regressors (as considered here)
an alternative is to use the weighted Eicker-test statistic $T_{E,\mathsf{W}%
} $ that is based on a long-run covariance estimator for $\mathbf{u}_{t}$.
The test statistic $T_{E,\mathsf{W}}$ is of the same form as given in (\ref%
{tlr}), but with the estimator $\hat{\Psi}_{E,\mathsf{W}}(y)=n^{-1}X^{\prime
}(\hat{K}(y)\bullet \mathsf{W})X$ being used instead of $\hat{\Psi}_{w}$,
where $\hat{K}(y)$ is the symmetric $n\times n$ Toeplitz matrix with $(i,j)$%
-th element given by $\hat{K}_{ij}(y)=n^{-1}\sum_{l=\left\vert
i-j\right\vert +1}^{n}\hat{u}_{l}(y)\hat{u}_{l-\left\vert i-j\right\vert
}(y) $, where $\mathsf{W}$ is an $n\times n$ symmetric and nonnegative
definite Toeplitz matrix of weights with ones on the main diagonal, and
where $\bullet $ denotes the Hadamard product. The (unweighted) Eicker-test
statistic $T_{E}$ discussed on pp.~283-284 of \cite{PP2016} corresponds to
the case where $\mathsf{W}$ is the matrix of all ones. As discussed in that
reference, the corresponding matrix $\hat{\Omega}_{E}(y)$ is always
nonnegative definite and the set where $\hat{\Omega}_{E}(y)$ is singular is
given by $\func{span}(X)$, which is a $\lambda _{\mathbb{R}^{n}}$-null set
by our maintained assumption that $k<n$ holds. By a result of Schur (see
Theorem 3.1 in \cite{Horn}) it is easy to see that all this is also true in
the weighted case, i.e., for $\hat{\Omega}_{E,\mathsf{W}}(y)=nR(X^{\prime
}X)^{-1}\hat{\Psi}_{E,\mathsf{W}}(y)(X^{\prime }X)^{-1}R^{\prime }$.

\subsection{Results on size control\label{main}}

To state the main result of this section we need to introduce some further
notation: Let $\omega \in \lbrack 0,\pi ]$ and let $s\geq 0$ be an integer.
Define $E_{n,s}(\omega )$ as the $n\times 2$-dimensional matrix with $j$-th
row equal to $\left( j^{s}\cos (j\omega ),j^{s}\sin (j\omega )\right) $.

\begin{definition}
\label{subspaces} Given a linear subspace $\mathcal{L}$ of $\mathbb{R}^{n}$
with $\dim (\mathcal{L})<n$, define for every $\omega \in \lbrack 0,\pi ]$%
\begin{equation}
\rho (\omega ,\mathcal{L})=\min \{s\in \mathbb{N\cup \{}0\mathbb{\}}:%
\limfunc{span}(E_{n,s}(\omega ))\not\subseteq \mathcal{L}\}.  \label{def_rho}
\end{equation}
\end{definition}

We note that the set on the right-hand side of (\ref{def_rho}) is nonempty
for every $\omega \in \lbrack 0,\pi ]$, and hence $\rho $ is well-defined
and takes values in $\mathbb{N\cup \{}0\mathbb{\}}$. Furthermore, $\rho
(\omega ,\mathcal{L})>0$ holds at most for finitely many $\omega \in \lbrack
0,\pi ]$. See Appendix \ref{App A} for a proof of these claims. We denote by 
$\mathfrak{M}_{0}^{lin}$ the linear space parallel to the affine space $%
\mathfrak{M}_{0}$, i.e., $\mathfrak{M}_{0}^{lin}=\mathfrak{M}_{0}-\mu _{0}$
for $\mu _{0}\in \mathfrak{M}_{0}$ (clearly, this does not depend on the
choice of $\mu _{0}\in \mathfrak{M}_{0}$).

\subsubsection{Result for $T_{w}$\label{w}}

Our first result concerning size control is given next and is an immediate
consequence of Theorem \ref{Ctoep_corr} in Section \ref{szctrl_3}. This
result is given for the test statistic $T_{w}$. A similar result for the
weighted Eicker-test statistic $T_{E,\mathsf{W}}$ is given in Subsection \ref%
{eicker}.

\begin{theorem}
\label{HAR_F_all}\footnote{%
Condition (\ref{non-incl}) clearly implies that the set $\mathsf{B}$ is a
proper subset of $\mathbb{R}^{n}$ and thus implies Assumption \ref{R_and_X}.
Hence, we could have dropped this assumption from the formulation of the
theorem. For clarity of presentation we have, however, chosen to explicitly
mention Assumption \ref{R_and_X}$.$} Suppose Assumptions \ref{AW} and \ref%
{R_and_X} are satisfied and $T_{w}$ is defined by (\ref{tlr}). Then for
every $0<\alpha <1$ there exists a real number $C(\alpha )$ such that%
\begin{equation}
\sup_{\mu _{0}\in \mathfrak{M}_{0}}\sup_{0<\sigma ^{2}<\infty }\sup_{f\in 
\mathfrak{F}_{\mathrm{all}}}P_{\mu _{0},\sigma ^{2}\Sigma (f)}(T_{w}\geq
C(\alpha ))\leq \alpha  \label{size-control_1}
\end{equation}%
holds, provided that 
\begin{equation}
\limfunc{span}\left( E_{n,\rho (\gamma ,\mathfrak{M}_{0}^{lin})}(\gamma
)\right) \not\subseteq \mathsf{B}\text{ \ \ for every \ }\gamma \in \lbrack
0,\pi ].  \label{non-incl}
\end{equation}%
In case the set $\mathsf{B}$ coincides with $\func{span}(X)$, condition (\ref%
{non-incl}) can equivalently be expressed as 
\begin{equation}
\func{rank}\left( X,E_{n,\rho (\gamma ,\mathfrak{M}_{0}^{lin})}(\gamma
)\right) >k\text{ \ \ for every \ }\gamma \in \lbrack 0,\pi ].
\label{non-incl-rank}
\end{equation}%
Furthermore, under the same condition (\ref{non-incl}) even equality can be
achieved in (\ref{size-control_1}) by a proper choice of $C(\alpha )$,
provided $\alpha \in (0,\alpha ^{\ast }]\cap (0,1)$, where $\alpha ^{\ast }$
is defined in (\ref{alpha*}) (with $T$ replaced by $T_{w}$) further below.
\end{theorem}

It turns out, see Subsection \ref{gen_size_contol} below, that for many
combinations of design matrices $X$ and restriction matrices $R$ the set $%
\mathsf{B}$ coincides with $\func{span}(X)$, and hence (\ref{non-incl})
reduces to the simpler rank condition (\ref{non-incl-rank}). Furthermore,
although a trivial observation, it should be kept in mind that $C(\alpha )$
depends not only on $\alpha $ but also on the testing problem at hand (i.e.,
on $X$, $R$, $r$, and the covariance model (here $\mathfrak{C}(\mathfrak{F}_{%
\mathrm{all}})$) as well as on the choice of test statistic (here on the
choice of weights $w(j,n)$); but see Remark \ref{delta} further below. We
furthermore note that, under the conditions of Theorem \ref{HAR_F_all}, a
smallest critical value $C_{\Diamond }(\alpha )$ exists that satisfies (\ref%
{size-control_1}) and that this critical value achieves equality in (\ref%
{size-control_1}) provided $\alpha \in (0,\alpha ^{\ast }]\cap (0,1)$. This
follows from Remark \ref{exist} in Section \ref{szctrl_3}.

The preceding theorem provides simple sufficient conditions under which size
control of commonly used autocorrelation robust tests is possible over the
class $\mathfrak{F}_{\mathrm{all}}$, and hence a fortiori over \emph{any} $%
\mathfrak{F}\subseteq \mathfrak{F}_{\mathrm{all}}$.\footnote{%
In particular, $\mathfrak{F}=\mathfrak{F}_{\mathrm{AR(}p\mathrm{)}}$ or $%
\mathfrak{F}=\mathfrak{F}_{\mathrm{ARMA(}p,q\mathrm{)}}$ is covered for
arbitrary $p$ ($(p,q)$, respectively), where $\mathfrak{F}_{\mathrm{ARMA(}p,q%
\mathrm{)}}$ denotes the set of \emph{all }normalized spectral densities
corresponding to stationary autoregressive moving average processes with
autoregressive (moving average) order not larger than $p$ ($q$,
respectively). Note that the definition of $\mathfrak{F}_{\mathrm{AR(}p%
\mathrm{)}}$ and $\mathfrak{F}_{\mathrm{ARMA(}p,q\mathrm{)}}$ does \emph{not 
}require the zeros of the autoregressive or moving average polynomial to be
bounded away from the unit circle.} Of course, critical values achieving
size control over a subset $\mathfrak{F}$ (i.e., critical values satisfying (%
\ref{size-control_1}) with $\mathfrak{F}_{\mathrm{all}}$ replaced by $%
\mathfrak{F}$) may be smaller than critical values achieving size control
over $\mathfrak{F}_{\mathrm{all}}$.

As mentioned before, the theorem is a special case of Theorem \ref%
{Ctoep_corr} in Section \ref{szctrl_3}, which provides more refined
sufficient conditions for the possibility of size control at the expense of
more complicated conditions and notation. Theorem \ref{Ctoep_corr} is
especially of importance if one is interested in sufficient conditions for
the possibility of size control over classes $\mathfrak{F}$ that are much
smaller than $\mathfrak{F}_{\mathrm{all}}$, since then the conditions in
Theorem \ref{HAR_F_all} may be unnecessarily restrictive; cf. Remarks \ref%
{crudeness} and \ref{bounded away} given below, but see also Remark \ref%
{AR(2)}. We further note that Theorem \ref{Ctoep_corr} is in turn a
corollary to Theorem \ref{Ctoep} in Section \ref{szctrl_3}, which applies to
a much larger class of test statistics than the one considered in the
present section. Finally, we note that the sufficient conditions in Theorem %
\ref{HAR_F_all}, as well as the sufficient conditions provided in the
theorems in Section \ref{szctrl_3}, are also necessary for size control in a
sense made precise in \cite{PP4}, provided that $\mathfrak{F}$ is rich
enough to encompass $\mathfrak{F}_{\mathrm{AR(}2\mathrm{)}}$.

\begin{remark}
\label{delta} \emph{(Independence of value of }$r$\emph{) }(i) Since $%
\mathfrak{M}_{0}^{lin}$ does not depend on the value of $r$, the sufficient
conditions in Theorem \ref{HAR_F_all} -- while depending on $X$ and $R$ --
do not depend on the value of $r$.

(ii) For a large class of test statistics (including $T_{w}$ considered
here) the size of the corresponding tests, and hence the size-controlling
critical values $C(\alpha )$, do not depend on the value of $r$, see Lemma %
\ref{neu} in Section \ref{two}. This observation is of some importance, as
it allows one easily to obtain confidence sets for $R\beta $ by
\textquotedblleft inverting\textquotedblright\ the test without the need of
recomputing the critical value for every value of $r$. [Of course, it is
here understood that the weights $w$ are not related to the value of $r$.]
\end{remark}

In the subsequent examples we apply Theorem \ref{HAR_F_all} and discuss
simple sufficient conditions under which size control over $\mathfrak{F}_{%
\mathrm{all}}$ (and hence a fortiori over any $\mathfrak{F}\subseteq 
\mathfrak{F}_{\mathrm{all}}$) is possible. We also discuss that these
sufficient conditions are generically satisfied in the appropriate universe
of design matrices. For these examples we always assume that $T_{w}$ is
given by (\ref{tlr}) with Assumption \ref{AW} being satisfied; also recall
that the $n\times k$-dimensional design matrix $X$ always has rank $k$, that 
$1\leq $ $k<n$ holds, and that the restriction matrix $R$ is always $q\times
k$-dimensional and has rank equal to $q\geq 1$.

\begin{example}
\label{exa1}Assume that the design matrix $X$ satisfies $\func{rank}%
(X,E_{n,0}(\gamma ))>k$ for every $\gamma \in \lbrack 0,\pi ]$ (which is
tantamount to $\limfunc{span}(E_{n,0}(\gamma ))\not\subseteq \limfunc{span}%
(X)$ for every $\gamma \in \lbrack 0,\pi ]$ since $\func{rank}(X)=k$ by our
assumptions). It is then easy to see that $\rho (\gamma ,\mathfrak{M}%
_{0}^{lin})=0$ for every $\gamma \in \lbrack 0,\pi ]$. Hence, the conditions
enabling size control over $\mathfrak{F}_{\mathrm{all}}$ in Theorem \ref%
{HAR_F_all} are then all satisfied if $X$ and the restriction matrix $R$ are
such that $\mathsf{B}$ coincides with $\limfunc{span}(X)$. In particular, it
follows easily from Theorem \ref{genericity_cor} further below that the just
mentioned rank condition as well as the condition on $\mathsf{B}$ is -- for
given restriction matrix $R$ -- generically satisfied in the set of all $%
n\times k$ matrices. [For the last claim use Theorem \ref{genericity_cor}
with $F$ the empty matrix, and observe that the conditions of that theorem
are clearly satisfied. Also observe that $\rho _{F}(\gamma )=0$ holds for
every $\gamma \in \lbrack 0,\pi ]$.]
\end{example}

The next example treats the case where an intercept is required to be in the
model.

\begin{example}
\label{exa2}Assume that the design matrix $X$ contains an intercept, in the
sense that $X$ has $e_{+}$ as its first column ($e_{+}$ is defined in
Section \ref{frame}). We also assume $k\geq 2$ and write $X=(e_{+},\tilde{X}%
) $. We are interested in testing restrictions that do not involve the
intercept, i.e., $R=(0,\tilde{R}),$ with $\tilde{R}$ of dimension $q\times
(k-1)$. [Recall that autocorrelation robust testing of restrictions that
involve the intercept by means of $T_{w}$ is futile whenever $\mathfrak{F}%
\supseteq \mathfrak{F}_{\mathrm{AR(1)}}$, see \cite{PP2016}, Example 3.1.]
It is now obvious that the rank condition on $(X,E_{n,0}(\gamma ))$ in the
preceding example is violated for $\gamma =0$, and hence the conclusions of
the preceding example do not apply in the situation considered here.
However, assume now instead that $\func{rank}(X,E_{n,0}(\gamma ))>k$ for
every $\gamma \in (0,\pi ]$ and that $\func{rank}(X,E_{n,1}(\gamma ))>k$ for 
$\gamma =0$ hold. It then follows that $\rho (\gamma ,\mathfrak{M}%
_{0}^{lin})=0$ for every $\gamma \in (0,\pi ]$ and that $\rho (\gamma ,%
\mathfrak{M}_{0}^{lin})=1$ for $\gamma =0$. Then again the conditions for
size control over $\mathfrak{F}_{\mathrm{all}}$ in Theorem \ref{HAR_F_all}
are seen to be satisfied if $X$ and the restriction matrix $R$ are such that 
$\mathsf{B}$ coincides with $\limfunc{span}(X)$. Similarly as in the
preceding example, it follows from Theorem \ref{genericity_cor} further
below that the new rank conditions as well as the condition on $\mathsf{B}$
are -- for given restriction matrix $R=(0,\tilde{R})$ -- generically
satisfied in the set of all $n\times k$ matrices of the form $X=(e_{+},%
\tilde{X})$; for a proof of this claim see Appendix \ref{App A}.
\end{example}

A completely analogous discussion as in the preceding example can be given
for the case where $X=(e_{-},\tilde{X})$ and is omitted.

\begin{example}
\label{exa3}Assume that the design matrix $X$ contains $e_{+}$ as its first
and $e_{-}$ (defined in Section \ref{frame}) as its second column, i.e., $%
X=(e_{+},e_{-},\tilde{X})$, and assume $k\geq 3$. Further assume that $R=(0,%
\tilde{R}),$ where now $\tilde{R}$ is of dimension $q\times (k-2)$. [Recall
that testing restrictions that involve the intercept by means of $T_{w}$ is
futile whenever $\mathfrak{F}\supseteq \mathfrak{F}_{\mathrm{AR(1)}}$, see 
\cite{PP2016}, Example 3.1, and a similar remark applies to restrictions
involving the coefficient of $e_{-}$.] Similar as before, the rank condition
on $(X,E_{n,0}(\gamma ))$ in the preceding Example \ref{exa1} is now
violated for $\gamma \in \left\{ 0,\pi \right\} $, and so is the rank
condition in Example \ref{exa2} for $\gamma =\pi $. However, if we require
instead that (i) $\func{rank}(X,E_{n,0}(\gamma ))>k$ for every $\gamma \in
(0,\pi )$, (ii) $\func{rank}(X,E_{n,1}(\gamma ))>k$ \ for every $\gamma \in
\left\{ 0,\pi \right\} $, and that (iii) $X$ and $R$ are such that $\mathsf{B%
}$ coincides with $\limfunc{span}(X)$, then it is not difficult to see that $%
\rho (\gamma ,\mathfrak{M}_{0}^{lin})=0$ for $\gamma \in (0,\pi )$, $\rho
(\gamma ,\mathfrak{M}_{0}^{lin})=1$ for $\gamma \in \left\{ 0,\pi \right\} $%
, and that the sufficient conditions in Theorem \ref{HAR_F_all} are
satisfied, implying that size control over $\mathfrak{F}_{\mathrm{all}}$ is
possible. Similarly as in the preceding examples, it follows from Theorem %
\ref{genericity_cor} further below that conditions (i)-(iii) are -- for
given restriction matrix $R=(0,\tilde{R})$ -- generically satisfied in the
set of all $n\times k$ matrices of the form $X=(e_{+},e_{-},\tilde{X})$,
provided only that the mild condition $q\leq (n/2)-1$ is satisfied; for a
proof of this claim see Appendix \ref{App A}.
\end{example}

In the next remark we exemplarily discuss how the sufficient conditions in
the preceding examples obtained from Theorem \ref{HAR_F_all} can be
weakened, if one is concerned with size control only over a set of spectral
densities much smaller than $\mathfrak{F}_{\mathrm{all}}$. But see also
Remark \ref{AR(2)} further below.

\begin{remark}
\label{crudeness} Suppose we are interested in size control of tests of the
form (\ref{tlr}) but now only over the much smaller set $\mathfrak{F}_{%
\mathrm{AR(}1\mathrm{)}}$.

(i) Assume that the design matrix $X$ satisfies $\func{rank}%
(X,E_{n,0}(\gamma ))>k$ for every $\gamma \in \left\{ 0,\pi \right\} $
(which is tantamount to $e_{+}\notin \limfunc{span}(X)$ and $e_{-}\notin 
\limfunc{span}(X)$ as we always assume $\func{rank}(X)=k$). Assume also that 
$X$ and the restriction matrix $R$ are such that $\mathsf{B}$ coincides with 
$\limfunc{span}(X)$. Then Theorem \ref{Ctoep_corr} in Section \ref{szctrl_3}
implies that size control over $\mathfrak{F}_{\mathrm{AR(}1\mathrm{)}}$ is
possible. [This is so since the set $\mathbb{S}(\mathfrak{F}_{\mathrm{AR(}1%
\mathrm{)}},\mathfrak{M}_{0}^{lin})$ appearing in that theorem is $\left\{
\left\{ 0\right\} ,\left\{ \pi \right\} \right\} $ as shown in Example \ref%
{chsingularAR(1)}. Furthermore, it is easy to see that here $\rho (\gamma ,%
\mathfrak{M}_{0}^{lin})=0$ for every $\gamma \in \left\{ 0,\pi \right\} $
holds.]

(ii) Assume that $X=(e_{+},\tilde{X})$ with $k\geq 2$, that $R=(0,\tilde{R}%
), $ with $\tilde{R}$ of dimension $q\times (k-1)$, and that $\mathsf{B}$
coincides with $\limfunc{span}(X)$. If $\func{rank}(X,E_{n,0}(\pi ))>k$
(i.e., if $e_{-}\notin \limfunc{span}(X)$), then size control over $%
\mathfrak{F}_{\mathrm{AR(}1\mathrm{)}}$ is possible. [This follows from
Theorem \ref{Ctoep_corr} since $\mathbb{S}(\mathfrak{F}_{\mathrm{AR(}1%
\mathrm{)}},\mathfrak{M}_{0}^{lin})$ now equals $\left\{ \left\{ \pi
\right\} \right\} $, see Example \ref{chsingularAR(1)}, and since $\rho (\pi
,\mathfrak{M}_{0}^{lin})=0$ holds.] The case where $X=(e_{-},X)$ can be
treated analogously, and we do not provide the details.

(iii) Assume $X=(e_{+},e_{-},\tilde{X})$ with $k\geq 3$, that $R=(0,\tilde{R}%
),$ with $\tilde{R}$ of dimension $q\times (k-2)$. Then size control over $%
\mathfrak{F}_{\mathrm{AR(}1\mathrm{)}}$ is possible without any further
conditions. [Again this follows from Theorem \ref{Ctoep_corr} upon observing
that now $\mathbb{S}(\mathfrak{F}_{\mathrm{AR(}1\mathrm{)}},\mathfrak{M}%
_{0}^{lin})$ is empty in view of Example \ref{chsingularAR(1)}.] We note
that this result is also in line with Theorem 3.7 of \cite{PP2016}.
\end{remark}

We proceed to providing two more examples illustrating Theorem \ref%
{HAR_F_all}. The first one concerns the case where a linear trend is present
in the model.

\begin{example}
\label{exa4} Assume $k\geq 3$ and that the design matrix $X$ contains $e_{+}$
as its first column and the vector $v=(1,2,\ldots ,n)^{\prime }$ as its
second column, i.e., $X=(e_{+},v,\tilde{X})$. That is, the linear model
contains a linear trend. Suppose $R=(0,\tilde{R}),$ where $\tilde{R}$ is of
dimension $q\times (k-2)$, i.e., the restriction to be tested does not
involve the coefficients appearing in the trend. This is, of course, a
special case of the model considered in Example \ref{exa2}. However, it is
plain that the condition $\func{rank}(X,E_{n,1}(\gamma ))>k$ for $\gamma =0$%
, used in that example, is not satisfied in the present context (as $%
E_{n,1}(0)=(v:0)$). A simple set of sufficient conditions for size control
in the present example is now as follows: (i) $\func{rank}(X,E_{n,0}(\gamma
))>k$ for every $\gamma \in (0,\pi ]$, (ii) $\func{rank}(X,E_{n,2}(\gamma
))>k$ for $\gamma =0$, and that (iii) $X$ and $R$ are such that $\mathsf{B}$
coincides with $\limfunc{span}(X)$. It is then not difficult to see that $%
\rho (\gamma ,\mathfrak{M}_{0}^{lin})=0$ for $\gamma \in (0,\pi ]$, $\rho
(\gamma ,\mathfrak{M}_{0}^{lin})=2$ for $\gamma =0$, and that the sufficient
conditions in Theorem \ref{HAR_F_all} are satisfied, implying that size
control over $\mathfrak{F}_{\mathrm{all}}$ is possible. Again, it follows
from Theorem \ref{genericity_cor} further below that conditions (i)-(iii)
are -- for given restriction matrix $R=(0,\tilde{R})$ -- generically
satisfied in the set of all $n\times k$ matrices of the form $X=(e_{+},v,%
\tilde{X})$; for a proof of this claim see Appendix \ref{App A}.
\end{example}

The preceding example can easily be generalized to arbitrary polynomial
trends, but we abstain from providing the details. The last example
considers a model with a cyclical component.

\begin{example}
\label{exa5}Assume that the design matrix $X$ has the form $%
(e_{+},E_{n,0}(\gamma _{0}),\tilde{X})$ for some $\gamma _{0}\in (0,\pi )$
and that $k\geq 4$. That is, the model contains a cyclical component.
Suppose $R=(0,\tilde{R}),$ where $\tilde{R}$ is of dimension $q\times (k-3)$%
, i.e., the restriction to be tested does neither involve the intercept nor
the coefficients appearing in the cyclical component. While this is a
special case of the model considered in Example \ref{exa2}, it is also plain
that the conditions provided in that example do not work here (since $\func{%
rank}(X,E_{n,0}(\gamma ))>k$ is violated for $\gamma =\gamma _{0}$). A
simple set of sufficient conditions enabling size control in the present
example is now as follows: (i) $\func{rank}(X,E_{n,0}(\gamma ))>k$ for every 
$\gamma \in (0,\pi ]\backslash \left\{ \gamma _{0}\right\} $, (ii) $\func{%
rank}(X,E_{n,1}(\gamma ))>k$ for $\gamma =0$ as well as $\gamma =\gamma _{0}$%
, and that (iii) $X$ and $R$ are such that $\mathsf{B}$ coincides with $%
\limfunc{span}(X)$. It is then not difficult to see that $\rho (\gamma ,%
\mathfrak{M}_{0}^{lin})=0$ for $\gamma \in (0,\pi ]\backslash \left\{ \gamma
_{0}\right\} $, $\rho (\gamma ,\mathfrak{M}_{0}^{lin})=1$ for $\gamma =0$ as
well as $\gamma =\gamma _{0}$, and that the sufficient conditions in Theorem %
\ref{HAR_F_all} are satisfied, implying that size control over $\mathfrak{F}%
_{\mathrm{all}}$ is possible. Again, it follows from Theorem \ref%
{genericity_cor} further below that conditions (i)-(iii) are -- for given
restriction matrix $R=(0,\tilde{R})$ -- generically satisfied in the set of
all $n\times k$ matrices of the form $X=(e_{+},E_{n,0}(\gamma _{0}),\tilde{X}%
)$, provided only that the mild condition $q\leq (n/3)-1$ is satisfied; for
a proof of this claim see Appendix \ref{App A}.
\end{example}

\begin{remark}
\emph{(No size control) }The size control result given above, as well as the
more refined Theorem \ref{Ctoep_corr} in Section \ref{szctrl_3}, do -- for
example -- not apply to the following testing problems: (i) testing the
intercept if $\mathfrak{F}\supseteq \mathfrak{F}_{\mathrm{AR(}1\mathrm{)}}$,
(ii) testing the coefficient of the regressor $e_{-}$ if $\mathfrak{F}%
\supseteq \mathfrak{F}_{\mathrm{AR(}1\mathrm{)}}$, (iii) testing a
hypothesis regarding the coefficients of a linear trend appearing in the
model provided $\mathfrak{F}\supseteq \mathfrak{F}_{\mathrm{AR(}2\mathrm{)}}$%
, and (iv) testing a hypothesis regarding the coefficients of a cyclical
component appearing in the model provided $\mathfrak{F}\supseteq \mathfrak{F}%
_{\mathrm{AR(}2\mathrm{)}}$. In fact, it is known that the size of
autocorrelation robust tests based on $T_{w}$ is one in any of these testing
problems regardless of the choice of critical value. For testing problems
(i) and (ii) this follows from Theorem 3.3 of \cite{PP2016} (see also
Example 3.1 in that reference). For testing problem (iv) this follows from
Theorem 3.12 of \cite{PP2016}. For testing problem (iii), see Section 5 in 
\cite{PP4}. [These results are closely related to the fact that the
sufficient conditions for size control given in Theorem \ref{HAR_F_all} and
in Section \ref{szctrl_3} are in fact necessary in a certain sense; see \cite%
{PP4}.]
\end{remark}

\begin{remark}
\label{AR(2)} Suppose that in the context of Theorem \ref{HAR_F_all} we are
interested in size control over a set $\mathfrak{F}$ with $\mathfrak{F}%
\supseteq \mathfrak{F}_{\mathrm{AR(}2\mathrm{)}}$. It then follows from
Remark \ref{nec} (with $\mathcal{L}=\mathfrak{M}_{0}^{lin}$) and Remark \ref%
{nec_2} that the sufficient condition (\ref{non-incl}) given in Theorem \ref%
{HAR_F_all} is in fact equivalent to the more refined sufficient conditions
given in Part 1 of Theorem \ref{Ctoep_corr}.
\end{remark}

\begin{remark}
\label{bounded away} (i) Suppose $\mathfrak{F}\subseteq \mathfrak{F}_{%
\mathrm{ARMA(}p,q\mathrm{)}}$ consists only of normalized spectral densities
corresponding to stationary autoregressive moving average processes with the
property that the zeros of all the autoregressive polynomials are bounded
away from the unit circle in the complex plane by a fixed amount $\delta >0$%
, say. Then it is easy to see that $\mathfrak{F}\subseteq \mathfrak{F}_{%
\mathrm{all}}^{B}$\ for some $B<\infty $ holds, where $\mathfrak{F}_{\mathrm{%
all}}^{B}$ is defined in Example \ref{chsingularToepmin} in Section \ref%
{struct}. Theorem \ref{Ctoep_corr} together with Example \ref%
{chsingularToepmin} now shows that size control over $\mathfrak{F}$ is
possible even if the condition (\ref{non-incl}) in Theorem \ref{HAR_F_all}
is not satisfied. [In fact, this conclusion is true for any of the sets $%
\mathfrak{F}_{\mathrm{all}}^{B}$ themselves.] Note, however, that choosing $%
\delta $ small may nevertheless result in large critical values, especially,
if $X$ and $R$ are such that condition (\ref{non-incl}) is violated (and
thus we are not in general guaranteed that size control is possible over $%
\mathfrak{F}_{\mathrm{ARMA(}p,q\mathrm{)}}$).

(ii) More generally, size control in the setting of Theorem \ref{HAR_F_all}
is always possible (i.e., even when (\ref{non-incl}) is violated) if, e.g.,
the covariance model $\mathfrak{C}(\mathfrak{F})$ employed does not have any
singular limit points; cf. also Remarks \ref{star} and \ref{doublestar}
further below.
\end{remark}

\subsubsection{Result for $T_{E,\mathsf{W}}$\label{eicker}}

Here we give a result similar to Theorem \ref{HAR_F_all} but for the
weighted Eicker-test statistic $T_{E,\mathsf{W}}$. The theorem follows
immediately from Remark \ref{ext}(ii).

\begin{theorem}
\label{HAR_F_all_Eicker} Let $T_{E,\mathsf{W}}$ be the weighted Eicker-test
statistic where $\mathsf{W}$ is an $n\times n$ symmetric and nonnegative
definite Toeplitz matrix of weights with ones on the main diagonal. Then for
every $0<\alpha <1$ there exists a real number $C(\alpha )$ such that%
\begin{equation}
\sup_{\mu _{0}\in \mathfrak{M}_{0}}\sup_{0<\sigma ^{2}<\infty }\sup_{f\in 
\mathfrak{F}_{\mathrm{all}}}P_{\mu _{0},\sigma ^{2}\Sigma (f)}(T_{E,\mathsf{W%
}}\geq C(\alpha ))\leq \alpha  \label{size-control_1_Eicker}
\end{equation}%
holds, provided that 
\begin{equation}
\limfunc{span}\left( E_{n,\rho (\gamma ,\mathfrak{M}_{0}^{lin})}(\gamma
)\right) \not\subseteq \limfunc{span}(X)\text{ \ \ for every \ }\gamma \in
\lbrack 0,\pi ].  \label{non-incl_Eicker}
\end{equation}%
This can equivalently be expressed as 
\begin{equation}
\func{rank}\left( X,E_{n,\rho (\gamma ,\mathfrak{M}_{0}^{lin})}(\gamma
)\right) >k\text{ \ \ for every \ }\gamma \in \lbrack 0,\pi ].
\label{non-incl-rank_Eicker}
\end{equation}%
Furthermore, under the same condition (\ref{non-incl_Eicker}) even equality
can be achieved in (\ref{size-control_1_Eicker}) by a proper choice of $%
C(\alpha )$, provided $\alpha \in (0,\alpha ^{\ast }]\cap (0,1)$, where $%
\alpha ^{\ast }$ is defined in (\ref{alpha*}) (with $T$ replaced by $T_{E,%
\mathsf{W}}$) further below.
\end{theorem}

Mutatis mutandis, the entire discussion in Subsection \ref{w} following
Theorem \ref{HAR_F_all} also applies to Theorem \ref{HAR_F_all_Eicker}.

\subsection{Generic size control\label{gen_size_contol}}

In Theorem \ref{genericity_cor} below we now show that -- for given
restriction matrix $R$ -- the set of design matrices, for which the
conditions in Theorem \ref{HAR_F_all} (Theorem \ref{HAR_F_all_Eicker},
respectively) are satisfied and hence size control is possible, is generic.
We provide this genericity result for a variety of universes of design
matrices. For example, if $F$ in Theorem \ref{genericity_cor} is absent
(more precisely, corresponds to the empty matrix), the genericity holds in
the class of all $n\times k$ design matrices. If one is only interested in
regression models containing an intercept, then one sets $F=e_{+}$, and the
theorem delivers a genericity result in the subuniverse of all $n\times k$
design matrices that contain $e_{+}$ as its first column. In general, $F$
stands for that subset of columns of the design matrix that are a priori
held fixed in the genericity analysis. The subsequent theorem follows
immediately by combining Theorem \ref{HAR_F_all} (Theorem \ref%
{HAR_F_all_Eicker}, respectively) with Lemmata \ref{generic_lem_1} and \ref%
{generic_lem_2} in Appendix \ref{App A}. Recall that $\mathsf{B}$ as well as 
$\rho (\gamma ,\mathfrak{M}_{0}^{lin})$ depend on $X$ (and $R$), which,
however, is not shown in the notation; and that the relation $\mathsf{B}=%
\func{span}(X)$ implies that Assumption \ref{R_and_X} holds. Also recall
that $n>k$ is assumed throughout.

\begin{theorem}
\label{genericity_cor}Let $F$ be a given $n\times k_{F}$ matrix of rank $%
k_{F}$, where $0\leq k_{F}<k$ (with the convention that $F$ is the empty
matrix in case $k_{F}=0$, that the rank of the empty matrix is zero, and
that its span is $\left\{ 0\right\} $). Assume that the given $q\times k$
restriction matrix $R$ of rank $q$ has the form $(0,\tilde{R})$ with $\tilde{%
R}$ a $q\times (k-k_{F})$ matrix. Suppose the columns of $F$ and $%
e_{i_{1}}(n),\ldots ,e_{i_{q}}(n)$ are linearly independent for every choice
of $1\leq i_{1}<\ldots <i_{q}\leq n$. Furthermore, suppose that (i) $n>k+2$
holds, or (ii) $\func{rank}(F,E_{n,0}(\gamma ^{\ast }))=k_{F}+2$ holds for
some $\gamma ^{\ast }\in (0,\pi )$. Then the following holds generically for
design matrices $X$ of the form $(F,\tilde{X})$ (i.e., holds on the
complement of a $\lambda _{\mathbb{R}^{n\times (k-k_{F})}}$-null set of
matrices $\tilde{X}$):

\begin{enumerate}
\item $X=(F,\tilde{X})$ has rank $k$.

\item $\mathsf{B}=\func{span}(X)$.

\item Assumption \ref{R_and_X} is satisfied.

\item $\rho (\gamma ,\mathfrak{M}_{0}^{lin})=\rho _{F}(\gamma )$ holds for
every\ $\gamma \in \lbrack 0,\pi ]$ where $\rho _{F}(\gamma )=\rho (\gamma ,%
\func{span}(F))$.

\item Conditions (\ref{non-incl}), (\ref{non-incl-rank}), (\ref%
{non-incl_Eicker}), and (\ref{non-incl-rank_Eicker}) hold.

\item Suppose $T_{w}$ is defined by (\ref{tlr}) with Assumption \ref{AW}
being satisfied. Then for every $0<\alpha <1$ there exists a real number $%
C(\alpha )$ such that (\ref{size-control_1}) holds; and if $\alpha \in
(0,\alpha ^{\ast }]\cap (0,1)$ even equality can be achieved in (\ref%
{size-control_1}), where $\alpha ^{\ast }$ is as in Theorem \ref{HAR_F_all}.

\item Let $T_{E,\mathsf{W}}$ be the weighted Eicker-test statistic where $%
\mathsf{W}$ is an $n\times n$ symmetric and nonnegative definite Toeplitz
matrix of weights with ones on the main diagonal. Then for every $0<\alpha
<1 $ there exists a real number $C(\alpha )$ such that (\ref%
{size-control_1_Eicker}) holds; and if $\alpha \in (0,\alpha ^{\ast }]\cap
(0,1)$ even equality can be achieved in (\ref{size-control_1_Eicker}), where 
$\alpha ^{\ast }$ is as in Theorem \ref{HAR_F_all_Eicker}.
\end{enumerate}
\end{theorem}

Note that neither the assumptions nor the first five conclusions nor the $%
\lambda _{\mathbb{R}^{n\times (k-k_{F})}}$-null set depend on the value of $%
r $ at all (this is obvious upon noting that $\rho (\gamma ,\mathfrak{M}%
_{0}^{lin})$ depends on $\mathfrak{M}_{0}$ only via $\mathfrak{M}_{0}^{lin}$%
, which is independent of the value of $r$). For the last two conclusions
note that they hold whatever the value of $r$ is.

\begin{remark}
(i) Theorem \ref{genericity_cor} assumes that the restrictions to be tested
do not involve the coefficients of the regressors corresponding to the
columns of $F$. This assumption can be traded-off with an assumption
ensuring that no singular limit point of the covariance model $\mathfrak{C}(%
\mathfrak{F}_{\mathrm{all}})$ concentrates on the space spanned by the
columns of $F$. We abstain from providing such results.

(ii) In case $k_{F}<\lfloor n/2\rfloor $ the rank-condition in (ii) of the
theorem is always satisfied: Suppose not, then $E_{n,0}(\gamma )v(\gamma
)\in \func{span}(F)$ for some nonzero vector $v(\gamma )$ and for every $%
\gamma \in (0,\pi )$. Choose $\gamma _{i}\in (0,\pi )$, $i=1,\ldots ,\lfloor
n/2\rfloor $, all $\gamma _{i}$ being different. By Lemma \ref{fullrank} in
Appendix \ref{App C}, it follows that the corresponding collection of
vectors $E_{n,0}(\gamma _{i})v(\gamma _{i})$ is linearly independent,
implying that $k_{F}\geq \lfloor n/2\rfloor $, a contradiction. [In case $%
k_{F}>n/2$ examples can be given where this condition is not satisfied.]
\end{remark}

\subsection{Comments on power properties\label{power}}

Classical autocorrelation robust tests can have, in fact not infrequently
will have, infimal power equal to zero if the underlying set $\mathfrak{F}$
is sufficiently rich; cf. Theorem 3.3 and Corollary 5.17 in \cite{PP2016} as
well as Lemma \ref{L4} in Section \ref{general}. In the special case where $%
\mathfrak{F}=\mathfrak{F}_{\mathrm{AR(}1\mathrm{)}}$, it has been shown in 
\cite{PP2016} and \cite{Prein2014b} how adjusted tests can be constructed
that have correct size and at the same time do not suffer from infimal power
being zero. In a companion paper, which builds on the results of the present
paper, we investigate power properties in more detail and provide
adjustments to the test statistics $T_{w}$ and $T_{E,\mathsf{W}}$ that
typically lead to improvements in power properties, at least over certain
important subsets of $\mathfrak{F}_{\mathrm{all}}$, while retaining size
control over $\mathfrak{F}_{\mathrm{all}}$ as in Theorems \ref{HAR_F_all}
and \ref{HAR_F_all_Eicker}.

We also note here that despite what has just been said, one can show for
autocorrelation robust tests based on $T_{w}$ or $T_{E,\mathsf{W}}$ (size
corrected or not) that power goes to one as one moves away from the null
hypothesis along sequences of the following form: Let $(\mu _{l},\sigma
_{l}^{2},f_{l})$ be such that $\mu _{l}$ moves further and further away from 
$\mathfrak{M}_{0}$ (the affine space of means described by the restrictions $%
R\beta =r$) in an orthogonal direction, where $\sigma _{l}^{2}$ converges to
some finite and positive $\sigma ^{2}$, and $f_{l}$ is such that $\Sigma
(f_{l})$ converges to a positive definite matrix. Note, however, that this
result rules out sequences $f_{l}$ for which $\Sigma (f_{l})$ degenerates as 
$l\rightarrow \infty $.

\subsection{Extensions to other test statistics\label{extension}}

\textbf{A. }\emph{(Adjusted tests)} In \cite{PP2016} we have discussed
adjusted autocorrelation robust tests, which are nothing else than standard
autocorrelation robust tests but computed from an augmented regression model
that contains not only the regressors in $X$, but also strategically chosen
auxiliary regressors. The above results can easily accommodate adjusted
tests: Simply view the augmented model as the true model. Since the adjusted
test then is just a standard autocorrelation robust test in the augmented
model, the above results can be applied. Note that the null hypothesis in
the augmented model encompasses the null hypothesis in the original model,
hence size control over the null hypothesis in the augmented model certainly
implies size control in the originally given model. For more discussion see
Theorem 3.8, Proposition 5.23, and especially Remark 5.24(iii) in \cite%
{PP2016}.

\textbf{B.} \emph{(Tests based on general quadratic covariance estimators)}
The test statistics $T_{GQ}$ in this class are of the form (\ref{tlr}) but
-- instead of $\hat{\Psi}_{w}$ -- use the estimator%
\begin{equation}
\hat{\Psi}_{GQ}(y)=\sum_{t,s=1}^{n}w(t,s;n)\hat{v}_{t}(y)\hat{v}%
_{s}(y)^{\prime }  \label{GQ-estimator}
\end{equation}%
for $y\in \mathbb{R}^{n}$, where the $n\times n$ weighting matrix $\mathcal{W%
}_{n}^{\ast }=(w(t,s;n))_{t,s}$ is symmetric and data-independent. For some
background on this more general class of estimators see Section 3.2.1 of 
\cite{PP2016}. Note that $\hat{\Psi}_{GQ}(y)$, and thus the corresponding $%
\hat{\Omega}_{GQ}(y)=nR(X^{\prime }X)^{-1}\hat{\Psi}_{GQ}(y)(X^{\prime
}X)^{-1}R^{\prime }$, is nonnegative definite for every $y\in \mathbb{R}^{n}$
provided $\mathcal{W}_{n}^{\ast }$ is nonnegative definite. In the important
case where $\mathcal{W}_{n}^{\ast }$ is additionally positive definite,
inspection of the proofs (together with Lemma 3.11 of \cite{PP2016}) shows
that all results given above for the test statistic $T_{w}$ based on $\hat{%
\Psi}_{w}$ remain valid for $T_{GQ}$ (provided Assumption \ref{AW} is
replaced by the assumption on $\mathcal{W}_{n}^{\ast }$ made here including
positive definiteness of $\mathcal{W}_{n}^{\ast }$); see also Remark \ref%
{ext}(i). In the case where $\mathcal{W}_{n}^{\ast }$ is nonnegative
definite, but not positive definite, conditions under which size control is
possible can be derived from Theorem \ref{Ctoep_nonsphericity} in Section %
\ref{szctrl_3} (with the help of Lemma 3.11 of \cite{PP2016}); we do not
provide details. In fact, even cases where $\mathcal{W}_{n}^{\ast }$ is not
nonnegative definite can be accommodated by this result under appropriate
conditions.

\textbf{C.} A referee has suggested a test statistic $T_{ref}$ which is of
the form (\ref{tlr}), but where $\hat{\Omega}_{w}(y)$ is replaced by $\hat{%
\omega}_{w}(y)R(X^{\prime }X)^{-1}R^{\prime }$. Here $\hat{\omega}%
_{w}(y)=\sum\nolimits_{j=-(n-1)}^{n-1}w(j,n)\hat{K}_{j}(y)$ where $\hat{K}%
_{j}(y)=\hat{K}_{-j}(y)=n^{-1}\sum_{l=j+1}^{n}\hat{u}_{l}(y)\hat{u}_{l-j}(y)$
for $j\geq 0$. It is easy to see that a size control result can be
established for $T_{ref}$: If Assumption \ref{AW} holds, the conclusion of
Theorem \ref{HAR_F_all_Eicker} with $T_{E,\mathsf{W}}$ replaced by $T_{ref}$
still holds. However, the form of the long-run covariance estimator $\hat{%
\omega}_{w}(y)R(X^{\prime }X)^{-1}R^{\prime }$ used by this test statistic
takes its justification from a well-known result of \cite{grenander1954},
which holds under certain conditions on the regressors only. A leading case
where these conditions are satisfied is polynomial regression.
Unfortunately, precisely for such regressors it turns out that the
conditions for size control are violated and, in fact, it can be shown that
in this case the test based on $T_{ref}$ has size $1$ for every choice of
critical value; see \cite{PP4}.

\textbf{D.} \emph{(Random bandwidth, prewithening, flat-top kernels,
GLS-based tests, general nonsphericity-corrected F-type tests)} Tests based
on weighted autocovariance estimators $\hat{\Psi}_{w}$, but where the
weights are allowed to depend on the data (e.g., lag-window estimators with
data-driven bandwidth choice), or where prewithening is used, can be viewed
as special cases of nonsphericity-corrected F-type tests (under appropriate
conditions, see \cite{Prein2014b}). The same is true for tests using
long-run variance estimators based on flat-top kernels. Another example are
tests based on parametric long-run variance estimators or tests based on
feasible generalized least squares. A size control result for general
nonsphericity-corrected F-type tests, i.e., Wald-type tests that may use
estimators for $\beta $ other than ordinary least squares or may use
estimators for the long-run covariance matrix other than the ones mentioned
so far, is given in Theorem \ref{Ctoep_nonsphericity} in Section \ref%
{szctrl_3}. However, we abstain in this paper from making the conditions in
that theorem more concrete for the subclasses of tests just mentioned.
Theorem \ref{Ctoep} in Section \ref{szctrl_3} furthermore even applies to
classes of tests more general than the class of nonsphericity-corrected
F-type tests.

\section{Computational issues and numerical results\label{numerical}}

In the previous section we have obtained conditions under which the size of
autocorrelation robust tests can be controlled by a proper choice of the
critical value. In the present section we now turn to the question of how
critical values guaranteeing size control can be determined numerically
(provided they exist). Additionally, we also address the related question of
how to numerically compute the size of an autocorrelation robust test for a
given choice of the critical value (e.g., when choosing a critical value
suggested by asymptotic theory). We illustrate our algorithms by computing
size-controlling critical values for regression models using US
macroeconomic time series from the FRED-MD database of \cite{mccracken} as
regressors and where the correlation structure of the errors in the
regression is governed by a variety of families of spectral densities. In
these models we also compute the actual size of standard autocorrelation
robust tests that employ critical values suggested in the literature.

We emphasize that, although the algorithms we shall discuss below are
designed for regression models with stationary autocorrelated errors and for
the test statistics $T_{w}$ defined in (\ref{tlr}) or $T_{E,\mathsf{W}}$
defined in Section \ref{autocorr}, the basic ideas extend to other test
statistics (e.g., the ones discussed in Section \ref{extension}) and also to
other covariance models with mostly obvious modifications.

\subsection{Computation of critical values and size}

Consider testing the hypothesis given in (\ref{testing problem 0}) at the
significance level $\alpha \in (0,1)$ by means of the test statistic $T_{w}$
as defined in (\ref{tlr}), and suppose that Assumptions \ref{AW} and \ref%
{R_and_X} are satisfied. Furthermore, suppose that one knows, e.g., by
having checked the sufficient conditions given in Theorem \ref{HAR_F_all},
that existence of a critical value $C(\alpha )$ satisfying 
\begin{equation}
\sup_{\mu _{0}\in \mathfrak{M}_{0}}\sup_{0<\sigma ^{2}<\infty }\sup_{f\in 
\mathfrak{F}}P_{\mu _{0},\sigma ^{2}\Sigma (f)}\left( T_{w}\geq C(\alpha
)\right) \leq \alpha  \label{triplesup}
\end{equation}%
is guaranteed, where $\mathfrak{F}$ is a user-specified subset of $\mathfrak{%
F}_{\mathrm{all}}$. Because such a critical value is certainly not unique
and because of power considerations, it is reasonable to try to find the
\textquotedblleft smallest\textquotedblright\ critical value satisfying the
inequality in the previous display. From the discussion following Theorem %
\ref{HAR_F_all} and from Remark \ref{exist} we conclude that such a smallest
critical value $C_{\Diamond }(\alpha )$ indeed exists; furthermore, if
equality is achievable in the preceding display, $C_{\Diamond }(\alpha )$
then certainly also achieves it. We note the obvious facts that $C_{\Diamond
}(\alpha )$ depends on $\mathfrak{F}$, and that any critical value smaller
than $C_{\Diamond }(\alpha )$ will lead to a test that violates the size
constraint (\ref{triplesup}).

Now, because of $G(\mathfrak{M}_{0})$-invariance of $T_{w}$ (see Lemma \ref%
{Fact1} and Remark \ref{special cases} in Section \ref{suff_con}) and
because of Remark 5.5(iii) in \cite{PP2016}, the inequality (\ref{triplesup}%
) is equivalent to%
\begin{equation*}
\sup_{f\in \mathfrak{F}}P_{\mu _{0},\Sigma (f)}\left( T_{w}\geq C(\alpha
)\right) \leq \alpha
\end{equation*}%
where we have chosen $\mu _{0}$ as an arbitrary but fixed element of $%
\mathfrak{M}_{0}$ and have set $\sigma =1$. Exploiting the fact that $P_{\mu
_{0},\Sigma (f)}\left( T_{w}=C\right) =0$ for every $C\in \mathbb{R}$ (since 
$\lambda _{\mathbb{R}^{n}}(T_{w}=C)=0$ for every $C\in \mathbb{R}$, see
Lemma \ref{Fact1} and Remark \ref{special cases} in Section \ref{suff_con}),
the preceding display implies that%
\begin{equation}
C_{\Diamond }(\alpha )=\sup_{f\in \mathfrak{F}}F_{\Sigma (f)}^{-1}(1-\alpha
),  \label{eqn:def}
\end{equation}%
where $F_{\Sigma (f)}$ denotes the cumulative distribution function (cdf) of 
$P_{\mu _{0},\Sigma (f)}\circ T_{w}$ (since $\mu _{0}\in \mathfrak{M}_{0}$
is fixed we do not need to show dependence on $\mu _{0}$ in the notation).
As usual, for a cdf $F$ we denote by $F^{-1}$ the corresponding quantile
function $F^{-1}(x)=\inf \{z\in \mathbb{R}:F(z)\geq x\}$. In order to obtain 
$C_{\Diamond }(\alpha )$, one must hence solve the optimization problem in (%
\ref{eqn:def}).

We shall now provide an heuristic optimization algorithm to solve (\ref%
{eqn:def}) in the case where $\mathfrak{F}=\mathfrak{F}_{\mathrm{AR(}p%
\mathrm{)}}$ with $1\leq p\leq n-1$. We write $C_{\Diamond }(\alpha ,p)$ to
emphasize the dependence of the critical value on the autoregressive order;
apart from $p$, the critical value only depends on $X$, $R$, and the weights 
$w$, but not on the value of $r$ (cf. Remark \ref{delta}(ii)). [We do not
show the dependence on $X$, $R$, and $w$ in the notation.] Note that by
Remark \ref{Fall=FAR(n-1)} the families $\mathfrak{F}_{\mathrm{AR(}n-1%
\mathrm{)}}$ and $\mathfrak{F}_{\mathrm{all}}$ induce the same testing
problem, and hence the critical value $C_{\Diamond }(\alpha ,n-1)$ achieves
size control also over $\mathfrak{F}_{\mathrm{all}}$. Consequently, the
subsequent discussion covers testing problems where $\mathfrak{F}=\mathfrak{F%
}_{\mathrm{all}}$ as a special case. We start by reparameterizing the
optimization problem (\ref{eqn:def}), exploiting the fact that $\mathfrak{F}=%
\mathfrak{F}_{\mathrm{AR(}p\mathrm{)}}$ can be parameterized through the
partial autocorrelation coefficients (reflection coefficients), cf. \cite%
{barndorff}: To each $p$-vector of partial autocorrelation coefficients $%
\rho \in (-1,1)^{p}$ there corresponds a unique normalized $\mathrm{AR(}p%
\mathrm{)}$ spectral density $f_{\rho }$, say, and vice versa. Hence,
writing $F_{\rho }$ for $F_{\Sigma (f_{\rho })}$, it follows that 
\begin{equation}
C_{\Diamond }(\alpha ,p)=\sup_{f\in \mathfrak{F}_{\mathrm{AR(}p\mathrm{)}%
}}F_{\Sigma (f)}^{-1}(1-\alpha )=\sup_{\rho \in (-1,1)^{p}}F_{\rho
}^{-1}(1-\alpha ).  \label{eqn:oprob}
\end{equation}%
That is, $C_{\Diamond }(\alpha ,p)$ can be found by maximizing the objective
function $\rho \mapsto F_{\rho }^{-1}(1-\alpha )$ over $(-1,1)^{p}$.
Compared to other parameterizations of the set of all stationary $\mathrm{AR(%
}p\mathrm{)}$ spectral densities, e.g., through the autoregression
coefficients or the set of zeros of the AR polynomial, working with partial
autocorrelation coefficients has the clear advantage that no
cross-restrictions are present. One aspect that complicates the optimization
problem, in addition to being potentially high-dimensional, is that the
objective function $\rho \mapsto F_{\rho }^{-1}(1-\alpha )$ needs to be
approximated numerically, e.g., by a Monte Carlo algorithm, since an
analytical expression for $F_{\rho }^{-1}$ is unknown in general. Therefore,
an optimization algorithm for determining the supremum in the previous
display needs to determine a quantile via a Monte Carlo algorithm each time
a function evaluation is required, which can be computationally intensive,
but is amenable to parallelization.

The optimization algorithm we use for numerically determining $C_{\Diamond
}(\alpha ,p)$ is described in detail in Algorithm \ref{alg:AR} which can be
found in Appendix \ref{App F}. Roughly speaking, the algorithm starts with a
preliminary step that selects candidate values $\rho $ from $(-1,1)^{p}$,
which are then used as initial values in a (local) optimization step. This
step returns improved candidate values of $\rho $, the best of which are in
turn used as initial values in a second (local) optimization step that now
uses a larger number of replications in the Monte-Carlo evaluation of the
objective function $F_{\rho }^{-1}(1-\alpha )$ than was used in the previous
steps.

A related problem is to numerically determine the size of the test that
rejects if $T_{w}$ exceeds a certain \emph{given} critical value $C$, i.e.,
one wants to obtain 
\begin{equation*}
\sup_{\mu _{0}\in \mathfrak{M}_{0}}\sup_{0<\sigma ^{2}<\infty }\sup_{f\in 
\mathfrak{F}_{\mathrm{AR(}p\mathrm{)}}}P_{\mu _{0},\sigma ^{2}\Sigma
(f)}\left( T_{w}\geq C\right) .
\end{equation*}%
Similarly as above this can be reduced to determining 
\begin{equation}
\sup_{f\in \mathfrak{F}_{\mathrm{AR(}p\mathrm{)}}}P_{\mu _{0},\Sigma
(f)}\left( T_{w}\geq C\right)  \label{2nd problem}
\end{equation}%
for a fixed value of $\mu _{0}\in \mathfrak{M}_{0}$. One can then use a
variant of Algorithm \ref{alg:AR}, which is described in Algorithm \ref%
{alg:size} in Appendix \ref{App F}, to solve that problem.

Finally we emphasize that, as is typical for numerical optimization
problems, there is no guarantee that the algorithms mentioned above do
return the exact critical value $C_{\Diamond }(\alpha ,p)$ or the exact size
of a test given a critical value $C$. The algorithms are heuristics that
numerically approximate the quantities of interest.

\begin{remark}
In case $p=0$, i.e., when the errors are i.i.d., the algorithms simplify
considerably in an obvious way as no optimization over $\rho $ is then
necessary.
\end{remark}

\begin{remark}
\label{other}\emph{(Other test statistics) }The above development has been
given for the test statistic $T_{w}$. It applies to any other test statistic 
$T$ as long as (i) $T$ is $G(\mathfrak{M}_{0})$-invariant and (ii) satisfies 
$\lambda _{\mathbb{R}^{n}}(T=C)=0$ for every $C\in \mathbb{R}$, and (iii)
one can ensure that a size-controlling critical value exists. It thus, in
particular, applies to the weighted Eicker-test statistic $T_{E,\mathsf{W}}$
as defined in Section \ref{autocorr} (cf. Lemma \ref{Fact1} and Remark \ref%
{special cases} in Section \ref{suff_con}). Note that for problem (\ref{2nd
problem}) the just given condition (ii) on $T$ is actually not needed. [If $%
T $ does not belong to the class of nonsphericity-corrected $F$-type test
statistics (cf. \cite{PP2016}, Section 5.4), we can, however, no longer
conclude that the corresponding critical values $C_{\Diamond }(\alpha ,p)$
are independent of the value of $r$.] Furthermore, for covariance models $%
\mathfrak{C}$ not of the form $\mathfrak{C}(\mathfrak{F})$ the general
principles underlying the reduction of (\ref{triplesup}) to (\ref{eqn:def})
still apply, provided $T$ satisfies (i)-(iii) given above. Algorithms that
perform optimization of the so-obtained analogue of (\ref{eqn:def}) can then
be developed in a similar way by exploiting the structure of the given
covariance model $\mathfrak{C}$.
\end{remark}

\subsection{An illustration for regression models based on US macroeconomic
time series and autoregressive error processes\label{sec:macro}}

We now apply Algorithms \ref{alg:AR} and \ref{alg:size} introduced above to
regression models based on data from the FRED-MD database, which consists of
128 macroeconomic time series that have been subjected to
stationarity-inducing transformations (see \cite{mccracken} for detailed
information concerning the database). More specifically, we consider
regression models of the form 
\begin{equation*}
\mathbf{y}_{t}=\beta _{1}+\beta _{2}t+\beta _{3}x_{t}+\mathbf{u}_{t}\quad 
\text{ for }t=1,\ldots ,n,
\end{equation*}%
where $\mathbf{u}_{1},\ldots ,\mathbf{u}_{n}$ are distributed like
consecutive observations from a mean zero stationary Gaussian process with
spectral density $\sigma ^{2}f$, $f\in \mathfrak{F}_{\mathrm{AR(}p\mathrm{)}%
} $. Here $x_{t}$ is one of the 128 macroeconomic variables in the FRED-MD
database, and where we use the most recent $n=100$ observations from each
time series.\footnote{%
The database was downloaded on October 25, 2016 from
https://research.stlouisfed.org/econ/mccracken/fred-databases/} For each of
these 128 regressors and for every $p\in \{0,1,2,5,10,25,50,99\}$ we
consider the problem of testing a restriction on the coefficient $\beta _{3}$
at the $5\%$ level, i.e., we consider testing problem (\ref{testing problem
0}) with $\mathfrak{F}=\mathfrak{F}_{\mathrm{AR(}p\mathrm{)}}$, $R=(0,0,1)$,
and with $r$ arbitrary (the results presented below do not depend on $r$,
cf. Remark \ref{delta}(ii)). Recall that the case $p=99$ realizes the
testing problem for the case $\mathfrak{F}=\mathfrak{F}_{\mathrm{all}}$, and
that the case $p=0$ corresponds to i.i.d. disturbances. In each setting we
consider the test statistic $T_{w}$ as defined in (\ref{tlr}) as well as $%
T_{E,\mathsf{W}}$ as defined in Section \ref{autocorr}, with the design
matrix $X$ corresponding to the regression model in the previous display.
Bartlett weights $w(j,n)=(1-\left\vert j\right\vert /M_{n})\mathbf{1}%
_{(-1,1)}(j/M_{n})$ with $M_{n}=n/10$ (i.e., bandwidth equal to $1/10$) are
used for the test statistic $T_{w}$, and the same weights are used for the
matrix $\mathsf{W}$ appearing in $T_{E,\mathsf{W}}$. Since $q=1$, rejecting
for large values of $T_{w}$ is equivalent to rejecting if the t-type test
statistic corresponding to $T_{w}$, i.e., if 
\begin{equation*}
t_{w}(y)=%
\begin{cases}
(\hat{\beta}_{3}(y)-r)/\hat{\Omega}_{w}^{1/2}\left( y\right) & \text{if }%
\hat{\Omega}_{w}\left( y\right) \neq 0 \\[7pt] 
0 & \text{ else },%
\end{cases}%
\end{equation*}%
is large in absolute value. A similar observation applies to tests obtained
from $T_{E,\mathsf{W}}$, the corresponding t-type test statistic being
denoted by $t_{E,\mathsf{W}}$. For the test statistic $t_{w}$, critical
values that are based on fixed-bandwidth asymptotics are provided in \cite%
{KV2005}, p. 1146, and this critical value (for the bandwidth and kernel
chosen here) is given by $2.260568$. For the sake of comparability, and
because critical values for t-statistics are usually easier to interpret, we
shall present critical values for the t-type version of the test statistics
in what follows. Critical values for $T_{w}$ and $T_{E,\mathsf{W}}$ can
easily be obtained by taking the square. We also note that the critical
values obtained below can be used for the construction of confidence
intervals for $\beta _{3}$. We shall now apply Algorithm \ref{alg:AR} to
numerically compute the critical value that is needed to control size in
each scenario. Additionally, we also apply Algorithm \ref{alg:size} to
numerically compute the size of the test that rejects if $\left\vert
t_{w}\right\vert $ exceeds the above mentioned critical value provided by 
\cite{KV2005}.\footnote{%
A referee has questioned if using this critical value here is appropriate
given that the regressors are treated as nonrandom and that a linear trend
is included. However, note that the theory developed in \cite{KVB2000}, \cite%
{KiefVogl2002, KV2002, KV2005} is based on high-level assumptions that are
compatible with random as well as nonrandom regressors. Furthermore, linear
trends (as long as their coefficients are not subject to tests as is the
case here) can be accommodated in this framework in the random as well as in
the nonrandom regressor case by a reasoning based on the Frisch-Waugh-Lovell
theorem as in \cite{KVB2000} (under appropriate assumptions on the
regressors and errors).} The particular settings used in Algorithms \ref%
{alg:AR} and \ref{alg:size} for the computations in this section are
described in detail in Appendix \ref{App H}.

To ensure, for each of the 128 design matrices, existence of a critical
value for $T_{w}$ (and hence for $t_{w}$) that controls size for any
(nonempty) family of (normalized) spectral densities, we now check the
sufficient conditions of Theorem \ref{HAR_F_all}: That Assumption \ref{AW}
is satisfied follows from the discussion after that assumption since we use
the Bartlett kernel. Assumption \ref{R_and_X} is satisfied for all 128 cases
as none of these design matrices contains an element of the canonical basis
in its span, which is easily verified numerically, and which then implies
the assumption, since $\func{rank}(R(X^{\prime }X)^{-1}X^{\prime })=q$
always holds. It remains to check condition (\ref{non-incl}) for each of the
128 design matrices. This can successfully be done numerically and we
describe the details of this computation in Appendix \ref{App G}. Since $%
\limfunc{span}(X)\subseteq \mathsf{B}$ always holds, this then also implies
validity of condition (\ref{non-incl_Eicker}), and thus also implies
existence of a corresponding critical value for $t_{E,\mathsf{W}}$ (cf.
Theorem \ref{HAR_F_all_Eicker}).

\begin{figure}[tbp]
\centering
\includegraphics[scale = .65]{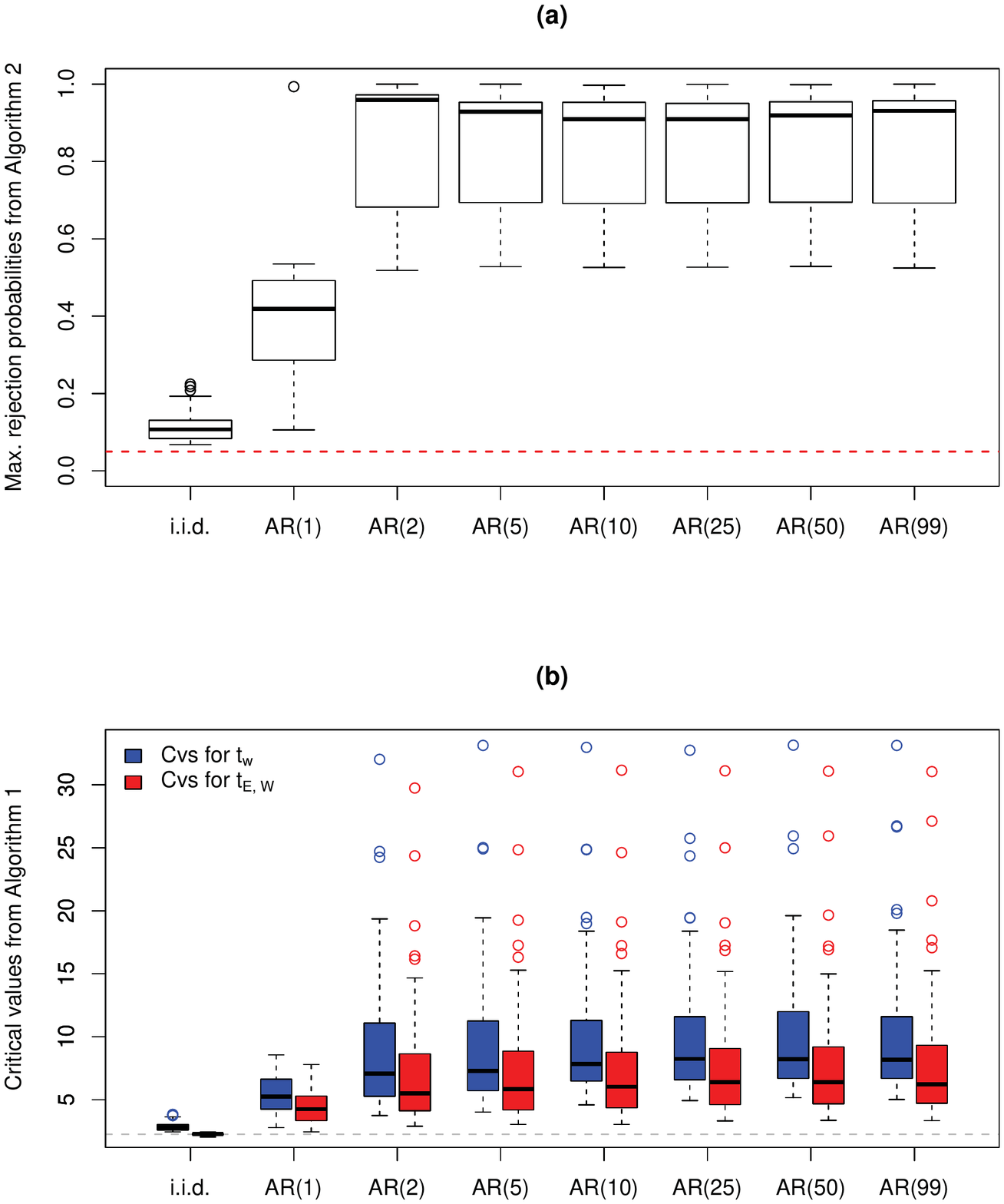}
\caption{(a) Sizes of the test which rejects if $|t_{w}|\geq 2.260568$
(Kiefer-Vogelsang critical value). The horizontal dashed red line
corresponds to 0.05. (b) Critical values guaranteeing a 5 \% level for the
t-type tests corresponding to $T_{w}$\ and $T_{E,\mathsf{W}}$. The
horizontal dashed grey line corresponds to the Kiefer-Vogelsang critical
value $2.260568$.}
\end{figure}


Figure 1(a) summarizes the numerical results for the size of the test that
rejects if $|t_{w}|$ exceeds the Kiefer-Vogelsang critical value $2.260568$:
For each autoregressive order $p\in \{1,2,5,10,25,50,99\}$ as well as for
the i.i.d. case (i.e., $p=0$) and for each of the 128 models considered we
obtained the size (i.e., the maximal rejection probability under the null)
by means of Algorithm \ref{alg:size} and summarize them in Figure 1(a) in
the form of boxplots (each boxplot representing the 128 sizes obtained for a
given order $p$). For a complete list of results see Table \ref{tab:size} in
Appendix \ref{app:tables}. As is apparent from Figure 1(a), the
Kiefer-Vogelsang critical value does not control size (not even in the
i.i.d. setting) at the desired $5\%$ level for neither one of the 128
regression models. This observation a fortiori extends to critical values
smaller than the Kiefer-Vogelsang critical value such as, e.g., the standard
normal critical value $1.96$, or the critical value obtained from
third-order asymptotic expansions (for the location model) in \cite{SPJ08}
which equals $2.242583$ (for the bandwidth and kernel chosen here). Figure
1(a) furthermore shows a large increase in the size when passing from the
i.i.d. to the $\mathrm{AR(}1\mathrm{)}$ case, and another large increase in
size when passing from the $\mathrm{AR(}1\mathrm{)}$ to the $\mathrm{AR(}2%
\mathrm{)}$ case.\footnote{%
The increase in size when passing form the i.i.d. to the $\mathrm{AR(}1%
\mathrm{)}$ case is connected to the fact that there are no concentration
spaces in the i.i.d. case, whereas in the $\mathrm{AR(}1\mathrm{)}$ case two
concentration spaces corresponding to angular frequencies $\gamma =0$ and $%
\gamma =\pi $ exist. The further increase in size when passing form $\mathrm{%
AR(}1\mathrm{)}$ to $\mathrm{AR(}2\mathrm{)}$ is related to the fact that $%
\mathrm{AR(}2\mathrm{)}$ models allow additional concentration spaces
corresponding to angular frequencies $\gamma \in (0,\pi )$.} In the $\mathrm{%
AR(}2\mathrm{)}$ case severe size distortions are present for all of the 128
regression models. Figure 1(a) also suggests that the sizes for the cases
with $p>2$ are comparable to the sizes in the $\mathrm{AR(}2\mathrm{)}$ case.

In Figure 1(b) we present the critical values which guarantee size control
at the $5\%$ level as computed by an application of Algorithm \ref{alg:AR}
for the test statistics $t_{w}$ as well as $t_{E,\mathsf{W}}$. Again we
present boxplots, and refer the reader to Tables \ref{tab:cv} and \ref%
{tab:cve} in Appendix \ref{app:tables} for a complete list of results.
Figure 1(b) suggests that the critical values required to control size
increase strongly when passing from the i.i.d. case to the $\mathrm{AR(}1%
\mathrm{)}$ model, and again when passing from the $\mathrm{AR(}1\mathrm{)}$
to the $\mathrm{AR(}2\mathrm{)}$ model. For larger $p$, the critical values,
while still increasing with $p$, seem to stabilize. Figure 1(b) also
illustrates the dependence of the critical value on the design matrix: For
some of the 128 regressors in the FRED-MD database, the critical values
needed to control size are very large, while for other regressors the
critical values are about 2-3 times as large as the Kiefer-Vogelsang
critical value (which, however, does \textit{not} provide size control).
Figure 1(b) further suggests that the critical values for $t_{w}$ needed to
control size at the $5\%$ level tend to be larger than the corresponding
critical values for $t_{E,\mathsf{W}}$.

While it is plain that the size-controlling critical values can never fall
when passing from an $\mathrm{AR(}p\mathrm{)}$ model to an $\mathrm{AR(}%
p^{\prime }\mathrm{)}$ model with $p<p^{\prime }$, this is not always
guaranteed for the \emph{numerically} determined critical values due to
numerical errors. We could have \textquotedblleft
monotonized\textquotedblright\ the results in Figure 1(b), but have decided
not to. A similar remark applies to Figure 1(a) as well as to Figure 2 given
further below.

A referee has suggested to examine also the critical value that is computed
under the presumption that the errors would follow a random walk. We discuss
this in Appendix \ref{App I}.

\begin{figure}[tbp]
\centering
\includegraphics[scale = .65]{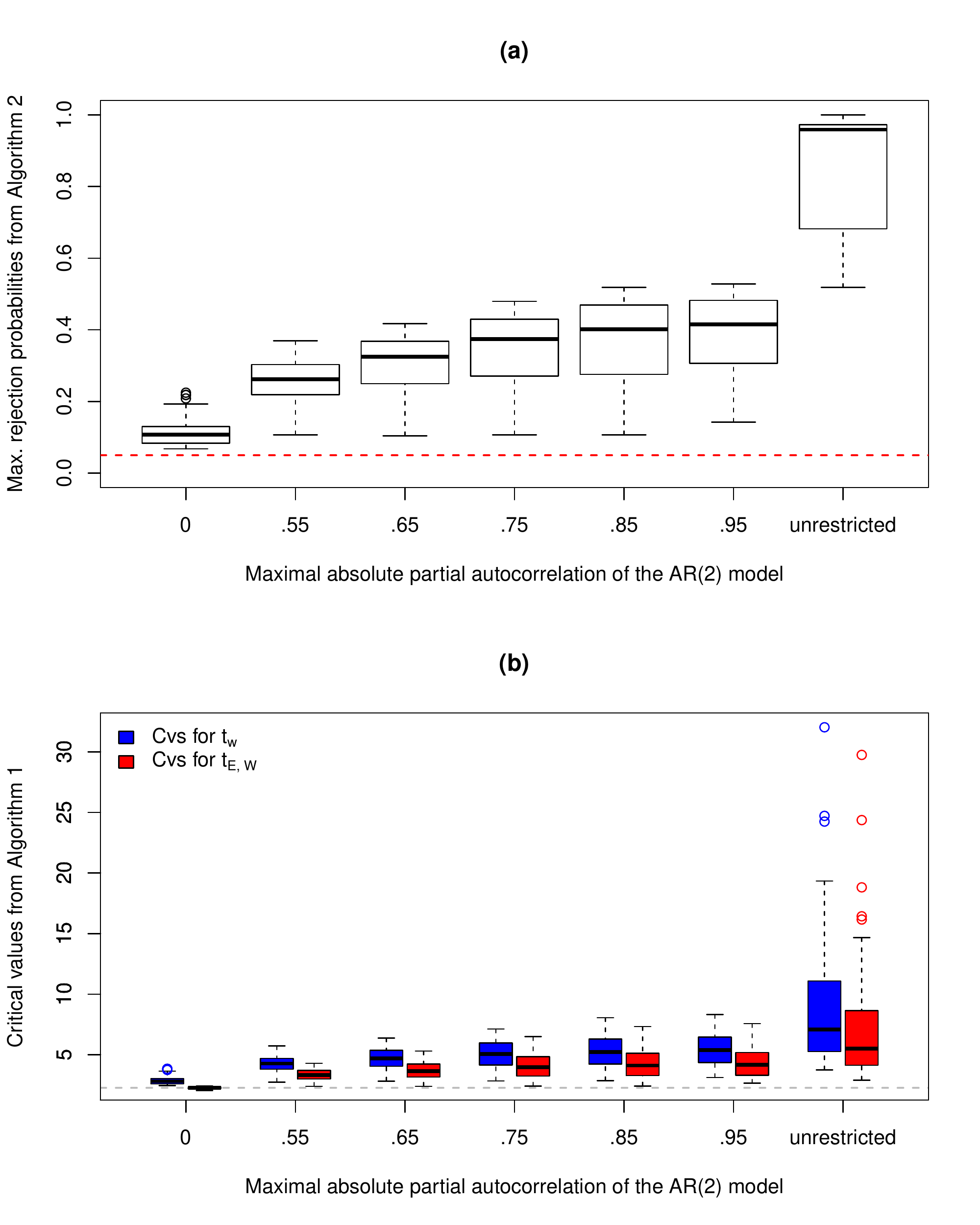}
\caption{(a) Sizes of the test which rejects if $|t_{w}|\geq 2.260568$
(Kiefer-Vogelsang critical value). The horizontal dashed red line
corresponds to 0.05. (b) Critical values guaranteeing a 5 \% level for the
t-type tests corresponding to $T_{w}$\ and $T_{E,\mathsf{W}}$. The
horizontal dashed grey line corresponds to the Kiefer-Vogelsang critical
value $2.260568$.}
\end{figure}


The final issue we shall investigate is how the size distortions of the test
using the Kiefer-Vogelsang critical values, documented in Figure 1(a) for
(unrestricted) $\mathrm{AR(}p\mathrm{)}$ models, are influenced if we
consider \emph{restricted} $\mathrm{AR(}p\mathrm{)}$ error processes where
the restrictions amount to placing a bound on the partial autocorrelations.
More precisely, we consider $\mathrm{AR(}p\mathrm{)}$ models where now the
partial autocorrelations are restricted to sets of the form $(-1+\varepsilon
,1-\varepsilon )^{p}\subseteq (-1,1)^{p}$ for some choice of $\varepsilon $, 
$0<\varepsilon <1$. Figure 1(a) suggests that, in order to obtain some
insight, we can focus on the $\mathrm{AR(}2\mathrm{)}$ case. We apply a
variant of Algorithm \ref{alg:size} to numerically compute the size of the
test based on $t_{w}$ together with the Kiefer-Vogelsang critical value,
where $\mathfrak{F}$ is now the set of all normalized $\mathrm{AR(}2\mathrm{)%
}$ spectral densities with maximal absolute partial autocorrelation
coefficient not exceeding a certain threshold in absolute value.
Furthermore, we apply a variant of Algorithm \ref{alg:AR} to numerically
determine the critical values needed to control the sizes of the tests based
on $t_{w}$ and $t_{E,\mathsf{W}}$ at the 5\% level over these sets of
spectral densities. As discussed in Remark \ref{rem:restrpart} in Appendix %
\ref{App F}, the algorithms now have to be modified in such a way that the
feasible set of the optimization problems in Stages 1 and 2 are restricted
sets of partial autocorrelation coefficients of the form $(-1+\varepsilon
,1-\varepsilon )^{2}$, and so that the starting values in Stage 0 fall
within this feasible set. The starting values are randomly generated as
described above, but in order to force them into $(-1+\varepsilon
,1-\varepsilon )^{2}$, they are all multiplied by $1-\varepsilon $.

The size of the test based on the Kiefer-Vogelsang critical value over the
so restricted $\mathrm{AR(}2\mathrm{)}$ models are summarized in Figure 2(a)
for a range of values for $1-\varepsilon $. From these results we see that
even if one is willing to impose the (questionable) assumption that partial
autocorrelation coefficients are known not to exceed $0.55$ in absolute
value, the size of the test rejecting whenever $\left\vert t_{w}\right\vert $
exceeds the Kiefer-Vogelsang critical value is considerable larger than $%
0.05 $ for most of the 128 regression models under consideration.
Unsurprisingly, the degree of size distortion increases steadily as the
bound for the maximal absolute partial autocorrelation increases, where we
observe a steep increase from $0.95$ to the unrestricted case. This shows
that even if the practitioner has good reasons to believe that each partial
autocorrelation is bounded away from one in modulus, critical values that
are typically used in practice, such as the Kiefer-Vogelsang critical value,
still fail to provide size control by a considerable margin. Numerical
computations of critical values that do provide size control are given in
Figure 2(b). Apart from the i.i.d. case, the corresponding boxplots do not
\textquotedblleft cover\textquotedblright\ the Kiefer-Vogelsang critical
value. In line with Figure 2(a), the size-controlling critical values
considerably exceed the Kiefer-Vogelsang critical value.

\section{Size control of tests of affine restrictions in regression models
with nonspherical disturbances: General theory \label{general}}

In this section we lay the foundation for all the size control results in
the paper. Other than in the preceding Sections \ref{autocorr} and \ref%
{numerical}, we here do not require that the disturbance vector in the
regression model (\ref{lm}) is induced by a stationary process, but we
revert to the more general framework specified in Section \ref{frame}. In
Subsection \ref{one} we obtain conditions under which the size of a
rejection region $W$ (satisfying certain invariance properties) is smaller
than one when testing (\ref{testing problem}). This result is then further
specialized to the important case when $W=\{T\geq C\}$ for a test statistic $%
T$ satisfying weak regularity conditions. In Subsection \ref{two} we
consider a \textit{family} of regions $W_{j}$ for $j\in \mathbb{N}$
(satisfying certain invariance properties), and we obtain conditions under
which for every $\alpha \in (0,1)$ there exists an element of this family $%
W_{j(\alpha )}$, say, the size of which does not exceed $\alpha $. These
results are then specialized to families of the form $W_{j}=\{T\geq C_{j}\}$
for $C_{j}\nearrow \infty $, and where the test statistic $T$ satisfies
certain weak regularity conditions. For such families we also obtain a lower
bound for the critical values that possibly can lead to size control, and we
study the problem under which conditions \emph{exact} size control can be
achieved, i.e., when for a given $\alpha \in (0,1)$ a critical value does
exist so that the size of the corresponding test equals $\alpha $. In
Subsection \ref{suff_con} we then show how some of the conditions arising in
the results in Subsections \ref{one} and \ref{two} can be implied from
lower-level conditions; see also Section \ref{suff_con_2} of Appendix \ref%
{App B}.

We start by defining a certain collection of linear subspaces of $\mathbb{R}%
^{n}$, where $n$ is sample size, that plays a central r\^{o}le in the size
control results. Loosely speaking, the linear spaces belonging to this
collection are either (nontrivial) projections of concentration spaces of
the covariance model $\mathfrak{C}$ (in the sense of \cite{PP2016}) on $%
\mathcal{L}^{\bot }$, where $\mathcal{L}$ is an appropriately chosen
subspace related to invariance properties of the tests under consideration,
or are what one could call \textquotedblleft higher-order\textquotedblright\
concentration spaces. For a more precise discussion see Section \ref%
{suff_con_2} in Appendix \ref{App B}. Since the tests we are interested in
are all at least $G(\mathfrak{M}_{0})$-invariant, a typical choice for $%
\mathcal{L}$ is $\mathfrak{M}_{0}^{lin}$, the linear space parallel to $%
\mathfrak{M}_{0}$. However, it proves useful to allow for the more general
case where $\mathcal{L}$ is an arbitrary linear space (typically containing $%
\mathfrak{M}_{0}^{lin}$). Recall from Section 5.1 of \cite{PP2016} that $G(%
\mathfrak{M}_{0})$ denotes the group of all maps of the form $y\mapsto
\delta (y-\mu _{0})+\mu _{0}^{\prime }$, where $\delta \in \mathbb{R}$, $%
\delta \neq 0$, and where $\mu _{0}$ as well as $\mu _{0}^{\prime }$ belong
to $\mathfrak{M}_{0}$.

\begin{definition}
Given a linear subspace $\mathcal{L}$ of $\mathbb{R}^{n}$ with $\dim (%
\mathcal{L})<n$ and a symmetric positive definite $n\times n$-dimensional
matrix $\Sigma $, we let 
\begin{equation}
\mathcal{L}(\Sigma )=\frac{\Pi _{\mathcal{L}^{\bot }}\Sigma \Pi _{\mathcal{L}%
^{\bot }}}{\Vert {\Pi _{\mathcal{L}^{\bot }}\Sigma \Pi _{\mathcal{L}^{\bot }}%
}\Vert }.  \label{L(Sigma)}
\end{equation}%
Given a covariance model $\mathfrak{C}$, we let $\mathcal{L}(\mathfrak{C}%
)=\left\{ \mathcal{L}(\Sigma ):\Sigma \in \mathfrak{C}\right\} $.
Furthermore, we define 
\begin{equation*}
\mathbb{J}(\mathcal{L},\mathfrak{C})=\left\{ \mathrm{\limfunc{span}}(\bar{%
\Sigma}):\bar{\Sigma}\in \limfunc{cl}(\mathcal{L}(\mathfrak{C})),\ \limfunc{%
rank}(\bar{\Sigma})<n-\dim (\mathcal{L})\right\} ,
\end{equation*}%
where the closure is here understood w.r.t. $\mathbb{R}^{n\times n}$.
\end{definition}

Note that the denominator in (\ref{L(Sigma)}) is always nonzero, that $%
\mathbb{J}(\mathcal{L},\mathfrak{C})$ can be empty, and that $\mathbb{J}(%
\mathcal{L},\mathfrak{C})$ neither contains $\mathcal{L}^{\bot }$ nor the
trivial space $\left\{ 0\right\} $. Also note that $\mathbb{J}(\mathcal{L},%
\mathfrak{C})$ is independent of the particular choice of norm used in the
above definition.

\begin{remark}
\label{rescale}(i) Even in the special case where $\mathcal{L}=\left\{
0\right\} $, the set $\mathbb{J}(\mathcal{L},\mathfrak{C})$ need not
coincide with the set of all concentration spaces in the sense of \cite%
{PP2016} (i.e., with the set, the union of which is $J(\mathfrak{C})$
defined in Section 5.3 of \cite{PP2016}). Both sets coincide in case $%
\mathcal{L}=\left\{ 0\right\} $ if and only if $\mathfrak{C}$ is bounded
away from zero, i.e., if there is no sequence $\Sigma _{j}\in \mathfrak{C}$
that converges to the zero matrix. Note that this latter condition on $%
\mathfrak{C}$ is satisfied for many covariance models, e.g., if $\mathfrak{C}
$ consists of correlation matrices or under similar normalization
assumptions.

(ii) In fact, as long as we are only concerned with $G(\mathfrak{M}_{0})$%
-invariant tests and their null-rejection probabilities, we could without
loss of generality always assume that the covariance model is bounded and is
bounded away from zero: note that, e.g., passing from $\mathfrak{C}$ to the
normalized covariance model $\left\{ \Sigma /\left\Vert \Sigma \right\Vert
:\Sigma \in \mathfrak{C}\right\} $ does not affect the null-rejection
probabilities of $G(\mathfrak{M}_{0})$-invariant tests, see Proposition 5.4
in \cite{PP2016}.\footnote{%
The effect of replacing $\Sigma $ by $\Sigma /\left\Vert \Sigma \right\Vert $
on non-null-rejection probabilities is that $\sigma ^{2}$ is replaced by $%
\sigma ^{2}\left\Vert \Sigma \right\Vert $. As a consequence, the
nuisance-minimal power at any $\mu \in \mathfrak{M}_{1}$ remains unaffected.}
Furthermore, note that $\mathbb{J}(\mathcal{L},\mathfrak{C})$ does not
change if $\mathfrak{C}$ is being rescaled.
\end{remark}

\subsection{Size less than one\label{one}}

The subsequent proposition gives simple sufficient conditions under which
the size of a test is less than one.

\begin{proposition}
\label{L1} Let $\mathfrak{C}$ be a covariance model, and let $W$ be a $G(%
\mathfrak{M}_{0})$-invariant rejection region that is also invariant w.r.t.
addition of elements of a linear subspace $\mathcal{V}$ of $\mathbb{R}^{n}$.
Define $\mathcal{L}=\mathrm{\limfunc{span}}(\mathfrak{M}_{0}^{lin}\cup 
\mathcal{V})$, i.e., $\mathcal{L}$ is the linear space generated by $%
\mathfrak{M}_{0}^{lin}\cup \mathcal{V}$, and assume that $\dim (\mathcal{L}%
)<n$. Suppose that

(a)$\ W^{c}$ is not a $\lambda _{\mathbb{R}^{n}}$-null set, and

(b) if $\mathcal{S}\in \mathbb{J}(\mathcal{L},\mathfrak{C})$ then $\left( 
\limfunc{cl}(W)\right) ^{c}\cap (\mu _{0}+\mathcal{S})\neq \emptyset $ for
some $\mu _{0}\in \mathfrak{M}_{0}$ (and hence for all $\mu _{0}\in 
\mathfrak{M}_{0}$ in view of $G(\mathfrak{M}_{0})$-invariance of $W$).

Then the size of the test given by the rejection region $W$ satisfies%
\begin{equation}
\sup_{\mu _{0}\in \mathfrak{M}_{0}}\sup_{0<\sigma ^{2}<\infty }\sup_{\Sigma
\in \mathfrak{C}}P_{\mu _{0},\sigma ^{2}\Sigma }(W)<1.  \label{size}
\end{equation}
\end{proposition}

The leading case in this proposition and in the results in the subsequent
section is the case where $\mathcal{V}=\left\{ 0\right\} $, and hence $%
\mathcal{L}=\mathfrak{M}_{0}^{lin}$.

\begin{remark}
(i) In case $W$ satisfies the invariance assumptions appearing in
Proposition \ref{L1} and $\dim (\mathcal{L})=n$ holds, it follows that $W$
is either empty or the entire space $\mathbb{R}^{n}$, both being trivial
cases. Similar remarks apply to the other results in this and the next
subsection and will not be stated.

(ii) If a rejection region $W^{\ast }$ differs from a rejection region $W$
that satisfies the assumptions of Proposition \ref{L1} only by a $\lambda _{%
\mathbb{R}^{n}}$-null set, then the conclusions of Proposition \ref{L1} also
hold for $W^{\ast }$ (even if it does not satisfy the assumptions of that
proposition), since $W$ and $W^{\ast }$ have the same rejection
probabilities. Similar remarks apply to the other results in this and the
next subsection and will not be stated.

(iii) It is not difficult to see that $G(\mathfrak{M}_{0})$-invariance
together with invariance w.r.t. addition of the elements of $\mathcal{V}$ is
equivalent to $G(\mu _{0}+\mathcal{L})$-invariance for some $\mu _{0}\in 
\mathfrak{M}_{0}$ (and hence every $\mu _{0}\in \mathfrak{M}_{0}$).
\end{remark}

\begin{remark}
\label{sufff}(i) Proposition \ref{L1} applies, in particular, to rejection
regions of the form $W=\{y\in \mathbb{R}^{n}:T(y)\geq C\}$ for some $-\infty
<C<\infty $, where $T:\mathbb{R}^{n}\rightarrow \mathbb{R}$ is
Borel-measurable, is $G(\mathfrak{M}_{0})$-invariant and is also invariant
w.r.t. addition of elements of a linear subspace $\mathcal{V}$ of $\mathbb{R}%
^{n}$.

(ii) If, additionally to the conditions in (i), $T$ is continuous on the
complement of a closed set $N^{\dag }$, then a sufficient condition for
condition (b) in Proposition \ref{L1} is as follows: if $\mathcal{S}\in 
\mathbb{J}(\mathcal{L},\mathfrak{C})$, then $W^{c}\cap (N^{\dag })^{c}\cap
(\mu _{0}+\mathcal{S})\neq \emptyset $ for some $\mu _{0}\in \mathfrak{M}%
_{0} $. This follows from Lemma \ref{L3} in Appendix \ref{App B}. [An
equivalent formulation of this sufficient condition is that whenever $%
\mathcal{S}\in \mathbb{J}(\mathcal{L},\mathfrak{C})$ then there exists an $%
s\in \mathcal{S}$ such that $T(\mu _{0}+s)<C$ and $\mu _{0}+s\notin N^{\dag
} $ for some $\mu _{0}\in \mathfrak{M}_{0}$.]

(iii) Suppose the conditions on $T$ in (ii) are satisfied. In case $\mathbb{J%
}(\mathcal{L},\mathfrak{C})$ is finite, i.e., $\mathbb{J}(\mathcal{L},%
\mathfrak{C})=\left\{ \mathcal{S}_{1},\ldots ,\mathcal{S}_{r}\right\} $, and
provided that for every $\mathcal{S}_{i}\in \mathbb{J}(\mathcal{L},\mathfrak{%
C})$ one can find an element $s_{i}\in \mathcal{S}_{i}$ with $\mu
_{0,i}+s_{i}\notin N^{\dag }$ for some $\mu _{0,i}\in \mathfrak{M}_{0}$, any 
$C$ satisfying $C>\max_{1\leq i\leq r}T(\mu _{0,i}+s_{i})$ gives rise to a
critical region $W$ that has size less than $1$. To see this, observe that
condition (b) in Proposition \ref{L1} is satisfied in view of (ii) and the
construction of $C$. Furthermore, observe that $\mu _{0,i}+s_{i}\in W^{c}$.
Since $T$ is continuous at $\mu _{0,i}+s_{i}$ (as $\mu _{0,i}+s_{i}\notin
N^{\dag }$), a sufficiently small open ball with center at $\mu _{0,i}+s_{i}$
also belongs to $W^{c}$, showing that $W^{c}$ is not a $\lambda _{\mathbb{R}%
^{n}}$-null set.
\end{remark}

\subsection{Size less than $\protect\alpha $\label{two}}

The next proposition is the basis for the size control results we want to
obtain.

\begin{proposition}
\label{L2}Let $\mathfrak{C}$ be a covariance model, and let $W_{j}$ be a
sequence of $G(\mathfrak{M}_{0})$-invariant rejection regions that are also
invariant w.r.t. addition of elements of a linear subspace $\mathcal{V}$ of $%
\mathbb{R}^{n}$. Define $\mathcal{L}=\mathrm{\limfunc{span}}(\mathfrak{M}%
_{0}^{lin}\cup \mathcal{V})$, and assume that $\dim (\mathcal{L})<n$. Assume
that the rejection regions satisfy $W_{j}\supseteq W_{j+1}$ for all $j\in 
\mathbb{N}$. Suppose that

(a) $\bigcap_{j\in \mathbb{N}}W_{j}$ is a $\lambda _{\mathbb{R}^{n}}$-null
set, and

(b) if $\mathcal{S}\in \mathbb{J}(\mathcal{L},\mathfrak{C})$ then $%
\bigcap_{j\in \mathbb{N}}\limfunc{cl}(W_{j})$ is a $\lambda _{\mu _{0}+%
\mathcal{S}}$-null set for some $\mu _{0}\in \mathfrak{M}_{0}$ (and hence
every $\mu _{0}\in \mathfrak{M}_{0}$).

Then the sizes of the tests given by the rejection regions $W_{j}$ satisfy%
\begin{equation}
\sup_{\mu _{0}\in \mathfrak{M}_{0}}\sup_{0<\sigma ^{2}<\infty }\sup_{\Sigma
\in \mathfrak{C}}P_{\mu _{0},\sigma ^{2}\Sigma }(W_{j})\rightarrow 0\quad 
\text{ as }\quad j\rightarrow \infty ;  \label{size_2}
\end{equation}%
in particular, for every $\alpha \in (0,1)$ there exists a $j(\alpha )\in 
\mathbb{N}$ so that 
\begin{equation*}
\sup_{\mu _{0}\in \mathfrak{M}_{0}}\sup_{0<\sigma ^{2}<\infty }\sup_{\Sigma
\in \mathfrak{C}}P_{\mu _{0},\sigma ^{2}\Sigma }(W_{j(\alpha )})\leq \alpha .
\end{equation*}
\end{proposition}

In the important special case, where $W_{j}=\{y\in \mathbb{R}^{n}:T(y)\geq
C_{j}\}$ for some real-valued test statistic $T$ and where $C_{j}\nearrow
\infty $ as $j\rightarrow \infty $, condition (a) in the preceding
proposition is clearly always satisfied since the intersection is empty in
this case. The subsequent corollary now provides sufficient conditions for
condition (b) in case of such rejection regions $W_{j}$. The conditions
imposed on $T$ in this corollary are widely satisfied (at least in the
leading case where $\mathcal{V}=\left\{ 0\right\} $), see Subsection \ref%
{suff_con}.

\begin{corollary}
\label{C5} Let $\mathfrak{C}$ be a covariance model, and assume that the
test statistic $T:\mathbb{R}^{n}\rightarrow \mathbb{R}$ is Borel-measurable
and is continuous on the complement of a closed set $N^{\dag }$.\footnote{%
While this condition is trivially satisfied for every Borel-measurable $T$
upon choosing $N^{\dag }$ equal to $\mathbb{R}^{n}$, satisfying the other
conditions in the corollary will rule out this case, except if $\mathbb{J}(%
\mathcal{L},\mathfrak{C})$ is empty. Also note that for typical test
statistics the set $N^{\dag }$ will turn out to be a `small' set, e.g., a $%
\lambda _{\mathbb{R}^{n}}$-null set.} Assume that $T$ and $N^{\dag }$ are $G(%
\mathfrak{M}_{0})$-invariant, and are also invariant w.r.t. addition of
elements of a linear subspace $\mathcal{V}$ of $\mathbb{R}^{n}$.\footnote{%
If $T$ is Borel-measurable, is continuous on the complement of a closed set $%
N^{\dag }$, and if $T$ satisfies the invariance requirements in the
corollary, then it is easy to see that one can always find a closed subset $%
N^{\dag \dag }$ of $N^{\dag }$ such that $T$ is continuous on the complement
of $N^{\dag \dag }$ and such that $N^{\dag \dag }$ satisfies the invariance
properties asked for in the corollary. Hence, requiring in the corollary
that the set $N^{\dag }$ satisfies the invariance conditions imposes no loss
of generality.} Define $\mathcal{L}=\mathrm{\limfunc{span}}(\mathfrak{M}%
_{0}^{lin}\cup \mathcal{V})$, and assume that $\dim (\mathcal{L})<n$. Let $%
W_{j}=\{y\in \mathbb{R}^{n}:T(y)\geq C_{j}\}$ for $-\infty <C_{j}<\infty $
with $C_{j}\nearrow \infty $ as $j\rightarrow \infty $. Then:

\begin{enumerate}
\item The conclusion of Proposition \ref{L2} holds, provided condition (b)
in that proposition is satisfied.

\item A sufficient condition for condition (b) in Proposition \ref{L2} to
hold is the following: if $\mathcal{S}\in \mathbb{J}(\mathcal{L},\mathfrak{C}%
)$, then the set $N^{\dag }$ is a $\lambda _{\mu _{0}+\mathcal{S}}$-null set
for some $\mu _{0}\in \mathfrak{M}_{0}$ (and hence for all $\mu _{0}\in 
\mathfrak{M}_{0}$).

\item In case $N^{\dag }$ is a finite or countable union of affine
subspaces, the sufficient condition given in 2. is equivalent to: if $%
\mathcal{S}\in \mathbb{J}(\mathcal{L},\mathfrak{C})$, then $\mu _{0}+%
\mathcal{S}\not\subseteq N^{\dag }$ for some $\mu _{0}\in \mathfrak{M}_{0}$
(and hence for all $\mu _{0}\in \mathfrak{M}_{0}$).
\end{enumerate}
\end{corollary}

\begin{remark}
\label{star} The corollary implies, in particular, that size control is
always possible in case $\mathbb{J}(\mathcal{L},\mathfrak{C})$ is empty.%
\footnote{%
Observe that in this case setting $N^{\dag }$ equal to $\mathbb{R}^{n}$ does
not restrict the applicability of the corollary and hence voids the
continuity requirement on $T$.} This is, e.g., the case if $\mathfrak{C}$
has no singular limit points and is norm bounded, or more generally if $%
\mathfrak{C}$ is such that $\left\{ \Sigma /\left\Vert \Sigma \right\Vert
:\Sigma \in \mathfrak{C}\right\} $ has no singular limit points. Of course,
this is in line with Theorems 5.10 and 5.12 in \cite{PP2016}.
\end{remark}

\begin{remark}
\label{simplif_1}If $\mathfrak{M}_{0}$ is a linear space (i.e., if $r=0$),
we can w.l.o.g. set $\mu _{0}=0$ in the sufficient conditions in Part 2 and
3 of Corollary \ref{C5}, leading to a simplification of the conditions. Even
if $\mathfrak{M}_{0}$ is not a linear space (i.e., if $r\neq 0$), the same
simplification can be made provided that $N^{\dag }$ is invariant under
addition of elements of $\mathfrak{M}_{0}$. This latter invariance property
is, in particular, satisfied whenever $N^{\dag }$ is $G(\mathfrak{M})$%
-invariant, which is the case for a large class of test statistics including
nonsphericity-corrected F-type test statistics that satisfy Assumption 5 in 
\cite{PP2016} (cf. Lemma \ref{Fact1} further below).
\end{remark}

\begin{remark}
\label{some comments}(i) The argument that establishes Part 3 of the
preceding corollary also shows that in condition (b) in Proposition \ref{L2}
we can replace \textquotedblleft $\bigcap_{j\in \mathbb{N}}\limfunc{cl}%
(W_{j})$ is a $\lambda _{\mu _{0}+\mathcal{S}}$-null set\textquotedblright\
equivalently by \textquotedblleft $\mu _{0}+\mathcal{S}\not\subseteq
\bigcap_{j\in \mathbb{N}}\limfunc{cl}(W_{j})$\textquotedblright\ provided
that $\bigcap_{j\in \mathbb{N}}\limfunc{cl}(W_{j})$ is a finite or countable
union of affine spaces.

(ii) The condition that $\bigcap_{j\in \mathbb{N}}\limfunc{cl}(W_{j})$ ($%
N^{\dag }$, respectively) is a $\lambda _{\mu _{0}+\mathcal{S}}$-null set in
Proposition \ref{L2} (Corollary \ref{C5}, respectively) is -- for $1$%
-dimensional $\mathcal{S}\in \mathbb{J}(\mathcal{L},\mathfrak{C})$ --
equivalent to $\bigcap_{j\in \mathbb{N}}\limfunc{cl}(W_{j})\cap (\mu _{0}+%
\mathcal{S})\subseteq \left\{ \mu _{0}\right\} $ ($N^{\dag }\cap (\mu _{0}+%
\mathcal{S})\subseteq \left\{ \mu _{0}\right\} $, respectively); i.e.,
except possibly for $\mu _{0}$, the entire set $\mu _{0}+\mathcal{S}$ lies
outside of $\bigcap_{j\in \mathbb{N}}\limfunc{cl}(W_{j})$ ($N^{\dag }$,
respectively). For a proof see Appendix \ref{App B}.
\end{remark}

Corollary \ref{C5} implies that, under its conditions, for every $\alpha \in
(0,1)$ there exists a real number $C(\alpha )$ such that%
\begin{equation}
\sup_{\mu _{0}\in \mathfrak{M}_{0}}\sup_{0<\sigma ^{2}<\infty }\sup_{\Sigma
\in \mathfrak{C}}P_{\mu _{0},\sigma ^{2}\Sigma }(T\geq C(\alpha ))\leq \alpha
\label{size_2.5}
\end{equation}%
holds. Under mild conditions a smallest such critical value exists as
discussed in the next remark.

\begin{remark}
\label{smallest C}\emph{(Existence of a smallest critical value guaranteeing
size control)}\textit{\ }(i) Let $\alpha \in (0,1)$ and suppose\textit{\ }$T$
is a test statistic\ such that (\ref{size_2.5}) holds. Let $\mathsf{CV}%
_{\leq }(\alpha )$ be the set of all real $C(\alpha )$ such that (\ref%
{size_2.5}) holds. Then the set $\mathsf{CV}_{\leq }(\alpha )$ is clearly
nonempty and is an interval either of the form $(C_{\Diamond }(\alpha
),\infty )$ or of the form $[C_{\Diamond }(\alpha ),\infty )$ for some real
number $C_{\Diamond }(\alpha )$. If the map $C\mapsto \sup_{\mu _{0}\in 
\mathfrak{M}_{0}}\sup_{0<\sigma ^{2}<\infty }\sup_{\Sigma \in \mathfrak{C}%
}P_{\mu _{0},\sigma ^{2}\Sigma }(T\geq C)$ is continuous from the right at $%
C=C_{\Diamond }(\alpha )$, then $\mathsf{CV}_{\leq }(\alpha )=[C_{\Diamond
}(\alpha ),\infty )$ must hold, i.e., a \emph{smallest }critical value
guaranteeing (\ref{size_2.5}) does exist and is given by $C_{\Diamond
}(\alpha )$. The just mentioned right-continuity property is easily seen to
be satisfied, whenever the test statistic $T$ has the property that $\lambda
_{\mathbb{R}^{n}}(T=C_{\Diamond }(\alpha ))=0$. This latter condition is
satisfied by a large class of test statistics (cf. Lemma \ref{Fact1} further
below).

(ii) Let $\alpha \in (0,1)$ and suppose $T$ is a test statistic\ such that
equality in (\ref{size_2.5}) holds for at least one real number $C(\alpha )$%
. Let $\mathsf{CV}_{=}(\alpha )$ denote the set of all such $C(\alpha )$,
which then clearly is a nonempty subinterval of $\mathsf{CV}_{\leq }(\alpha
) $ consisting of an initial (w.r.t. the order on the real line) segment of $%
\mathsf{CV}_{\leq }(\alpha )$.\footnote{%
That is, if $c\in \mathsf{CV}_{=}(\alpha )$, then every $c^{\prime }\in 
\mathsf{CV}_{\leq }(\alpha )$ with $c^{\prime }\leq c$ also belongs to $%
\mathsf{CV}_{=}(\alpha )$.} If $\mathsf{CV}_{\leq }(\alpha )=[C_{\Diamond
}(\alpha ),\infty )$, as is, e.g., the case under the condition discussed in
(i), then $C_{\Diamond }(\alpha )$ is also the smallest element of $\mathsf{%
CV}_{=}(\alpha )$, i.e., a \emph{smallest }critical value guaranteeing
equality in (\ref{size_2.5}) exists and is given by $C_{\Diamond }(\alpha )$.
\end{remark}

The next lemma provides a lower bound $C^{\ast }$ for the critical values $C$
that possibly can lead to size control and is a building block for the next
proposition.

\begin{lemma}
\label{L4} Let $\mathfrak{C}$, $T$, $N^{\dag }$, $\mathcal{V}$, and $%
\mathcal{L}$ be as in Corollary \ref{C5}. Denote by $\mathbb{H}$ the set of
all $\mathcal{S\in }\mathbb{J}(\mathcal{L},\mathfrak{C})$ such that $T$
restricted to $\mu _{0}+\mathcal{S}$ is equal to a constant $C(\mathcal{S})$%
, say, $\lambda _{\mu _{0}+\mathcal{S}}$-almost everywhere for some $\mu
_{0}\in \mathfrak{M}_{0}$ (and hence for all $\mu _{0}\in \mathfrak{M}_{0}$%
). Define $C_{\ast }=\inf_{\mathcal{S}\in \mathbb{H}}C(\mathcal{S})$ and $%
C^{\ast }=\sup_{\mathcal{S}\in \mathbb{H}}C(\mathcal{S})$, with the
convention that $C_{\ast }=\infty $ and $C^{\ast }=-\infty $ if $\mathbb{H}$
is empty. Then the following hold:

\begin{enumerate}
\item Any one-dimensional $\mathcal{S}\in \mathbb{J}(\mathcal{L},\mathfrak{C}%
)$ belongs to $\mathbb{H}$. [In particular, $\mathbb{H}$ is nonempty if a
one-dimensional $\mathcal{S}\in \mathbb{J}(\mathcal{L},\mathfrak{C})$
exists.]

\item Suppose that for every $\mathcal{S}\in \mathbb{H}$ the set $N^{\dag }$
is a $\lambda _{\mu _{0}+\mathcal{S}}$-null set for some $\mu _{0}\in 
\mathfrak{M}_{0}$ (and hence for all $\mu _{0}\in \mathfrak{M}_{0}$). Then
for $C\in (-\infty ,C^{\ast })$ the size of the test with critical region $%
\left\{ T\geq C\right\} $ satisfies%
\begin{equation}
\sup_{\mu _{0}\in \mathfrak{M}_{0}}\sup_{0<\sigma ^{2}<\infty }\sup_{\Sigma
\in \mathfrak{C}}P_{\mu _{0},\sigma ^{2}\Sigma }(T\geq C)=1.  \label{size_3}
\end{equation}%
If, additionally, $N^{\dag }$ is a $\lambda _{\mu _{0}+\mathcal{S}}$-null
set for some $\mu _{0}\in \mathfrak{M}_{0}$ (and hence for all $\mu _{0}\in 
\mathfrak{M}_{0}$) not only for $\mathcal{S}\in \mathbb{H}$ but for every $%
\mathcal{S\in }\mathbb{J}(\mathcal{L},\mathfrak{C})$, then the l.h.s. in (%
\ref{size_3}) converges to zero for $C\rightarrow \infty $ (implying that
then $C^{\ast }<\infty $ holds).

\item Suppose that for every $\mathcal{S}\in \mathbb{H}$ the set $N^{\dag }$
is a $\lambda _{\mu _{0}+\mathcal{S}}$-null set for some $\mu _{0}\in 
\mathfrak{M}_{0}$ (and hence for all $\mu _{0}\in \mathfrak{M}_{0}$). Then
for $C\in (C_{\ast },\infty )$%
\begin{equation}
\inf_{\mu _{0}\in \mathfrak{M}_{0}}\inf_{0<\sigma ^{2}<\infty }\inf_{\Sigma
\in \mathfrak{C}}P_{\mu _{0},\sigma ^{2}\Sigma }(T\geq C)=0.  \label{bias}
\end{equation}
\end{enumerate}
\end{lemma}

Part 3 of the preceding lemma also provides some information on the infimal
rejection probability under the null hypothesis (and thus on power
properties \textquotedblleft near\textquotedblright\ the null) of the test
in case $\mathbb{H}$ is nonempty: In this case, we clearly have $C_{\ast
}\leq C^{\ast }$. It follows that, under the assumptions of the lemma, a
test based on $T$ which has size less than or equal to $\alpha <1$, \emph{%
must} have a critical value $C$ larger than or equal to $C^{\ast }$($\geq
C_{\ast }$). Part 3 then implies severe biasedness of the test (except
possibly if $C=C_{\ast }=C^{\ast }$ holds), which typically entails bad
power properties in certain regions of the alternative hypothesis (in view
of Part 3 of Theorem 5.7 and Remark 5.5(iii) in \cite{PP2016}). However,
note that in case $\mathbb{H}$ is empty, we have $C_{\ast }=\infty $ and
consequently Part 3 of the lemma does not convey any information about
rejection probabilities. Since we concentrate exclusively on size properties
in this paper, we postpone a detailed discussion of power properties to a
companion paper. Furthermore, in view of Remark 5.5(iii) in \cite{PP2016},
relation (\ref{size_3}) remains valid even after one removes the suprema
over $\mu _{0}$ and $\sigma ^{2}$ in (\ref{size_3}). A similar remark
applies to (\ref{bias}).

The size control result in Corollary \ref{C5} can be sharpened to an exact
size control result under some additional assumptions.

\begin{proposition}
\label{P6} Let $\mathfrak{C}$, $T$, $N^{\dag }$, $\mathcal{V}$, and $%
\mathcal{L}$ be as in Corollary \ref{C5}. Suppose that for every $\mathcal{S}%
\in \mathbb{J}(\mathcal{L},\mathfrak{C})$ the set $N^{\dag }$ is a $\lambda
_{\mu _{0}+\mathcal{S}}$-null set for some $\mu _{0}\in \mathfrak{M}_{0}$
(and hence for all $\mu _{0}\in \mathfrak{M}_{0}$). Let $\mathbb{H}$, $C(%
\mathcal{S})$, and $C^{\ast }$ be as in Lemma \ref{L4} (note that $C^{\ast
}<\infty $ holds in view of Part 2 of that lemma).

\begin{description}
\item[A.] Suppose that for every $C\in (C^{\ast },\infty )$

(a) $\lambda _{\mathbb{R}^{n}}(T=C)=0$, and

(b) if $\mathcal{S}\in \mathbb{J}(\mathcal{L},\mathfrak{C})\backslash 
\mathbb{H}$ then $\lambda _{\mu _{0}+\mathcal{S}}(T=C)=0$ for some $\mu
_{0}\in \mathfrak{M}_{0}$ (and hence for every $\mu _{0}\in \mathfrak{M}_{0}$%
).

Then the following holds:
\end{description}

\begin{enumerate}
\item The function%
\begin{equation}
C\mapsto \sup_{\mu _{0}\in \mathfrak{M}_{0}}\sup_{0<\sigma ^{2}<\infty
}\sup_{\Sigma \in \mathfrak{C}}P_{\mu _{0},\sigma ^{2}\Sigma }(T\geq C)
\label{func}
\end{equation}%
is nonincreasing, equals one for $C\in (-\infty ,C^{\ast })$, is continuous
on $(C^{\ast },\infty )$, and converges to $0$ as $C\rightarrow \infty $.

\item Define%
\begin{equation}
\alpha ^{\ast }=\sup_{C\in (C^{\ast },\infty )}\sup_{\mu _{0}\in \mathfrak{M}%
_{0}}\sup_{0<\sigma ^{2}<\infty }\sup_{\Sigma \in \mathfrak{C}}P_{\mu
_{0},\sigma ^{2}\Sigma }(T\geq C).  \label{alpha*}
\end{equation}%
Then, for every $\alpha \in (0,1)$ there exists a $C(\alpha )\in (C^{\ast
},\infty )$ such that 
\begin{equation}
\sup_{\mu _{0}\in \mathfrak{M}_{0}}\sup_{0<\sigma ^{2}<\infty }\sup_{\Sigma
\in \mathfrak{C}}P_{\mu _{0},\sigma ^{2}\Sigma }(T\geq C(\alpha ))\leq \min
(\alpha ,\alpha ^{\ast })\leq \alpha .  \label{alpha*-ineq}
\end{equation}%
[Note that the l.h.s.~of (\ref{alpha*-ineq}) necessarily is less than or
equal to $\alpha ^{\ast }$ for every choice of $C(\alpha )\in (C^{\ast
},\infty )$.] Furthermore, for every $\alpha \in (0,\alpha ^{\ast })$ the
constant $C(\alpha )\in (C^{\ast },\infty )$ can be chosen such that
(\textquotedblleft exact size control\textquotedblright ) 
\begin{equation}
\sup_{\mu _{0}\in \mathfrak{M}_{0}}\sup_{0<\sigma ^{2}<\infty }\sup_{\Sigma
\in \mathfrak{C}}P_{\mu _{0},\sigma ^{2}\Sigma }(T\geq C(\alpha ))=\alpha .
\label{alpha-equal}
\end{equation}
\end{enumerate}

\begin{description}
\item[B.] Suppose $C^{\ast \ast }\geq C^{\ast }$ is a real number and
suppose that the conditions (a) and (b) given above are only known to hold
for every $C\in (C^{\ast \ast },\infty )$. Define $\alpha ^{\ast \ast }$ in
the same way as $\alpha ^{\ast }$, but with $C^{\ast \ast }$ replacing $%
C^{\ast }$. Then the claims in Part A.1 continue to hold as they stand
except for the fact that continuity is now only guaranteed on $(C^{\ast \ast
},\infty )$. As a consequence, all the conclusions of Part A.2 continue to
hold provided $C^{\ast }$ in that part is everywhere replaced by $C^{\ast
\ast }$ and $\alpha ^{\ast }$ by $\alpha ^{\ast \ast }$. [Since $C^{\ast
\ast }\geq C^{\ast }$ and hence $\alpha ^{\ast \ast }\leq \alpha ^{\ast }$,
also the claim in (\ref{alpha*-ineq}) continues to hold as it stands.]
\end{description}
\end{proposition}

Again by Remark 5.5(iii) in \cite{PP2016}, the suprema w.r.t. $\mu _{0}$ and 
$\sigma ^{2}$ can be removed from (\ref{func})-(\ref{alpha-equal}), and the
resulting expressions then do not depend on the particular choice for $\mu
_{0}\in \mathfrak{M}_{0}$ and $\sigma ^{2}\in (0,\infty )$.

\begin{remark}
\label{smallest C_2}\emph{(Continuity properties and existence of smallest
critical values) }Suppose the assumptions for Part A of Proposition \ref{P6}
are satisfied.

(i) If $C^{\ast }=-\infty $ then the function given by (\ref{func}) is
continuous on $\mathbb{R}$, whereas if $C^{\ast }>-\infty $ it is only
guaranteed to be continuous on $\mathbb{R}\backslash \left\{ C^{\ast
}\right\} $. If $C^{\ast }>-\infty $ and if, additionally, $\lambda _{%
\mathbb{R}^{n}}(T=C^{\ast })=0$ holds, it is easy to see that (\ref{func})
is then at least continuous from the right at $C^{\ast }$.

(ii) If $C^{\ast }=-\infty $ then clearly $\alpha ^{\ast }=1$ and hence
exact size control (\ref{alpha-equal}) is possible for every $\alpha \in
(0,1)$, whereas in case $C^{\ast }>-\infty $ we only can conclude that $%
\alpha ^{\ast }\leq 1$ and that (\ref{alpha-equal}) holds for $\alpha \in
(0,\alpha ^{\ast })$. If $C^{\ast }>-\infty $ and if, additionally, the map (%
\ref{func}) is continuous from the right at $C^{\ast }$, it follows from (i)
that (\ref{alpha-equal}) can also be achieved for $\alpha =\alpha ^{\ast }$
(with a $C(\alpha )$ belonging to $[C^{\ast },\infty )$).

(iii) An upper bound for $\alpha ^{\ast }$ is obviously given by $\sup_{\mu
_{0}\in \mathfrak{M}_{0}}\sup_{0<\sigma ^{2}<\infty }\sup_{\Sigma \in 
\mathfrak{C}}P_{\mu _{0},\sigma ^{2}\Sigma }(T\geq C^{\ast })$. If $C^{\ast
}=-\infty $ holds, or if $C^{\ast }>-\infty $ and the map (\ref{func}) is
continuous from the right at $C^{\ast }$, then this upper bound coincides
with $\alpha ^{\ast }$.

(iv) If $C^{\ast }=-\infty $ holds, or if $C^{\ast }>-\infty $ and the map (%
\ref{func}) is continuous from the right at $C^{\ast }$, then for every $%
\alpha \in (0,1)$ a smallest critical value $C_{\Diamond }(\alpha )\in 
\mathbb{R}$ satisfying (\ref{alpha*-ineq}) exists (i.e., $\mathsf{CV}_{\leq
}(\alpha )=[C_{\Diamond }(\alpha ),\infty )$ holds) in view of (i) and
Remark \ref{smallest C}. [Note that for $\alpha \in (0,1)$ here $C_{\Diamond
}(\alpha )>C^{\ast }$ must hold if $\alpha <\alpha ^{\ast }$, and that $%
C_{\Diamond }(\alpha )=C^{\ast }$ must hold if $\alpha =\alpha ^{\ast }$.]
Under the same conditions, for every $\alpha \in (0,\alpha ^{\ast }]\cap
(0,1)$ there is a smallest critical value satisfying (\ref{alpha-equal})
which is again given by $C_{\Diamond }(\alpha )$ (in fact, $\mathsf{CV}%
_{=}(\alpha )=[C_{\Diamond }(\alpha ),b]$ for some real $b\geq C_{\Diamond
}(\alpha )$ holds). This follows from (i), (ii) above and Remark \ref%
{smallest C} (and the fact that (\ref{func}) goes to zero for $C\rightarrow
\infty $).

(v) For every $\alpha \in (0,\alpha ^{\ast })$ a smallest critical value $%
C_{\Diamond }(\alpha )\in \mathbb{R}$ satisfying (\ref{alpha*-ineq}) and (%
\ref{alpha-equal}) always exists, even without the right-continuity
condition in case $C^{\ast }>-\infty $. Necessarily, $C_{\Diamond }(\alpha
)>C^{\ast }$ has to hold.

(vi) The case $\alpha ^{\ast }=0$ can occur, e.g., if $T$ is identically
equal to a constant. However, for large classes of test statistics such as
nonsphericity-corrected $F$-type test statistic as defined in (28) in
Section 5.4 of \cite{PP2016} and satisfying Assumptions 5 and 6 in that
paper we always have $\alpha ^{\ast }>0$. This follows from Part 5 of Lemma
5.15 in that reference.
\end{remark}

\begin{remark}
\label{dimL=n-1}(i) Suppose $\dim (\mathcal{L})=n-1$. Then clearly $\mathbb{J%
}(\mathcal{L},\mathfrak{C})$ is empty. Furthermore, any $T$ that satisfies
the invariance properties mentioned in Corollary \ref{C5} is then constant $%
\lambda _{\mathbb{R}^{n}}$-almost everywhere, hence size control is
trivially possible. Similarly, $W$ in Proposition \ref{L1} then is an $%
\lambda _{\mathbb{R}^{n}}$-null set and the proposition holds trivially.
Similarly, the sets $W_{j}$ in Proposition \ref{L2} are then $\lambda _{%
\mathbb{R}^{n}}$-null sets from a certain $j$ onwards.

(ii) In case $\mathcal{L}=\mathfrak{M}_{0}^{lin}$, the case $\dim (\mathcal{L%
})=n-1$ can not arise, since $k<n$ and $q\geq 1$ are always assumed.
\end{remark}

The following observation applies to a large class of test statistics and
is, e.g., useful when constructing confidence sets by \textquotedblleft
inverting\textquotedblright\ the corresponding test.

\begin{lemma}
\label{neu}Suppose $T$ is a nonsphericity-corrected $F$-type test statistic
as defined in (28) in Section 5.4 of \cite{PP2016} and that Assumption 5 in
that paper is satisfied. Then the rejection probabilities under the null
hypothesis described by (\ref{testing problem}) do not depend on the value
of $r$. As a consequence, the size-controlling critical values $C(\alpha )$
(if they exist) do not depend on the value of $r$ either. [It is understood
here that the estimators $\check{\beta}$ and $\check{\Omega}$ used to define
the test statistic $T$ have been chosen independently of the value of $r$.]
\end{lemma}

\subsection{Some sufficient conditions\label{suff_con}}

We collect here sufficient conditions for some of the assumptions on $T$ in
the preceding results. For sufficient conditions relating to $\mathbb{J}(%
\mathcal{L},\mathfrak{C})$ see Section \ref{suff_con_2} in Appendix \ref{App
B}.

\begin{lemma}
\label{Fact1}Suppose $T$ is a nonsphericity-corrected $F$-type test
statistic as defined in (28) in Section 5.4 of \cite{PP2016} and that
Assumption 5 in that paper is satisfied. Then:

\begin{enumerate}
\item $T$ is Borel-measurable and is continuous on the complement of a
closed $\lambda _{\mathbb{R}^{n}}$-null set $N^{\ast }$ (with $N^{\ast }$
given in (27) of \cite{PP2016}). Furthermore, $T$ and $N^{\ast }$ are $G(%
\mathfrak{M}_{0})$-invariant (in fact, $N^{\ast }$ is even $G(\mathfrak{M})$%
-invariant).

\item $\lambda _{\mathbb{R}^{n}}(T=C)=0$ holds for $-\infty <C<\infty $.

\item The complement of the rejection region $\left\{ T\geq C\right\} $ is
not a $\lambda _{\mathbb{R}^{n}}$-null set for every $C>0$.
\end{enumerate}
\end{lemma}

\begin{remark}
\label{special cases}\emph{(Special cases)} Lemma \ref{Fact1} applies, in
particular, to the commonly used autocorrelation robust test statistic $%
T_{w} $ given in (\ref{tlr}) provided Assumptions \ref{AW} and \ref{R_and_X}
are satisfied, since such test statistics then are nonsphericity-corrected $%
F $-type test statistics and the above mentioned Assumption 5 is satisfied,
cf. Lemma A.1 in \cite{PP2016}. The same is true, more generally, for the
test statistics $T_{GQ}$ defined in Section \ref{extension} whenever the
weighting matrix $\mathcal{W}_{n}^{\ast }$ (also defined in that section) is
positive definite and Assumption \ref{R_and_X} holds (this is proved in the
same way as Lemma A.1 in \cite{PP2016} using Lemma 3.11 instead of Lemma 3.1
in that reference). Furthermore, the weighted Eicker-test statistic $T_{E,%
\mathsf{W}}$ defined in Section \ref{autocorr} (with $\mathsf{W}$ a
symmetric and nonnegative definite $n\times n$ Toeplitz matrix with ones on
the main diagonal) is a nonsphericity-corrected $F$-type test statistic with
the above mentioned Assumption 5 being always satisfied; hence Lemma \ref%
{Fact1} also applies to the weighted Eicker-test statistic. [In fact, also
Assumptions 6 and 7 in \cite{PP2016} are satisfied for these three test
statistics (under the respective assumptions mentioned above), since any one
of $\hat{\Omega}_{w}(y)$, $\hat{\Omega}_{GQ}(y)$, and $\hat{\Omega}_{E,%
\mathsf{W}}(y)$ is then nonnegative definite for every $y\in \mathbb{R}^{n}$%
.]
\end{remark}

For the test statistics mentioned in the preceding remark more can be said
about the set $N^{\ast }$. For the weighted Eicker-test statistic (with $%
\mathsf{W}$ as in Remark \ref{special cases}) we always have $N^{\ast }=%
\limfunc{span}(X)$, thus it is a proper linear subspace of $\mathbb{R}^{n}$.
For autocorrelation robust test statistics of the form (\ref{tlr}) and under
Assumptions \ref{AW} and \ref{R_and_X} the set $N^{\ast }$ turns out to be
the set $\mathsf{B}$ defined in Section \ref{autocorr}, cf. Lemmata 3.1,
A.1, and 5.15 of \cite{PP2016}; and this is more generally true for the test
statistics $T_{GQ}$ provided the weighting matrix $\mathcal{W}_{n}^{\ast }$
is positive definite and Assumption \ref{R_and_X} holds, cf. Lemma 3.11 of 
\cite{PP2016}. The following is now true for the set $\mathsf{B}$.

\begin{lemma}
\label{Fact2}Suppose Assumption \ref{R_and_X} holds. Then $\mathsf{B}$ is a
finite union of proper linear subspaces of $\mathbb{R}^{n}$; in case $q=1$, $%
\mathsf{B}$ is a proper linear subspace. [Consequently, the same applies to
the set $N^{\ast }$ associated with the autocorrelation robust test
statistic $T_{w}$ defined in (\ref{tlr}) (or with $T_{GQ}$, respectively) if
also Assumption \ref{AW} (or the condition $\mathcal{W}_{n}^{\ast }$
positive definite, respectively) applies.]
\end{lemma}

Simple examples show, that in case $q>1$, the set $\mathsf{B}$ in Lemma \ref%
{Fact2} need not be a linear space itself. However, generically $\mathsf{B}=%
\limfunc{span}(X)$ holds under typical assumptions as is shown in Lemma \ref%
{generic_lem_2} in Appendix \ref{App A}, cf. also Theorem \ref%
{genericity_cor}. The next lemma verifies that the condition (b) in
Proposition \ref{P6} is often satisfied.

\begin{lemma}
\label{Fact3}(i) Let Assumptions \ref{AW} and \ref{R_and_X} hold and let $%
T_{w}$ be defined as in (\ref{tlr}). Suppose $\mathcal{S}$ is a linear
subspace of $\mathbb{R}^{n}$ and $\mu $ is an element of $\mathbb{R}^{n}$
such that $T_{w}$ restricted to $\mu +\mathcal{S}$ is not equal to a
constant $\lambda _{\mu +\mathcal{S}}$-almost everywhere. Then, for every
real number $C$ it holds that $\lambda _{\mu +\mathcal{S}}(T_{w}=C)=0$.

(ii) Let Assumption \ref{R_and_X} hold and let $T_{GQ}$ be the test
statistic as defined in Section \ref{extension} with a positive definite
weighting matrix $\mathcal{W}_{n}^{\ast }$. Then the same conclusion as in
(i) holds with $T_{w}$ replaced by $T_{GQ}$.

(iii) Let $\mathsf{W}$ be an $n\times n$ symmetric and nonnegative definite
Toeplitz matrix of weights with ones on the main diagonal. Then the same
conclusion as in (i) holds with $T_{w}$ replaced by $T_{E,\mathsf{W}}$.
\end{lemma}

Lemmata \ref{Fact1}-\ref{Fact3} can obviously be used to provide streamlined
versions of Propositions \ref{L1}, \ref{L2}, \ref{P6}, Lemma \ref{L4}, as
well as Corollary \ref{C5} in the case of nonsphericity-corrected $F$-type
tests, and, in particular, in the case of autocorrelation robust tests based
on $T_{w}$, $T_{GQ}$, or $T_{E,\mathsf{W}}$. We abstain from presenting such
results.

\section{Size control of tests of affine restrictions in regression models
with stationary autocorrelated errors: General results\label{structure under
Topelitz}}

It transpires from the results in Section \ref{general} that characterizing
the elements of the collection $\mathbb{J}(\mathcal{L},\mathfrak{C})$ is
central to achieving explicit conditions for size control. In this section
we undertake such a characterization for the important case where $\mathfrak{%
C}=\mathfrak{C}(\mathfrak{F})$, i.e., when the errors in the regression
model come from a stationary process. In Subsection \ref{struct} we present
the characterization result, which in turn forms the basis for the size
control results in Subsection \ref{szctrl_3}.

\subsection{The structure of $\mathbb{J}(\mathcal{L},\mathfrak{C}(\mathfrak{F%
}))$ \label{struct}}

Before we can state the main results of this subsection we need to introduce
some more notation.

\begin{definition}
For $\omega \in \lbrack 0,\pi ]$ and for $d\in \mathbb{N}$ we define $\kappa
(\omega ,d)=d$ if $\omega \in \{0,\pi \}$, and $\kappa (\omega ,d)=2d$ if $%
\omega \in (0,\pi )$. For a positive integer $p$, for $\underline{\omega }%
=(\omega _{1},\ldots ,\omega _{p})\in \lbrack 0,\pi ]^{p}$, and for $%
\underline{d}=(d_{1},\ldots ,d_{p})\in \mathbb{N}^{p}$ we define 
\begin{equation}
\kappa (\underline{\omega },\underline{d})=\sum_{i=1}^{p}\kappa (\omega
_{i},d_{i}).  \label{deltadef}
\end{equation}%
In case $p=1$, we shall often simply write $\omega $ for $\underline{\omega }
$ and $d$ for $\underline{d}$.
\end{definition}

It proves useful to introduce the convention that $\underline{\omega }$ and $%
\underline{d}$ are the $0$-tupels for $p=0$ and to set $\kappa (\underline{%
\omega },\underline{d})=0$ in this case. The following notation will be
helpful.

\begin{definition}
\label{subspaces_2} Let $\mathcal{L}$ be a linear subspace of $\mathbb{R}%
^{n} $ with $\dim (\mathcal{L})<n$. We write $\underline{\omega }(\mathcal{L}%
)$ for the vector obtained by ordering the elements of 
\begin{equation}
\{\omega \in \lbrack 0,\pi ]:\rho (\omega ,\mathcal{L})>0\}  \label{subs}
\end{equation}%
from smallest to largest, provided the set in (\ref{subs}) is nonempty, and
we denote by $p(\mathcal{L})$ the dimension of this vector (clearly $p(%
\mathcal{L})>0$ then holds); furthermore, we set $d_{i}(\mathcal{L})=\rho
(\omega _{i}(\mathcal{L}),\mathcal{L})$ for every $i=1,\ldots ,p(\mathcal{L}%
) $ (where $\omega _{i}(\mathcal{L})$ denotes the $i$-th coordinate of $%
\underline{\omega }(\mathcal{L})$), and we write $\underline{d}(\mathcal{L})$
for the vector with $i$-th coordinate equal to $d_{i}(\mathcal{L})$. If the
set in (\ref{subs}) is empty, we take $\underline{\omega }(\mathcal{L})$ as
well as $\underline{d}(\mathcal{L})$ as the $0$-tupel (which we may identify
with the empty set) and we set $p(\mathcal{L})=0$.
\end{definition}

Recall that the set (\ref{subs}) is always a finite set as discussed
subsequent to Definition \ref{subspaces} and hence $\underline{\omega }(%
\mathcal{L})$ is well-defined. Clearly, the set (\ref{subs}) coincides with
the set $\{\omega \in \lbrack 0,\pi ]:\limfunc{span}(E_{n,0}(\omega
))\subseteq \mathcal{L}\}$. Furthermore, $\kappa (\underline{\omega }(%
\mathcal{L}),\underline{d}(\mathcal{L}))=0$ if $p(\mathcal{L})=0$ in view of
the above conventions.

\begin{definition}
\label{difffactor}For $\omega \in \lbrack 0,\pi ]$ define polynomials in the
variable $z$ as $\Delta _{\omega }(z)=1-\cos (\omega )z$ if $\omega \in
\{0,\pi \}$ and as $\Delta _{\omega }(z)=1-2\cos (\omega )z+z^{2}$ if $%
\omega \in (0,\pi )$. For $p$ a positive integer, for $\underline{\omega }%
=(\omega _{1},\ldots ,\omega _{p})\in \lbrack 0,\pi ]^{p}$, and for $%
\underline{d}=(d_{1},\ldots ,d_{p})\in \mathbb{N}^{p}$ let the polynomial $%
\Delta _{\underline{\omega },\underline{d}}(z)$ be defined as the product%
\begin{equation*}
\Delta _{\underline{\omega },\underline{d}}(z)=\prod_{i=1}^{p}\Delta
_{\omega _{i}}^{d_{i}}(z).
\end{equation*}%
In case $p=0$ (and thus $\underline{\omega }$ and $\underline{d}$ are the $0$%
-tupels) we define $\Delta _{\underline{\omega },\underline{d}}$ as the
constant polynomial $1$.
\end{definition}

Note that the degree of $\Delta _{\omega _{i}}^{d_{i}}(z)$ is $\kappa
(\omega _{i},d_{i})$ and that of $\Delta _{\underline{\omega },\underline{d}%
}(z)$ is $\kappa (\underline{\omega },\underline{d})$.

A finite Borel measure $\mathsf{m}$ on $[-\pi ,\pi ]$ is said to be
symmetric, if $\mathsf{m}(A)=\mathsf{m}(-A)$ for every Borel subset $A$ of $%
[-\pi ,\pi ]$. Recall that the finite and symmetric Borel measures on $[-\pi
,\pi ]$ are precisely the spectral measures of real weakly stationary
processes. For a spectral density $g$, we denote by $\mathsf{m}_{g}$ the
Borel measure on $[-\pi ,\pi ]$ with density $g$ (w.r.t. Lebesgue measure $%
\lambda _{\lbrack -\pi ,\pi ]}$ on $[-\pi ,\pi ]$).

\begin{definition}
\label{singularToep} Let $\mathfrak{F}\subseteq \mathfrak{F}_{\mathrm{all}}$
be nonempty, and let $\mathcal{L}$ be a linear subspace of $\mathbb{R}^{n}$
with $\dim (\mathcal{L})<n$. We define $\mathbb{M}(\mathfrak{F},\mathcal{L})$
to be the set of all finite and symmetric Borel measures $\mathsf{m}$ on $%
[-\pi ,\pi ]$ with finite support, such that (i) $\mathsf{m}$ is the weak
limit of a sequence $\mathsf{m}_{g_{j}}$, $j\in \mathbb{N}$, where 
\begin{equation}
g_{j}(\nu )=\left\vert \Delta _{\underline{\omega }(\mathcal{L}),\underline{d%
}(\mathcal{L})}(e^{\iota \nu })\right\vert ^{2}f_{j}(\nu )/\int_{-\pi }^{\pi
}\left\vert \Delta _{\underline{\omega }(\mathcal{L}),\underline{d}(\mathcal{%
L})}(e^{\iota \nu })\right\vert ^{2}f_{j}(\nu )d\nu  \label{g_j}
\end{equation}%
for some sequence $f_{j}\in \mathfrak{F}$, and (ii) $\sum_{\gamma \in 
\mathrm{\limfunc{supp}}(\mathsf{m})\cap \lbrack 0,\pi ]}\kappa (\gamma
,1)<n-\kappa (\underline{\omega }(\mathcal{L}),\underline{d}(\mathcal{L}))$
holds. We furthermore define $\mathbb{S}(\mathfrak{F},\mathcal{L})=\left\{ 
\mathrm{\limfunc{supp}}(\mathsf{m})\cap \lbrack 0,\pi ]:\mathsf{m}\in 
\mathbb{M}(\mathfrak{F},\mathcal{L})\right\} $.
\end{definition}

We note that by construction the elements of $\mathbb{M}(\mathfrak{F},%
\mathcal{L})$ all have total mass equal to one, and have a nonempty and
finite support. Furthermore, $\mathbb{M}(\mathfrak{F},\mathcal{L})$, and
thus $\mathbb{S}(\mathfrak{F},\mathcal{L})$, can be empty. Also recall that $%
\kappa (\underline{\omega }(\mathcal{L}),\underline{d}(\mathcal{L}))\leq
\dim (\mathcal{L})<n$ holds by Lemma \ref{dimsubspaces} in Appendix \ref{App
D}.

To illustrate the concepts introduced in the preceding definition we shall
now determine the set $\mathbb{S}(\mathfrak{F},\mathcal{L})$ for a few
choices of $\mathfrak{F}$. The proofs of the claims made in the next two
examples can be found in Appendix \ref{App D}. Recall that $\mathfrak{F}_{%
\mathrm{all}}$ denotes the set of all normalized spectral densities.

\begin{example}
\label{chsingularToepmax} Let $\mathcal{L}$ be a linear subspace of $\mathbb{%
R}^{n}$ with $\dim (\mathcal{L})<n$. Then%
\begin{equation}
\mathbb{S}(\mathfrak{F}_{\mathrm{all}},\mathcal{L})=\left\{ \Gamma \subseteq
\lbrack 0,\pi ]:\func{card}(\Gamma )\in \mathbb{N},\sum_{\gamma \in \Gamma
}\kappa (\gamma ,1)<n-\kappa (\underline{\omega }(\mathcal{L}),\underline{d}(%
\mathcal{L}))\right\} .  \label{S_of_F_all}
\end{equation}%
Note that this set is empty if $n=\kappa (\underline{\omega }(\mathcal{L}),%
\underline{d}(\mathcal{L}))+1$, and is equal to $\{\left\{ 0\right\}
,\left\{ \pi \right\} \}$ if $n=\kappa (\underline{\omega }(\mathcal{L}),%
\underline{d}(\mathcal{L}))+2$. In case $n>\kappa (\underline{\omega }(%
\mathcal{L}),\underline{d}(\mathcal{L}))+2$, it is an infinite set with the
property that $\bigcup \mathbb{S}(\mathfrak{F}_{\mathrm{all}},\mathcal{L})$
is equal to $[0,\pi ]$; in fact, even $\left\{ \gamma \right\} \in \mathbb{S}%
(\mathfrak{F}_{\mathrm{all}},\mathcal{L})$ holds for every $\gamma \in
\lbrack 0,\pi ]$. [Note that any $\Gamma \in \mathbb{S}(\mathfrak{F}_{%
\mathrm{all}},\mathcal{L})$ in particular satisfies $2\func{card}(\Gamma
)-2<n-\kappa (\underline{\omega }(\mathcal{L}),\underline{d}(\mathcal{L}))$.]
\end{example}

\begin{example}
\label{chsingularToepmin} Let $\mathcal{L}$ be a linear subspace of $\mathbb{%
R}^{n}$ with $\dim (\mathcal{L})<n$. For $B$ such that $0<B<\infty $ let $%
\mathfrak{F}_{\mathrm{all}}^{B}$ denote the subset of elements of $\mathfrak{%
F}_{\mathrm{all}}$ that are $\lambda _{\lbrack -\pi ,\pi ]}$-essentially
bounded by $B$. Then $\mathbb{S}(\mathfrak{F}_{\mathrm{all}}^{B},\mathcal{L}%
)=\emptyset $ for every $0<B<\infty $.
\end{example}

The result in the next example is easy to derive and we leave the proof to
the reader.

\begin{example}
\label{chsingularAR(1)} Let $\mathcal{L}$ be a linear subspace of $\mathbb{R}%
^{n}$ with $\dim (\mathcal{L})<n$ and consider $\mathfrak{F}_{\mathrm{AR(}1%
\mathrm{)}}$. Then we have the following: If $n=\kappa (\underline{\omega }(%
\mathcal{L}),\underline{d}(\mathcal{L}))+1$ then $\mathbb{S}(\mathfrak{F}_{%
\mathrm{AR(}1\mathrm{)}},\mathcal{L})$ is empty. Otherwise, we have four
cases: (i) If neither $e_{+}$ nor $e_{-}$ belong to $\mathcal{L}$, then $%
\mathbb{S}(\mathfrak{F}_{\mathrm{AR(}1\mathrm{)}},\mathcal{L})=\left\{
\left\{ 0\right\} ,\left\{ \pi \right\} \right\} $. (ii) If $e_{+}$ belongs
to $\mathcal{L}$, but $e_{-}$ does not, then $\mathbb{S}(\mathfrak{F}_{%
\mathrm{AR(}1\mathrm{)}},\mathcal{L})=\left\{ \left\{ \pi \right\} \right\} $%
. (iii) If $e_{+}$ does not belong to $\mathcal{L}$, but $e_{-}$ does, then $%
\mathbb{S}(\mathfrak{F}_{\mathrm{AR(}1\mathrm{)}},\mathcal{L})=\left\{
\left\{ 0\right\} \right\} $. (iv) If both $e_{+}$ and $e_{-}$ belong to $%
\mathcal{L}$, then $\mathbb{S}(\mathfrak{F}_{\mathrm{AR(}1\mathrm{)}},%
\mathcal{L})$ is empty.
\end{example}

\begin{proposition}
\label{limitsToep} Let $\mathfrak{F}\subseteq \mathfrak{F}_{\mathrm{all}}$
be nonempty. Let $\mathcal{L}$ be a linear subspace of $\mathbb{R}^{n}$ with 
$\dim (\mathcal{L})<n$.

\begin{enumerate}
\item For every $\bar{\Sigma}\in \limfunc{cl}(\mathcal{L}(\mathfrak{C}(%
\mathfrak{F})))$ such that $\limfunc{rank}(\bar{\Sigma})<n-\dim (\mathcal{L}%
) $ there exists a set $\Gamma \in \mathbb{S}(\mathfrak{F},\mathcal{L})$ and
positive real numbers $c(\gamma )$ for $\gamma \in \Gamma $, such that 
\begin{equation}
\bar{\Sigma}=\frac{\Pi _{\mathcal{L}^{\bot }}\sum_{\gamma \in \Gamma
}c(\gamma )E_{n,\rho (\gamma ,\mathcal{L})}(\gamma )E_{n,\rho (\gamma ,%
\mathcal{L})}^{\prime }(\gamma )\Pi _{\mathcal{L}^{\bot }}}{\Vert {\Pi _{%
\mathcal{L}^{\bot }}\sum_{\gamma \in \Gamma }c(\gamma )E_{n,\rho (\gamma ,%
\mathcal{L})}(\gamma )E_{n,\rho (\gamma ,\mathcal{L})}^{\prime }(\gamma )\Pi
_{\mathcal{L}^{\bot }}}\Vert }  \label{represent}
\end{equation}%
holds. Furthermore, for every $\Gamma \in \mathbb{S}(\mathfrak{F},\mathcal{L}%
)$ there exists $\bar{\Sigma}\in \limfunc{cl}(\mathcal{L}(\mathfrak{C}(%
\mathfrak{F})))$ and positive real numbers $c(\gamma )$ for $\gamma \in
\Gamma $, such that (\ref{represent}) holds (and clearly $\limfunc{rank}(%
\bar{\Sigma})\leq n-\dim (\mathcal{L})$ holds).

\item The set $\mathbb{J}(\mathcal{L},\mathfrak{C}(\mathfrak{F}))$ coincides
with the set of all linear subspaces of $\mathbb{R}^{n}$ that (i) have
dimension smaller than $n-\dim (\mathcal{L})$, and (ii) that can be
expressed as 
\begin{equation*}
\limfunc{span}\left( \Pi _{\mathcal{L}^{\bot }}\left( E_{n,\rho (\gamma _{1},%
\mathcal{L})}(\gamma _{1}),\ldots ,E_{n,\rho (\gamma _{p},\mathcal{L}%
)}(\gamma _{p})\right) \right)
\end{equation*}%
for some $\Gamma \in \mathbb{S}(\mathfrak{F},\mathcal{L})$, where the $%
\gamma _{i}$'s denote the elements of $\Gamma $ and $p$ denotes $\func{card}%
(\Gamma )$.

\item Every element of $\mathbb{J}(\mathcal{L},\mathfrak{C}(\mathfrak{F}))$
contains a subspace of the form $\limfunc{span}(\Pi _{\mathcal{L}^{\bot
}}E_{n,\rho (\gamma ,\mathcal{L})}(\gamma ))$ for some $\gamma \in \bigcup 
\mathbb{S}(\mathfrak{F},\mathcal{L})$.
\end{enumerate}
\end{proposition}

We illustrate the proposition by two examples.

\begin{example}
\emph{(Example \ref{chsingularAR(1)} continued)} If $n=\dim (\mathcal{L})+1$
(which includes the case where $n=\kappa (\underline{\omega }(\mathcal{L}),%
\underline{d}(\mathcal{L}))+1$), the set $\mathbb{J}(\mathcal{L},\mathfrak{C}%
(\mathfrak{F}_{\mathrm{AR(}1\mathrm{)}}))$ is empty. Otherwise, we have the
following cases: (i) If neither $e_{+}$ nor $e_{-}$ belongs to $\mathcal{L}$%
, then $\mathbb{J}(\mathcal{L},\mathfrak{C}(\mathfrak{F}_{\mathrm{AR(}1%
\mathrm{)}}))=\left\{ \limfunc{span}(\Pi _{\mathcal{L}^{\bot }}e_{+}),%
\limfunc{span}(\Pi _{\mathcal{L}^{\bot }}e_{-})\right\} $. (ii) If $e_{+}$
belongs to $\mathcal{L}$, but $e_{-}$ does not, then $\mathbb{J}(\mathcal{L},%
\mathfrak{C}(\mathfrak{F}_{\mathrm{AR(}1\mathrm{)}}))=\left\{ \limfunc{span}%
(\Pi _{\mathcal{L}^{\bot }}e_{-})\right\} $. (iii) If $e_{+}$ does not
belong to $\mathcal{L}$, but $e_{-}$ does, then $\mathbb{J}(\mathcal{L},%
\mathfrak{C}(\mathfrak{F}_{\mathrm{AR(}1\mathrm{)}}))=\left\{ \limfunc{span}%
(\Pi _{\mathcal{L}^{\bot }}e_{+})\right\} $. (iv) If both $e_{+}$ and $e_{-}$
belong to $\mathcal{L}$, then $\mathbb{J}(\mathcal{L},\mathfrak{C}(\mathfrak{%
F}_{\mathrm{AR(}1\mathrm{)}}))$ is empty. [In this simple example these
results can alternatively be obtained from the fact that the concentration
spaces of $\mathfrak{C}(\mathfrak{F}_{\mathrm{AR(}1\mathrm{)}})$ in the
sense of \cite{PP2016} are given by $\limfunc{span}(e_{+})$ and $\limfunc{%
span}(e_{-})$, combined with Lemma \ref{suff_con_for J} in Section \ref%
{suff_con_2} of Appendix \ref{App B} as well as with Lemma G.1 of \cite%
{PP2016}.]
\end{example}

\begin{example}
\label{chsingularAR(2)}Let $\mathcal{L}$ be a linear subspace of $\mathbb{R}%
^{n}$ with $\dim (\mathcal{L})<n$ and consider $\mathfrak{F}\subseteq 
\mathfrak{F}_{\mathrm{all}}$ such that $\mathfrak{F\supseteq }$ $\mathfrak{F}%
_{\mathrm{AR(}2\mathrm{)}}$ (which, in particular, covers the cases $%
\mathfrak{F}=\mathfrak{F}_{\mathrm{all}}$ as well as $\mathfrak{F}=\mathfrak{%
F}_{\mathrm{AR(}p\mathrm{)}}$ for $p\geq 2$). The following is shown in
Section 3 of \cite{PP4}:

(i) Every $\gamma \in \bigcup \mathbb{S}(\mathfrak{F},\mathcal{L})$
satisfies $\left\{ \gamma \right\} \in \mathbb{S}(\mathfrak{F},\mathcal{L})$.

(ii) Suppose $\dim (\mathcal{L})+2<n$. Then $\left\{ \gamma \right\} \in 
\mathbb{S}(\mathfrak{F},\mathcal{L})$ holds for every $\gamma \in \lbrack
0,\pi ]$. And it easily follows that $\limfunc{span}(\Pi _{\mathcal{L}^{\bot
}}E_{n,\rho (\gamma ,\mathcal{L})}(\gamma ))$ belongs to $\mathbb{J}(%
\mathcal{L},\mathfrak{C}(\mathfrak{F}))$ for every $\gamma \in \lbrack 0,\pi
]$.

(iii) Suppose $\dim (\mathcal{L})+2\geq n$. Then $\left\{ \gamma \right\}
\in \mathbb{S}(\mathfrak{F},\mathcal{L})$ holds for $\gamma \in \lbrack
0,\pi ]$ precisely when $\kappa (\gamma ,1)<n-\kappa (\underline{\omega }(%
\mathcal{L}),\underline{d}(\mathcal{L}))$. Furthermore, $\limfunc{span}(\Pi
_{\mathcal{L}^{\bot }}E_{n,\rho (\gamma ,\mathcal{L})}(\gamma ))$ belongs to 
$\mathbb{J}(\mathcal{L},\mathfrak{C}(\mathfrak{F}))$ for every $\gamma \in
\lbrack 0,\pi ]$ that satisfies (a) $\kappa (\gamma ,1)<n-\kappa (\underline{%
\omega }(\mathcal{L}),\underline{d}(\mathcal{L}))$ and (b) $\dim \limfunc{%
span}(\Pi _{\mathcal{L}^{\bot }}E_{n,\rho (\gamma ,\mathcal{L})}(\gamma
))<n-\dim (\mathcal{L})$.
\end{example}

\subsection{Results on size control \label{szctrl_3}}

\begin{theorem}
\label{Ctoep} Let $\mathfrak{F}\subseteq \mathfrak{F}_{\mathrm{all}}$ be
nonempty, and assume that the test statistic $T:\mathbb{R}^{n}\rightarrow 
\mathbb{R}$ is Borel-measurable and is continuous on the complement of a
closed set $N^{\dag }$. Assume that $T$ and $N^{\dag }$ are $G(\mathfrak{M}%
_{0})$-invariant, and are also invariant w.r.t. addition of elements of a
linear subspace $\mathcal{V}$ of $\mathbb{R}^{n}$. Define $\mathcal{L}=%
\mathrm{\limfunc{span}}(\mathfrak{M}_{0}^{lin}\cup \mathcal{V})$, and assume
that $\dim (\mathcal{L})<n$.

\begin{enumerate}
\item Then for every $0<\alpha <1$ there exists a real number $C(\alpha )$
such that%
\begin{equation}
\sup_{\mu _{0}\in \mathfrak{M}_{0}}\sup_{0<\sigma ^{2}<\infty }\sup_{f\in 
\mathfrak{F}}P_{\mu _{0},\sigma ^{2}\Sigma (f)}(T\geq C(\alpha ))\leq \alpha
\label{0limit}
\end{equation}%
holds, provided every linear subspace $\mathcal{S}$, say, of $\mathbb{R}^{n}$
that (i) has dimension smaller than $n-\dim (\mathcal{L})$, and that (ii)
can be written as%
\begin{equation}
\mathcal{S}=\limfunc{span}\left( \Pi _{\mathcal{L}^{\bot }}\left( E_{n,\rho
(\gamma _{1},\mathcal{L})}(\gamma _{1}),\ldots ,E_{n,\rho (\gamma _{p},%
\mathcal{L})}(\gamma _{p})\right) \right) \quad \text{ for some }\quad
\Gamma \in \mathbb{S}(\mathfrak{F},\mathcal{L}),  \label{form_of_S}
\end{equation}%
satisfies $\lambda _{\mu _{0}+\mathcal{S}}(N^{\dag })=0$ for some $\mu
_{0}\in \mathfrak{M}_{0}$ (and hence all $\mu _{0}\in \mathfrak{M}_{0}$).
Here the $\gamma _{i}$'s denote the elements of $\Gamma $ and $p=\func{card}%
(\Gamma )$. [In case $N^{\dag }$ is a finite or countable union of affine
subspaces, we may replace $\lambda _{\mu _{0}+\mathcal{S}}(N^{\dag })=0$ by $%
\mu _{0}+\mathcal{S}\not\subseteq N^{\dag }$ in that condition.]

\item Suppose $N^{\dag }$ is a finite or countable union of affine
subspaces. Then a sufficient condition for the condition given in Part 1 is
that for all $\gamma \in \bigcup \mathbb{S}(\mathfrak{F},\mathcal{L})$ it
holds that $\mu _{0}+\limfunc{span}(\Pi _{\mathcal{L}^{\bot }}E_{n,\rho
(\gamma ,\mathcal{L})}(\gamma ))\not\subseteq N^{\dag }$ for some $\mu
_{0}\in \mathfrak{M}_{0}$ (and hence all $\mu _{0}\in \mathfrak{M}_{0}$). [A
stricter sufficient condition is obtained if $\gamma $ is required to range
over $[0,\pi ]$ rather than only $\bigcup \mathbb{S}(\mathfrak{F},\mathcal{L}%
)$.]
\end{enumerate}
\end{theorem}

\begin{remark}
\label{special_cases}\emph{(Special cases) }(i) The measurability and
continuity assumption on $T$ as well as the $G(\mathfrak{M}_{0})$-invariance
assumption on $T$ and $N^{\dag }$ are automatically satisfied for
nonsphericity-corrected F-type test statistics as defined in (28) in Section
5.4 of \cite{PP2016} provided that Assumption 5 in that paper is satisfied,
see Lemma \ref{Fact1} in Section \ref{suff_con} above.

(ii) Specializing further (cf. Remark \ref{special cases}) to the case where 
$T_{GQ}$ is the test statistic defined in Section \ref{extension} (with
Assumption \ref{R_and_X} being satisfied and with a positive definite
weighting matrix $\mathcal{W}_{n}^{\ast }$), which, in particular, includes
the case of the test statistic $T_{w}$ defined in (\ref{tlr}) (with
Assumptions \ref{AW} and \ref{R_and_X} being satisfied), or to the case of
the weighted Eicker-test statistic $T_{E,\mathsf{W}}$ defined in Section \ref%
{autocorr} (with $\mathsf{W}$ a symmetric and nonnegative definite $n\times
n $ Toeplitz matrix with ones on the main diagonal), we then have that the
set $N^{\dag }$($=N^{\ast }$) is additionally always guaranteed to be a
proper linear subspace or a finite union of proper linear subspaces, see
Lemma \ref{Fact2} in Section \ref{suff_con} and the discussion preceding
this lemma; cf. also Theorem \ref{Ctoep_corr} below and the results in
Section \ref{autocorr}.
\end{remark}

\begin{remark}
\label{simplifications} \emph{(Some simplifications)} (i) The invariance
assumptions on $N^{\dag }$ in Theorem \ref{Ctoep} imply, in particular, that 
$N^{\dag }$ is invariant under addition of elements in $\mathcal{L}$. Hence,
the sufficient condition in Part 1 of Theorem \ref{Ctoep}, namely that $%
\lambda _{\mu _{0}+\mathcal{S}}(N^{\dag })=0$ holds for every linear space $%
\mathcal{S}$ satisfying (i) and (ii) in that theorem, can equivalently be
expressed as the condition that $\lambda _{\mu _{0}+\mathcal{T}(\mathcal{S}%
)}(N^{\dag })=0$ holds for every linear space $\mathcal{S}$ satisfying (i)
and (ii), where $\mathcal{T}(\mathcal{S})$ is shorthand for $\limfunc{span}%
((E_{n,\rho (\gamma _{1},\mathcal{L})}(\gamma _{1}),\ldots ,E_{n,\rho
(\gamma _{p},\mathcal{L})}(\gamma _{p})))$. Similarly, the relation $\mu
_{0}+\mathcal{S}\not\subseteq N^{\dag }$ given in the sentence in
parenthesis in Part 1 of the theorem can equivalently be expressed as $\mu
_{0}+\mathcal{T}(\mathcal{S})\not\subseteq N^{\dag }$. Finally, the
sufficient conditions given in Part 2 of the theorem can equivalently be
expressed as $\mu _{0}+\func{span}(E_{n,\rho (\gamma ,\mathcal{L})}(\gamma
))\not\subseteq N^{\dag }$ for every $\gamma \in \bigcup \mathbb{S}(%
\mathfrak{F},\mathcal{L})$ (for every $\gamma \in \lbrack 0,\pi ]$,
respectively).

(ii) If $r=0$, or if $r\neq 0$ but $N^{\dag }$ is invariant under addition
of elements of $\mathfrak{M}_{0}$ (which is, e.g., the case if $N^{\dag }$
is $G(\mathfrak{M})$-invariant), we may set $\mu _{0}=0$ in any of the
sufficient conditions in Parts 1 and 2 of Theorem \ref{Ctoep} that involve $%
N^{\dag }$, cf. Remark \ref{simplif_1} in Section \ref{two}; and the same
applies to the equivalent formulations of these conditions discussed in (i)
above. [Recall that $N^{\dag }$ is $G(\mathfrak{M})$-invariant for a large
class of test statistics including nonsphericity-corrected F-type test
statistics that satisfy Assumption 5 in \cite{PP2016}; in particular, this
applies to the test statistics mentioned in Remark \ref{special_cases}(ii)
above.]

(iii) Imposing that $\lambda _{\mu _{0}+\mathcal{S}}(N^{\dag })=0$ holds for
every linear space $\mathcal{S}$ satisfying (ii) in Part 1 of Theorem \ref%
{Ctoep} (but not necessarily (i)) leads to a potentially stricter sufficient
condition. However, if $N^{\dag }$ is an $\lambda _{\mathbb{R}^{n}}$-null
set, this potentially stricter condition is in fact equivalent to the
condition given in Part 1 of Theorem \ref{Ctoep}. For a proof see Appendix %
\ref{App D}.

(iv) In case $\dim (\mathcal{L})=n-1$, the test statistic $T$ is $\lambda _{%
\mathbb{R}^{n}}$-almost everywhere constant and size control is hence
trivially possible. Note also that the sufficient condition in Part 1 of
Theorem \ref{Ctoep} is trivially satisfied since $\mathbb{S}(\mathfrak{F},%
\mathcal{L})$ is then empty and hence no $\mathcal{S}$ satisfies (\ref%
{form_of_S}); cf. Remark \ref{dimL=n-1}.
\end{remark}

Specializing to nonsphericity-corrected F-type test statistics as defined in
Section 5.4 of \cite{PP2016} and to the important case $\mathcal{L}=%
\mathfrak{M}_{0}^{lin}$ gives the following result.

\begin{theorem}
\label{Ctoep_nonsphericity} Let $\mathfrak{F}\subseteq \mathfrak{F}_{\mathrm{%
all}}$ be nonempty, and suppose that $T$ is a nonsphericity-corrected $F$%
-type test statistic as defined in (28) in Section 5.4 of \cite{PP2016} and
that Assumption 5 in that paper is satisfied.

\begin{enumerate}
\item Then for every $0<\alpha <1$ there exists a real number $C(\alpha )$
such that (\ref{0limit}) holds, provided%
\begin{equation}
\lambda _{\limfunc{span}((E_{n,\rho (\gamma _{1},\mathfrak{M}%
_{0}^{lin})}(\gamma _{1}),\ldots ,E_{n,\rho (\gamma _{p},\mathfrak{M}%
_{0}^{lin})}(\gamma _{p})))}(N^{\ast })=0  \label{suff_cond_nonspher}
\end{equation}%
holds for every $\Gamma \in \mathbb{S}(\mathfrak{F},\mathfrak{M}_{0}^{lin})$%
. Here the $\gamma _{i}$'s denote the elements of $\Gamma $ and $p=\func{card%
}(\Gamma )$.

\item Suppose that $N^{\ast }$ is a finite or countable union of affine
subspaces. Then, for every $\Gamma \in \mathbb{S}(\mathfrak{F},\mathfrak{M}%
_{0}^{lin})$, we may rewrite (\ref{suff_cond_nonspher}) equivalently as 
\begin{equation}
\limfunc{span}\left( (E_{n,\rho (\gamma _{1},\mathfrak{M}_{0}^{lin})}(\gamma
_{1}),\ldots ,E_{n,\rho (\gamma _{p},\mathfrak{M}_{0}^{lin})}(\gamma
_{p}))\right) \not\subseteq N^{\ast }.  \label{non-incl-infty}
\end{equation}%
Furthermore, a sufficient condition for this is given by $\limfunc{span}%
(E_{n,\rho (\gamma ,\mathfrak{M}_{0}^{lin})}(\gamma ))\not\subseteq N^{\ast
} $ for every $\gamma \in \bigcup \mathbb{S}(\mathfrak{F},\mathfrak{M}%
_{0}^{lin})$ (and an even stricter sufficient condition is obtained if $%
\gamma $ is here required to range over $[0,\pi ]$ instead of only over $%
\bigcup \mathbb{S}(\mathfrak{F},\mathfrak{M}_{0}^{lin})$).
\end{enumerate}
\end{theorem}

Specializing Theorem \ref{Ctoep_nonsphericity} to commonly used
autocorrelation robust tests based on $T_{w}$ we obtain the following
result. Since for these tests statistics also all assumptions of Proposition %
\ref{P6} can be shown to be satisfied (cf. Section \ref{suff_con}), the size
control result can furthermore be sharpened to an \emph{exact} size control
result.

\begin{theorem}
\label{Ctoep_corr} Let $\mathfrak{F}\subseteq \mathfrak{F}_{\mathrm{all}}$
be nonempty. Suppose Assumptions \ref{AW} and \ref{R_and_X} are satisfied
and $T_{w}$ is defined by (\ref{tlr}).

\begin{enumerate}
\item Then for every $0<\alpha <1$ there exists a real number $C(\alpha )$
such that (\ref{0limit}) holds (with $T$ replaced by $T_{w}$), provided%
\begin{equation}
\limfunc{span}\left( (E_{n,\rho (\gamma _{1},\mathfrak{M}_{0}^{lin})}(\gamma
_{1}),\ldots ,E_{n,\rho (\gamma _{p},\mathfrak{M}_{0}^{lin})}(\gamma
_{p}))\right) \not\subseteq \mathsf{B}  \label{non-incl_1}
\end{equation}%
holds for every $\Gamma \in \mathbb{S}(\mathfrak{F},\mathfrak{M}_{0}^{lin})$%
. Here the $\gamma _{i}$'s denote the elements of $\Gamma $ and $p=\func{card%
}(\Gamma )$. Furthermore, under the same condition even equality can be
achieved in (\ref{0limit}) (with $T$ replaced by $T_{w}$) by a proper choice
of $C(\alpha )$, provided $\alpha \in (0,\alpha ^{\ast }]\cap (0,1)$, where $%
\alpha ^{\ast }$ is defined as in (\ref{alpha*}) (with $T$ replaced by $%
T_{w} $).

\item A sufficient condition for (\ref{non-incl_1}) to hold for every $%
\Gamma \in \mathbb{S}(\mathfrak{F},\mathfrak{M}_{0}^{lin})$ is given by%
\begin{equation}
\limfunc{span}\left( E_{n,\rho (\gamma ,\mathfrak{M}_{0}^{lin})}(\gamma
)\right) \not\subseteq \mathsf{B}  \label{non-incl_2}
\end{equation}%
for every $\gamma \in \bigcup \mathbb{S}(\mathfrak{F},\mathfrak{M}%
_{0}^{lin}) $ (and an even stricter sufficient condition is obtained if we
require (\ref{non-incl_2}) to hold for every $\gamma \in \lbrack 0,\pi ]$).

\item In case the set $\mathsf{B}$ coincides with $\func{span}(X)$,
condition (\ref{non-incl_1}) ((\ref{non-incl_2}), respectively) can
equivalently be expressed as $\func{rank}(X,E_{n,\rho (\gamma _{1},\mathfrak{%
M}_{0}^{lin})}(\gamma _{1}),\ldots ,E_{n,\rho (\gamma _{p},\mathfrak{M}%
_{0}^{lin})}(\gamma _{p}))>k$ ($\func{rank}(X,E_{n,\rho (\gamma ,\mathfrak{M}%
_{0}^{lin})}(\gamma ))>k$, respectively).
\end{enumerate}
\end{theorem}

For the following remark recall that $\mathsf{CV}_{\leq }(\alpha )$ and $%
\mathsf{CV}_{=}(\alpha )$ have been defined in Remark \ref{smallest C}.

\begin{remark}
\label{exist}(i) Under the assumptions of Theorem \ref{Ctoep_nonsphericity}
or Theorem \ref{Ctoep_corr} we have that $\mathsf{CV}_{\leq }(\alpha
)=[C_{\Diamond }(\alpha ),\infty )$ holds, i.e., a \emph{smallest }critical
value guaranteeing size control exists and is given by $C_{\Diamond }(\alpha
)$. This follows from Remark \ref{smallest C}, Lemma \ref{Fact1}, and Remark %
\ref{special cases}.

(ii) Under the assumptions of Theorem \ref{Ctoep_corr} and if $\alpha \in
(0,\alpha ^{\ast }]\cap (0,1)$, we have that $\mathsf{CV}_{=}(\alpha )$ has $%
C_{\Diamond }(\alpha )$ as its smallest element, i.e., a \emph{smallest }%
critical value guaranteeing exact size control exists and is given by $%
C_{\Diamond }(\alpha )$.

(iii) Under the assumptions of Theorems \ref{Ctoep_nonsphericity} or \ref%
{Ctoep_corr} the size of the test, and hence the size-controlling critical
values $C(\alpha )$, do not depend on the value of $r$; cf. Lemma \ref{neu}.
Also the sufficient conditions in both theorems do not depend on the value
of $r$. [It is understood here that the estimator for $\beta $ as well as
the covariance matrix estimator used to define the test statistic have been
chosen independently of the value of $r$.]
\end{remark}

\begin{remark}
\label{ext}(i) Theorem \ref{Ctoep_corr} carries over to the test statistics $%
T_{GQ}$ (see Section \ref{extension}) if Assumption \ref{AW} is replaced by
the assumption that the weighting matrix $\mathcal{W}_{n}^{\ast }$ is
positive definite. This follows from an inspection of the proof of Theorem %
\ref{Ctoep_corr} and the fact that $T_{GQ}$ is a special case of a
nonsphericity-corrected $F$-type test statistic with $N^{\ast }=\mathsf{B}$,
see Section \ref{suff_con}. A statement analogous to Remark \ref{exist} also
applies here.

(ii) A result for the weighted Eicker-test similar to Theorem \ref%
{Ctoep_corr} is obtained by replacing Assumptions \ref{AW} and \ref{R_and_X}
in that theorem by the assumption that $\mathsf{W}$ is a symmetric and
nonnegative definite $n\times n$ Toeplitz matrix with ones on the main
diagonal and by replacing the set $\mathsf{B}$ by $\func{span}(X)$. This
follows again from an inspection of the proof of Theorem \ref{Ctoep_corr}
and the fact that $T_{E,\mathsf{W}}$ is a special case of a
nonsphericity-corrected $F$-type test statistic with $N^{\ast }=\func{span}%
(X)$, see Section \ref{suff_con}. A statement analogous to Remark \ref{exist}
also applies here.
\end{remark}

\begin{remark}
\label{nec}Let $\mathcal{L}$ be a linear subspace of $\mathbb{R}^{n}$ with $%
\dim (\mathcal{L})<n$ and suppose $\mathfrak{F}\subseteq \mathfrak{F}_{%
\mathrm{all}}$ has the property that $\gamma \in \bigcup \mathbb{S}(%
\mathfrak{F},\mathcal{L})$ implies $\left\{ \gamma \right\} \in \mathbb{S}(%
\mathfrak{F},\mathcal{L})$. [Note that this is always the case if $\mathfrak{%
F}$ contains $\mathfrak{F}_{\mathrm{AR(}2\mathrm{)}}$ as discussed in
Example \ref{chsingularAR(2)}.]

(i) Suppose $\mathcal{L}=\mathfrak{M}_{0}^{lin}$. Then in the context of
Theorem \ref{Ctoep_corr} it is obvious that the first sufficient condition
in Part 2 of that theorem (i.e., the condition that (\ref{non-incl_2}) holds
for every $\gamma \in \bigcup \mathbb{S}(\mathfrak{F},\mathfrak{M}%
_{0}^{lin}) $) is actually equivalent to the sufficient condition in Part 1
(i.e., to the condition that (\ref{non-incl_1}) holds for every $\Gamma \in 
\mathbb{S}(\mathfrak{F},\mathfrak{M}_{0}^{lin})$). A similar remark applies
to the versions of Theorem \ref{Ctoep_corr} for $T_{GQ}$ and $T_{E,\mathsf{W}%
}$ outlined in Remark \ref{ext}.

(ii) Suppose $\mathcal{L}=\mathfrak{M}_{0}^{lin}$. Then similarly in the
context of Part 2 of Theorem \ref{Ctoep_nonsphericity} the condition that (%
\ref{non-incl-infty}) holds for every $\Gamma \in \mathbb{S}(\mathfrak{F},%
\mathfrak{M}_{0}^{lin})$ is equivalent to $\limfunc{span}(E_{n,\rho (\gamma ,%
\mathfrak{M}_{0}^{lin})}(\gamma ))\not\subseteq N^{\ast }$ for every $\gamma
\in \bigcup \mathbb{S}(\mathfrak{F},\mathfrak{M}_{0}^{lin})$.

(iii) Under the assumptions of Part 2 of Theorem \ref{Ctoep}, the first
sufficient condition given in Part 2 is in fact equivalent to the sufficient
condition given in Part 1 of that theorem provided $N^{\dag }\varsubsetneqq 
\mathbb{R}^{n}$. This follows from Remark \ref{simplifications}(iii), noting
that $N^{\dag }$, as a finite or countable union of affine subspaces, is
then a $\lambda _{\mathbb{R}^{n}}$-null set. [In case $N^{\dag }=\mathbb{R}%
^{n}$, the claim is clearly also true provided $\limfunc{rank}(\Pi _{%
\mathcal{L}^{\bot }}E_{n,\rho (\gamma ,\mathcal{L})}(\gamma ))<n-\dim (%
\mathcal{L})$ holds for some $\gamma \in \bigcup \mathbb{S}(\mathfrak{F},%
\mathcal{L})$ or if $\bigcup \mathbb{S}(\mathfrak{F},\mathcal{L})$ is empty.
A sufficient condition for the rank condition just mentioned is that $\kappa
(\gamma ,1)<n-\dim (\mathcal{L})$ holds for such a $\gamma $; and this
latter condition certainly holds if $\dim (\mathcal{L})<n-2$.]\footnote{%
Suppose that $\limfunc{rank}(\Pi _{\mathcal{L}^{\bot }}E_{n,\rho (\gamma ,%
\mathcal{L})}(\gamma ))<n-\dim (\mathcal{L})$ holds whenever $\gamma \in
\bigcup \mathbb{S}(\mathfrak{F},\mathcal{L})$. Then actually a much simpler
argument shows that the claimed equivalence is true, regardless of whether $%
N^{\dag }\varsubsetneqq \mathbb{R}^{n}$ or $N^{\dag }=\mathbb{R}^{n}$. A
sufficient condition for the just mentioned rank condition is that $\kappa
(\gamma ,1)<n-\dim (\mathcal{L})$ holds whenever $\gamma \in \bigcup \mathbb{%
S}(\mathfrak{F},\mathcal{L})$, which is in turn implied by $\dim (\mathcal{L}%
)<n-2$.}
\end{remark}

\begin{remark}
\label{nec_2}Suppose $\mathfrak{F}\subseteq \mathfrak{F}_{\mathrm{all}}$,
but $\mathfrak{F}\supseteq \mathfrak{F}_{\mathrm{AR(}2\mathrm{)}}$. Under
the assumptions of Part 2 of Theorem \ref{Ctoep}, the sufficient condition
given in Part 2 as well as the stricter sufficient condition given in the
final sentence in brackets are in fact equivalent provided $N^{\dag
}\varsubsetneqq \mathbb{R}^{n}$.\footnote{\label{fn_1000}If $\mathbb{S}(%
\mathfrak{F},\mathcal{L})$ is nonempty, then the claimed equivalence is also
true if $N^{\dag }=\mathbb{R}^{n}$ as is easily seen. Note that $\mathbb{S}(%
\mathfrak{F},\mathcal{L})$ is nonempty if and only if $\kappa (\underline{%
\omega }(\mathcal{L}),\underline{d}(\mathcal{L}))<n-1$ in view of Example %
\ref{chsingularAR(2)}. A sufficient condition for the latter inequality is
that $\dim (\mathcal{L})<n-1$, cf. Lemma \ref{dimsubspaces}.} Furthermore,
in the context of Part 2 of Theorem \ref{Ctoep_nonsphericity}, the
sufficient condition $\limfunc{span}(E_{n,\rho (\gamma ,\mathfrak{M}%
_{0}^{lin})}(\gamma ))\not\subseteq N^{\ast }$ for every $\gamma \in \bigcup 
\mathbb{S}(\mathfrak{F},\mathfrak{M}_{0}^{lin})$ is in fact equivalent to
the stricter condition $\limfunc{span}(E_{n,\rho (\gamma ,\mathfrak{M}%
_{0}^{lin})}(\gamma ))\not\subseteq N^{\ast }$ for every $\gamma \in \lbrack
0,\pi ]$. Similarly, in the context of Part 2 of Theorem \ref{Ctoep_corr},
the sufficient condition that (\ref{non-incl_2}) holds for every $\gamma \in
\bigcup \mathbb{S}(\mathfrak{F},\mathfrak{M}_{0}^{lin})$ is in fact
equivalent to the stricter condition that (\ref{non-incl_2}) holds for every 
$\gamma \in \lbrack 0,\pi ]$. A similar remark applies to the versions of
Theorem \ref{Ctoep_corr} for $T_{GQ}$ and $T_{E,\mathsf{W}}$ outlined in
Remark \ref{ext}. A proof of these claims can be found in Appendix \ref{App
D}.\footnote{%
In case $\dim (\mathcal{L})+2<n$, these claims as well as the claim in
Footnote \ref{fn_1000} are actually a simple consequence of Example \ref%
{chsingularAR(2)}(ii) as then $\bigcup \mathbb{S}(\mathfrak{F},\mathcal{L})$
coincides with $[0,\pi ]$.}
\end{remark}

\begin{remark}
\label{doublestar}If $\mathbb{J}(\mathcal{L},\mathfrak{C}(\mathfrak{F}))$ is
empty (which is, e.g., the case if $\mathfrak{C}(\mathfrak{F})$ has no
singular limit points), then size control in the contexts of the theorems in
this subsection is always possible; cf. Remark \ref{star}.
\end{remark}

\appendix

\section{Appendix\label{App A}: Proofs for Section \protect\ref{autocorr}}

The quantities $V_{n}^{(0)}\left( \underline{\omega },\underline{d}\right) $
and $\kappa (\underline{\omega },\underline{d})$ used in the subsequent
proofs are defined in Appendix \ref{App C} and Section \ref{struct},
respectively.

\textbf{Proof of claims regarding Definition \ref{subspaces}:} Suppose $%
p\geq 1$ and suppose $\underline{\omega }$ is a $p\times 1$ vector with
distinct coordinates $\omega _{i}\in \lbrack 0,\pi ]$, such that $\limfunc{%
span}(E_{n,j}(\omega _{i}))\subseteq \mathcal{L}$ for $j=0,\ldots ,d_{i}-1$,
holds where $d_{i}\in \mathbb{N}$. Set $\underline{d}=(d_{1},\ldots ,d_{p})$%
. Since $\mathcal{L}$ is a linear space, it follows that the span of the
matrix $V_{n}^{(0)}\left( \underline{\omega },\underline{d}\right) $ is
contained in $\mathcal{L}$, and thus $\limfunc{rank}(V_{n}^{(0)}\left( 
\underline{\omega },\underline{d}\right) )\leq \dim (\mathcal{L})<n$ must
hold. But as shown in Lemma \ref{fullrank} in Appendix \ref{App C}, the rank
of $V_{n}^{(0)}\left( \underline{\omega },\underline{d}\right) $ equals $%
\min (n,\kappa (\underline{\omega },\underline{d}))$. Consequently, $\kappa (%
\underline{\omega },\underline{d})<n$ follows. Inspection of the definition
of $\kappa (\underline{\omega },\underline{d})$ now shows that this
inequality implies an upper bound on $p$ and the coordinates of $\underline{d%
}$. This obviously then implies that the set on the right-hand side of (\ref%
{def_rho}) contains $0$ for every $\omega \in \lbrack 0,\pi ]$ except
possibly for at most finitely many $\omega $'s; it also implies that the set
on the right-hand side of (\ref{def_rho}) is nonempty for every $\omega \in
\lbrack 0,\pi ]$. But this establishes the claims. $\blacksquare $

\textbf{Proof of claim in Example \ref{exa2}: }To prove this claim we use
Theorem \ref{genericity_cor} with $F=(e_{+})$. First, observe that $F$ is
clearly linearly independent of every set of $q$ standard basis vectors,
since $q\leq k-1<n$. Second, condition (ii) in that theorem is satisfied
since the rank of $(F,E_{n,0}(\gamma ^{\ast }))$, $\gamma ^{\ast }\in (0,\pi
)$, coincides with the rank of the matrix $V_{n}^{(0)}\left( \underline{%
\omega },\underline{d}\right) $, where $\underline{\omega }=(0,\gamma ^{\ast
})$ and $\underline{d}=(1,1)$, and since $V_{n}^{(0)}\left( \underline{%
\omega },\underline{d}\right) $ has rank $\min (n,3)=3$ by Lemma \ref%
{fullrank}. Recall that $n\geq 3$ holds since $n>k\geq 2$. Also note that $%
\rho _{F}(\gamma )=1$ for $\gamma =0$ and $\rho _{F}(\gamma )=0$ for $\gamma
\in (0,\pi ]$, again by Lemma \ref{fullrank}. $\blacksquare $

\textbf{Proof of claim in Example \ref{exa3}:} To prove this claim we use
Theorem \ref{genericity_cor} with $F=(e_{+},e_{-})$. First, observe that $F$
is linearly independent of every set of $q$ standard basis vectors provided $%
q\leq (n/2)-1$ in view of of Part 2 of Lemma \ref%
{independence_F_standard_basis} given below. Second, condition (ii) in
Theorem \ref{genericity_cor} is satisfied since the rank of $%
(F,E_{n,0}(\gamma ^{\ast }))$, $\gamma ^{\ast }\in (0,\pi )$, coincides with
the rank of the matrix $V_{n}^{(0)}\left( \underline{\omega },\underline{d}%
\right) $, where $\underline{\omega }=(0,\pi ,\gamma ^{\ast })$ and $%
\underline{d}=(1,1,1)$, and since $V_{n}^{(0)}\left( \underline{\omega },%
\underline{d}\right) $ has rank $\min (n,4)=4$ by Lemma \ref{fullrank}.
Recall that $n\geq 4$ holds since $n>k\geq 3$. Also note that $\rho
_{F}(\gamma )=1$ for $\gamma =0,\pi $ and $\rho _{F}(\gamma )=0$ for $\gamma
\in (0,\pi )$, again by Lemma \ref{fullrank}. $\blacksquare $

\textbf{Proof of claim in Example \ref{exa4}:} To prove this claim we use
Theorem \ref{genericity_cor} with $F=(e_{+},v)$. First, observe that $F$ is
linearly independent of every set of $q$ standard basis vectors in view of
Part 1 of Lemma \ref{independence_F_standard_basis} below, since $q\leq
k-2<n-2$ and since every nonzero element of $\limfunc{span}(F)$ can have at
most one zero coordinate. Second, condition (ii) in Theorem \ref%
{genericity_cor} is satisfied since the rank of $(F,E_{n,0}(\gamma ^{\ast
})) $, $\gamma ^{\ast }\in (0,\pi )$, coincides with the rank of the matrix $%
V_{n}^{(0)}\left( \underline{\omega },\underline{d}\right) $, where $%
\underline{\omega }=(0,\gamma ^{\ast })$ and $\underline{d}=(2,1)$, and
since $V_{n}^{(0)}\left( \underline{\omega },\underline{d}\right) $ has rank 
$\min (n,4)=4$ by Lemma \ref{fullrank}. Recall that $n\geq 4$ holds since $%
n>k\geq 3$. Also note that $\rho _{F}(\gamma )=2$ for $\gamma =0$ and $\rho
_{F}(\gamma )=0$ for $\gamma \in (0,\pi ]$, again by Lemma \ref{fullrank}. $%
\blacksquare $

\textbf{Proof of claim in Example \ref{exa5}:} To prove this claim we use
Theorem \ref{genericity_cor} with $F=(e_{+},E_{n,0}(\gamma _{0}))$. First,
observe that $F$ is linearly independent of every set of $q$ standard basis
vectors provided $q\leq (n/3)-1$ in view of Part 2 of Lemma \ref%
{independence_F_standard_basis} below. Second, condition (ii) in Theorem \ref%
{genericity_cor} is satisfied since the rank of $(F,E_{n,0}(\gamma ^{\ast
})) $, for $\gamma ^{\ast }\in (0,\pi )$ with $\gamma ^{\ast }\neq \gamma
_{0}$, coincides with the rank of the matrix $V_{n}^{(0)}\left( \underline{%
\omega },\underline{d}\right) $ defined in Appendix \ref{App C}, where $%
\underline{\omega }=(0,\gamma _{0},\gamma ^{\ast })$ and $\underline{d}%
=(1,1,1)$, and since $V_{n}^{(0)}\left( \underline{\omega },\underline{d}%
\right) $ has rank $\min (n,5)=5$ by Lemma \ref{fullrank}. Recall that $%
n\geq 5$ holds since $n>k\geq 4$. Also note that $\rho _{F}(\gamma )=1$ for $%
\gamma =0,\gamma _{0}$ and $\rho _{F}(\gamma )=0$ for $\gamma \in (0,\pi
]\backslash \{\gamma _{0}\} $, again by Lemma \ref{fullrank}. $\blacksquare $

\begin{lemma}
\label{independence_F_standard_basis}Let $F$ be an $n\times k_{F}$%
-dimensional matrix with $\limfunc{rank}(F)=k_{F}\geq 1$. Let $s\in \mathbb{N%
}$ and $s<n$.

1. If the maximum of the number of zero-coordinates of nonzero elements of $%
\limfunc{span}(F)$ is smaller than $n-s$, then the columns of $F$ and $%
e_{i_{1}}(n),\ldots ,e_{i_{s}}(n)$ are linearly independent for every choice
of $1\leq i_{1}<\ldots <i_{s}\leq n$.

2. Suppose $F$ is such that $\limfunc{span}(F)=\limfunc{span}(V_{n}^{(0)}(%
\underline{\omega },\underline{d}))$ holds for some $\underline{\omega }\in
\lbrack 0,\pi ]^{p}$ with distinct coordinates, for some $\underline{d}\in 
\mathbb{N}^{p}$, and for some positive integer $p$. If $k_{F}\leq n/(s+1)$,
then the columns of $F$ and $e_{i_{1}}(n),\ldots ,e_{i_{s}}(n)$ are linearly
independent for every choice of $1\leq i_{1}<\ldots <i_{s}\leq n$.
\end{lemma}

\textbf{Proof:} The first claim is trivial. For the second claim we argue by
contradiction: Suppose there exist indices $1\leq i_{1}<\ldots <i_{s}\leq n$
such that $F$ and $e_{i_{1}}(n),\ldots ,e_{i_{s}}(n)$ are linearly
dependent. Then we can find an $s$-dimensional vector $b\neq 0$ so that $%
v:=\sum_{j=1}^{s}b_{j}e_{i_{j}}(n)\in \limfunc{span}(V_{n}^{(0)}(\underline{%
\omega },\underline{d}))$ holds. The finite sequence $v_{1},\ldots ,v_{n}$
of components of $v$ must then contain a string of consecutive zeros of
length at least $\lfloor n/(s+1)\rfloor $. Now, $v$ is obviously nonzero and
thus $v$ must have a nonzero coordinate $v_{i^{\ast }}$, say, that is
preceded or succeeded by at least $\lfloor n/(s+1)\rfloor $ consecutive
zeros (note that $\lfloor n/(s+1)\rfloor \geq 1$ since $s<n$). Observe that $%
\kappa (\underline{\omega },\underline{d})=k_{F}$ must hold, since $k_{F}=%
\limfunc{rank}(F)=\limfunc{rank}(V_{n}^{(0)}(\underline{\omega },\underline{d%
}))=\min (n,\kappa (\underline{\omega },\underline{d}))$ (by Lemma \ref%
{fullrank} in Appendix \ref{App C}) and since $k_{F}\leq n/(s+1)<n$. It then
follows from Lemma \ref{basis} in Appendix \ref{App C} that $D_{n}(\Delta _{%
\underline{\omega },\underline{d}})v=0$, where $D_{n}(\Delta _{\underline{%
\omega },\underline{d}})$ is a $(n-k_{F})\times n$ matrix as defined in
Appendix \ref{App C} and where $\Delta _{\underline{\omega },\underline{d}}$
is defined in Section \ref{struct}. Note that the coefficient of the highest
power occurring in $\Delta _{\underline{\omega },\underline{d}}$ is always $%
\pm 1$, that the constant term is $1$, and that the degree of $\Delta _{%
\underline{\omega },\underline{d}}$ is $\kappa (\underline{\omega },%
\underline{d})=k_{F}$. Inspection of the equation system $D_{n}(\Delta _{%
\underline{\omega },\underline{d}})v=0$ and noting that $k_{F}\leq n/(s+1)$
is equivalent to $k_{F}\leq \lfloor n/(s+1)\rfloor $, now reveals that one
of the equations is of the form%
\begin{equation*}
\pm v_{i^{\ast }}+c_{1}v_{i^{\ast }-1}+\ldots +c_{k_{F}-1}v_{i^{\ast
}-(k_{F}-1)}+v_{i^{\ast }-k_{F}}=0
\end{equation*}%
where $v_{i^{\ast }-1}=\ldots =v_{i^{\ast }-k_{F}}=0$, or of the form%
\begin{equation*}
\pm v_{i^{\ast }+k_{F}}+c_{1}v_{i^{\ast }+k_{F}-1}+\ldots
+c_{k_{F}-1}v_{i^{\ast }+1}+v_{i^{\ast }}=0
\end{equation*}%
where $v_{i^{\ast }+1}=\ldots =v_{i^{\ast }+k_{F}}=0$. But this implies $%
v_{i^{\ast }}=0$, a contradiction. $\blacksquare $

\begin{lemma}
\label{generic_lem_1}Let $F$ be an $n\times k_{F}$ matrix with $\func{rank}%
(F)=k_{F}$ where $0\leq k_{F}<k$ (with the convention that $F$ is the empty
matrix in case $k_{F}=0$, that the rank of the empty matrix is zero, and
that its span is $\left\{ 0\right\} $). Define $\rho _{F}(\gamma )=\rho
(\gamma ,\func{span}(F))$. Let%
\begin{equation}
\mathfrak{\tilde{X}}_{1}=\left\{ \tilde{X}\in \mathbb{R}^{n\times (k-k_{F})}:%
\func{rank}(F,\tilde{X})=k,\func{rank}(F,\tilde{X},E_{n,\rho _{F}(\gamma
)}(\gamma ))>k\text{ \ for all\ }\gamma \in \lbrack 0,\pi ]\right\} .
\label{X_1_def}
\end{equation}

\begin{enumerate}
\item If the $q\times k$ restriction matrix $R$ of rank $q$ is of the form $%
(0,\tilde{R})$ where $\tilde{R}$ is $q\times (k-k_{F})$, then for every $%
X=(F,\tilde{X})$ with $\tilde{X}\in \mathfrak{\tilde{X}}_{1}$ we have $\rho
(\gamma ,\mathfrak{M}_{0}^{lin})=\rho _{F}(\gamma )$ for every\ $\gamma \in
\lbrack 0,\pi ]$. [Note that $\rho (\gamma ,\mathfrak{M}_{0}^{lin})$ depends
on $X$ and $R$, but this is not expressed in the notation.]

\item The set $\mathbb{R}^{n\times (k-k_{F})}\backslash \mathfrak{\tilde{X}}%
_{1}$ is contained in a $\lambda _{\mathbb{R}^{n\times (k-k_{F})}}$-null set
if $n>k+2$ holds.

\item The set $\mathbb{R}^{n\times (k-k_{F})}\backslash \mathfrak{\tilde{X}}%
_{1}$ is contained in a $\lambda _{\mathbb{R}^{n\times (k-k_{F})}}$-null set
if $\func{rank}(F,E_{n,0}(\gamma ^{\ast }))=k_{F}+2$ holds for some $\gamma
^{\ast }\in (0,\pi )$.
\end{enumerate}
\end{lemma}

\textbf{Proof: }1. For every $X=(F,\tilde{X})$ with $\tilde{X}\in \mathbb{R}%
^{n\times (k-k_{F})}$ and for $R$ of the form assumed in the lemma we
immediately see that the associated space $\mathfrak{M}_{0}^{lin}$
(depending on $X$ and $R$, but the dependence not being shown in the
notation) contains $\func{span}(F)$. Consequently, $\rho (\gamma ,\mathfrak{M%
}_{0}^{lin})\geq \rho _{F}(\gamma )$ must hold for every\ $\gamma \in
\lbrack 0,\pi ]$. To prove the converse, note that for $\tilde{X}\in 
\mathfrak{\tilde{X}}_{1}$ the second rank condition in (\ref{X_1_def})
implies that $E_{n,\rho _{F}(\gamma )}(\gamma )$ is not contained in $\func{%
span}(X)=\func{span}((F,\tilde{X}))$, and hence a fortiori not in $\mathfrak{%
M}_{0}^{lin}$. This immediately implies $\rho (\gamma ,\mathfrak{M}%
_{0}^{lin})\leq \rho _{F}(\gamma )$ for every\ $\gamma \in \lbrack 0,\pi ]$.

2. Since the set of all $\tilde{X}$ such that $\func{rank}(F,\tilde{X})<k$
is obviously a $\lambda _{\mathbb{R}^{n\times (k-k_{F})}}$-null set, it
suffices to show that the set $A=\bigcup_{\gamma \in \lbrack 0,\pi
]}A_{\gamma }$ where%
\begin{equation*}
A_{\gamma }=\left\{ \tilde{X}\in \mathbb{R}^{n\times (k-k_{F})}:\func{rank}%
(F,\tilde{X})=k,\func{rank}(F,\tilde{X},E_{n,\rho _{F}(\gamma )}(\gamma
))=k\right\}
\end{equation*}%
is contained in a $\lambda _{\mathbb{R}^{n\times (k-k_{F})}}$-null set. We
first show that $A_{\gamma }$ is a $\lambda _{\mathbb{R}^{n\times
(k-k_{F})}} $-null set for each fixed $\gamma $: By definition of $\rho
_{F}(\gamma )$ at least one of the columns of $E_{n,\rho _{F}(\gamma
)}(\gamma )$ does not belong to $\func{span}(F)$. Choose one such column and
denote it by $h$. Then $A_{\gamma }$ is contained in the set 
\begin{equation*}
B_{\gamma }=\left\{ \tilde{X}\in \mathbb{R}^{n\times (k-k_{F})}:\func{rank}%
(F,\tilde{X},h)\leq k\right\} .
\end{equation*}%
Since $(F,\tilde{X},h)$ has dimension $n\times (k+1)$ the set $B_{\gamma }$
is given by the zero-set of the polynomial $\det ((F,\tilde{X},h)^{\prime
}(F,\tilde{X},h))$. We next construct a matrix $\tilde{X}\in \mathbb{R}%
^{n\times (k-k_{F})}$ which does not belong to this set. Observe that $F$
and $h$ together span a linear space of dimension $k_{F}+1$ and that $%
k_{F}+1\leq k<n$ holds since we have assumed $k_{F}<k$ and since we always
maintain $k<n$. Hence we can find $k-k_{F}$ linearly independent vectors in $%
(\func{span}((F,h)))^{\bot },$ which we use as the columns of $\tilde{X}$.
Clearly, this $\tilde{X}$ does not belong to $B_{\gamma }$. Consequently, $%
B_{\gamma }$ is a $\lambda _{\mathbb{R}^{n\times (k-k_{F})}}$-null set, and
a fortiori the same is true for $A_{\gamma }$ (Borel-measurability of $%
A_{\gamma }$ being trivial). Let now $U$ be the finite set $\left\{ \gamma
\in \lbrack 0,\pi ]:\rho _{F}(\gamma )>0\right\} \cup \left\{ 0,\pi \right\} 
$ (cf. the discussion following Definition \ref{subspaces}). Then $%
\bigcup_{\gamma \in U}A_{\gamma }$ is a $\lambda _{\mathbb{R}^{n\times
(k-k_{F})}}$-null set. It remains to show that $\bigcup_{\gamma \in \lbrack
0,\pi ]\backslash U}A_{\gamma }$ is contained in a $\lambda _{\mathbb{R}%
^{n\times (k-k_{F})}}$-null set. Note that $[0,\pi ]\backslash U\subseteq
(0,\pi )$ is an open set and that $\rho _{F}(\gamma )=0$ holds for every $%
\gamma \in \lbrack 0,\pi ]\backslash U$. Hence, $\func{span}(E_{n,0}(\gamma
))$ is not contained in $\func{span}(F)$ for every $\gamma \in \lbrack 0,\pi
]\backslash U$. Thus, if $\tilde{X}\in A_{\gamma }$ with $\gamma \in \lbrack
0,\pi ]\backslash U$ we can then find an index $i(\gamma )$, $1\leq i(\gamma
)\leq k-k_{F}$ and a $(k+1)\times 1$ vector $v(\gamma )$, such that%
\begin{equation*}
\tilde{X}_{.i(\gamma )}=(F,\tilde{X}(\lnot i(\gamma )),E_{n,0}(\gamma
))v(\gamma )
\end{equation*}%
holds, where $\tilde{X}_{.i(\gamma )}$ denotes the $i(\gamma )$-th column of 
$\tilde{X}$, and $\tilde{X}(\lnot i(\gamma ))$ denotes the matrix $\tilde{X}$
after the $i(\gamma )$-th column has been deleted. In other words, $\tilde{X}
$ coincides with 
\begin{equation*}
((F,\tilde{X}(\lnot i(\gamma )),E_{n,0}(\gamma ))v(\gamma ),\tilde{X}(\lnot
i(\gamma ))),
\end{equation*}%
up to a permutation of columns. [Note that $A_{\gamma }$ may be empty,
namely if $\func{rank}(F,E_{n,0}(\gamma ))=k_{F}+2$ and $k-k_{F}=1$.] Now,
for every $(k-k_{F})\times (k-k_{F})$ permutation matrix $Per$ define the
map $\Xi _{Per}:\mathbb{R}^{n\times (k-k_{F}-1)}\times \mathbb{R}%
^{k+1}\times \lbrack 0,\pi ]\backslash U\rightarrow \mathbb{R}^{n\times
(k-k_{F})}$ via $(\bar{X},\bar{v},\gamma )\mapsto ((F,\bar{X},E_{n,0}(\gamma
))\bar{v},\bar{X})Per$. [In case $k-k_{F}=1$ the symbol $\mathbb{R}^{n\times
(k-k_{F}-1)}$ is to be interpreted as $\left\{ 0\right\} $ and hence the map 
$\Xi _{Per}$ is effectively defined on $\mathbb{R}^{k+1}\times \lbrack 0,\pi
]\backslash U$.] Because of what has been said before, we now see that $%
\bigcup_{\gamma \in \lbrack 0,\pi ]\backslash U}A_{\gamma }$ is contained in
the union of the images of all the maps $\Xi _{Per}$ when $Per$ varies in
the indicated set of permutation matrices. Note that this union is a finite
union. It hence suffices to show that the image of each $\Xi _{Per}$ is
contained in a $\lambda _{\mathbb{R}^{n\times (k-k_{F})}}$-null set.
Clearly, the domain of definition of each $\Xi _{Per}$ is an open set in
Euclidean space of dimension $n(k-k_{F}-1)+k+2$ and each $\Xi _{Per}$ is a
smooth map. Sard's Theorem (see, e.g., \cite{milnor1997}) now implies the
desired conclusion provided $n(k-k_{F}-1)+k+2$ is smaller than the dimension 
$n(k-k_{F})$ of the range space. But this is guaranteed by the assumption
that $n>k+2$.

3. We proceed as in 2.~up to the point where the set $U$ has been defined.
Now define $U^{\ast }$ as the union of $U$ and the set $\left\{ \gamma \in
(0,\pi ):\func{rank}(F,E_{n,0}(\gamma ))<k_{F}+2\right\} $. The latter set
is clearly contained in the zero-set of the function $\gamma \mapsto \det
((F,E_{n,0}(\gamma ))^{\prime }(F,E_{n,0}(\gamma )))$. Obviously this
function can be expressed as a rational function in $\exp (\iota \gamma )$
and thus only has finitely many zeros, except if it is identically zero. But
the latter cannot happen because of the assumption made for Part 3. We have
now established that $U^{\ast }$ is a finite set. It follows that $%
\bigcup_{\gamma \in U^{\ast }}A_{\gamma }$ is a $\lambda _{\mathbb{R}%
^{n\times (k-k_{F})}}$-null set. It remains to show that $\bigcup_{\gamma
\in \lbrack 0,\pi ]\backslash U^{\ast }}A_{\gamma }$ is contained in a $%
\lambda _{\mathbb{R}^{n\times (k-k_{F})}}$-null set. Note that $[0,\pi
]\backslash U^{\ast }\subseteq (0,\pi )$ is an open set and that $\rho
_{F}(\gamma )=0$ holds for every $\gamma \in \lbrack 0,\pi ]\backslash
U^{\ast }$. In case $k-k_{F}=1$, it follows that $\bigcup_{\gamma \in
\lbrack 0,\pi ]\backslash U^{\ast }}A_{\gamma }$ is the empty set (since $%
\func{rank}(F,E_{n,0}(\gamma ))=k_{F}+2$ for $\gamma \in \lbrack 0,\pi
]\backslash U^{\ast }$) and we are done. Hence assume $k-k_{F}\geq 2$. Then,
if $\tilde{X}\in A_{\gamma }$ with $\gamma \in \lbrack 0,\pi ]\backslash
U^{\ast }$ we can find indices $i_{1}(\gamma )$, $i_{2}(\gamma )$ with $%
1\leq i_{1}(\gamma )<i_{2}(\gamma )\leq k-k_{F}$, and two $k\times 1$
vectors $v_{1}(\gamma )$ and $v_{2}(\gamma )$ such that%
\begin{equation*}
(\tilde{X}_{.i_{1}(\gamma )},\tilde{X}_{.i_{2}(\gamma )})=(F,\tilde{X}(\lnot
i_{1}(\gamma ),\lnot i_{1}(\gamma )),E_{n,0}(\gamma ))(v_{1}(\gamma
),v_{2}(\gamma ))
\end{equation*}%
holds, where $\tilde{X}(\lnot i_{1}(\gamma ),\lnot i_{2}(\gamma ))$ denotes
the matrix $\tilde{X}$ after the columns $i_{1}(\gamma )$ and $i_{2}(\gamma
) $ have been deleted. In other words, $\tilde{X}$ coincides with 
\begin{equation*}
((F,\tilde{X}(\lnot i_{1}(\gamma ),\lnot i_{1}(\gamma )),E_{n,0}(\gamma
))(v_{1}(\gamma ),v_{2}(\gamma )),\tilde{X}(\lnot i_{1}(\gamma ),\lnot
i_{1}(\gamma ))),
\end{equation*}%
up to a permutation of columns. Now, for every $(k-k_{F})\times (k-k_{F})$
permutation matrix $Per$ define the map $\Xi _{Per}:\mathbb{R}^{n\times
(k-k_{F}-2)}\times \mathbb{R}^{k}\times \mathbb{R}^{k}\times \lbrack 0,\pi
]\backslash U^{\ast }\rightarrow \mathbb{R}^{n\times (k-k_{F})}$ via $(\bar{X%
},\bar{v}_{1},\bar{v}_{2},\gamma )\mapsto ((F,\bar{X},E_{n,0}(\gamma ))(\bar{%
v}_{1},\bar{v}_{2}),\bar{X})Per$. [In case $k-k_{F}=2$ the symbol $\mathbb{R}%
^{n\times (k-k_{F}-2)}$ is to be interpreted as $\left\{ 0\right\} $ and
hence the map $\Xi _{Per}$ is effectively defined on $\mathbb{R}^{k}\times 
\mathbb{R}^{k}\times \lbrack 0,\pi ]\backslash U$.] Because of what has been
said before, we now see that $\bigcup_{\gamma \in \lbrack 0,\pi ]\backslash
U^{\ast }}A_{\gamma }$ is contained in the union of the images of all the
maps $\Xi _{Per}$ when $Per$ varies in the indicated set of permutation
matrices. This union is again a finite union. It hence suffices to show that
the image of each $\Xi _{Per}$ is contained in a $\lambda _{\mathbb{R}%
^{n\times (k-k_{F})}}$-null set. Clearly, the domain of definition of each $%
\Xi _{Per}$\ is an open set in Euclidean space of dimension $%
n(k-k_{F}-2)+2k+1$ and each $\Xi _{Per}$ is a smooth map. Sard's Theorem
again now implies the desired conclusion provided $n(k-k_{F}-2)+2k+1$ is
smaller than the dimension $n(k-k_{F})$ of the range space. But this is
guaranteed by our general assumption that $n>k$. $\blacksquare $

\begin{lemma}
\label{generic_lem_2}Let $F$ be an $n\times k_{F}$ matrix with $\func{rank}%
(F)=k_{F}$ where $0\leq k_{F}<k$ (with the same convention as before if $%
k_{F}=0$). Assume that the $q\times k$ restriction matrix $R$ of rank $q$ is
of the form $(0,\tilde{R})$ where $\tilde{R}$ is $q\times (k-k_{F})$.
Furthermore, assume that the columns of $F$ and $e_{i_{1}}(n),\ldots
,e_{i_{q}}(n)$ are linearly independent for every choice of $1\leq
i_{1}<\ldots <i_{q}\leq n$. Then the complement of the set 
\begin{equation*}
\mathfrak{\tilde{X}}_{2}=\left\{ \tilde{X}\in \mathbb{R}^{n\times (k-k_{F})}:%
\func{rank}(X)=k,\text{ }\mathsf{B}=\func{span}(X)\right\}
\end{equation*}%
in $\mathbb{R}^{n\times (k-k_{F})}$ is contained in a $\lambda _{\mathbb{R}%
^{n\times (k-k_{F})}}$-null set, where $X=(F,\tilde{X})$. [Recall that $%
\mathsf{B}$ depends on $X$, which, however, is not shown in the notation.]
\end{lemma}

\textbf{Proof: Step 1:} Let $1\leq i_{1}<\ldots <i_{q}\leq n$ be given and
define the set%
\begin{equation*}
\mathfrak{\tilde{X}}(i_{1},\ldots ,i_{q})=\left\{ \tilde{X}\in \mathbb{R}%
^{n\times (k-k_{F})}:\func{rank}(X)=k,\limfunc{rank}(R(X^{\prime
}X)^{-1}X^{\prime }\left( e_{i_{1}}(n),\ldots ,e_{i_{q}}(n)\right)
)=q\right\} ,
\end{equation*}%
where $X=(F,\tilde{X})$. First, we show that $\mathfrak{\tilde{X}}%
(i_{1},\ldots ,i_{q})$ is nonempty: Since $\func{span}((F,e_{i_{1}}(n),%
\ldots ,e_{i_{q}}(n)))$ has dimension $k_{F}+q$ by the assumptions on $F$
and since $k_{F}+q\leq k<n$ in view of the assumptions on $R$, we can find
orthonormal $n\times 1$ vectors $a_{1},\ldots ,a_{k-(k_{F}+q)}$ in the
orthogonal complement of $\func{span}((F,e_{i_{1}}(n),\ldots ,e_{i_{q}}(n)))$%
. Define $\tilde{X}^{\ast }=(a_{1},\ldots
,a_{k-(k_{F}+q)},e_{i_{1}}(n),\ldots ,e_{i_{q}}(n))$, with the convention
that $\tilde{X}^{\ast }=(e_{i_{1}}(n),\ldots ,e_{i_{q}}(n))$ in case $%
k_{F}+q=k$, and set $X^{\ast }=(F,\tilde{X}^{\ast }$). Obviously $X^{\ast }$
has rank equal to $k$. Let $\tilde{V}$ be a nonsingular $(k-k_{F})\times
(k-k_{F})$ matrix such that $\tilde{R}\tilde{V}^{-1}=(0,I_{q})$ and define $%
V $ as the $k\times k$ block-diagonal matrix with first diagonal block $%
I_{k_{F}}$ and second diagonal block $\tilde{V}$. Clearly, $%
RV^{-1}=(0,I_{q}) $ holds. Set $X=X^{\ast }V$ and note that $X=(F,\tilde{X})$
with $\tilde{X}=\tilde{X}^{\ast }\tilde{V}$ and that $\func{rank}(X)=k$.
Furthermore, we have%
\begin{equation*}
R(X^{\prime }X)^{-1}X^{\prime }\left( e_{i_{1}}(n),\ldots
,e_{i_{q}}(n)\right) =(0,I_{q})(X^{\ast \prime }X^{\ast })^{-1}X^{\ast
\prime }\left( e_{i_{1}}(n),\ldots ,e_{i_{q}}(n)\right) =I_{q},
\end{equation*}%
showing that the so-constructed $\tilde{X}$ belongs to $\mathfrak{\tilde{X}}%
(i_{1},\ldots ,i_{q})$. Second, observe that $\tilde{X}\notin \mathfrak{%
\tilde{X}}(i_{1},\ldots ,i_{q})$ is equivalent to 
\begin{equation*}
\det (X^{\prime }X)=0\text{ or }(\det (X^{\prime }X)\neq 0\text{ and }\det
(R(X^{\prime }X)^{-1}X^{\prime }(e_{i_{1}}(n),\ldots ,e_{i_{q}}(n)))=0),
\end{equation*}%
which in turn is equivalent to%
\begin{equation*}
\det (X^{\prime }X)\det (R\limfunc{adj}(X^{\prime }X)X^{\prime
}(e_{i_{1}}(n),\ldots ,e_{i_{q}}(n)))=0\text{.}
\end{equation*}%
This is a polynomial in the entries of $\tilde{X}$ and does not vanish
identically, because we have shown that $\mathfrak{\tilde{X}}(i_{1},\ldots
,i_{q})$ is nonempty. Consequently, the complement of $\mathfrak{\tilde{X}}%
(i_{1},\ldots ,i_{q})$ is a $\lambda _{\mathbb{R}^{n\times (k-k_{F})}}$-null
set.

\textbf{Step 2:} It follows that the set $\mathfrak{\tilde{X}}_{20}$ defined
as the intersection of all sets of the form $\mathfrak{\tilde{X}}%
(i_{1},\ldots ,i_{q})$, where we vary over all possible combinations of
indices satisfying $1\leq i_{1}<\ldots <i_{q}\leq n$, is the complement of a 
$\lambda _{\mathbb{R}^{n\times (k-k_{F})}}$-null set, since this is an
intersection of finitely many sets.

\textbf{Step 3:} Let $1\leq j_{1}<\ldots <j_{n-q+1}\leq n$ be given and
define the set%
\begin{equation*}
\mathfrak{\tilde{X}}^{\ast }(j_{1},\ldots ,j_{n-q+1})=\left\{ \tilde{X}\in 
\mathbb{R}^{n\times (k-k_{F})}:z_{j_{1}}=\ldots =z_{j_{n-q+1}}=0\text{ \ for
some }z\in (\func{span}(X))^{\bot }\backslash \{0\}\right\} .
\end{equation*}%
We show that this set is a $\lambda _{\mathbb{R}^{n\times (k-k_{F})}}$-null
set: If $q=1$, then $\mathfrak{\tilde{X}}^{\ast }(j_{1},\ldots ,j_{n-q+1})$
is obviously empty. Hence consider the case $q>1$. Observe that $\tilde{X}%
\in \mathfrak{\tilde{X}}^{\ast }(j_{1},\ldots ,j_{n-q+1})$ is equivalent to
the equation system $A(\tilde{X})z=0$ having a nonzero solution, where%
\begin{equation*}
A(\tilde{X})=\left( e_{j_{1}}(n),\ldots ,e_{j_{n-q+1}}(n),X\right) ^{\prime }%
\text{.}
\end{equation*}%
Observe that $A(\tilde{X})$ is of dimension $(n-q+1+k)\times n$ and that the
row-dimension is larger than $n$, since $q\leq k$. Consequently, $\tilde{X}%
\in \mathfrak{\tilde{X}}^{\ast }(j_{1},\ldots ,j_{n-q+1})$ is equivalent to 
\begin{equation*}
\det (A(\tilde{X})^{\prime }A(\tilde{X}))=0,
\end{equation*}%
a polynomial equation in the elements of $\tilde{X}$. The solution set is
thus a $\lambda _{\mathbb{R}^{n\times (k-k_{F})}}$-null set if we can
exhibit an element $\tilde{X}\notin \mathfrak{\tilde{X}}^{\ast
}(j_{1},\ldots ,j_{n-q+1})$. We now construct such an element as follows:
Let $i_{1},\ldots ,i_{q-1}$ be the elements of $\left\{ 1,\ldots ,n\right\} $
not appearing in the list $j_{1},\ldots ,j_{n-q+1}$. Choose linearly
independent elements $a_{1},\ldots ,a_{k-(k_{F}+q-1)}$ of the orthogonal
complement of $\func{span}((F,e_{i_{1}}(n),\ldots ,e_{i_{q-1}}(n)))$; note
that such a choice is possible since the dimension of this span is $%
k_{F}+q-1<k<n$ in view of our assumptions. Now define%
\begin{equation*}
\tilde{X}=(e_{i_{1}}(n),\ldots ,e_{i_{q-1}}(n),a_{1},\ldots
,a_{k-(k_{F}+q-1)})
\end{equation*}%
and observe that clearly $\tilde{X}\notin \mathfrak{\tilde{X}}^{\ast
}(j_{1},\ldots ,j_{n-q+1})$ is satisfied.

\textbf{Step 4:} Define $\mathfrak{\tilde{X}}_{21}$ as the complement (in $%
\mathbb{R}^{n\times (k-k_{F})}$) of the set 
\begin{equation*}
\left( \bigcup \mathfrak{\tilde{X}}^{\ast }(j_{1},\ldots ,j_{n-q+1})\right)
\cup \left\{ \tilde{X}\in \mathbb{R}^{n\times (k-k_{F})}:\func{rank}%
(X)<k\right\} ,
\end{equation*}%
where the union extends over all tuples $(j_{1},\ldots ,j_{n-q+1})$
satisfying $1\leq j_{1}<\ldots <j_{n-q+1}\leq n$. Then $\mathfrak{\tilde{X}}%
_{21}$ is the complement of a $\lambda _{\mathbb{R}^{n\times (k-k_{F})}}$%
-null set: This is obvious because of Step 3, the fact that the union is a
union of finitely many sets, and because the set $\left\{ \tilde{X}\in 
\mathbb{R}^{n\times (k-k_{F})}:\func{rank}(X)<k\right\} $ is clearly a $%
\lambda _{\mathbb{R}^{n\times (k-k_{F})}}$-null set.

\textbf{Step 5:} It follows from the preceding steps that $\mathfrak{\tilde{X%
}}_{20}\cap \mathfrak{\tilde{X}}_{21}$ is the complement of a $\lambda _{%
\mathbb{R}^{n\times (k-k_{F})}}$-null set. We now show that $\mathfrak{%
\tilde{X}}_{20}\cap \mathfrak{\tilde{X}}_{21}$ is contained in $\mathfrak{%
\tilde{X}}_{2}$: First, if $\tilde{X}\in \mathfrak{\tilde{X}}_{20}\cap 
\mathfrak{\tilde{X}}_{21}$, then $X=(F,\tilde{X})$ clearly has rank equal to 
$k$. Second, we need to show that $\tilde{X}\in \mathfrak{\tilde{X}}%
_{20}\cap \mathfrak{\tilde{X}}_{21}$ implies $\mathsf{B}\subseteq \func{span}%
(X)$, the other inclusion being trivial: Now, suppose that $\tilde{X}\in 
\mathfrak{\tilde{X}}_{20}\cap \mathfrak{\tilde{X}}_{21}$, and that $y\in 
\mathsf{B}$ but $y\notin \func{span}(X)$ would be possible. Then $z=\Pi _{(%
\func{span}(X))^{\bot }}y\neq 0$ would follow. Hence, there would have to
exist indices $1\leq i_{1}<\ldots <i_{q}\leq n$ such that 
\begin{equation}
e_{i_{1}}^{\prime }(n)z\neq 0,\ldots ,e_{i_{q}}^{\prime }(n)z\neq 0,
\label{nonzero}
\end{equation}%
because of $\tilde{X}\in \mathfrak{\tilde{X}}_{21}$. But then we could
conclude that%
\begin{eqnarray*}
q &\geq &\func{rank}(B(y))\geq \func{rank}\left( R(X^{\prime
}X)^{-1}X^{\prime }\left( e_{i_{1}}(n)e_{i_{1}}^{\prime }(n)z,\ldots
,e_{i_{q}}(n)e_{i_{q}}^{\prime }(n)z\right) \right) \\
&=&\func{rank}\left( R(X^{\prime }X)^{-1}X^{\prime }\left(
e_{i_{1}}(n),\ldots ,e_{i_{q}}(n)\right) \right) =q,
\end{eqnarray*}%
the last but one equality holding because of (\ref{nonzero}), and the last
one holding in view of $\tilde{X}\in \mathfrak{\tilde{X}}_{20}$. This would
entail $y\notin \mathsf{B}$, a contradiction. This now completes the proof. $%
\blacksquare $

\section{Appendix\label{App B}: Proofs for Section \protect\ref{general} and
some sufficient conditions relating to $\mathbb{J}(\mathcal{L},\mathfrak{C})$%
}

\textbf{Proof of Proposition \ref{L1}:} Fix an arbitrary $\mu _{0}\in 
\mathfrak{M}_{0}$. Due to $G(\mathfrak{M}_{0})$-invariance of $W$ the l.h.s.
of (\ref{size}) coincides with $\sup_{\Sigma \in \mathfrak{C}}P_{\mu
_{0},\Sigma }(W)$, cf. Remark 5.5(iii) in \cite{PP2016}. Now, let $\Sigma
_{j}$ be a sequence in $\mathfrak{C}$ such that $P_{\mu _{0},\Sigma
_{j}}(W)\rightarrow \sup_{\Sigma \in \mathfrak{C}}P_{\mu _{0},\Sigma }(W)$
as $j\rightarrow \infty $. By relative compactness of $\mathcal{L}(\mathfrak{%
C})$, we may assume that $\mathcal{L}(\Sigma _{j})\rightarrow \bar{\Sigma}$
(possibly after passing to a subsequence). For notational convenience,
define the sequence $L_{j}=\mathcal{L}(\Sigma _{j})+\Pi _{\mathcal{L}}$ and
denote its limit by $L=\bar{\Sigma}+\Pi _{\mathcal{L}}$. We first claim that%
\begin{equation}
P_{\mu _{0},\Sigma _{j}}(W)=P_{\mu _{0},\mathcal{L}(\Sigma _{j})}(W)=P_{\mu
_{0},L_{j}}(W).  \label{invcons}
\end{equation}%
To see this, let $\mathbf{Z}$ be a standard Gaussian $n\times 1$ vector.
Then 
\begin{equation*}
P_{\mu _{0},\Sigma _{j}}(W)=\Pr \left( \mu _{0}+\Sigma _{j}^{1/2}\mathbf{Z}%
\in W\right) =\Pr \left( \mu _{0}+\Pi _{\mathcal{L}^{\bot }}\Sigma _{j}^{1/2}%
\mathbf{Z}\in W\right)
\end{equation*}%
because $W$ is easily seen to be invariant w.r.t. addition of elements of $%
\mathcal{L}$ and since $\Pi _{\mathcal{L}}\Sigma _{m}^{1/2}\mathbf{Z}$
clearly belongs to $\mathcal{L}$. The latter probability equals 
\begin{eqnarray*}
&&\Pr \left( \mu _{0}+\Pi _{\mathcal{L}^{\bot }}\Sigma _{j}^{1/2}\mathbf{Z/}%
\left\Vert \Pi _{\mathcal{L}^{\bot }}\Sigma _{j}\Pi _{\mathcal{L}^{\bot
}}\right\Vert ^{1/2}\in \mu _{0}+\left( W-\mu _{0}\right) \mathbf{/}%
\left\Vert \Pi _{\mathcal{L}^{\bot }}\Sigma _{j}\Pi _{\mathcal{L}^{\bot
}}\right\Vert ^{1/2}\right) \\
&=&\Pr \left( \mu _{0}+\Pi _{\mathcal{L}^{\bot }}\Sigma _{j}^{1/2}\mathbf{Z/}%
\left\Vert \Pi _{\mathcal{L}^{\bot }}\Sigma _{j}\Pi _{\mathcal{L}^{\bot
}}\right\Vert ^{1/2}\in W\right) =P_{\mu _{0},\mathcal{L}(\Sigma _{j})}(W)
\end{eqnarray*}%
where the first equality follows from $G(\mathfrak{M}_{0})$-invariance of $W$%
. Furthermore, using the invariance of $W$ w.r.t. to addition of elements of 
$\mathcal{L}$ again we obtain%
\begin{equation*}
P_{\mu _{0},\mathcal{L}(\Sigma _{j})}(W)=\Pr \left( \mu _{0}+\mathcal{L}%
(\Sigma _{j})^{1/2}\mathbf{Z}\in W\right) =\Pr \left( \mu _{0}+\left( 
\mathcal{L}(\Sigma _{j})^{1/2}+\Pi _{\mathcal{L}}\right) \mathbf{Z}\in
W\right) =P_{\mu _{0},L_{j}}(W),
\end{equation*}%
the last equality holding because $\mathcal{L}(\Sigma _{j})^{1/2}+\Pi _{%
\mathcal{L}}$ is a square-root of $L_{j}$. This establishes (\ref{invcons}).

We now distinguish two cases.

\textit{Case 1:} Suppose $L$ is positive definite. Then $P_{\mu _{0},\Sigma
_{j}}(W)\rightarrow P_{\mu _{0},L}(W)$ follows from (\ref{invcons}) in view
of total variation convergence of $P_{\mu _{0},L_{j}}$ to $P_{\mu _{0},L}$.
But note that $P_{\mu _{0},L}(W)<1$, since $W^{c}$ is not a $\lambda _{%
\mathbb{R}^{n}}$-null set by assumption and since $\lambda _{\mathbb{R}^{n}}$
is equivalent to $P_{\mu _{0},L}$.

\textit{Case 2:} Suppose $L$ is singular. We note that 
\begin{equation*}
\limsup_{j\rightarrow \infty }P_{\mu _{0},\Sigma _{j}}(W)\leq
\limsup_{j\rightarrow \infty }P_{\mu _{0},\mathcal{L}(\Sigma _{j})}(\limfunc{%
cl}(W))\leq P_{\mu _{0},\bar{\Sigma}}(\limfunc{cl}(W)),
\end{equation*}%
where we have used Equation (\ref{invcons}), the inclusion $W\subseteq 
\limfunc{cl}(W)$, weak convergence of $P_{\mu _{0},\mathcal{L}(\Sigma _{j})}$
to $P_{\mu _{0},\bar{\Sigma}}$ (cf. Lemma E.1 in \cite{PP2016}), and the
Portmanteau theorem. Define $\mathcal{S}=\mathrm{\limfunc{span}}(\bar{\Sigma}%
)$ and observe that 
\begin{equation*}
P_{\mu _{0},\bar{\Sigma}}(\limfunc{cl}(W))=1-P_{\mu _{0},\bar{\Sigma}}\left(
\left( \limfunc{cl}(W)\right) ^{c}\cap (\mu _{0}+\mathcal{S})\right) ,
\end{equation*}%
since $P_{\mu _{0},\bar{\Sigma}}$ is concentrated on $\mu _{0}+\mathcal{S}$.
Furthermore, singularity of $L$ implies $\limfunc{rank}(\bar{\Sigma})<n-\dim
(\mathcal{L})$, and thus $\mathcal{S}\in \mathbb{J}(\mathcal{L},\mathfrak{C}%
) $. Hence, by assumption, $\left( \limfunc{cl}(W)\right) ^{c}\cap (\mu _{0}+%
\mathcal{S})\neq \emptyset $, and thus there exists a nonempty open subset $%
U\subseteq \mu _{0}+\mathcal{S}$ (w.r.t. the topology induced from $\mathbb{R%
}^{n}$) such that $U\subseteq \left( \limfunc{cl}(W)\right) ^{c}\cap (\mu
_{0}+\mathcal{S})$. Equivalence of the measures $P_{\mu _{0},\bar{\Sigma}}$
and $\lambda _{\mu _{0}+\mathcal{S}}$ then implies $P_{\mu _{0},\bar{\Sigma}%
}(U)>0$. Together with the preceding display this then gives $P_{\mu _{0},%
\bar{\Sigma}}(\limfunc{cl}(W))\leq 1-P_{\mu _{0},\bar{\Sigma}}(U)<1$. $%
\blacksquare $

\textbf{Proof of Proposition \ref{L2}: }Fix an arbitrary $\mu _{0}\in 
\mathfrak{M}_{0}$. Due to $G(\mathfrak{M}_{0})$-invariance of $W_{j}$ the
l.h.s. of (\ref{size_2}) coincides with $\sup_{\Sigma \in \mathfrak{C}%
}P_{\mu _{0},\Sigma }(W_{j})$, cf. Remark 5.5(iii) in \cite{PP2016}. Let $%
\Sigma _{j}$ be a sequence in $\mathfrak{C}$ such that%
\begin{equation}
\left\vert P_{\mu _{0},\Sigma _{j}}(W_{j})-\sup_{\Sigma \in \mathfrak{C}%
}P_{\mu _{0},\Sigma }(W_{j})\right\vert \rightarrow 0\quad \text{ as }\quad
j\rightarrow \infty .  \label{supapprox}
\end{equation}%
By relative compactness of $\mathcal{L}(\mathfrak{C})$, we may assume that $%
\mathcal{L}(\Sigma _{j})\rightarrow \bar{\Sigma}$ (possibly after passing to
a subsequence). Again we define the sequence $L_{j}=\mathcal{L}(\Sigma
_{j})+\Pi _{\mathcal{L}}$ and denote its limit by $L=\bar{\Sigma}+\Pi _{%
\mathcal{L}}$. Using the assumed invariance properties of $W_{j}$ we obtain
as before that%
\begin{equation}
P_{\mu _{0},\Sigma _{j}}(W_{j})=P_{\mu _{0},\mathcal{L}(\Sigma
_{j})}(W_{j})=P_{\mu _{0},L_{j}}(W_{j}).  \label{invconsm}
\end{equation}

\textit{Case 1:} Suppose $L$ is positive definite. Let $\varepsilon >0$ and
note that $P_{\mu _{0},L}$ is equivalent to $\lambda _{\mathbb{R}^{n}}$.
Since $W_{j}\supseteq W_{j+1}$, and since $\bigcap_{j\in \mathbb{N}}W_{j}$
is a $\lambda _{\mathbb{R}^{n}}$-null set by assumption, we can find $%
j_{0}\in \mathbb{N}$ such that $P_{\mu _{0},L}(W_{j_{0}})<\varepsilon $.
Since $P_{\mu _{0},L_{j}}$ converges to $P_{\mu _{0},L}$ in total variation,
we arrive at 
\begin{equation*}
\lim_{j\rightarrow \infty }P_{\mu _{0},L_{j}}(W_{j_{0}})<\varepsilon .
\end{equation*}%
From $W_{j}\supseteq W_{j+1}$ we obtain 
\begin{equation*}
\limsup_{j\rightarrow \infty }P_{\mu _{0},L_{j}}(W_{j})\leq
\limsup_{j\rightarrow \infty }P_{\mu _{0},L_{j}}(W_{j_{0}})<\varepsilon ,
\end{equation*}%
which proves the claim in view of (\ref{invconsm}) and (\ref{supapprox}).

\textit{Case 2:} Suppose $L$ is singular. Let $\varepsilon >0$. Singularity
of $L$ implies $\limfunc{rank}(\bar{\Sigma})<n-\dim (\mathcal{L})$, and
therefore $\mathcal{S}=\limfunc{span}(\bar{\Sigma})\in \mathbb{J}(\mathcal{L}%
,\mathfrak{C})$ must hold. By assumption, $\bigcap_{j\in \mathbb{N}}\limfunc{%
cl}(W_{j})$ is hence a $\lambda _{\mu _{0}+\mathcal{S}}$-null set. Since
furthermore $\limfunc{cl}(W_{j})\supseteq \limfunc{cl}(W_{j+1})$ holds, and
since $P_{\mu _{0},\bar{\Sigma}}$ is equivalent to $\lambda _{\mu _{0}+%
\mathcal{S}}$, we can hence find a $j_{0}\in \mathbb{N}$ such that 
\begin{equation*}
P_{\mu _{0},\bar{\Sigma}}(\limfunc{cl}(W_{j_{0}}))<\varepsilon .
\end{equation*}%
But $\limfunc{cl}(W_{j})\supseteq \limfunc{cl}(W_{j+1})$ implies 
\begin{align*}
\limsup_{j\rightarrow \infty }P_{\mu _{0},\mathcal{L}(\Sigma _{j})}(W_{j})&
\leq \limsup_{j\rightarrow \infty }P_{\mu _{0},\mathcal{L}(\Sigma _{j})}(%
\limfunc{cl}(W_{j})) \\
& \leq \limsup_{j\rightarrow \infty }P_{\mu _{0},\mathcal{L}(\Sigma _{j})}(%
\limfunc{cl}(W_{j_{0}})).
\end{align*}%
Since $P_{\mu _{0},\mathcal{L}(\Sigma _{j})}$ converges to $P_{\mu _{0},\bar{%
\Sigma}}$ weakly (cf. Lemma E.1 in \cite{PP2016}), the Portmanteau theorem
gives 
\begin{equation*}
\limsup_{j\rightarrow \infty }P_{\mu _{0},\mathcal{L}(\Sigma _{j})}(\limfunc{%
cl}(W_{j_{0}}))\leq P_{\mu _{0},\bar{\Sigma}}(\limfunc{cl}%
(W_{j_{0}}))<\varepsilon ,
\end{equation*}%
which then proves the claim via (\ref{invconsm}) and (\ref{supapprox}). $%
\blacksquare $

\begin{lemma}
\label{L3} Suppose $T:\mathbb{R}^{n}\rightarrow \mathbb{R}$ is continuous on
the complement of a closed set $N^{\dag }\subseteq \mathbb{R}^{n}$. Then for
every $-\infty <C<\infty $ the set $\{y\in \mathbb{R}^{n}:T(y)\geq C\}\cup
N^{\dag }$ is closed, and hence contains $\limfunc{cl}(\{y\in \mathbb{R}%
^{n}:T(y)\geq C\})$. Furthermore, $\limfunc{bd}(\{y\in \mathbb{R}%
^{n}:T(y)\geq C\})$ is contained in the closed set $\{y\in \mathbb{R}%
^{n}:T(y)=C\}\cup N^{\dag }$. [Analogous statements hold for the set $\{y\in 
\mathbb{R}^{n}:T(y)\leq C\}$.]
\end{lemma}

\textbf{Proof:} Suppose $y_{0}$ is an accumulation point of $A:=\{y\in 
\mathbb{R}^{n}:T(y)\geq C\}\cup N^{\dag }$. If $y_{0}\in N^{\dag }\subseteq
A $ we are done. Suppose next that $y_{0}\notin N^{\dag }$. Let $y_{j}\in A$
be such that $y_{j}\rightarrow y_{0}$. Then $y_{j}\notin N^{\dag }$
eventually, since $N^{\dag }$ is closed. By assumption $T$ is continuous on $%
\mathbb{R}^{n}\backslash N^{\dag }$. Hence $T(y_{j})\rightarrow T(y_{0})$.
Now, since $y_{j}\in A$ and since $y_{j}\notin N^{\dag }$ eventually, we
have $T(y_{j})\geq C$ eventually. But then $T(y_{0})\geq C$ follows,
implying that $y_{0}\in A$. We turn to the second claim. That $B:=\{y\in 
\mathbb{R}^{n}:T(y)=C\}\cup N^{\dag }$ is closed, is proved in a completely
analogous way. Next we establish the claimed inclusion: If $y_{\ast }$ is an
element of $\limfunc{bd}(\{y\in \mathbb{R}^{n}:T(y)\geq C\})$ and $y_{\ast }$
belongs to $N^{\dag }$ we are done. Assume $y_{\ast }$ does not belong to $%
N^{\dag }$. Since $y_{\ast }$ must also be an element of $\limfunc{cl}%
(\{y\in \mathbb{R}^{n}:T(y)\geq C\})$, and hence of $A$ by what has already
been shown, it follows that $T(y_{\ast })\geq C$ must hold. If $T(y_{\ast
})>C$ would be true, then $T$ would have to be larger than $C$ on an open
neighborhood of $y_{\ast }$, since $y_{\ast }\notin N^{\dag }$. This would
lead to the contradiction that $y_{\ast }$ belongs to the interior of $%
\{y\in \mathbb{R}^{n}:T(y)\geq C\}$. The claim in parenthesis is proved
completely analogously. $\blacksquare $

\textbf{Proof of Corollary \ref{C5}:} For the first statement we check the
conditions of Proposition \ref{L2}. The invariance properties of $W_{j}$
follow from invariance of $T$. That $W_{j}\supseteq W_{j+1}$ holds is
obvious. Condition (a) of Proposition \ref{L2} is satisfied because $%
\bigcap_{j\in \mathbb{N}}W_{j}$ is empty as $T$ is real-valued and $%
C_{j}\nearrow \infty $ as $j\rightarrow \infty $. Hence, Part 1 follows.
Next observe that%
\begin{equation*}
\bigcap_{j\in \mathbb{N}}\limfunc{cl}(W_{j})\subseteq \bigcap_{j\in \mathbb{N%
}}(W_{j}\cup N^{\dag })=\bigcap_{j\in \mathbb{N}}W_{j}\cup N^{\dag
}=\emptyset \cup N^{\dag }=N^{\dag },
\end{equation*}%
where we have used Lemma \ref{L3} to obtain the first inclusion. Part 2 then
follows immediately. For Part 3 note that $\lambda _{\mu _{0}+\mathcal{S}%
}(N^{\dag })=0$ obviously implies $\mu _{0}+\mathcal{S}\not\subseteq N^{\dag
}$, since $\lambda _{\mu _{0}+\mathcal{S}}$ is supported by $\mu _{0}+%
\mathcal{S}$. The converse is seen as follows: Clearly, $(\mu _{0}+\mathcal{%
S)}\cap N^{\dag }$ is a finite or countable union of sets of the form $(\mu
_{0}+\mathcal{S})\cap \mathcal{A}_{j}$, where the $\mathcal{A}_{j}$'s are
affine spaces, the union of which is $N^{\dag }$. Since $\mu _{0}+\mathcal{S}%
\not\subseteq N^{\dag }$, the sets $(\mu _{0}+\mathcal{S})\cap \mathcal{A}%
_{j}$ must be proper affine subspaces of $\mu _{0}+\mathcal{S}$ or must be
empty, entailing $\lambda _{\mu _{0}+\mathcal{S}}((\mu _{0}+\mathcal{S})\cap 
\mathcal{A}_{j})=0$ for every $j$. Since $\lambda _{\mu _{0}+\mathcal{S}%
}(N^{\dag })=\lambda _{\mu _{0}+\mathcal{S}}((\mu _{0}+\mathcal{S)}\cap
N^{\dag })$, we conclude that $\lambda _{\mu _{0}+\mathcal{S}}(N^{\dag })=0$%
. $\blacksquare $

\textbf{Proof of claim in Remark \ref{some comments}(ii):} First note that $%
\bigcap_{j\in \mathbb{N}}\limfunc{cl}(W_{j})$ and $N^{\dag }$ are $G(%
\mathfrak{M}_{0})$-invariant. The reverse implication is trivial since $%
\lambda _{\mu _{0}+\mathcal{S}}$ is supported on $\mu _{0}+\mathcal{S}$ and
has a density there, and thus assigns zero mass to $\left\{ \mu _{0}\right\} 
$ as $\dim \mathcal{S}>0$. To prove the other implication suppose that $x\in
\bigcap_{j\in \mathbb{N}}\limfunc{cl}(W_{j})\cap (\mu _{0}+\mathcal{S})$ but
is different from $\mu _{0}$. Then $\delta (x-\mu _{0})+\mu _{0}$ also
belongs to this set for every $\delta \neq 0$ in view of $G(\mathfrak{M}%
_{0}) $-invariance of $\bigcap_{j\in \mathbb{N}}\limfunc{cl}(W_{j})$ and
because of the particular form of $\mu _{0}+\mathcal{S}$. Note that $s=x-\mu
_{0}$ is nonzero and belongs to $\mathcal{S}$, and hence spans $\mathcal{S}$%
. We conclude that $\left( \mu _{0}+\mathcal{S}\right) \backslash \left\{
\mu _{0}\right\} $ is a subset of $\bigcap_{j\in \mathbb{N}}\limfunc{cl}%
(W_{j})$. But this contradicts the assumption that $\bigcap_{j\in \mathbb{N}}%
\limfunc{cl}(W_{j})$ is a $\lambda _{\mu _{0}+\mathcal{S}}$-null set. The
proof for $N^{\dag }$ is completely analogous. $\blacksquare $

\textbf{Proof of Lemma \ref{L4}:} Since $T$ is $G(\mathfrak{M}_{0})$%
-invariant and since for a one-dimensional $\mathcal{S}$ every element of $%
\mu _{0}+\mathcal{S}$ has the form $\mu _{0}+\gamma s_{0}$ for a fixed $%
s_{0}\in \mathcal{S}$, we conclude that 
\begin{equation*}
T\left( \mu _{0}+\gamma s_{0}\right) =T\left( \gamma ^{-1}\left( \left( \mu
_{0}+\gamma s_{0}\right) -\mu _{0}\right) +\mu _{0}\right) =T\left( \mu
_{0}+s_{0}\right)
\end{equation*}%
for every $\gamma \neq 0$. Hence, any one-dimensional $\mathcal{S}\in 
\mathbb{J}(\mathcal{L},\mathfrak{C})$ satisfies $\mathcal{S}\in \mathbb{H}$.
This proves Part 1. In view of the assumed $G(\mathfrak{M}_{0})$-invariance
of $T$ we may for the rest of the proof fix an arbitrary $\mu _{0}\in 
\mathfrak{M}_{0}$, set $\sigma ^{2}=1$, and drop the suprema (infima) w.r.t. 
$\mu _{0}$ and $\sigma ^{2}$ from (\ref{size_3}) and (\ref{bias}), cf.
Remark 5.5(iii) in \cite{PP2016}. We turn to Part 2: Concerning (\ref{size_3}%
), note that in case $\mathbb{H}$ is empty there is nothing to prove, hence
assume that $\mathbb{H}$ is nonempty and choose a $C\in (-\infty ,C^{\ast })$%
. Since $C<C^{\ast }$, there exists an $\mathcal{S}\in \mathbb{H}$ such that
the corresponding constant $C(\mathcal{S})$ satisfies $C<C(\mathcal{S})\leq
C^{\ast }$. By assumption, $T$ is continuous on the complement of $N^{\dag }$%
, a closed $\lambda _{\mu _{0}+\mathcal{S}}$-null set. By definition of $C(%
\mathcal{S})$ we have that $\lambda _{\mu _{0}+\mathcal{S}}(T\neq C(\mathcal{%
S}))=0$ which together with $C(\mathcal{S})>C$ and Lemma \ref{L3} implies
that%
\begin{eqnarray*}
\lambda _{\mu _{0}+\mathcal{S}}([\limfunc{int}(\{y &\in &\mathbb{R}%
^{n}:T(y)\geq C\})]^{c})=\lambda _{\mu _{0}+\mathcal{S}}(\limfunc{cl}(\{y\in 
\mathbb{R}^{n}:T(y)<C\})) \\
&\leq &\lambda _{\mu _{0}+\mathcal{S}}(\limfunc{cl}(\{y\in \mathbb{R}%
^{n}:T(y)\leq C\}))\leq \lambda _{\mu _{0}+\mathcal{S}}(\{y\in \mathbb{R}%
^{n}:T(y)\leq C\})+\lambda _{\mu _{0}+\mathcal{S}}(N^{\dag }) \\
&\leq &\lambda _{\mu _{0}+\mathcal{S}}(T\neq C(\mathcal{S}))+\lambda _{\mu
_{0}+\mathcal{S}}(N^{\dag })=0.
\end{eqnarray*}%
Now, let $\Sigma _{j}$ be a sequence in $\mathfrak{C}$ such that $\mathcal{L}%
(\Sigma _{j})\rightarrow \bar{\Sigma}$ with $\limfunc{span}(\bar{\Sigma})=%
\mathcal{S}$. Then by $G(\mathfrak{M}_{0})$-invariance of $T$ and its
invariance w.r.t. addition of elements of $\mathcal{V}$ we have (cf. the
proof of (\ref{invcons}))%
\begin{align*}
\sup_{\Sigma \in \mathfrak{C}}P_{\mu _{0},\Sigma }(T\geq C)& \geq
\liminf_{j\rightarrow \infty }P_{\mu _{0},\Sigma _{j}}(T\geq
C)=\liminf_{j\rightarrow \infty }P_{\mu _{0},\mathcal{L}(\Sigma _{j})}(T\geq
C) \\
& \geq \liminf_{j\rightarrow \infty }P_{\mu _{0},\mathcal{L}(\Sigma
_{j})}\left( \limfunc{int}\left( \left\{ y\in \mathbb{R}^{n}:T(y)\geq
C\right\} \right) \right) \\
& \geq P_{\mu _{0},\bar{\Sigma}}\left( \limfunc{int}\left( \left\{ y\in 
\mathbb{R}^{n}:T(y)\geq C\right\} \right) \right) ,
\end{align*}%
where the last inequality follows from weak convergence of $P_{\mu _{0},%
\mathcal{L}(\Sigma _{j})}$ to $P_{\mu _{0},\bar{\Sigma}}$ (cf. Lemma E.1 in 
\cite{PP2016}) and the Portmanteau theorem. But absolute continuity of $%
P_{\mu _{0},\bar{\Sigma}}$ w.r.t. $\lambda _{\mu _{0}+\mathcal{S}}$ then
implies%
\begin{equation*}
P_{\mu _{0},\bar{\Sigma}}\left( \limfunc{int}\left( \left\{ y\in \mathbb{R}%
^{n}:T(y)\geq C\right\} \right) \right) =1-P_{\mu _{0},\bar{\Sigma}}(\left[ 
\limfunc{int}\left( \left\{ y\in \mathbb{R}^{n}:T(y)\geq C\right\} \right) %
\right] ^{c})=1.
\end{equation*}%
This proves (\ref{size_3}). The statement concerning the convergence
behavior of the size as $C\rightarrow \infty $ is a consequence of Corollary %
\ref{C5}.

Finally consider Part 3. Again, if $\mathbb{H}$ is empty there is nothing to
prove. Hence assume $\mathbb{H\neq \varnothing }$. Let $C\in (C_{\ast
},\infty )$. Then we can find $\mathcal{S}\in \mathbb{H}$ with $C(\mathcal{S}%
)<C$. Furthermore, by definition we can find a sequence $\Sigma _{j}\in 
\mathfrak{C}$ such that $\bar{\Sigma}_{j}=\mathcal{L}\left( \Sigma
_{j}\right) $ converges to $\bar{\Sigma}$ with $\limfunc{span}(\bar{\Sigma})=%
\mathcal{S}$. It follows that $T=C(\mathcal{S})<C$ holds $\lambda _{\mu _{0}+%
\limfunc{span}(\bar{\Sigma})}$-a.e., since $\mathcal{S}\in \mathbb{H}$.
Hence, $\lambda _{\mu _{0}+\limfunc{span}(\bar{\Sigma})}(T\geq C)=0$, which
entails $P_{\mu _{0},\bar{\Sigma}}(T\geq C)=0$ by equivalence of the
measures involved. Using (\ref{invcons}), the Portmanteau theorem, and Lemma %
\ref{L3}, we obtain%
\begin{eqnarray*}
\inf_{\Sigma \in \mathfrak{C}}P_{\mu _{0},\Sigma }(T\geq C) &\leq
&\limsup_{j\rightarrow \infty }P_{\mu _{0},\Sigma _{j}}(T\geq
C)=\limsup_{j\rightarrow \infty }P_{\mu _{0},\bar{\Sigma}_{j}}(T\geq C) \\
&\leq &\limsup_{j\rightarrow \infty }P_{\mu _{0},\bar{\Sigma}_{j}}\left( 
\limfunc{cl}\left( \left\{ y\in \mathbb{R}^{n}:T(y)\geq C\right\} \right)
\right) \\
&\leq &P_{\mu _{0},\bar{\Sigma}}\left( \limfunc{cl}\left( \left\{ y\in 
\mathbb{R}^{n}:T(y)\geq C\right\} \right) \right) \\
&\leq &P_{\mu _{0},\bar{\Sigma}}\left( T\geq C\right) +P_{\mu _{0},\bar{%
\Sigma}}(N^{\dag })=P_{\mu _{0},\bar{\Sigma}}(N^{\dag })=0,
\end{eqnarray*}%
the last equality following from equivalence of $P_{\mu _{0},\bar{\Sigma}}$
with $\lambda _{\mu _{0}+\limfunc{span}(\bar{\Sigma})}$, from $\limfunc{span}%
(\bar{\Sigma})\in \mathbb{H}$, and the assumptions made on $N^{\dag }$. $%
\blacksquare $

\textbf{Proof of Proposition \ref{P6}:} We start with the observation that
the statements \textquotedblleft and hence for all $\mu _{0}\in \mathfrak{M}%
_{0}$\textquotedblright\ in parentheses follow from the corresponding
statements involving \textquotedblleft and hence for some $\mu _{0}\in 
\mathfrak{M}_{0}$\textquotedblright\ in view of the assumed $G(\mathfrak{M}%
_{0})$-invariance of $T$ and $N^{\dag }$. In view of the assumed $G(%
\mathfrak{M}_{0})$-invariance of $T$ we may for the rest of the proof fix an
arbitrary $\mu _{0}\in \mathfrak{M}_{0}$, set $\sigma ^{2}=1$, and drop the
suprema w.r.t. $\mu _{0}$ and $\sigma ^{2}$ from the displayed expressions
shown in the lemma, cf. Remark 5.5(iii) in \cite{PP2016}.

We now establish Part A.1 of the lemma: The statement concerning
nonincreasingness is obvious. The constancy property of the size as well as
its convergence to zero for $C\rightarrow \infty $ has already been
established in Lemma \ref{L4}. \ Therefore, it remains to verify that the
function $C\mapsto \sup_{\Sigma \in \mathfrak{C}}P_{\mu _{0},\Sigma }(T\geq
C)$ is continuous on $(C^{\ast },\infty )$, where we note that $C^{\ast
}<\infty $ holds. In order to achieve this we proceed in two steps:

\textit{Step 1:}\textbf{\ }We show that the map $(C,\bar{\Sigma})\mapsto
P_{\mu _{0},\bar{\Sigma}}(T\geq C)$ is continuous on $(C^{\ast },\infty
)\times \limfunc{cl}(\mathcal{L}\left( \mathfrak{C}\right) )$ and that $%
\limfunc{cl}(\mathcal{L}\left( \mathfrak{C}\right) )$ is compact.
Compactness is obvious, since $\mathcal{L}\left( \mathfrak{C}\right) $ is
norm bounded by construction. In order to establish continuity, let $%
C_{j}\in (C^{\ast },\infty )$ be a sequence such that $C_{j}\rightarrow C\in
(C^{\ast },\infty )$ and let $\bar{\Sigma}_{j}\in \limfunc{cl}(\mathcal{L}%
\left( \mathfrak{C}\right) )$ converge to $\bar{\Sigma}\in \limfunc{cl}(%
\mathcal{L}\left( \mathfrak{C}\right) )$. In view of the assumed invariance
properties of $T$ we have $P_{\mu _{0},\bar{\Sigma}_{j}}(T\geq \cdot
)=P_{\mu _{0},\Omega _{j}}(T\geq \cdot )$ and $P_{\mu _{0},\bar{\Sigma}%
}(T\geq \cdot )=P_{\mu _{0},\Omega }(T\geq \cdot )$ where $\Omega _{j}=\bar{%
\Sigma}_{j}+\Pi _{\mathcal{L}}$ and $\Omega =\bar{\Sigma}+\Pi _{\mathcal{L}}$
(cf. the proof of (\ref{invcons})). Note that $\Omega _{j}$ converges to $%
\Omega $, and thus $P_{\mu _{0},\Omega _{j}}$ converges to $P_{\mu
_{0},\Omega }$ weakly (cf. Lemma E.1 in \cite{PP2016}), and in fact in total
variation distance if $\Omega $ is nonsingular.

\textit{Case 1:} Assume that\textit{\ }$\Omega $ is nonsingular. Then
convergence of $P_{\mu _{0},\Omega _{j}}$ to $P_{\mu _{0},\Omega }$ in total
variation distance implies%
\begin{equation*}
P_{\mu _{0},\bar{\Sigma}_{j}}(T\geq C_{j})-P_{\mu _{0},\bar{\Sigma}}(T\geq
C_{j})=P_{\mu _{0},\Omega _{j}}(T\geq C_{j})-P_{\mu _{0},\Omega }(T\geq
C_{j})\rightarrow 0\quad \text{for}\quad j\rightarrow \infty .
\end{equation*}%
Furthermore, 
\begin{equation*}
P_{\mu _{0},\bar{\Sigma}}(T\geq C_{j})=P_{\mu _{0},\Omega }(T\geq
C_{j})\rightarrow P_{\mu _{0},\Omega }(T\geq C)=P_{\mu _{0},\bar{\Sigma}%
}(T\geq C)\quad \text{for}\quad j\rightarrow \infty ,
\end{equation*}%
since $P_{\mu _{0},\Omega }(T=C)=0$ holds in view of the assumption $\lambda
_{\mathbb{R}^{n}}(T=C)=0$ for $C>C^{\ast }$ and equivalence of $P_{\mu
_{0},\Omega }$ and $\lambda _{\mathbb{R}^{n}}$. Together this implies the
desired convergence.

\textit{Case 2: }Assume that $\Omega $ is singular. Then $\limfunc{span}(%
\bar{\Sigma})\in \mathbb{J}(\mathcal{L},\mathfrak{C})$ follows. We
distinguish two subcases:

\textit{Case 2a: }Assume that $\limfunc{span}(\bar{\Sigma})\in \mathbb{H}$.
Choose $\varepsilon >0$ small enough such that $C-\varepsilon >C^{\ast }$
holds. Then $T\leq C^{\ast }<C-\varepsilon $ holds $\lambda _{\mu _{0}+%
\limfunc{span}(\bar{\Sigma})}$-a.e., since $\limfunc{span}(\bar{\Sigma})\in 
\mathbb{H}$. Consequently, $\lambda _{\mu _{0}+\limfunc{span}(\bar{\Sigma}%
)}(T\geq C-\varepsilon )=0$, which entails $P_{\mu _{0},\bar{\Sigma}}(T\geq
C-\varepsilon )=0$ by equivalence of the measures involved. Consequently,
also $P_{\mu _{0},\bar{\Sigma}}(T\geq C)=0$ holds. Using the just
established identities, the Portmanteau theorem, and Lemma \ref{L3}, we
obtain%
\begin{eqnarray*}
\limsup_{j\rightarrow \infty }P_{\mu _{0},\bar{\Sigma}_{j}}(T\geq C_{j})
&\leq &\limsup_{j\rightarrow \infty }P_{\mu _{0},\bar{\Sigma}_{j}}\left( 
\limfunc{cl}\left( \left\{ y\in \mathbb{R}^{n}:T(y)\geq C-\varepsilon
\right\} \right) \right) \\
&\leq &P_{\mu _{0},\bar{\Sigma}}\left( \limfunc{cl}\left( \left\{ y\in 
\mathbb{R}^{n}:T(y)\geq C-\varepsilon \right\} \right) \right) \\
&\leq &P_{\mu _{0},\bar{\Sigma}}\left( T\geq C-\varepsilon \right) +P_{\mu
_{0},\bar{\Sigma}}(N^{\dag })=P_{\mu _{0},\bar{\Sigma}}(N^{\dag })=0,
\end{eqnarray*}%
the last equality following from equivalence of $P_{\mu _{0},\bar{\Sigma}}$
with $\lambda _{\mu _{0}+\limfunc{span}(\bar{\Sigma})}$, from $\limfunc{span}%
(\bar{\Sigma})\in \mathbb{J}(\mathcal{L},\mathfrak{C})$, and the assumptions
made on $N^{\dag }$. This establishes the desired convergence.

\textit{Case 2b: }Assume that $\limfunc{span}(\bar{\Sigma})\in \mathbb{J}(%
\mathcal{L},\mathfrak{C})\backslash \mathbb{H}$. Choose $\varepsilon >0$ as
before. By Lemma \ref{L3} and our assumptions on $N^{\dag }$, and since $%
\lambda _{\mu _{0}+\limfunc{span}(\bar{\Sigma})}(T=C-\varepsilon )=0$ in
view of $\limfunc{span}(\bar{\Sigma})\in \mathbb{J}(\mathcal{L},\mathfrak{C}%
)\backslash \mathbb{H}$ and of $C-\varepsilon >C^{\ast }$, it follows that 
\begin{equation*}
\lambda _{\mu _{0}+\limfunc{span}(\bar{\Sigma})}\left( \limfunc{bd}\left(
\left\{ y\in \mathbb{R}^{n}:T(y)\geq C-\varepsilon \right\} \right) \right)
=0.
\end{equation*}%
Using equivalence of $P_{\mu _{0},\bar{\Sigma}}$ with $\lambda _{\mu _{0}+%
\limfunc{span}(\bar{\Sigma})}$ we also have 
\begin{equation*}
P_{\mu _{0},\bar{\Sigma}}\left( \limfunc{bd}\left( \left\{ y\in \mathbb{R}%
^{n}:T(y)\geq C-\varepsilon \right\} \right) \right) =0,
\end{equation*}%
and the same is true if $C-\varepsilon $ is replaced by $C+\varepsilon $.
From weak convergence we thus obtain%
\begin{equation*}
P_{\mu _{0},\bar{\Sigma}_{j}}(T\geq C\pm \varepsilon )\rightarrow P_{\mu
_{0},\bar{\Sigma}}(T\geq C\pm \varepsilon )\quad \text{for}\quad
j\rightarrow \infty \text{.}
\end{equation*}%
This implies%
\begin{equation*}
P_{\mu _{0},\bar{\Sigma}}(T\geq C-\varepsilon )\geq \limsup_{j\rightarrow
\infty }P_{\mu _{0},\bar{\Sigma}_{j}}(T\geq C_{j})\geq \liminf_{j\rightarrow
\infty }P_{\mu _{0},\bar{\Sigma}_{j}}(T\geq C_{j})\geq P_{\mu _{0},\bar{%
\Sigma}}(T\geq C+\varepsilon ).
\end{equation*}%
Observe that $P_{\mu _{0},\bar{\Sigma}}(T=C)=0$ holds, since $\lambda _{\mu
_{0}+\limfunc{span}(\bar{\Sigma})}(T=C)=0$ by our assumptions (note that $%
C>C^{\ast }$ and that $\limfunc{span}(\bar{\Sigma})\in \mathbb{J}(\mathcal{L}%
,\mathfrak{C})\backslash \mathbb{H}$). Letting $\varepsilon $ go to zero in
the above display then gives the desired convergence of $P_{\mu _{0},\bar{%
\Sigma}_{j}}(T\geq C_{j})$ to $P_{\mu _{0},\bar{\Sigma}}(T\geq C)$.

\textit{Step 2:} Note that by the assumed invariance properties and by the
definition of $\mathcal{L}\left( \mathfrak{C}\right) $ we have for every $%
C>C^{\ast }$ that%
\begin{equation*}
\sup_{\Sigma \in \mathfrak{C}}P_{\mu _{0},\Sigma }(T\geq C)=\sup_{\Sigma \in 
\mathfrak{C}}P_{\mu _{0},\mathcal{L}(\Sigma )}(T\geq C)=\sup_{\bar{\Sigma}%
\in \mathcal{L}\left( \mathfrak{C}\right) }P_{\mu _{0},\bar{\Sigma}}(T\geq
C)=\sup_{\bar{\Sigma}\in \limfunc{cl}(\mathcal{L}\left( \mathfrak{C}\right)
)}P_{\mu _{0},\bar{\Sigma}}(T\geq C),
\end{equation*}%
the last equality following from the continuity established in Step 1. But
the right-most supremum is continuous on $(C^{\ast },\infty )$ as a
consequence of the claim established in Step 1 and Lemma \ref{max-theorem}
given below. [Note that $(C^{\ast },\infty )$ as well as $\limfunc{cl}(%
\mathcal{L}\left( \mathfrak{C}\right) )$ are not empty, since $C^{\ast
}<\infty $ has been established before and since $\mathfrak{C}\neq
\varnothing $ by assumption.] This completes the proof of Part A.1.

The claims in Part A.2 are now immediate consequences of the already
established Part A.1. Inspection of the proof of Part A.1 shows that under
the assumptions of Part B continuity of the size on $(C^{\ast \ast },\infty
) $ follows. Everything else in Part B is then proved similarly as the
corresponding claims in Part A. $\blacksquare $

The following lemma is a special case of Berge's maximum theorem, see \cite%
{berge1963}, Chapter VI, Section 3.

\begin{lemma}
\label{max-theorem}Let $f:A\times B\rightarrow \mathbb{R}$ be a continuous
map, where $A$ is a (nonempty) topological space and $B$ is a (nonempty)
compact topological space. Then $g(a)=\sup_{b\in B}f(a,b)$ for $a\in A$
defines a continuous map $g:A\rightarrow \mathbb{R}$.
\end{lemma}

\textbf{Proof of Lemma \ref{neu}:} Since $T$ is $G(\mathfrak{M}_{0})$%
-invariant (see Lemma 5.15 and Proposition 5.4 in \cite{PP2016}) it follows
that $P_{\mu _{0},\sigma ^{2}\Sigma }(T\geq C)$ for $\mu _{0}\in \mathfrak{M}%
_{0}$ does neither depend on the choice of $\mu _{0}\in \mathfrak{M}_{0}$
nor on $\sigma ^{2}$. We hence we may fix $\mu _{0}\in \mathfrak{M}_{0}$ and
set $\sigma ^{2}=1$. Furthermore, $T(y)=T_{0}(y-\mu _{0})$ with%
\begin{equation*}
T_{0}(y)=\left\{ 
\begin{array}{cc}
(R\check{\beta}(y))^{\prime }\check{\Omega}^{-1}(y)(R\check{\beta}(y)) & 
y\in \mathbb{R}^{n}\backslash N^{\ast } \\ 
0 & y\in N^{\ast }%
\end{array}%
\right.
\end{equation*}
follows, because $N^{\ast }$ is $G(\mathfrak{M})$-invariant by Lemma 5.15 in 
\cite{PP2016} and because of the equivariance (invariance) requirements on $%
\check{\beta}$ ($\check{\Omega}$, respectively) made in Assumption 5 of \cite%
{PP2016}. Consequently,%
\begin{equation*}
P_{\mu _{0},\Sigma }(T(y)\geq C)=P_{\mu _{0},\Sigma }(T_{0}(y-\mu _{0})\geq
C)=P_{0,\Sigma }(T_{0}(y)\geq C)
\end{equation*}%
holds, where the most right-hand expression does not depend on the value of $%
r$. $\blacksquare $

\textbf{Proof of Lemma \ref{Fact1}:} Parts 1 and 3 have been established in
Lemma 5.15 of \cite{PP2016}, Borel-measurability being trivial. Consider
next Part 2: Lemma 5.15 of \cite{PP2016} shows that $\lambda _{\mathbb{R}%
^{n}}(T=C)=0$ holds for $C>0$. It follows immediately, that this also holds
for $C<0$ (by passing from $T$ to $-T$, absorbing the sign into $\check{%
\Omega}$, and by applying Lemma 5.15 in that reference to $-T$). That $%
\lambda _{\mathbb{R}^{n}}(T=C)=0$ also holds for $C=0$ is seen as follows:
Write the set $O=\left\{ y\in \mathbb{R}^{n}:T(y)=0\right\} $ as $O_{\ast
}\cup N^{\ast }$ where $O_{\ast }=\left\{ y\in \mathbb{R}^{n}\backslash
N^{\ast }:T(y)=0\right\} $. Certainly $O_{\ast }\subseteq \mathbb{R}%
^{n}\backslash N^{\ast }$ by construction. It suffices to show that $O_{\ast
}$ is a $\lambda _{\mathbb{R}^{n}}$-null set: But this follows from
repeating the arguments given in the proof of Part 4 of Lemma 5.15 in \cite%
{PP2016} for the set $O_{\ast }$ (instead of for $O$), with the only change
that the argument that the set $O(y_{2})$ constructed in the proof is empty
if $y_{2}\in N^{\ast }\cap \mathfrak{M}^{\bot }$ now has to be deduced from
the observation that $y=y_{1}+y_{2}\in \mathbb{R}^{n}\backslash N^{\ast }$
is not possible if $y_{2}\in N^{\ast }\cap \mathfrak{M}^{\bot }$, since $%
y_{1}\in \mathfrak{M}$ and since $N^{\ast }$ is $G(\mathfrak{M})$-invariant. 
$\blacksquare $

\textbf{Proof of Lemma \ref{Fact2}:} By definition $\mathsf{B}$ is the set
where 
\begin{eqnarray*}
B(y) &=&R(X^{\prime }X)^{-1}X^{\prime }\limfunc{diag}\left( e_{1}^{\prime
}(n)\Pi _{\limfunc{span}(X)^{\bot }}y,\ldots ,e_{n}^{\prime }(n)\Pi _{%
\limfunc{span}(X)^{\bot }}y\right) \\
&=&R(X^{\prime }X)^{-1}X^{\prime }\left[ e_{1}(n)e_{1}^{\prime }(n)\Pi _{%
\limfunc{span}(X)^{\bot }}y,\ldots ,e_{n}(n)e_{n}^{\prime }(n)\Pi _{\limfunc{%
span}(X)^{\bot }}y\right]
\end{eqnarray*}%
has rank less than $q$. Define the set 
\begin{equation*}
D=\left\{ (j_{1},\ldots ,j_{s}):1\leq s\leq n,1\leq j_{1}<\ldots <j_{s}\leq
n,\limfunc{rank}\left( R(X^{\prime }X)^{-1}X^{\prime }\left[
e_{j_{1}}(n),\ldots ,e_{j_{s}}(n)\right] \right) <q\right\} ,
\end{equation*}%
which may be empty in case $q=1$. Consider first the case where $D$ is
nonempty: Since $R(X^{\prime }X)^{-1}X^{\prime }$ has rank $q$, it is then
easy to see that we have $y\in \mathsf{B}$ if and only if there exists $%
(j_{1},\ldots ,j_{s})\in D$ such that $e_{j}^{\prime }(n)\Pi _{\limfunc{span}%
(X)^{\bot }}y=0$ for $j\neq j_{i}$ for $i=1,\ldots ,s$. This shows, that $%
\mathsf{B}$ is a finite union of (not necessarily distinct) linear
subspaces. In case $D$ is empty, $\limfunc{rank}(R(X^{\prime
}X)^{-1}X^{\prime })=q$ implies that $y\in \mathsf{B}$ if and only if $%
e_{j}^{\prime }(n)\Pi _{\limfunc{span}(X)^{\bot }}y=0$ for all $1\leq j\leq
n $, i.e., if and only if $y\in \limfunc{span}(X)$. That the linear
subspaces making up $\mathsf{B}$ are proper, follows since otherwise $%
\mathsf{B}$ would be all of $\mathbb{R}^{n}$, which is impossible under
Assumptions \ref{R_and_X} as transpires from an inspection of Lemma 3.1 (and
its proof) in \cite{PP2016}. To prove the second claim, observe that in case 
$q=1$ the condition that $\limfunc{rank}(B(y))$ is less than $q$ is
equivalent to $B(y)=0$. Since the expressions $R(X^{\prime }X)^{-1}X^{\prime
}e_{j}(n)$ are now scalar, we may thus write the condition $B(y)=0$
equivalently as%
\begin{equation*}
\left[ R(X^{\prime }X)^{-1}X^{\prime }e_{1}(n)e_{1}(n),\ldots ,R(X^{\prime
}X)^{-1}X^{\prime }e_{n}(n)e_{n}(n)\right] \Pi _{\limfunc{span}(X)^{\bot
}}y=0.
\end{equation*}%
But this shows that $\mathsf{B}$ is a linear space, namely the kernel of the
matrix appearing on the l.h.s. of the preceding display. $\blacksquare $

\textbf{Proof of Lemma \ref{Fact3}:} We first prove (ii): Note that the set $%
N^{\ast }$ on which $\hat{\Omega}_{GQ}(y)$ is singular coincides with $%
\mathsf{B}$, and hence is a finite union of proper linear subspaces of $%
\mathbb{R}^{n}$ by Lemma \ref{Fact2}. Since $T_{GQ}$ is constant on $N^{\ast
}$ by definition, it follows that $\mu +\mathcal{S}\not\subseteq N^{\ast }$
must hold. An argument like the one discussed in Remark \ref{some comments}%
(i) then shows that $N^{\ast }$ is a $\lambda _{\mu +\mathcal{S}}$ null set.
Consequently, $T_{GQ}$ restricted to $(\mu +\mathcal{S})\backslash N^{\ast }$
is not constant. Suppose now there exists a $C$ so that $\lambda _{\mu +%
\mathcal{S}}(\{y\in \mathbb{R}^{n}:T(y)=C\})>0$. Then, since $N^{\ast }$ is
a $\lambda _{\mu +\mathcal{S}}$-null set, it follows that even $\lambda
_{\mu +\mathcal{S}}(\{y\in \mathbb{R}^{n}\backslash N^{\ast }:T(y)=C\})>0$,
which can be written as $\lambda _{\mu +\mathcal{S}}(\{y\in \mathbb{R}%
^{n}\backslash N^{\ast }:p(y)=0\})>0$, with $p(y)=(R\hat{\beta}%
(y)-r)^{\prime }\limfunc{adj}(\hat{\Omega}_{GQ}(y))(R\hat{\beta}(y)-r)-\det (%
\hat{\Omega}_{GQ}(y))C$, a polynomial in $y$. This implies that $p$
restricted to $\mu +\mathcal{S}$ vanishes on a set of positive $\lambda
_{\mu +\mathcal{S}}$-measure. Since $p$ can clearly be expressed as a
polynomial in coordinates parameterizing the affine space $\mu +\mathcal{S}$%
, it follows that $p$ vanishes identically on $\mu +\mathcal{S}$. But this
implies that $T_{GQ}$ restricted to $(\mu +\mathcal{S})\backslash N^{\ast }$
is constant equal to $C$, a contradiction. Part (i) follows as a special
case of Part (ii). The proof of Part (iii) is completely analogous, noting
that for the weighted Eicker-test statistic the set $N^{\ast }$ is always $%
\limfunc{span}(X)$. $\blacksquare $

\subsection{Sufficient conditions relating to $\mathbb{J}(\mathcal{L},%
\mathfrak{C})$ \label{suff_con_2}}

The subsequent lemma sheds light on the relation between the collection $%
\mathbb{J}(\mathcal{L},\mathfrak{C})$ and the set of concentration spaces of 
$\mathfrak{C}$ in an important case and leads to sufficient conditions
discussed in the remark given below. Recall from Definition 2.1 in \cite%
{PP2016} that a linear subspace $\mathcal{Z}$ of $\mathbb{R}^{n}$ is said to
be a concentration space of $\mathfrak{C}$, if $\dim (\mathcal{Z})<n$ and if
there exists a sequence $\Sigma _{m}\in \mathfrak{C}$ such that $\Sigma
_{m}\rightarrow \Sigma ^{\ast }$ with $\limfunc{span}(\Sigma ^{\ast })=%
\mathcal{Z}$.

\begin{lemma}
\label{suff_con_for J} Let $\mathfrak{C}$ be a covariance model and let $%
\mathcal{L}$ be a linear subspace of $\mathbb{R}^{n}$ with $\dim (\mathcal{L}%
)<n$. Then the following hold:

\begin{enumerate}
\item If $\mathcal{Z}$ is a concentration space of $\mathfrak{C}$, then
either $\Pi _{\mathcal{L}^{\bot }}\mathcal{Z}$ is an element of $\mathbb{J}(%
\mathcal{L},\mathfrak{C})$ or $\Pi _{\mathcal{L}^{\bot }}\mathcal{Z}=%
\mathcal{L}^{\bot }$ or $\Pi _{\mathcal{L}^{\bot }}\mathcal{Z}=\left\{
0\right\} $.

\item Suppose that $\mathfrak{C}$ is bounded. If $\mathcal{S}\in \mathbb{J}(%
\mathcal{L},\mathfrak{C})$ then either (i) there is a concentration space $%
\mathcal{Z}$ of $\mathfrak{C}$ with the property that $\mathcal{S}=\Pi _{%
\mathcal{L}^{\bot }}\mathcal{Z}$ or (ii) there exists a concentration space $%
\mathcal{Z}$ of $\mathfrak{C}$ satisfying $\Pi _{\mathcal{L}^{\bot }}%
\mathcal{Z}=\left\{ 0\right\} $ (i.e., $\mathcal{Z}\subseteq \mathcal{L}$)
and a sequence $\Sigma _{j}\in \mathfrak{C}$ such that $\Sigma _{j}$
converges to $\Sigma ^{\ast }$ satisfying $\mathcal{Z}=\limfunc{span}(\Sigma
^{\ast })$ and such that $\Pi _{\mathcal{L}^{\bot }}\Sigma _{j}\Pi _{%
\mathcal{L}^{\bot }}/\left\Vert \Pi _{\mathcal{L}^{\bot }}\Sigma _{j}\Pi _{%
\mathcal{L}^{\bot }}\right\Vert $ converges to $\bar{\Sigma}$ satisfying $%
\mathcal{S}=\limfunc{span}(\bar{\Sigma})$.
\end{enumerate}
\end{lemma}

\textbf{Proof:} 1. If $\mathcal{Z}$ is a concentration space, we can find a
sequence $\Sigma _{j}\in \mathfrak{C}$ such that $\Sigma _{j}$ converges to
a singular matrix $\tilde{\Sigma}$ with $\mathcal{Z}=\limfunc{span}(\tilde{%
\Sigma})$. If $\Pi _{\mathcal{L}^{\bot }}\mathcal{Z}=\mathcal{L}^{\bot }$ or 
$\Pi _{\mathcal{L}^{\bot }}\mathcal{Z}=\left\{ 0\right\} $ holds, we are
done, since neither $\mathcal{L}^{\bot }$ nor $\left\{ 0\right\} $ can
belong to $\mathbb{J}(\mathcal{L},\mathfrak{C})$.\footnote{%
If $\mathcal{L}$ $=\left\{ 0\right\} $, (i) the case $\Pi _{\mathcal{L}%
^{\bot }}\mathcal{Z}=\mathcal{L}^{\bot }$ is impossible, and (ii) the case $%
\Pi _{\mathcal{L}^{\bot }}\mathcal{Z}=\left\{ 0\right\} $ cannot arise if $%
\mathfrak{C}$ is bounded away from the zero matrix.} Hence assume that $%
\left\{ 0\right\} \subsetneq \Pi _{\mathcal{L}^{\bot }}\mathcal{Z}\subsetneq 
\mathcal{L}^{\bot }$. Then it is easy to see that $\Pi _{\mathcal{L}^{\bot }}%
\tilde{\Sigma}\Pi _{\mathcal{L}^{\bot }}\neq 0$ must hold. But then $\Pi _{%
\mathcal{L}^{\bot }}\Sigma _{j}\Pi _{\mathcal{L}^{\bot }}/\left\Vert \Pi _{%
\mathcal{L}^{\bot }}\Sigma _{j}\Pi _{\mathcal{L}^{\bot }}\right\Vert $
converges to $\bar{\Sigma}:=\Pi _{\mathcal{L}^{\bot }}\tilde{\Sigma}\Pi _{%
\mathcal{L}^{\bot }}/\left\Vert \Pi _{\mathcal{L}^{\bot }}\tilde{\Sigma}\Pi
_{\mathcal{L}^{\bot }}\right\Vert $. Because of $\Pi _{\mathcal{L}^{\bot }}%
\mathcal{Z}\subsetneq \mathcal{L}^{\bot }$, it follows that $\limfunc{rank}(%
\bar{\Sigma})<n-\dim (\mathcal{L})$ must hold, showing that $\limfunc{span}(%
\bar{\Sigma})\in \mathbb{J}(\mathcal{L},\mathfrak{C})$. It remains to show
that $\limfunc{span}(\bar{\Sigma})=\Pi _{\mathcal{L}^{\bot }}\mathcal{Z}$,
i.e., that $\limfunc{span}(\Pi _{\mathcal{L}^{\bot }}\tilde{\Sigma}\Pi _{%
\mathcal{L}^{\bot }})=\Pi _{\mathcal{L}^{\bot }}\limfunc{span}(\tilde{\Sigma}%
)$. But this follows since $\limfunc{span}(\Pi _{\mathcal{L}^{\bot }}\tilde{%
\Sigma}\Pi _{\mathcal{L}^{\bot }})=\limfunc{span}(\Pi _{\mathcal{L}^{\bot }}%
\tilde{\Sigma}^{1/2})=\Pi _{\mathcal{L}^{\bot }}\limfunc{span}(\tilde{\Sigma}%
^{1/2})=\Pi _{\mathcal{L}^{\bot }}\limfunc{span}(\tilde{\Sigma})$.

2. Suppose $\mathcal{S}$ is as in Part 2 of the lemma. Since $\mathcal{S}\in 
\mathbb{J}(\mathcal{L},\mathfrak{C})$ there exists a sequence $\Sigma
_{j}\in \mathfrak{C}$ such that $\Pi _{\mathcal{L}^{\bot }}\Sigma _{j}\Pi _{%
\mathcal{L}^{\bot }}/\left\Vert \Pi _{\mathcal{L}^{\bot }}\Sigma _{j}\Pi _{%
\mathcal{L}^{\bot }}\right\Vert $ converges to a singular matrix $\bar{\Sigma%
}$ with $\limfunc{rank}(\bar{\Sigma})<n-\dim (\mathcal{L})$ and such that $%
\mathcal{S}=\limfunc{span}(\bar{\Sigma})$ holds. By the assumption on $%
\mathfrak{C}$, we can find a subsequence $j_{i}$ along which $\Sigma
_{j_{i}} $ converges to a matrix $\Sigma ^{\ast }$. Note that $\Sigma ^{\ast
}$ must be singular, since otherwise $\Pi _{\mathcal{L}^{\bot }}\Sigma
_{j_{i}}\Pi _{\mathcal{L}^{\bot }}/\left\Vert \Pi _{\mathcal{L}^{\bot
}}\Sigma _{j_{i}}\Pi _{\mathcal{L}^{\bot }}\right\Vert $ would converge to
the matrix $\Pi _{\mathcal{L}^{\bot }}\Sigma ^{\ast }\Pi _{\mathcal{L}^{\bot
}}/\left\Vert \Pi _{\mathcal{L}^{\bot }}\Sigma ^{\ast }\Pi _{\mathcal{L}%
^{\bot }}\right\Vert $ which would have rank equal to $n-\dim (\mathcal{L})$%
, but at the same time would have to be equal to $\bar{\Sigma}$, which has
smaller rank. Hence, $\mathcal{Z}:=\limfunc{span}(\Sigma ^{\ast })$ is a
concentration space of $\mathfrak{C}$. Consider first the case where $%
\left\Vert \Pi _{\mathcal{L}^{\bot }}\Sigma ^{\ast }\Pi _{\mathcal{L}^{\bot
}}\right\Vert \neq 0$. Then we can conclude that $\bar{\Sigma}$ and $\Pi _{%
\mathcal{L}^{\bot }}\Sigma ^{\ast }\Pi _{\mathcal{L}^{\bot }}$ coincide up
to a positive proportionality factor. By construction of $\bar{\Sigma}$ we
have $\limfunc{span}(\bar{\Sigma})=\mathcal{S}$. Hence, $\limfunc{span}(\Pi
_{\mathcal{L}^{\bot }}\Sigma ^{\ast }\Pi _{\mathcal{L}^{\bot }})=\mathcal{S}$
holds. But the same argument as in the proof of Part 1 shows that $\limfunc{%
span}(\Pi _{\mathcal{L}^{\bot }}\Sigma ^{\ast }\Pi _{\mathcal{L}^{\bot
}})=\Pi _{\mathcal{L}^{\bot }}\limfunc{span}(\Sigma ^{\ast })$, which leads
to $\Pi _{\mathcal{L}^{\bot }}\mathcal{Z}=\Pi _{\mathcal{L}^{\bot }}\limfunc{%
span}(\Sigma ^{\ast })=\mathcal{S}$. Next consider the case where $%
\left\Vert \Pi _{\mathcal{L}^{\bot }}\Sigma ^{\ast }\Pi _{\mathcal{L}^{\bot
}}\right\Vert =0 $.\footnote{%
Note that this case cannot arise in case $\mathcal{L}=\left\{ 0\right\} $
and $\mathfrak{C}$ is bounded away from the zero matrix.} In this case we
can clearly conclude that $\mathcal{Z}=\limfunc{span}(\Sigma ^{\ast })$ is a
concentration space satisfying $\Pi _{\mathcal{L}^{\bot }}\mathcal{Z}%
=\left\{ 0\right\} $. The remaining claims follow for the sequence $\Sigma
_{j_{i}}$ just constructed. $\blacksquare $

The assumption that $\mathfrak{C}$ is bounded used in Part 2 of the lemma
can be made without much loss of generality if the tests one is interested
in are invariant under $G(\mathfrak{M}_{0})$, see Remark \ref{rescale}(ii).
The preceding lemma allows one in certain circumstances to reduce checking
the conditions on $\mathcal{S}\in \mathbb{J}(\mathcal{L},\mathfrak{C})$
postulated in Proposition \ref{L1} or in Proposition \ref{L2} (or the
corresponding sufficient conditions appearing in Remark \ref{sufff}(ii) or
Corollary \ref{C5}) to checking similar conditions that are expressed in
terms of the concentration spaces $\mathcal{Z}$. This can be advantageous
since the concentration spaces are sometimes easier to obtain than the
spaces $\mathcal{S}\in \mathbb{J}(\mathcal{L},\mathfrak{C})$. We illustrate
this in the following remark. However, it is important to note that for many
covariance models of interest this \textquotedblleft
reduction\textquotedblright\ trick does \emph{not} work and $\mathbb{J}(%
\mathcal{L},\mathfrak{C})$ has to be determined. For example, this is the
case for $\mathfrak{C}(\mathfrak{F}_{\mathrm{all}})$ and related covariance
models, necessitating the developments in Section \ref{struct}.

\begin{remark}
Let $\mathfrak{C}$ be a covariance model that is bounded and is bounded away
from the zero matrix. Assume that $T$, $N^{\dag }$, and $W$ are as in
Corollary \ref{C5} and assume furthermore that $\mathcal{L}\subseteq 
\limfunc{span}(X)$ (which is, in particular, the case if $\mathcal{V}%
=\left\{ 0\right\} $).

(i) Assume that $N^{\dag }$ is a finite or countable union of affine
subspaces. Suppose furthermore that we are in a scenario where we can show
that no concentration space $\mathcal{Z}$ is entirely contained in $\limfunc{%
span}(X)$, and thus no $\mathcal{Z}$ is entirely contained in $\mathcal{L}$.
Then the sufficient condition given at the end of Corollary \ref{C5}, namely
that if $\mathcal{S}\in \mathbb{J}(\mathcal{L},\mathfrak{C})$ then $\mu _{0}+%
\mathcal{S}\not\subseteq N^{\dag }$ for some $\mu _{0}\in \mathfrak{M}_{0}$
(and hence for all $\mu _{0}\in \mathfrak{M}_{0}$) must hold, is satisfied
whenever $\mu _{0}+\mathcal{Z}\not\subseteq N^{\dag }$ for some $\mu _{0}\in 
\mathfrak{M}_{0}$ (and hence for all $\mu _{0}\in \mathfrak{M}_{0}$) holds
for every concentration space $\mathcal{Z}$ (this is so since $N^{\dag }$ is
invariant under addition of elements from $\mathcal{L}$ and since every $%
\mathcal{S}\in \mathbb{J}(\mathcal{L},\mathfrak{C})$ must be of the form $%
\Pi _{\mathcal{L}^{\bot }}\mathcal{Z}$ in view of the preceding lemma and
our assumption on $\mathcal{Z}$). Hence, we can check the sufficient
condition in Corollary \ref{C5} without explicitly computing the spaces $%
\mathcal{S}\in \mathbb{J}(\mathcal{L},\mathfrak{C})$. Furthermore, in many
cases of interest we have $N^{\dag }=\limfunc{span}(X)$ (cf. the discussion
surrounding Lemma \ref{Fact2}), in which case the condition $\mu _{0}+%
\mathcal{Z}\not\subseteq N^{\dag }$ is then an automatic consequence of the
assumption that no concentration space $\mathcal{Z}$ is entirely contained
in $\limfunc{span}(X)$ (since $\limfunc{span}(X)$ is a linear space
containing $\mathfrak{M}_{0}$).

(ii) The example in (i) can be generalized a bit: Suppose now that for every
concentration space $\mathcal{Z}$ we either have that (a) it is not
contained in $\limfunc{span}(X)$, or (b) that $\mathcal{Z}$ is contained in $%
\mathcal{L}$ but for every sequence $\Sigma _{j}$ converging to some $\Sigma
^{\ast }$ satisfying $\mathcal{Z}=\limfunc{span}(\Sigma ^{\ast })$ the limit
points of $\Pi _{\mathcal{L}^{\bot }}\Sigma _{j}\Pi _{\mathcal{L}^{\bot
}}/\left\Vert \Pi _{\mathcal{L}^{\bot }}\Sigma _{j}\Pi _{\mathcal{L}^{\bot
}}\right\Vert $ are regular on $\mathcal{L}^{\bot }$. Then it suffices to
check that $\mu _{0}+\mathcal{Z}\not\subseteq N^{\dag }$ for some $\mu
_{0}\in \mathfrak{M}_{0}$ (and hence for all $\mu _{0}\in \mathfrak{M}_{0}$)
for every concentration space $\mathcal{Z}$ that is not contained in $%
\limfunc{span}(X)$ [this follows from the discussion in (i), since by
property (b) and the preceding lemma any $\mathcal{S}\in \mathbb{J}(\mathcal{%
L},\mathfrak{C})$ must be of the form $\Pi _{\mathcal{L}^{\bot }}\mathcal{Z}$%
]. Again if $N^{\dag }=\limfunc{span}(X)$, then this latter condition is
automatically satisfied. Of course, a sufficient condition for the
aforementioned limit points to be regular on $\mathcal{L}^{\bot }$ is the
following condition: for each relevant $\mathcal{Z}$, $\Sigma $, and
sequence $\Sigma _{j}$ the limit points of the sequence $\Pi _{\mathcal{L}%
^{\bot }}\Sigma _{j}\Pi _{\mathcal{L}^{\bot }}/\left\Vert \Pi _{\mathcal{L}%
^{\bot }}\Sigma _{j}\Pi _{\mathcal{L}^{\bot }}\right\Vert $ are regular on $%
\limfunc{span}(\Sigma ^{\ast })^{\bot }$.

(iii) Suppose now that we are in the situation of (ii) except that now $%
N^{\dag }$ is not a finite or countable union of affine subspaces. Then the
relevant sufficient condition in Corollary \ref{C5} can be shown to be
implied by the condition that if $\mathcal{Z}$ is a concentration space of
the covariance model $\mathfrak{C}$ that is not contained in $\limfunc{span}%
(X)$, then the set $N^{\dag }$ is a $\lambda _{\mu _{0}+\mathcal{Z}}$-null
set for some $\mu _{0}\in \mathfrak{M}_{0}$ (and hence for all $\mu _{0}\in 
\mathfrak{M}_{0}$). This follows since any $\mathcal{Z}$ with $\Pi _{%
\mathcal{L}^{\bot }}\mathcal{Z}=\mathcal{S}$ can be shown to be of the form
of a direct sum $\mathcal{A}\oplus \mathcal{B}$, where $\mathcal{A}$ is a
linear space that is linearly isomorphic to $\mathcal{S}$ and $\mathcal{B}$
is a linear subspace of $\mathcal{L}$, and since $N^{\dag }-\mu _{0}$ is of
the form $N^{\dag \dag }\oplus \mathcal{L}$ for an appropriate Borel-set $%
N^{\dag \dag }\subseteq \mathcal{L}^{\bot }$ and where the direct sum is in
fact an orthogonal sum. We omit the details. Similar arguments can be
applied if condition (b) in Proposition \ref{L2} is to be verified instead
of the sufficient condition just considered.

(iv) Similar arguments apply to the sufficient condition given in Remark \ref%
{sufff}(ii) or to condition (b) of Proposition \ref{L1}.
\end{remark}

\section{Appendix\label{App C}: Auxiliary results for Section \protect\ref%
{structure under Topelitz}}

In this appendix we provide results that will be used in the proofs of the
results of Section \ref{structure under Topelitz} that are provided in
Appendix \ref{App D}.

\begin{definition}
\label{Vdef}Let $\omega \in \lbrack 0,\pi ]$, $l\in \mathbb{Z}$, $m\in 
\mathbb{N}$, and let $s\geq 0$ be an integer. Define $E_{m,s}^{(l)}(\omega )$
as the $m\times 2$-dimensional matrix with $j$-th row equal to 
\begin{equation*}
\left( (j+l)^{s}\cos ((j+l)\omega ),(j+l)^{s}\sin ((j+l)\omega )\right) ,
\end{equation*}%
where we shall often drop the superscript $l$ in case $l=0$. For a positive
integer $p$, for $\underline{\omega }=(\omega _{1},\ldots ,\omega _{p})\in
\lbrack 0,\pi ]^{p}$, and for $\underline{d}=(d_{1},\ldots ,d_{p})\in 
\mathbb{N}^{p}$ we define the $m\times 2\sum_{i=1}^{p}d_{i}$-dimensional
matrix 
\begin{equation*}
V_{m}^{(l)}\left( \underline{\omega },\underline{d}\right) =\left(
E_{m,0}^{(l)}(\omega _{1}),\ldots ,E_{m,d_{1}-1}^{(l)}(\omega _{1}),\ldots
,E_{m,0}^{(l)}(\omega _{p}),\ldots ,E_{m,d_{p}-1}^{(l)}(\omega _{p})\right) .
\end{equation*}%
In case $p=1$, we shall often simply write $\omega $ for $\underline{\omega }
$ and $d$ for $\underline{d}$.
\end{definition}

For the following recall that $\kappa (\underline{\omega },\underline{d})$
has been defined in Section \ref{struct}, and that we use the convention
that $\underline{\omega }$ and $\underline{d}$ are the $0$-tupels for $p=0$.

\begin{lemma}
\label{fullrank}Let $p$ be a positive integer, let $\underline{\omega }\in
\lbrack 0,\pi ]^{p}$ have distinct coordinates, and let $\underline{d}\in 
\mathbb{N}^{p}$. Then for every positive integer $m$ and for every integer $%
l $ it holds that 
\begin{equation*}
\limfunc{rank}(V_{m}^{(l)}(\underline{\omega },\underline{d}))=\min
(m,\kappa (\underline{\omega },\underline{d})).
\end{equation*}
\end{lemma}

\textbf{Proof of Lemma \ref{fullrank}:} Standard results concerning linear
difference equations (e.g., \cite{Kelley_and_peterson}, Chp. 3) can be used
to verify that the collection made up of the functions $j\mapsto j^{s}\cos
(j\omega _{i})$ and $j\mapsto j^{s}\sin (j\omega _{i})$ (defined for $j\in 
\mathbb{Z}$) for $\omega _{i}\in (0,\pi )$ and $s=0,\ldots ,d_{i}-1$ as well
as of $j\mapsto j^{s}\cos (j\omega _{i})$ for $\omega _{i}\in \{0,\pi \}$
and $s=0,\ldots ,d_{i}-1$ forms a fundamental set of solutions of the
difference equation $\Delta _{\underline{\omega },\underline{d}}(z)w_{t}=0$
with $t\in \mathbb{Z}$, where we abuse the symbol $z$ also to denote the
backshift-operator. Hence any Casorati-matrix associated with this
fundamental set is non-singular. Observe that striking the columns
corresponding to $(j+l)^{s}\sin ((j+l)\omega _{i})$ for $\omega _{i}\in
\{0,\pi \}$ from the matrix $V_{\kappa (\underline{\omega },\underline{d}%
)}^{(l)}(\underline{\omega },\underline{d})$, results in a $\kappa (%
\underline{\omega },\underline{d})\times \kappa (\underline{\omega },%
\underline{d})$ matrix which is precisely a Casorati-matrix. This shows that 
$\limfunc{rank}(V_{\kappa (\underline{\omega },\underline{d})}^{(l)}(%
\underline{\omega },\underline{d}))=\kappa (\underline{\omega },\underline{d}%
)$. The claim is then an immediate consequence. $\blacksquare $

\begin{definition}
For a polynomial $\Theta (z)=1+\theta _{1}z+\ldots +\theta _{a}z^{a}$ of
degree $a$ we define for $m>a$ the $(m-a)\times m$ matrix $D_{m}(\Theta )$
via%
\begin{equation*}
\begin{pmatrix}
\theta _{a} & \ldots & \theta _{1} & 1 & 0 & \ldots & 0 \\ 
0 & \theta _{a} & \ldots & \theta _{1} & 1 & \ldots & 0 \\ 
\vdots &  &  &  &  &  & \vdots \\ 
0 & \ldots & 0 & \theta _{a} & \ldots & \theta _{1} & 1%
\end{pmatrix}%
\end{equation*}%
with the convention that $D_{m}(\Theta )=I_{m}$ in case $a=0$.
\end{definition}

\begin{remark}
\label{D-matrix}(i) Obviously, $D_{m}(\Theta )$ has full row-rank, i.e., $%
\limfunc{rank}(D_{m}(\Theta ))=m-a$.

(ii) Let $\Theta _{1}(z)$ and $\Theta _{2}(z)$ be polynomials of degree $%
a_{1}$ and $a_{2}$, respectively, satisfying $\Theta _{1}(0)=\Theta
_{2}(0)=1 $. Then for $m>a_{1}+a_{2}$ we have $D_{m}(\Theta _{1}\Theta
_{2})=D_{m-a_{1}}(\Theta _{2})D_{m}(\Theta _{1})=D_{m-a_{2}}(\Theta
_{1})D_{m}(\Theta _{2})$ where $D_{m}(\Theta _{1}\Theta _{2})$ denotes the
matrix associated with the polynomial $\Theta _{1}(z)\Theta _{2}(z)$.
\end{remark}

\begin{lemma}
\label{basis} Let $p$ be a positive integer, let $\underline{\omega }\in
\lbrack 0,\pi ]^{p}$ have distinct coordinates, and let $\underline{d}\in 
\mathbb{N}^{p}$. Then for every positive integer $m$ satisfying $m>\kappa (%
\underline{\omega },\underline{d})$, and for every integer $l$ we have that
the transposes of the row vectors of $D_{m}(\Delta _{\underline{\omega },%
\underline{d}})$ constitute a basis of $\limfunc{span}(V_{m}^{(l)}(%
\underline{\omega },\underline{d}))^{\bot }$. In particular, $\limfunc{span}%
(V_{m}^{(l)}(\underline{\omega },\underline{d}))$ does not depend on $l$.
\end{lemma}

\textbf{Proof of Lemma \ref{basis}:} Since the columns of $V_{m}^{(l)}(%
\underline{\omega },\underline{d})$ are either zero or segments of length $m$
(with $m>\kappa (\underline{\omega },\underline{d})$) from the fundamental
set of solutions to the difference equation $\Delta _{\underline{\omega },%
\underline{d}}(z)w_{t}=0$, we obviously have $D_{m}(\Delta _{\underline{%
\omega },\underline{d}})V_{m}^{(l)}(\underline{\omega },\underline{d})=0$.
This implies $\limfunc{span}(D_{m}^{\prime }(\Delta _{\underline{\omega },%
\underline{d}}))\subseteq \limfunc{span}(V_{m}^{(l)}(\underline{\omega },%
\underline{d}))^{\bot }$. Since the $m-\kappa (\underline{\omega },%
\underline{d})$ rows of $D_{m}(\Delta _{\underline{\omega },\underline{d}})$
are linearly independent, cf. Remark \ref{D-matrix}, and since $\limfunc{rank%
}(V_{m}^{(l)}(\underline{\omega },\underline{d}))=\kappa (\underline{\omega }%
,\underline{d})$ by Lemma \ref{fullrank}, the result follows. $\blacksquare $

\begin{remark}
\label{basis_1}In case $p=1$, and hence $\underline{\omega }=\omega \in
\lbrack 0,\pi ]$, and $\underline{d}=d\in \mathbb{N}$, the result in Lemma %
\ref{basis} reduces to the fact that the transposes of the row vectors of $%
D_{m}(\Delta _{\omega }^{d})$ constitute a basis of $\limfunc{span}%
((E_{m,0}^{(l)}(\omega ),\ldots ,E_{m,d-1}^{(l)}(\omega )))^{\bot }$. In
particular, $D_{m}(\Delta _{\omega }^{d})(E_{m,0}^{(l)}(\omega ),\ldots
,E_{m,d-1}^{(l)}(\omega ))=0$ holds.
\end{remark}

\begin{lemma}
\label{recrel} Let $\omega \in \lbrack 0,\pi ]$, $d\in \mathbb{N}$, and $%
l\in \mathbb{Z}$. If $m\in \mathbb{N}$ satisfies $m>\kappa (\omega ,1)$ then
we have%
\begin{equation}
D_{m}(\Delta _{\omega })E_{m,d}^{(l)}(\omega )=\left\{ 
\begin{array}{cc}
2\sum_{i=0}^{d-1}\binom{d}{i}\cos ^{(d-i)}(\omega )E_{m-\kappa (\omega
,1),i}^{(l+1)}(\omega )P^{i-d} & \text{for \ }\omega \in (0,\pi ) \\ 
\cos (\omega )\sum_{i=0}^{d-1}\binom{d}{i}E_{m-\kappa (\omega
,1),i}^{(l)}(\omega ) & \text{for \ }\omega \in \{0,\pi \}%
\end{array}%
\right. ,  \label{diff_1a}
\end{equation}%
where $\cos ^{(i)}(\omega )$ denotes the $i$-th order derivative of the
cosine function and where $P$ is the $2\times 2$-dimensional orthogonal
matrix with first row $(0,1)$ and second row $(-1,0)$. If $m\in \mathbb{N}$
satisfies $m>\kappa (\omega ,d)$ then we have%
\begin{equation}
D_{m}(\Delta _{\omega }^{d})E_{m,d}^{(l)}(\omega )=\left\{ 
\begin{array}{cc}
2^{d}d!(-\sin (\omega ))^{d}E_{m-\kappa (\omega ,d),0}^{(l+d)}(\omega )P^{-d}
& \text{for \ }\omega \in (0,\pi ) \\ 
d!(\cos (\omega ))^{d}E_{m-\kappa (\omega ,d),0}^{(l)}(\omega ) & \text{for
\ }\omega \in \{0,\pi \}%
\end{array}%
\right. .  \label{diff_1b}
\end{equation}
\end{lemma}

\textbf{Proof of Lemma \ref{recrel}:} From Lemma \ref{basis} and Remark \ref%
{basis_1} we obtain that the identity $D_{m}(\Delta _{\omega
})E_{m,0}^{(l)}(\omega )=0$ for every $\omega \in \lbrack 0,\pi ]$ holds.
Consider first the case where $\omega \in (0,\pi )$: Since the left-hand
side of this identity is a smooth function of $\omega $ in this range, we
may differentiate the identity $d$ times leading to 
\begin{equation*}
\sum_{i=0}^{d}\binom{d}{i}\left( \frac{d^{d-i}}{d\omega ^{d-i}}D_{m}(\Delta
_{\omega })\right) \left( \frac{d^{i}}{d\omega ^{i}}E_{m,0}^{(l)}(\omega
)\right) =0.
\end{equation*}%
Rearranging terms and computing the derivatives gives 
\begin{equation*}
D_{m}(\Delta _{\omega })E_{m,d}^{(l)}(\omega )P^{d}=2\sum_{i=0}^{d-1}\binom{d%
}{i}\cos ^{(d-i)}(\omega )[0,I_{m-\kappa (\omega ,1)},0]E_{m,i}^{(l)}(\omega
)P^{i}.
\end{equation*}%
But clearly $[0,I_{m-\kappa (\omega ,1)},0]E_{m,i}^{(l)}(\omega
)=E_{m-\kappa (\omega ,1),i}^{(l+1)}(\omega )$ holds, from which we obtain (%
\ref{diff_1a}) in case $\omega \in (0,\pi )$. Next, consider the case where $%
\omega \in \{0,\pi \}$ holds: Using the Binomial formula, together with $%
\omega \in \{0,\pi \}$, we have 
\begin{equation*}
(r+1)^{d}\cos ((r+1)\omega )-\cos (\omega )r^{d}\cos (r\omega )=\cos (\omega
)\sum_{i=0}^{d-1}\binom{d}{i}r^{i}\cos (r\omega ).
\end{equation*}%
Since the second column of $E_{m,i}^{(l)}(\omega )$ is $0$ for any $i$ and $%
m $, the claim (\ref{diff_1a}) then also follows in case $\omega \in \{0,\pi
\} $. We next prove (\ref{diff_1b}) by induction over $d$. In case $d=1$,
the claim holds as it reduces to (\ref{diff_1a}). Suppose that the induction
hypothesis now holds for some $d\geq 1$ (and any $m>\kappa (\omega ,d)$).
Then for any $m\in \mathbb{N}$ satisfying $m>\kappa (\omega ,d+1)$ we have,
using Remark \ref{D-matrix}(ii) and (\ref{diff_1a}) with $d$ replaced by $%
d+1 $,%
\begin{eqnarray*}
&&D_{m}(\Delta _{\omega }^{d+1})E_{m,d+1}^{(l)}(\omega )=D_{m-\kappa (\omega
,1)}(\Delta _{\omega }^{d})D_{m}(\Delta _{\omega })E_{m,d+1}^{(l)}(\omega )
\\
&=&\left\{ 
\begin{array}{cc}
2\sum_{i=0}^{d}\binom{d+1}{i}\cos ^{(d+1-i)}(\omega )D_{m-\kappa (\omega
,1)}(\Delta _{\omega }^{d})E_{m-\kappa (\omega ,1),i}^{(l+1)}(\omega
)P^{i-d-1} & \text{for \ }\omega \in (0,\pi ) \\ 
\cos (\omega )\sum_{i=0}^{d}\binom{d+1}{i}D_{m-\kappa (\omega ,1)}(\Delta
_{\omega }^{d})E_{m-\kappa (\omega ,1),i}^{(l)}(\omega ) & \text{for \ }%
\omega \in \{0,\pi \}%
\end{array}%
\right. .
\end{eqnarray*}%
Observe that $m-\kappa (\omega ,1)>\kappa (\omega ,d)$ holds, since $%
m>\kappa (\omega ,d+1)$ and $\kappa (\omega ,d+1)=\kappa (\omega ,1)+\kappa
(\omega ,d)$. Hence we may apply Lemma \ref{basis} and Remark \ref{basis_1}
to obtain $D_{m-\kappa (\omega ,1)}(\Delta _{\omega }^{d})E_{m-\kappa
(\omega ,1),i}^{(l^{\prime })}(\omega )=0$ for $i<d$ and for $l^{\prime }=l$
or $l^{\prime }=l+1$. We thus obtain from the preceding display that 
\begin{equation*}
D_{m}(\Delta _{\omega }^{d+1})E_{m,d+1}^{(l)}(\omega )=\left\{ 
\begin{array}{cc}
2(d+1)\cos ^{(1)}(\omega )D_{m-\kappa (\omega ,1)}(\Delta _{\omega
}^{d})E_{m-\kappa (\omega ,1),d}^{(l+1)}(\omega )P^{-1} & \text{for \ }%
\omega \in (0,\pi ) \\ 
(d+1)\cos (\omega )D_{m-\kappa (\omega ,1)}(\Delta _{\omega
}^{d})E_{m-\kappa (\omega ,1),d}^{(l)}(\omega ) & \text{for \ }\omega \in
\{0,\pi \}%
\end{array}%
\right. .
\end{equation*}%
Together with the induction hypothesis (applied with $m$ replaced by $%
m-\kappa (\omega ,1)$) this establishes (\ref{diff_1b}) for $d+1$. $%
\blacksquare $

\begin{lemma}
\label{recrel2} Let $(\omega _{1},\omega _{2})\in \lbrack 0,\pi ]^{2}$, $%
l\in \mathbb{Z}$, $m\in \mathbb{N}$, and $d\in \mathbb{N}$. Assume $m>\kappa
(\omega _{1},d)$ holds. Then%
\begin{equation*}
D_{m}(\Delta _{\omega _{1}}^{d})E_{m,0}^{(l)}(\omega _{2})=E_{m-\kappa
(\omega _{1},d),0}^{(l)}(\omega _{2})(A(\omega _{1},\omega _{2}))^{d},
\end{equation*}%
where $A(\omega _{1},\omega _{2})=2(\cos (\omega _{2})-\cos (\omega
_{1}))P(\omega _{2})$ when $\omega _{1}\in (0,\pi )$, and $A(\omega
_{1},\omega _{2})=P(\omega _{2})-\cos (\omega _{1})I_{2}$ when $\omega
_{1}\in \{0,\pi \}$. The matrices $A(\omega _{1},\omega _{2})$ are multiples
of orthogonal matrices; they are nonsingular if $\omega _{1}\neq \omega _{2}$
and equal the zero matrix otherwise. [Here $P(\omega )$ denotes the $2\times
2$-dimensional orthogonal matrix with first row $(\cos (\omega ),\sin
(\omega ))$ and second row $(-\sin (\omega ),\cos (\omega ))$.]
\end{lemma}

\textbf{Proof of Lemma \ref{recrel2}:} Consider first the case where $d=1$.
We start with the following standard trigonometric identities for $j\in 
\mathbb{Z}$, which easily follow from the angle addition formulas,%
\begin{align}
\cos ((j+2)\omega _{2})-2\cos (\omega _{2})\cos ((j+1)\omega _{2})+\cos
(j\omega _{2})& =0  \label{trigid1} \\
\sin ((j+2)\omega _{2})-2\cos (\omega _{2})\sin ((j+1)\omega _{2})+\sin
(j\omega _{2})& =0.  \label{trigid2}
\end{align}%
Consider first the case where $\omega _{1}\in (0,\pi )$ holds. From (\ref%
{trigid1}) and (\ref{trigid2}) it follows that 
\begin{align}
& \cos ((j+2)\omega _{2})-2\cos (\omega _{1})\cos ((j+1)\omega _{2})+\cos
(j\omega _{2})=2\cos ((j+1)\omega _{2})(\cos (\omega _{2})-\cos (\omega
_{1}))  \label{rec} \\
& \sin ((j+2)\omega _{2})-2\cos (\omega _{1})\sin ((j+1)\omega _{2})+\sin
(j\omega _{2})=2\sin ((j+1)\omega _{2})(\cos (\omega _{2})-\cos (\omega
_{1})),  \label{rec*}
\end{align}%
and thus (by the angle addition formulas) 
\begin{equation*}
D_{m}(\Delta _{\omega _{1}})E_{m,0}^{(l)}(\omega _{2})=2(\cos (\omega
_{2})-\cos (\omega _{1}))E_{m-\kappa (\omega _{1},1),0}^{(l)}(\omega
_{2})P(\omega _{2})=E_{m-\kappa (\omega _{1},1),0}^{(l)}(\omega
_{2})A(\omega _{1},\omega _{2}).
\end{equation*}%
Next consider the case where $\omega _{1}\in \{0,\pi \}$. It is then easy to
see, using the angle addition formulas, that 
\begin{equation*}
D_{m}(\Delta _{\omega _{1}})E_{m,0}^{(l)}(\omega _{2})=E_{m-\kappa (\omega
_{1},1),0}^{(l)}(\omega _{2})(P(\omega _{2})-\cos (\omega
_{1})I_{2})=E_{m-\kappa (\omega _{1},1),0}^{(l)}(\omega _{2})A(\omega
_{1},\omega _{2}).
\end{equation*}%
The case $d>1$ now follows from (cf. Remark \ref{D-matrix}) 
\begin{equation*}
D_{m}(\Delta _{\omega _{1}}^{d})=D_{m-\kappa (\omega _{1},d-1)}(\Delta
_{\omega _{1}})\cdots D_{m-\kappa (\omega _{1},2)}(\Delta _{\omega
_{1}})D_{m-\kappa (\omega _{1},1)}(\Delta _{\omega _{1}})D_{m}(\Delta
_{\omega _{1}}),
\end{equation*}%
together with a repeated application of the already established result for $%
d=1$. $\blacksquare $

Recall that the finite and symmetric Borel measures on $[-\pi ,\pi ]$ are
precisely the spectral measures of real weakly stationary processes.

\begin{definition}
\label{spectral}For $\mathsf{m}$ a finite and symmetric Borel measure on $%
[-\pi ,\pi ]$ (symmetry here meaning $\mathsf{m}(A)=\mathsf{m}(-A)$ for
every Borel subset $A$ of $[-\pi ,\pi ]$) and for $m\in \mathbb{N}$ we
define the $m\times m$ matrix%
\begin{equation*}
\Sigma (\mathsf{m},m)=\left[ \int_{-\pi }^{\pi }e^{-\iota \nu (j-j^{\prime
})}d\mathsf{m}(\nu )\right] _{j,j^{\prime }=1}^{m}.
\end{equation*}%
For a spectral density $f$, i.e., an even $\lambda _{\lbrack -\pi ,\pi ]}$%
-integrable function from $[-\pi ,\pi ]$ to $[0,\infty )$, we denote by $%
\mathsf{m}_{f}$ the (finite and symmetric) Borel measure on $[-\pi ,\pi ]$
with density $f$ (w.r.t. Lebesgue measure $\lambda _{\lbrack -\pi ,\pi ]}$
on $[-\pi ,\pi ]$) and we abbreviate $\Sigma (\mathsf{m}_{f},m)$ to $\Sigma
(f,m)$. Finally, for $\mathsf{m}$ a finite and symmetric Borel measure on $%
[-\pi ,\pi ]$ and $\Theta $ a polynomial, we denote by $\Theta \odot \mathsf{%
m}$ the (finite and symmetric) Borel measure given by $(\Theta \odot \mathsf{%
m})(A)=\int_{A}\left\vert \Theta (e^{\iota \nu })\right\vert ^{2}d\mathsf{m}%
(\nu )$ for Borel sets $A\subseteq \lbrack -\pi ,\pi ]$.
\end{definition}

Given a nonempty set of spectral densities $\mathfrak{F}$ and $m\in \mathbb{N%
}$, we shall write $\mathfrak{C}(\mathfrak{F},m)=\{\Sigma (f,m):f\in 
\mathfrak{F}\}$. In case $m=n$, where $n$ is the sample size, we shall -- in
line with the notation introduced in Section \ref{autocorr} -- often simply
write $\Sigma (f)$ for $\Sigma (f,n)$ and $\mathfrak{C}(\mathfrak{F})$ for $%
\mathfrak{C}(\mathfrak{F},n)$. The next lemma is an immediate consequence of
a standard result concerning linear filters (e.g., \cite{rozanov}, Chp I.8)
applied to the linear filter $\Delta _{\underline{\omega },\underline{d}}(z)$
operating on a stationary process with spectral measure $\mathsf{m}$.

\begin{lemma}
\label{Msdgeneral}Let $p$ be a positive integer, let $\underline{\omega }\in
\lbrack 0,\pi ]^{p}$, let $\underline{d}\in \mathbb{N}^{p}$, and let $m\in 
\mathbb{N}$ satisfy $m>\kappa (\omega ,d)$. Let $\mathsf{m}$ be a finite and
symmetric Borel measure on $[-\pi ,\pi ]$. Then%
\begin{equation*}
D_{m}(\Delta _{\underline{\omega },\underline{d}})\Sigma (\mathsf{m}%
,m)D_{m}^{\prime }(\Delta _{\underline{\omega },\underline{d}})=\Sigma
(\Delta _{\underline{\omega },\underline{d}}\odot \mathsf{m},m-\kappa (%
\underline{\omega },\underline{d})).
\end{equation*}
\end{lemma}

The following result is well-known, but difficult to pinpoint in the
literature in this form. The first claim of the lemma can, e.g., be found in
Theorem 2.6 of \cite{krein}. This is closely related to a theorem of Carath%
\'{e}odory (Section 4.1 of \cite{gre_and_szeg}), on which we have chosen to
base the proof of the lemma.

\begin{lemma}
\label{Toeplitz}Let $\Phi $ be a real nonnegative definite symmetric
Toeplitz matrix of dimension $m\times m$ with $m\in \mathbb{N}$. Then there
exists a finite and symmetric Borel measure $\mathsf{m}$ on $[-\pi ,\pi ]$
such that $\Phi =\Sigma (\mathsf{m},m)$. If $\Phi $ is singular, the measure 
$\mathsf{m}$ is unique. Furthermore, if $\Phi $ is singular and $\Phi \neq 0$%
, $\mathsf{m}$ is of the form $\sum_{i=1}^{p}c_{i}(\delta _{-\omega
_{i}}+\delta _{\omega _{i}})$ for some $p\in \mathbb{N}$, for some positive $%
c_{i}$, $1\leq i\leq p$, and for some $\underline{\omega }=(\omega
_{1},\ldots ,\omega _{p})\in \lbrack 0,\pi ]^{p}$ with $\omega _{1}<\ldots
<\omega _{p}$ such that $1\leq \sum_{i=1}^{p}\kappa (\omega _{i},1)<m$; if $%
\Phi =0$, the measure $\mathsf{m}$ is the zero measure.
\end{lemma}

\textbf{Proof of Lemma \ref{Toeplitz}:} We start with a preparatory remark:
Recall that for any finite and symmetric Borel measure $\mathsf{m}$ the
matrix $\Sigma (\mathsf{m},m)$ is nonnegative definite, since $x^{\prime
}\Sigma (\mathsf{m},m)x=\int_{-\pi }^{\pi }\left\vert x(e^{\iota \nu
})\right\vert ^{2}d\mathsf{m}(\nu )\geq 0$ for every $x\in \mathbb{R}^{m}$,
where $x(e^{\iota \nu })=\sum_{j=1}^{m}x_{j}e^{\iota j\nu }$. If $\Sigma (%
\mathsf{m},m)$ is singular for some $\mathsf{m}$, then $0=x_{\ast }^{\prime
}\Sigma (\mathsf{m},m)x_{\ast }=\int_{-\pi }^{\pi }\left\vert x_{\ast
}(e^{\iota \nu })\right\vert ^{2}d\mathsf{m}(\nu )$ for some nonzero $%
x_{\ast }\in \mathbb{R}^{m}$ must hold. Consequently, the support of any
such measure $\mathsf{m}$ must be contained in the zero-set of the
trigonometric polynomial $\left\vert x_{\ast }(e^{\iota \nu })\right\vert
^{2}$ (which has degree at most $m-1$). Hence, $\mathsf{m}$ must be the zero
measure or must be of the form $\sum_{i=1}^{p}c_{i}(\delta _{-\omega
_{i}}+\delta _{\omega _{i}})$ for some $p\in \mathbb{N}$, for positive $%
c_{i} $'s, $1\leq i\leq p$, and for some $\underline{\omega }=(\omega
_{1},\ldots ,\omega _{p})\in \lbrack 0,\pi ]^{p}$ with $\omega _{1}<\ldots
<\omega _{p}$; obviously, $\sum_{i=1}^{p}\kappa (\omega _{i},1)<m$ must also
hold.

Now, if $\Phi =0$, the zero measure satisfies $\Phi =\Sigma (\mathsf{m},m)$
and obviously this is the only possible choice. Next consider the case where 
$\Phi $ is singular, but $\Phi \neq 0$. Then $m\geq 2$ must hold. In the
following let $\Phi (\left\vert i-j\right\vert )$ denote the $(i,j)$-th
element of $\Phi $. Consider first the case $m=2$: Then $\left\vert \Phi
(1)\right\vert =\Phi (0)$ has to hold. Consequently, $\Phi =\Sigma (\mathsf{m%
},m)$ holds for $\mathsf{m}=(\Phi (0)/2)(\delta _{-\omega _{1}}+\delta
_{\omega _{1}})$ where $\omega _{1}=0$ or $\pi $, depending on whether $\Phi
(1)$ is positive or negative (here $p=1$ and $\kappa (\omega _{1},1)=1<m=2$
is satisfied). If $\mathsf{m}^{\prime }$ is another finite and symmetric
Borel measure satisfying $\Phi =\Sigma (\mathsf{m}^{\prime },m)$, then the
preparatory remark shows that $\mathsf{m}^{\prime }$ must be a discrete
measure of the form as in the preparatory remark with $\sum_{i=1}^{p}\kappa
(\omega _{i},1)<m=2$. But this shows that $p=1$ and that $\omega _{1}=0$ or $%
\pi $. Uniqueness then follows immediately. We next turn to the case $m>2$.
Observe that not all $\Phi (r)$ for $r=1,\ldots ,m-1$ can be zero, since $%
\Phi \neq 0$, $\Phi $ is singular, and is Toeplitz. It now follows from a
theorem of Carath\'{e}odory (Section 4.1 of \cite{gre_and_szeg}) that%
\begin{equation}
\Phi (r)=\dint_{-\pi }^{\pi }e^{-\iota \nu r}d\mathsf{m}(\nu )\text{ \ \ \
for every }r=1,\ldots ,m-1  \label{cara}
\end{equation}%
holds for some measure $\mathsf{m}$ of the form $\sum_{i=1}^{p}c_{i}(\delta
_{-\omega _{i}}+\delta _{\omega _{i}})$ with $p\geq 1$, $c_{i}>0$, $%
\underline{\omega }=(\omega _{1},\ldots ,\omega _{p})\in \lbrack 0,\pi ]^{p}$
with $\omega _{1}<\ldots <\omega _{p}$, and with $\sum_{i=1}^{p}\kappa
(\omega _{i},1)<m$; it furthermore follows from that theorem that $p$ and
the constants $c_{i}$, $\omega _{i}$ are uniquely determined. [The theorem
in Section 4.1 of \cite{gre_and_szeg} is given for complex $\Phi (r)$ and
shows that (\ref{cara}) holds for a Borel measure of the form $%
\sum_{j=1}^{q}a_{j}\delta _{\varpi _{j}}$, where $1\leq q<m$, $a_{j}>0$, and
where $\varpi _{j}$ are distinct elements in the half-open interval $(-\pi
,\pi ]$. Furthermore, $q$ and the constants $a_{j}$, $\varpi _{j}$ are
uniquely determined. Exploiting that $\Phi (r)=\overline{\Phi (r)}$ in our
context, the just mentioned uniqueness immediately shows that those $\varpi
_{j}$'s, which are different from $0$ or $\pi $, must appear in pairs
symmetrically located around zero, and that $a_{j}=a_{j^{\prime }}$ must
hold if $\varpi _{j}=-\varpi _{j^{\prime }}$. It is now not difficult to see
that one can replace the Borel measure $\sum_{j=1}^{q}a_{j}\delta _{\varpi
_{j}}$ by an appropriate $\mathsf{m}$ of the form as given above (note that $%
e^{-\iota \nu r}=e^{\iota \nu r}$ for $\nu =0,\pi $) and that uniqueness of $%
q$, the $a_{j}$'s and the $\varpi _{j}$'s translates into uniqueness of $p$,
the $c_{i}$'s, and $\underline{\omega }\in \lbrack 0,\pi ]^{p}$. Also note
that $q<m$ translates into $\sum_{i=1}^{p}\kappa (\omega _{i},1)<m$.] It
thus suffices to show that $\Phi (0)=\int_{-\pi }^{\pi }d\mathsf{m}(\nu )$,
uniqueness of $\mathsf{m}$ then already following from the preparatory
remark, together with the uniqueness part of the cited theorem. Now, because
of $\sum_{i=1}^{p}\kappa (\omega _{i},1)<m$, we can find a vector $x_{0}\in 
\mathbb{R}^{m}$ such that $x_{0}(e^{\iota \nu })$ vanishes at $\nu =\omega
_{i}$ and at $\nu =-\omega _{i}$ for every coordinate $\omega _{i}$ of $%
\underline{\omega }$. [Just set $x_{0}(e^{\iota \nu })$ equal to $%
\prod_{i:\omega _{i}\in (0,\pi )}(1-e^{-\iota \omega _{i}}e^{\iota \nu
})(1-e^{\iota \omega _{i}}e^{\iota \nu })\prod_{i:\omega _{i}\in \{0,\pi
\}}(1-e^{-\iota \omega _{i}}e^{\iota \nu })$ and observe that the
coefficients are real.] This implies $x_{0}^{\prime }\Sigma (\mathsf{m}%
,m)x_{0}=0$, and hence that $\Sigma (\mathsf{m},m)$ is singular. Obviously, $%
\Sigma (\mathsf{m},m)=\Phi +(\int_{-\pi }^{\pi }d\mathsf{m}(\nu )-\Phi
(0))I_{m}$ holds. We thus obtain%
\begin{equation*}
0=\inf_{x^{\prime }x=1}x^{\prime }\Sigma (\mathsf{m},m)x=\inf_{x^{\prime
}x=1}x^{\prime }\Phi x+\dint_{-\pi }^{\pi }d\mathsf{m}(\nu )-\Phi
(0)=\dint_{-\pi }^{\pi }d\mathsf{m}(\nu )-\Phi (0),
\end{equation*}%
since $\inf_{x^{\prime }x=1}x^{\prime }\Phi x=0$ in view of singularity of $%
\Phi $. This completes the proof for singular $\Phi $. Finally, if $\Phi $
is positive definite, write $\Phi $ as $\Phi _{\ast }+cI_{m}$ where $c>0$
and $\Phi _{\ast }$ is nonnegative definite and singular. Obviously, $\Phi
_{\ast }$ is symmetric and Toeplitz. Hence $\Phi _{\ast }=\Sigma (\mathsf{m}%
_{\ast },m)$ for some measure $\mathsf{m}_{\ast }$ satisfying the conditions
in the theorem. Setting $\mathsf{m}=\mathsf{m}_{\ast }+(c/2\pi )\lambda
_{\lbrack -\pi ,\pi ]}$ completes the proof. $\blacksquare $

\begin{remark}
\label{Toeplitz_rem}(i) While the measure $\mathsf{m}$ in Lemma \ref%
{Toeplitz} is unique in case $\Phi $ is singular, it is never unique if $%
\Phi $ is positive definite. As shown in the proof, if $\Phi $ is positive
definite, one such measure is given by the sum of a discrete measure and the
spectral measure of white noise. As is well-known (see, e.g., Section 3.9.2
of \cite{stoica}), for positive definite $\Phi $, another representation $%
\Phi =\Sigma (\mathsf{m},m)$ can be found where $\mathsf{m}$ is the spectral
measure of an appropriate stationary autoregressive process of order at most 
$m-1$ (and thus is absolutely continuous w.r.t. Lebesgue measure).

(ii) An alternative proof of Lemma \ref{Toeplitz} can be based on the just
mentioned autoregressive representation result by applying a limiting
argument to cover also the case of singular $\Phi $.

(iii) As a converse to the last claim in Lemma \ref{Toeplitz} we have: If $%
\mathsf{m}=0$ or if $\mathsf{m}$ is of the form $\sum_{i=1}^{p}c_{i}(\delta
_{-\omega _{i}}+\delta _{\omega _{i}})$ for some $p\in \mathbb{N}$, for some
positive $c_{i}$, $1\leq i\leq p$, and for some $\underline{\omega }=(\omega
_{1},\ldots ,\omega _{p})\in \lbrack 0,\pi ]^{p}$ with $\omega _{1}<\ldots
<\omega _{p}$ such that $1\leq \sum_{i=1}^{p}\kappa (\omega _{i},1)<m$
holds, then $\Sigma (\mathsf{m},m)$ is singular. This follows since $%
x_{0}^{\prime }\Sigma (\mathsf{m},m)x_{0}=0$ for $x_{0}$ as defined in the
proof.

(iv) The lemma is somewhat more general than what is needed later in the
paper, but its extra generality is useful in other contexts.
\end{remark}

\section{Appendix\label{App D}: Proofs for Sections \protect\ref{struct} and 
\protect\ref{szctrl_3}}

Recall that by our conventions $\kappa (\underline{\omega }(\mathcal{L}),%
\underline{d}(\mathcal{L}))=0$ if $p(\mathcal{L})=0$.

\begin{lemma}
\label{dimsubspaces}Let $\mathcal{L}$ be a linear subspace of $\mathbb{R}%
^{n} $ with $\dim (\mathcal{L})<n$. Then $\kappa (\underline{\omega }(%
\mathcal{L}),\underline{d}(\mathcal{L}))\leq \dim (\mathcal{L})$ holds.
\end{lemma}

\textbf{Proof of Lemma \ref{dimsubspaces}:} If $p(\mathcal{L})=0$, there is
nothing to prove in view of the above convention. Suppose now that $p(%
\mathcal{L})>0$. By construction, $\limfunc{span}(V_{n}^{(0)}(\underline{%
\omega }(\mathcal{L}),\underline{d}(\mathcal{L})))\subseteq \mathcal{L}$.
From Lemma \ref{fullrank} we have $\limfunc{rank}(V_{n}^{(0)}(\underline{%
\omega }(\mathcal{L}),\underline{d}(\mathcal{L})))=\min (n,\kappa (%
\underline{\omega }(\mathcal{L}),\underline{d}(\mathcal{L})))$.
Consequently, $\min (n,\kappa (\underline{\omega }(\mathcal{L}),\underline{d}%
(\mathcal{L})))\leq \dim (\mathcal{L})<n$ must hold, which obviously proves
the desired result. $\blacksquare $

\begin{definition}
\label{Hdef} For a positive integer $p$, $\underline{\omega }\in \lbrack
0,\pi ]^{p}$, $\underline{d}\in \mathbb{N}^{p}$, and for $n$ satisfying $%
n>\kappa (\underline{\omega },\underline{d})$ we define the $n\times
(n-\kappa (\underline{\omega },\underline{d}))$-dimensional matrix 
\begin{equation*}
H_{n}(\underline{\omega },\underline{d}):=D_{n}^{\prime }(\Delta _{%
\underline{\omega },\underline{d}})\left( D_{n}(\Delta _{\underline{\omega },%
\underline{d}})D_{n}^{\prime }(\Delta _{\underline{\omega },\underline{d}%
})\right) ^{-1}.
\end{equation*}
\end{definition}

Observe that the inverse in the preceding display exists in view of Remark %
\ref{D-matrix}.

\begin{lemma}
\label{Mspectral} Let $\mathcal{L}$ be a linear subspace of $\mathbb{R}^{n}$
with $\dim (\mathcal{L})<n$. Furthermore, let $\mathsf{m}$ be a finite
symmetric Borel measure on $[-\pi ,\pi ]$. Then 
\begin{equation}
\Pi _{\mathcal{L}^{\bot }}\Sigma (\mathsf{m},n)\Pi _{\mathcal{L}^{\bot
}}=\Pi _{\mathcal{L}^{\bot }}H_{n}(\underline{\omega }(\mathcal{L}),%
\underline{d}(\mathcal{L}))\Sigma (\Delta _{\underline{\omega }(\mathcal{L}),%
\underline{d}(\mathcal{L})}\odot \mathsf{m},n-\kappa (\underline{\omega }(%
\mathcal{L}),\underline{d}(\mathcal{L})))H_{n}^{\prime }(\underline{\omega }(%
\mathcal{L}),\underline{d}(\mathcal{L}))\Pi _{\mathcal{L}^{\bot }}.
\label{project}
\end{equation}
\end{lemma}

\textbf{Proof of Lemma \ref{Mspectral}:} If $p(\mathcal{L})=0$, we have $%
\Delta _{\underline{\omega }(\mathcal{L}),\underline{d}(\mathcal{L})}=1$ and 
$\kappa (\underline{\omega }(\mathcal{L}),\underline{d}(\mathcal{L}))=0$ by
our conventions. Since then clearly $H_{n}(\underline{\omega }(\mathcal{L}),%
\underline{d}(\mathcal{L}))=I_{n}$ holds in view of the definition of $%
D_{n}(\Delta _{\underline{\omega },\underline{d}})$, there is nothing to
prove in this case. Assume that $p(\mathcal{L})>0$. By Lemma \ref%
{dimsubspaces} we have that $n>\kappa (\underline{\omega }(\mathcal{L}),%
\underline{d}(\mathcal{L}))$ holds. Lemma \ref{basis} shows that the columns
of $D_{n}^{\prime }(\Delta _{\underline{\omega }(\mathcal{L}),\underline{d}(%
\mathcal{L})})$ constitute a basis of $\limfunc{span}(V_{n}^{(0)}(\underline{%
\omega }(\mathcal{L}),\underline{d}(\mathcal{L})))^{\bot }$, which contains $%
\mathcal{L}^{\bot }$ in view of Definition \ref{subspaces_2}. Hence, $\Pi _{%
\mathcal{L}^{\bot }}=\Pi _{\mathcal{L}^{\bot }}H_{n}(\underline{\omega }(%
\mathcal{L}),\underline{d}(\mathcal{L}))D_{n}(\Delta _{\underline{\omega }(%
\mathcal{L}),\underline{d}(\mathcal{L})})$ holds. Inserting this into the
l.h.s. of (\ref{project}) and applying Lemma \ref{Msdgeneral} completes the
proof. $\blacksquare $

\begin{lemma}
\label{proj} Let $\mathcal{L}$ be a linear subspace of $\mathbb{R}^{n}$ with 
$\dim (\mathcal{L})<n$. Let $\omega \in \lbrack 0,\pi ]$. Then there exists
a $2\times 2$-dimensional regular matrix $B(\mathcal{L},\omega )$ that is
proportional to an orthogonal matrix such that 
\begin{equation}
\Pi _{\mathcal{L}^{\bot }}H_{n}(\underline{\omega }(\mathcal{L}),\underline{d%
}(\mathcal{L}))E_{n-\kappa (\underline{\omega }(\mathcal{L}),\underline{d}(%
\mathcal{L})),0}(\omega )=%
\begin{cases}
\Pi _{\mathcal{L}^{\bot }}E_{n,d_{i}(\mathcal{L})}(\omega )B(\mathcal{L}%
,\omega ) & \text{ if }\omega =\omega _{i}(\mathcal{L})\text{ for some }i \\ 
\Pi _{\mathcal{L}^{\bot }}E_{n,0}(\omega )B(\mathcal{L},\omega ) & \text{
else}.%
\end{cases}
\label{stat0}
\end{equation}%
Furthermore, $\Pi _{\mathcal{L}^{\bot }}E_{n,d_{i}(\mathcal{L})}(\omega _{i}(%
\mathcal{L}))B(\mathcal{L},\omega _{i}(\mathcal{L}))\neq 0$ for every $%
i=1,\ldots ,p(\mathcal{L})$, and $\Pi _{\mathcal{L}^{\bot }}E_{n,0}(\omega
)B(\mathcal{L},\omega )\neq 0$ for every $\omega $ not equal to some
coordinate $\omega _{i}(\mathcal{L})$ of $\underline{\omega }(\mathcal{L})$.
\end{lemma}

\textbf{Proof of Lemma \ref{proj}:} If $p(\mathcal{L})=0$, the result is
trivial with $B(\mathcal{L},\omega )$ the identity matrix. Hence, suppose $p(%
\mathcal{L})>0$ holds. We start with the case where $\omega $ does not
coincide with a coordinate of $\underline{\omega }(\mathcal{L})$. Because of
Remark \ref{D-matrix} we may write%
\begin{equation*}
D_{n}(\Delta _{\underline{\omega }(\mathcal{L}),\underline{d}(\mathcal{L}%
)})=\dprod\limits_{i=1}^{p(\mathcal{L})}D_{n-g_{i}}(\Delta _{\omega _{i}(%
\mathcal{L})}^{d_{i}(\mathcal{L})})
\end{equation*}%
where $g_{i}=\sum_{j=i+1}^{p(\mathcal{L})}\kappa _{j}$ and $\kappa _{j}$ is
shorthand for $\kappa (\omega _{j}(\mathcal{L}),d_{j}(\mathcal{L}))$. Now we
may repeatedly apply Lemma \ref{recrel2} (with $l=0$) to obtain 
\begin{equation}
D_{n}(\Delta _{\underline{\omega }(\mathcal{L}),\underline{d}(\mathcal{L}%
)})E_{n,0}(\omega )=E_{n-\kappa (\underline{\omega }(\mathcal{L}),\underline{%
d}(\mathcal{L})),0}(\omega )A(\mathcal{L},\omega ),  \label{stat1}
\end{equation}%
where $A(\mathcal{L},\omega )$ is proportional to an orthogonal matrix and
is nonsingular. Next assume that $\omega $ coincides with a coordinate of $%
\underline{\omega }(\mathcal{L})$, say $\omega _{i^{\ast }}(\mathcal{L})$.
We may write%
\begin{equation*}
D_{n}(\Delta _{\underline{\omega }(\mathcal{L}),\underline{d}(\mathcal{L}%
)})=\left( \dprod\limits_{i=1,i\neq i^{\ast }}^{p(\mathcal{L}%
)}D_{n-g_{i}^{\ast }}(\Delta _{\omega _{i}(\mathcal{L})}^{d_{i}(\mathcal{L}%
)})\right) D_{n}(\Delta _{\omega _{i^{\ast }}(\mathcal{L})}^{d_{i^{\ast }}(%
\mathcal{L})})
\end{equation*}%
where $g_{i}^{\ast }=\kappa _{i^{\ast }}+\sum_{j=i+1,j\neq i^{\ast }}^{p(%
\mathcal{L})}\kappa _{j}$. Applying first Lemma \ref{recrel} (with $l=0$)
and using the elementary fact that $E_{\cdot ,0}^{(s)}(\omega )=E_{\cdot
,0}(\omega )P(\omega )^{s}$, where $P(\omega )$ has been defined in Lemma %
\ref{recrel2}, and then repeatedly applying Lemma \ref{recrel2} gives%
\begin{equation}
D_{n}(\Delta _{\underline{\omega }(\mathcal{L}),\underline{d}(\mathcal{L}%
)})E_{n,d_{i^{\ast }}(\mathcal{L})}(\omega )=E_{n-\kappa (\underline{\omega }%
(\mathcal{L}),\underline{d}(\mathcal{L})),0}(\omega )A(\mathcal{L},\omega ),
\label{stat2}
\end{equation}%
for some $A(\mathcal{L},\omega )$ that is proportional to an orthogonal
matrix and is nonsingular. Multiplying (\ref{stat1}) and (\ref{stat2}) by $%
\Pi _{\mathcal{L}^{\bot }}H_{n}(\underline{\omega }(\mathcal{L}),\underline{d%
}(\mathcal{L}))$ from the left as well as by the inverse of $A(\mathcal{L}%
,\omega )$ from the right and using that $\Pi _{\mathcal{L}^{\bot }}=\Pi _{%
\mathcal{L}^{\bot }}H_{n}(\underline{\omega }(\mathcal{L}),\underline{d}(%
\mathcal{L}))D_{n}(\Delta _{\underline{\omega }(\mathcal{L}),\underline{d}(%
\mathcal{L})})$ holds, as noted in the proof of Lemma \ref{Mspectral}, now
completes the proof of (\ref{stat0}). The final claim in the lemma follows
immediately from the definition of $d_{i}(\mathcal{L})$ and $\omega _{i}(%
\mathcal{L})$, respectively, see Definition \ref{subspaces_2}. $\blacksquare 
$

\textbf{Proof of (\ref{S_of_F_all}) in Example \ref{chsingularToepmax}:} By
its definition, $\mathbb{S}(\mathfrak{F}_{\mathrm{all}},\mathcal{L})$ is
certainly contained in the set on the r.h.s. of (\ref{S_of_F_all}). To prove
the reverse inclusion, it suffices to show that any $\Gamma \subseteq
\lbrack 0,\pi ]$ with $\func{card}(\Gamma )\in \mathbb{N}$ is the
intersection of $[0,\pi ]$ with the support of a finite and symmetric Borel
measure $\mathsf{m}$ on $[-\pi ,\pi ]$ that arises as the weak limit of a
sequence $\mathsf{m}_{g_{j}}$ with $g_{j}$ as in (i) of Definition \ref%
{singularToep} (and with $\mathfrak{F}=\mathfrak{F}_{\mathrm{all}}$). Let
now $\Gamma \subseteq \lbrack 0,\pi ]$ with $\func{card}(\Gamma )\in \mathbb{%
N}$ be given. For every $\gamma \in \Gamma $ let $U_{j}(\gamma )$ be the
intersection of the open interval $(\gamma -1/j,\gamma +1/j)$ with $[0,\pi ]$%
. Define $V_{j}(\gamma )$ as $U_{j}(\gamma )\cup (-U_{j}(\gamma ))$. For $j$
large enough the elements of the collection $\left\{ V_{j}(\gamma ):\gamma
\in \Gamma \right\} $ are pairwise disjoint. For each $\gamma \in \Gamma $
one can then easily find functions $h_{\gamma ,j}$ on $[-\pi ,\pi ]$ such
that (i) $h_{\gamma ,j}$ vanishes outside of $V_{j}(\gamma )$, (ii) $%
h_{\gamma ,j}$ is positive on $V_{j}(\gamma )$, (iii) $h_{\gamma ,j}$ is
symmetric, and (iv) $h_{\gamma ,j}$ is Borel-measurable and satisfies $%
\int_{-\pi }^{\pi }h_{\gamma ,j}(\nu )d\nu =1$. [E.g., let $h_{\gamma ,j}$
be the indicator function of $V_{j}(\gamma )$, suitably normalized.] Define 
\begin{equation*}
h_{\gamma ,j}^{\ast }(\nu )=\left\vert \Delta _{\underline{\omega }(\mathcal{%
L}),\underline{d}(\mathcal{L})}(e^{\iota \nu })\right\vert ^{2}h_{\gamma
,j}(\nu )\left/ \int_{-\pi }^{\pi }\left\vert \Delta _{\underline{\omega }(%
\mathcal{L}),\underline{d}(\mathcal{L})}(e^{\iota \nu })\right\vert
^{2}h_{\gamma ,j}(\nu )d\nu \right. ,
\end{equation*}%
where we note that the integral in the denominator is certainly positive. It
is now obvious that for each $\gamma \in \Gamma $ the measure with density $%
h_{\gamma ,j}^{\ast }(\nu )$ converges weakly (as $j\rightarrow \infty $) to
the convex combination of pointmass at $\gamma $ and $-\gamma $, each with
weight $1/2$, if $\gamma \neq 0$, and to unit pointmass at zero if $\gamma
=0 $. Next define 
\begin{equation*}
h_{j}^{\ast \ast }=\sum_{\gamma \in \Gamma }h_{\gamma ,j}^{\ast }.
\end{equation*}%
It follows that the measure on $[-\pi ,\pi ]$ with density $h_{j}^{\ast \ast
}$ converges weakly to the discrete measure $\mathsf{m}^{\ast \ast }$ that
assigns mass $1/2$ to the points $\gamma $ and $-\gamma $ with $\gamma \in
\Gamma $, $\gamma \neq 0$, and assigns mass $1$ to $\gamma =0$ in case $0\in
\Gamma $. Define 
\begin{equation*}
f_{j}=c_{j}^{-1}\sum_{\gamma \in \Gamma }\left[ h_{\gamma ,j}(\nu )\left/
\int_{-\pi }^{\pi }\left\vert \Delta _{\underline{\omega }(\mathcal{L}),%
\underline{d}(\mathcal{L})}(e^{\iota \nu })\right\vert ^{2}h_{\gamma ,j}(\nu
)d\nu \right. \right] ,
\end{equation*}%
where the normalization constant $c_{j}$ is given by 
\begin{equation*}
c_{j}=\sum_{\gamma \in \Gamma }\left( \int_{-\pi }^{\pi }\left\vert \Delta _{%
\underline{\omega }(\mathcal{L}),\underline{d}(\mathcal{L})}(e^{\iota \nu
})\right\vert ^{2}h_{\gamma ,j}(\nu )d\nu \right) ^{-1}\text{.}
\end{equation*}%
Obviously, $f_{j}$ belongs to $\mathfrak{F}_{\mathrm{all}}$. Using this
sequence $f_{j}$, construct the sequence $g_{j}$ as in (\ref{g_j}). Clearly, 
$g_{j}$ coincides with $h_{j}^{\ast \ast }/\func{card}(\Gamma )$, and hence $%
\mathsf{m}_{g_{j}}$ converges to the finite and symmetric Borel measure $%
\mathsf{m}=\mathsf{m}^{\ast \ast }/\func{card}(\Gamma )$ which obviously is
exactly supported on $\Gamma \cup (-\Gamma )$. This completes the proof. $%
\blacksquare $

\textbf{Proof of the claims in Example \ref{chsingularToepmin}:} Suppose $%
\mathbb{S}(\mathfrak{F}_{\mathrm{all}}^{B},\mathcal{L})$ were nonempty. Then
there has to exist a $\mathsf{m}\in \mathbb{M}(\mathfrak{F},\mathcal{L})$.
Let $\Gamma =\mathrm{\limfunc{supp}}(\mathsf{m})\cap \lbrack 0,\pi ]$. Then $%
\func{card}(\Gamma )\in \mathbb{N}$ must hold. Let $g_{j}$ and $f_{j}$ be as
in (\ref{g_j}) with $\mathsf{m}_{g_{j}}$ converging weakly to $\mathsf{m}$.
Choose $\delta >0$ such that $B\delta <1$. Let $U$ be the union of the
finitely many intervals $(\eta -\varepsilon ,\eta +\varepsilon )\cap \lbrack
-\pi ,\pi ]$ where $\eta $ runs through the union of $\Gamma \cup (-\Gamma )$
with the set of zeros of the trigonometric polynomial $\left\vert \Delta _{%
\underline{\omega }(\mathcal{L}),\underline{d}(\mathcal{L})}(e^{\iota \nu
})\right\vert ^{2}$. Here $\varepsilon >0$ is chosen so small such that the
intervals are disjoint and such that the Lebesgue measure of $U$ is smaller
than $\delta $. Since the boundary of $V=[-\pi ,\pi ]\backslash U$ w.r.t. $%
[-\pi ,\pi ]$ has measure zero under $\mathsf{m}$ we have that $%
\int_{V}g_{j}(\nu )d\nu $ converges to $\mathsf{m}(V)=0$. Furthermore, $%
\left\vert \Delta _{\underline{\omega }(\mathcal{L}),\underline{d}(\mathcal{L%
})}(e^{\iota \nu })\right\vert ^{2}$ is bounded from below on $V$ by a
positive constant, $c$ say; and it is bounded from above by a finite
constant, $C$ say, on $[-\pi ,\pi ]$. We conclude that%
\begin{equation*}
\int_{V}g_{j}(\nu )d\nu \geq \left( 2\pi CB\right) ^{-1}c\int_{V}f_{j}(\nu
)d\nu \geq \left( 2\pi CB\right) ^{-1}c(1-B\delta )>0,
\end{equation*}%
a contradiction. $\blacksquare $

\textbf{Proof of Proposition \ref{limitsToep}:} 1. If $\bar{\Sigma}\in 
\limfunc{cl}(\mathcal{L}(\mathfrak{C}(\mathfrak{F})))$, then there must
exist a sequence $f_{j}\in \mathfrak{F}$ so that $\mathcal{L}(\Sigma
(f_{j}))\rightarrow \bar{\Sigma}$ holds as $j\rightarrow \infty $. Lemma \ref%
{Mspectral} (together with homogeneity of $\Sigma (\cdot ,\cdot )$ in its
first argument) now shows that 
\begin{equation}
\mathcal{L}(\Sigma (f_{j}))=\frac{\Pi _{\mathcal{L}^{\bot }}H_{n}(\underline{%
\omega }(\mathcal{L}),\underline{d}(\mathcal{L}))\Sigma (g_{j},n-\kappa (%
\underline{\omega }(\mathcal{L}),\underline{d}(\mathcal{L})))H_{n}^{\prime }(%
\underline{\omega }(\mathcal{L}),\underline{d}(\mathcal{L}))\Pi _{\mathcal{L}%
^{\bot }}}{\Vert \Pi _{\mathcal{L}^{\bot }}H_{n}(\underline{\omega }(%
\mathcal{L}),\underline{d}(\mathcal{L}))\Sigma (g_{j},n-\kappa (\underline{%
\omega }(\mathcal{L}),\underline{d}(\mathcal{L})))H_{n}^{\prime }(\underline{%
\omega }(\mathcal{L}),\underline{d}(\mathcal{L}))\Pi _{\mathcal{L}^{\bot
}}\Vert },  \label{Qseq}
\end{equation}%
where $g_{j}$ is as in Definition \ref{singularToep}, observing that the
denominator in the preceding display, i.e., $\Pi _{\mathcal{L}^{\bot
}}\Sigma (f_{j})\Pi _{\mathcal{L}^{\bot }}$, is nonzero because of positive
definiteness of $\Sigma (f_{j})$ and the assumption $\dim (\mathcal{L})<n$.
Since the sequence $\mathsf{m}_{g_{j}}$ is tight (as each $\mathsf{m}%
_{g_{j}} $ is a probability measure on the compact set $[-\pi ,\pi ]$), it
converges weakly, at least along a subsequence, to a finite and
(necessarily) symmetric Borel probability measure $\mathsf{m}$, say. We now
show that $\Sigma (\mathsf{m},n-\kappa (\underline{\omega }(\mathcal{L}),%
\underline{d}(\mathcal{L})))$ must be singular if $\limfunc{rank}(\bar{\Sigma%
})<n-\dim (\mathcal{L})$ holds: assume not, then $H_{n}(\underline{\omega }(%
\mathcal{L}),\underline{d}(\mathcal{L}))\Sigma (\mathsf{m},n-\kappa (%
\underline{\omega }(\mathcal{L}),\underline{d}(\mathcal{L})))H_{n}^{\prime }(%
\underline{\omega }(\mathcal{L}),\underline{d}(\mathcal{L}))$ defines a
bijection from $\limfunc{span}(V_{n}^{(0)}(\underline{\omega }(\mathcal{L}),%
\underline{d}(\mathcal{L})))^{\bot }$ onto itself (with the convention that
this latter space is $\mathbb{R}^{n}$ if $p(\mathcal{L})=0$). Consequently, 
\begin{equation*}
A:=\Pi _{\mathcal{L}^{\bot }}H_{n}(\underline{\omega }(\mathcal{L}),%
\underline{d}(\mathcal{L}))\Sigma (\mathsf{m},n-\kappa (\underline{\omega }(%
\mathcal{L}),\underline{d}(\mathcal{L})))H_{n}^{\prime }(\underline{\omega }(%
\mathcal{L}),\underline{d}(\mathcal{L}))\Pi _{\mathcal{L}^{\bot }}
\end{equation*}%
would have rank $n-\dim (\mathcal{L})>0$.\footnote{%
By construction of $A$, its rank can not exceed $n-\limfunc{dim}(\mathcal{L})
$. If its rank were less than $n-\limfunc{dim}(\mathcal{L})$, we could find
a nonzero $z \notin \mathcal{L}$ such that $A z = 0$ holds. Since $A$
vanishes on $\mathcal{L}$ by its construction, we may assume w.l.o.g.~that $%
z \in \mathcal{L}^{\bot} \backslash \{0\}$. But then $z^{\prime}Az = 0$ and
hence $z^{\prime}H_{n}(\underline{\omega }(\mathcal{L}),\underline{d}(%
\mathcal{L}))\Sigma (\mathsf{m},n-\kappa (\underline{\omega }(\mathcal{L}),%
\underline{d}(\mathcal{L})))H_{n}^{\prime }(\underline{\omega }(\mathcal{L}),%
\underline{d}(\mathcal{L}))z = 0$ follows as $\Pi_{\mathcal{L}^{\bot}} z = z$%
. By symmetry and nonnegative definiteness we arrive at $H_{n}(\underline{%
\omega }(\mathcal{L}),\underline{d}(\mathcal{L}))\Sigma (\mathsf{m},n-\kappa
(\underline{\omega }(\mathcal{L}),\underline{d}(\mathcal{L})))H_{n}^{\prime
}(\underline{\omega }(\mathcal{L}),\underline{d}(\mathcal{L}))z = 0$. Since $%
z \in \mathcal{L}^{\bot} \subseteq \limfunc{span}(V_{n}^{(0)}(\underline{%
\omega }(\mathcal{L}),\underline{d}(\mathcal{L})))^{\bot }$ in view of the
definition of $V_{n}^{(0)}(\underline{\omega }(\mathcal{L}),\underline{d}(%
\mathcal{L}))$, this would contradict the bijectivity established before.}
But passing to the limit in (\ref{Qseq}) along the above mentioned
subsequence would then entail $\bar{\Sigma}=A/\left\Vert A\right\Vert $,
contradicting $\func{rank}(\bar{\Sigma})<n-\dim (\mathcal{L})$. Having
established singularity of $\Sigma (\mathsf{m},n-\kappa (\underline{\omega }(%
\mathcal{L}),\underline{d}(\mathcal{L})))$ and noting that $\mathsf{m}$ can
not be the zero measure (as it must have total mass $1$), we can conclude
from Lemma \ref{Toeplitz} that $\mathsf{m}=\sum_{i=1}^{p}\bar{c}_{i}(\delta
_{-\gamma _{i}}+\delta _{\gamma _{i}})$ for some $p\in \mathbb{N}$ and
positive real numbers $\bar{c}_{i}$, $i=1,\ldots ,p$, and for $0\leq \gamma
_{1}<\gamma _{2}<\ldots <\gamma _{p}\leq \pi $, such that $%
\sum_{i=1}^{p}\kappa (\gamma _{i},1)<n-\kappa (\underline{\omega }(\mathcal{L%
}),\underline{d}(\mathcal{L}))$ holds. In particular, $\Gamma :=\left\{
\gamma _{1},\ldots ,\gamma _{p}\right\} \in \mathbb{S}(\mathfrak{F},\mathcal{%
L})$ holds. We may now conclude that 
\begin{equation}
\Sigma (\mathsf{m},n-\kappa (\underline{\omega }(\mathcal{L}),\underline{d}(%
\mathcal{L})))=\sum_{i=1}^{p}2\bar{c}_{i}E_{n-\kappa (\underline{\omega }(%
\mathcal{L}),\underline{d}(\mathcal{L})),0}(\gamma _{i})E_{n-\kappa (%
\underline{\omega }(\mathcal{L}),\underline{d}(\mathcal{L})),0}^{\prime
}(\gamma _{i}).  \label{Qseq_2}
\end{equation}%
As a consequence of Lemma \ref{proj} the matrix $A$ is nonzero and coincides
with the numerator in (\ref{represent}) for some positive $c(\gamma _{i})$%
's. Passing to the limit in (\ref{Qseq}) along the above mentioned
subsequence then establishes (\ref{represent}). We turn to the second claim
in Part 1 next: If $\Gamma \in \mathbb{S}(\mathfrak{F},\mathcal{L})$, we can
find a measure $\mathsf{m}$ and a sequence $f_{j}\in \mathfrak{F}$
satisfying all the requirements in Definition \ref{singularToep}. In
particular, (\ref{Qseq_2}) again holds for some positive constants $c_{i}$
and with $\gamma _{i}$ enumerating the elements of $\Gamma $. Consider the
sequence $\Sigma (f_{j})$. Again Lemma \ref{Mspectral} shows that (\ref{Qseq}%
) holds. It follows that the numerator of (\ref{Qseq}) converges to $A$
defined above. Now, $A\neq 0$ follows again from Lemma \ref{proj}. But then
we can conclude that $\mathcal{L}(\Sigma (f_{j}))$ converges to $\bar{\Sigma}%
:=A/\left\Vert A\right\Vert $, implying $\bar{\Sigma}\in \limfunc{cl}(%
\mathcal{L}(\mathfrak{C}(\mathfrak{F})))$. But Lemma \ref{proj} now implies
that (\ref{represent}) holds. The claim in parentheses is obvious since $%
\bar{\Sigma}$ vanishes on $\mathcal{L}$ in view of (\ref{represent}).

2.\&3. Part 2 is a simple consequence of Part 1 since $\limfunc{span}%
(\sum_{i=1}^{l}A_{i}A_{i}^{\prime })=\limfunc{span}((A_{1},\ldots ,A_{l}))$
holds for matrices $A_{i}$ of the same row-dimension. Part 3 follows
immediately from Part 2. $\blacksquare $

\textbf{Proof of Theorem \ref{Ctoep}:} Follows from Corollary \ref{C5} (with 
$\mathfrak{C}=\mathfrak{C}(\mathfrak{F})$) together with Proposition \ref%
{limitsToep}. $\blacksquare $

\textbf{Proof of the claim in Remark \ref{simplifications}(iii):} Note that
any $\mathcal{S}$ satisfying (ii), but not (i), coincides with $\mathcal{L}%
^{\bot }$. Since $N^{\dag }$ is invariant under addition of elements from $%
\mathcal{L}$, we can write $N^{\dag }$ as the direct sum $M\oplus \mathcal{L}
$, where $M=N^{\dag }\cap \mathcal{L}^{\bot }$. But then we have%
\begin{equation*}
0=\lambda _{\mathbb{R}^{n}}(N^{\dag })=\lambda _{\mathcal{L}^{\bot }\oplus 
\mathcal{L}}(M\oplus \mathcal{L}),
\end{equation*}%
from which we can conclude that $\lambda _{\mathcal{L}^{\bot }}(M)=0$, i.e., 
$\lambda _{\mathcal{S}}(N^{\dag })=0$ holds (note that $\lambda _{\mathcal{S}%
}(N^{\dag })=\lambda _{\mathcal{S}}(N^{\dag }\cap \mathcal{S})$). Let $\mu
_{0}\in \mathfrak{M}_{0}$ be arbitrary, and write $\mu _{0}=l+l^{\bot }$,
where $l\in \mathcal{L}$ and $l^{\bot }\in \mathcal{L}^{\bot }$. Since $%
l^{\bot }+\mathcal{S}=\mathcal{S}$ and since $N^{\dag }$ is invariant under
addition of elements from $\mathcal{L}$ we arrive at%
\begin{equation*}
\lambda _{\mu _{0}+\mathcal{S}}(N^{\dag })=\lambda _{l+l^{\bot }+\mathcal{S}%
}(N^{\dag })=\lambda _{l+\mathcal{S}}(N^{\dag })=\lambda _{\mathcal{S}%
}(N^{\dag }-l)=\lambda _{\mathcal{S}}(N^{\dag })=0.
\end{equation*}%
But this shows that any $\mathcal{S}$ satisfying (ii), but not (i),
automatically satisfies the condition $\lambda _{\mu _{0}+\mathcal{S}%
}(N^{\dag })=0$ for every $\mu _{0}\in \mathfrak{M}_{0}$. $\blacksquare $

\textbf{Proof of Theorem \ref{Ctoep_nonsphericity}:} Observe that here $%
\mathcal{L}=\mathfrak{M}_{0}^{lin}$ and that $\dim (\mathcal{L})=k-q<n$
holds, since we always assume $k<n$ and $q\geq 1$. The claims in Part 1 then
follow immediately from Part 1 of Theorem \ref{Ctoep} and Remarks \ref%
{special_cases}(i) and \ref{simplifications}(i)-(iii) (recalling that $%
N^{\dag }=N^{\ast }$ is an $G(\mathfrak{M})$-invariant $\lambda _{\mathbb{R}%
^{n}}$-null set, cf. Lemma \ref{Fact1}). Part 2 follows from Part 1, noting
that $\lambda _{\mathcal{T}}(N^{\ast })=0$ is equivalent to $\mathcal{T}%
\nsubseteq N^{\ast }$ if $\mathcal{T}$ is an affine space (cf. the proof of
Corollary \ref{C5}). $\blacksquare $

\textbf{Proof of Theorem \ref{Ctoep_corr}:} Recall from Section \ref%
{suff_con} that, under Assumptions \ref{AW} and \ref{R_and_X},
autocorrelation robust tests based on $T_{w}$ are a special case of
nonsphericity-corrected F-type tests and that Assumption 5 (as well as
Assumptions 6 and 7) of \cite{PP2016} are then satisfied. Furthermore,
recall from Section \ref{suff_con} that the set $N^{\ast }$ here is given by 
$\mathsf{B}$, which is a finite union of proper linear subspaces as a result
of Lemma \ref{Fact2}. Then, except for the last claim, Part 1 follows
immediately from Theorem \ref{Ctoep_nonsphericity}. To prove the last claim
in Part 1, note that the assumptions for Proposition \ref{P6} (with $%
\mathfrak{C}=\mathfrak{C}(\mathfrak{F})$ and $\mathcal{L}=\mathfrak{M}%
_{0}^{lin}$) are satisfied in view of Lemmata \ref{Fact1} and \ref{Fact3},
keeping in mind the characterization of $\mathbb{J}(\mathcal{L},\mathfrak{C}(%
\mathfrak{F}))$ provided in Proposition \ref{limitsToep} and Remark \ref%
{simplifications}. The result then follows from Proposition \ref{P6} in case 
$\alpha <\alpha ^{\ast }$. In case $\alpha =\alpha ^{\ast }<1$ it follows
from Remark \ref{smallest C_2}(i),(ii) together with Lemma \ref{Fact1}. Part
2 is now a trivial consequence of Part 1. Part 3 is obvious. $\blacksquare $

\textbf{Proof of the claims in Remark \ref{nec_2}:} Suppose the claimed
equivalence does not hold. Then we can find a $\gamma \in \lbrack 0,\pi ]$
such that $\gamma \notin \bigcup \mathbb{S}(\mathfrak{F},\mathcal{L})$ and $%
\mu _{0}+\limfunc{span}(\Pi _{\mathcal{L}^{\bot }}E_{n,\rho (\gamma ,%
\mathcal{L})}(\gamma ))\subseteq N^{\dag }$ holds for some $\mu _{0}\in 
\mathfrak{M}_{0}$. Observe that $N^{\dag }+\mathcal{L}\subseteq N^{\dag }$
holds by the assumed invariance properties. Consequently, 
\begin{eqnarray*}
\mathcal{A} &:&=\limfunc{span}(E_{n,\rho (\gamma ,\mathcal{L})}(\gamma ))+%
\mathcal{L}=\limfunc{span}(\Pi _{\mathcal{L}^{\bot }}E_{n,\rho (\gamma ,%
\mathcal{L})}(\gamma )+\Pi _{\mathcal{L}}E_{n,\rho (\gamma ,\mathcal{L}%
)}(\gamma ))+\mathcal{L} \\
&\subseteq &\limfunc{span}(\Pi _{\mathcal{L}^{\bot }}E_{n,\rho (\gamma ,%
\mathcal{L})}(\gamma ))+\mathcal{L}\subseteq N^{\dag }-\mu _{0}\neq \mathbb{R%
}^{n}
\end{eqnarray*}%
since $N^{\dag }\neq \mathbb{R}^{n}$. It follows that $\dim (\mathcal{A})<n$%
. Trivially, $\limfunc{span}(E_{n,\rho (\gamma ,\mathcal{L})}(\gamma
))\subseteq \mathcal{A}$ holds, and $\limfunc{span}(E_{n,\rho (\gamma ,%
\mathcal{L})}(\gamma ))\nsubseteqq \mathcal{L}$ in view of the definition of 
$\rho (\gamma ,\mathcal{L})$. Consequently, 
\begin{equation*}
\kappa (\underline{\omega }(\mathcal{L}),\underline{d}(\mathcal{L}))+\kappa
(\gamma ,1)\leq \kappa (\underline{\omega }(\mathcal{A}),\underline{d}(%
\mathcal{A}))\leq \dim (\mathcal{A})<n,
\end{equation*}%
where we have made use of Lemma \ref{dimsubspaces}. But in view of Example %
\ref{chsingularAR(2)}(ii),(iii) this shows that $\{\gamma \}\in \mathbb{S}(%
\mathfrak{F},\mathcal{L})$ and thus $\gamma \in \bigcup \mathbb{S}(\mathfrak{%
F},\mathcal{L})$, a contradiction. This establishes the first claim. The
remaining claims are proved analogously, observing that the relevant sets $%
N^{\ast }$, $\mathsf{B}$, etc are $\lambda _{\mathbb{R}^{n}}$-null sets and
thus are proper subsets of $\mathbb{R}^{n}$. $\blacksquare $

\section{Appendix\label{App E}: Extensions to non-Gaussian models}

As mentioned in the Introduction and in Section \ref{frame}, Gaussianity is
not essential to the results in the paper. Here we discuss various ways of
substantially weakening the Gaussianity assumption.

\subsection{Elliptically symmetric and related distributions}

Consider again the linear model (\ref{lm}) with all the assumptions made in
Section \ref{frame}, except that the disturbance vector $\mathbf{U}$ now
follows an elliptically symmetric distribution. More precisely, let $Z_{%
\mathrm{spher}}$ denote the set of spherically symmetric distributions $%
\zeta $ on $\mathbb{R}^{n}$ that have no atom at the origin. The vector $%
\mathbf{U}$ is assumed to be distributed as $\sigma \Sigma ^{1/2}\mathbf{z}$%
, where $\mathbf{z}$ has a distribution $\zeta \in Z_{\mathrm{spher}}$ and $%
\Sigma \in \mathfrak{C}$. For background information on elliptically
symmetric distributions see \cite{cambanis}. Three remarks are in order:
First, we do not assume that $\zeta $ is absolutely continuous. Second, if $%
\zeta $ has a finite first moment (which we, however, do not assume), then $%
\mathbf{U}$ has mean zero; otherwise, we can only say that the origin is the
(uniquely determined) center of elliptical symmetry of $\mathbf{U}$. Third,
if $\zeta $ has finite second moments, then $\mathbf{U}$ has covariance
matrix $c_{\zeta }\sigma ^{2}\Sigma $, where $c_{\zeta }$ is the variance of
the first component of $\mathbf{z}$; otherwise, $\sigma ^{2}\Sigma $ is a
parameter describing the ellipticity of the distribution of $\mathbf{U}$
(which is unique if we consider $\zeta $ as fixed and which is only unique
up to a scale factor if $\zeta $ can freely vary in $Z_{\mathrm{spher}}$).
Nevertheless, we shall -- in abuse of terminology -- continue to address the
set $\mathfrak{C}$ as the covariance model. Let $Q_{\mu ,\sigma ^{2}\Sigma
,\zeta }$ denote the distribution of $\mathbf{Y}$ resulting from model (\ref%
{lm}) under the preceding assumptions and where $\mu =X\beta $. Then, for 
\emph{every} distribution $\zeta \in Z_{spher}$, for \emph{any} $G(\mathfrak{%
M}_{0})$-invariant rejection region $W$, and for \emph{any }$\mu _{0}\in 
\mathfrak{M}_{0}$ we have that%
\begin{eqnarray*}
Q_{\mu _{0},\sigma ^{2}\Sigma ,\zeta }(W) &=&\Pr \left( \mu _{0}+\sigma
\Sigma ^{1/2}\mathbf{z}\in W\right) =\Pr \left( \sigma \Sigma ^{1/2}\mathbf{z%
}/\left\Vert \mathbf{z}\right\Vert \in (W-\mu _{0})/\left\Vert \mathbf{z}%
\right\Vert +\mu _{0}-\mu _{0}\right) \\
&=&\Pr \left( \sigma \Sigma ^{1/2}\mathbf{z}/\left\Vert \mathbf{z}%
\right\Vert \in W-\mu _{0}\right) =\Pr \left( \mu _{0}+\sigma \Sigma ^{1/2}%
\mathbf{z}/\left\Vert \mathbf{z}\right\Vert \in W\right) ,
\end{eqnarray*}%
where we have used $G(\mathfrak{M}_{0})$-invariance and the fact that $%
\left\Vert \mathbf{z}\right\Vert \neq 0$ almost surely in the last but one
step. Note that the distribution of $\mathbf{z}/\left\Vert \mathbf{z}%
\right\Vert $ is the uniform distribution on the unit sphere in $\mathbb{R}%
^{n}$ and hence does not depend on $\zeta $ at all. Since the Gaussian case
is a special case with $\zeta $ the $n$-dimensional standard normal
distribution, we have%
\begin{equation}
P_{\mu _{0},\sigma ^{2}\Sigma }(W)=Q_{\mu _{0},\sigma ^{2}\Sigma ,\zeta }(W)
\label{equality}
\end{equation}%
for every distribution $\zeta \in Z_{\mathrm{spher}}$ (and for every $\mu
_{0}\in \mathfrak{M}_{0}$, $0<\sigma ^{2}<\infty $, $\Sigma \in \mathfrak{C}$%
). That is, the rejection probabilities of any $G(\mathfrak{M}_{0})$%
-invariant rejection region $W$ under the null hypothesis are the same
whether we assume a Gaussian linear model or a linear model with
elliptically symmetric errors (satisfying the above made assumptions). In
particular, for the size we have%
\begin{equation}
\sup_{\zeta \in Z_{\mathrm{spher}}}\sup_{\mu _{0}\in \mathfrak{M}%
_{0}}\sup_{0<\sigma ^{2}<\infty }\sup_{\Sigma \in \mathfrak{C}}Q_{\mu
_{0},\sigma ^{2}\Sigma ,\zeta }(W)=\sup_{\mu _{0}\in \mathfrak{M}%
_{0}}\sup_{0<\sigma ^{2}<\infty }\sup_{\Sigma \in \mathfrak{C}}P_{\mu
_{0},\sigma ^{2}\Sigma }(W).  \label{size equality}
\end{equation}%
This shows that all results in the paper carry over to the case of
elliptically distributed errors as they stand. In particular, note that the
critical values $C(\alpha )$ computed under the assumption of Gaussianity
automatically also deliver size control under the more general assumption of
elliptical symmetry and thus are robust under this sort of deviations from
Gaussianity. Note that the model discussed here allows for heavy-tailed
disturbances.

The above discussion in fact shows that all results in the paper hold for a
class of distributions even wider than the class of elliptically symmetric
distributions: Let $Z_{\mathrm{ua}}$ denote the class of distributions $%
\zeta $ on $\mathbb{R}^{n}$ that (i) do not have any mass at the origin and
(ii) have the property that the distribution of $\mathbf{z}/\left\Vert 
\mathbf{z}\right\Vert $ under $\zeta $ is the uniform distribution on the
unit sphere in $\mathbb{R}^{n}$. [Clearly, $Z_{\mathrm{ua}}$ contains $Z_{%
\mathrm{spher}}$, but also other distributions under which the radial
component $\left\Vert \mathbf{z}\right\Vert $ may be dependent on the
uniformly distributed angular component $\mathbf{z}/\left\Vert \mathbf{z}%
\right\Vert $.] Then it is plain that (\ref{equality}) and hence (\ref{size
equality}) continue to hold with $Z_{\mathrm{spher}}$ replaced by $Z_{%
\mathrm{ua}}$. In particular, the critical values computed under Gaussianity
are valid in this much wider context. Note, however, that now in general $%
X\beta $ and $\sigma ^{2}\Sigma $ no longer have the same interpretation as
in the Gaussian or elliptically symmetric case.

\begin{remark}
\label{cross} The above discussion obviously also applies in case that (i) $%
\zeta $ is restricted to a subset of $Z_{\mathrm{spher}}$ (or $Z_{\mathrm{ua}%
}$, respectively) or (ii) there are cross-restrictions between $(\beta
,\sigma ^{2},\Sigma )$ and $\zeta $ in the sense that depending on $(\beta
,\sigma ^{2},\Sigma )$ the distribution $\zeta $ is restricted to a subset $%
Z(\beta ,\sigma ^{2},\Sigma )$ of $Z_{\mathrm{spher}}$ (or $Z_{\mathrm{ua}}$%
, respectively).
\end{remark}

\subsection{Other distributions}

Again we consider the linear model (\ref{lm}) with all the assumptions as in
Section \ref{frame}, except that $\mathbf{U}$ is now assumed to be
distributed as $\sigma \Sigma ^{1/2}U\mathbf{z}$ where $\Sigma \in \mathfrak{%
C}$, where $U$ is an orthogonal matrix and where $\mathbf{z}$ has a density $%
h$. It is assumed that the pair $(h,U)$ belongs to a given subset $\mathcal{D%
}$ of $H\times \mathcal{U}$ where $H$ is a given set of density functions
and $\mathcal{U}$ denotes the set of orthogonal $n\times n$ matrices.
[Important special cases are (i) the case where $U=I_{n}$ holds for every $%
(h,U)\in \mathcal{D}$, or (ii) when $H$ is a singleton.] We assume that $%
\mathbf{z}$ has mean zero and unit covariance matrix under each $h\in H$.%
\footnote{%
We make this assumption only in order for $X\beta $ and $\sigma ^{2}\Sigma $
to have the same interpretation as in the Gaussian case. Lemma \ref{Lcont2}
also holds without this assumption. However, note that then the
interpretation of $X\beta $ and $\sigma ^{2}\Sigma $ becomes somewhat
obscure and there is no guarantee that the parameters are identified.} Let $%
Q_{\mu ,\sigma ^{2}\Sigma ,h,U}$ denote the distribution of $\mathbf{Y}$
resulting from the preceding assumptions. Then we have the following result
which will allow us to easily extend the size control results from the
Gaussian case to the present setting. Observe that the condition on $H$ in
the subsequent lemma is trivially satisfied if $H=\left\{ h\right\} $ and
thus imposes no further condition on $h$ in this case (and the same is true
if $H$ is a finite set).

\begin{lemma}
\label{Lcont2} Suppose the maintained assumptions of this subsection hold
and that there is a $\lambda _{\mathbb{R}^{n}}$-integrable envelope $h^{\ast
}$ for $H$ (i.e., $h(z)\leq h^{\ast }(z)$ holds $\lambda _{\mathbb{R}^{n}}$%
-a.e. for every $h\in H$ and $\int h^{\ast }(z)d\lambda _{\mathbb{R}%
^{n}}(z)<\infty $). Let $W_{m}$ be a sequence of rejection regions such that%
\begin{equation}
\sup_{\mu _{0}\in \mathfrak{M}_{0}}\sup_{0<\sigma ^{2}<\infty }\sup_{\Sigma
\in \mathfrak{C}}P_{\mu _{0},\sigma ^{2}\Sigma }(W_{m})\rightarrow 0\text{ \
\ for \ }m\rightarrow \infty .  \label{supP}
\end{equation}%
Then%
\begin{equation}
\sup_{(h,U)\in \mathcal{D}}\sup_{\mu _{0}\in \mathfrak{M}_{0}}\sup_{0<\sigma
^{2}<\infty }\sup_{\Sigma \in \mathfrak{C}}Q_{\mu _{0},\sigma ^{2}\Sigma
,h,U}(W_{m})\rightarrow 0\text{ \ \ for \ }m\rightarrow \infty .
\label{supQ}
\end{equation}
\end{lemma}

\textbf{Proof:} Let $Q_{\mu _{0,m},\sigma _{m}^{2}\Sigma _{m},h_{m},U_{m}}$
with $\mu _{0,m}\in \mathfrak{M}_{0}$, $0<\sigma _{m}^{2}<\infty $, $\Sigma
_{m}\in \mathfrak{C}$, and $(h_{m},U_{m})\in \mathcal{D}$ be a sequence such
that $Q_{\mu _{0,m},\sigma _{m}^{2}\Sigma _{m},h_{m},U_{m}}(W_{m})$ differs
from the multiple supremum in (\ref{supQ}) only by a null sequence. Then we
have $P_{\mu _{0,m},\sigma _{m}^{2}\Sigma _{m}}(W_{m})\rightarrow 0$ as $%
m\rightarrow \infty $ as a consequence of (\ref{supP}). Furthermore,%
\begin{equation*}
P_{\mu _{0,m},\sigma _{m}^{2}\Sigma _{m}}(W_{m})=\Pr (\mu _{0,m}+\sigma
_{m}\Sigma _{m}^{1/2}\mathbf{G}\in W_{m})=\Pr (\mathbf{G}\in U_{m}B_{m})=\Pr
(\mathbf{G}\in B_{m})=P_{0,I_{n}}(B_{m})
\end{equation*}%
and%
\begin{equation*}
Q_{\mu _{0,m},\sigma _{m}^{2}\Sigma _{m},h_{m},U_{m}}(W_{m})=\Pr (\mu
_{0,m}+\sigma _{m}\Sigma _{m}^{1/2}U_{m}\mathbf{z}_{m}\in W_{m})=\Pr (%
\mathbf{z}_{m}\in B_{m}),
\end{equation*}%
where $\mathbf{G}$ is a standard Gaussian $n\times 1$ random vector, $%
\mathbf{z}_{m}$ is a random vector with density $h_{m}$, and $%
B_{m}=U_{m}^{\prime }\sigma _{m}^{-1}\Sigma _{m}^{-1/2}(W_{m}-\mu _{0,m})$.
It suffices to show that $\Pr (\mathbf{z}_{m}\in B_{m})$ converges to zero
as $m\rightarrow \infty $. Let $m^{\prime }$ be an arbitrary subsequence of $%
m$. Since $P_{0,I_{n}}(B_{m^{\prime }})\rightarrow 0$, there exists a
subsequence $m^{\prime \prime }$ of $m^{\prime }$ and a $\lambda _{\mathbb{R}%
^{n}}$-null set $A_{1}$ such that $\mathbf{1}_{B_{m^{\prime \prime
}}}(z)\rightarrow 0$ for every $z\notin A_{1}$ (e.g., Theorem 3.12 in \cite%
{rudin}). Since the envelope $h^{\ast }$ is $\lambda _{\mathbb{R}^{n}}$%
-integrable, it is finite everywhere, except possibly on a $\lambda _{%
\mathbb{R}^{n}}$-null set $A_{2}$. By assumption, $0\leq h_{m^{\prime \prime
}}(z)\leq h^{\ast }(z)$ holds for all $m^{\prime \prime }$ outside of a $%
\lambda _{\mathbb{R}^{n}}$-null set $A_{3}$. Set $A=A_{1}\cup A_{2}\cup
A_{3} $. But then $\mathbf{1}_{B_{m^{\prime \prime }}}(z)h_{m^{\prime \prime
}}(z)\rightarrow 0$ holds for every $z$ outside the $\lambda _{\mathbb{R}%
^{n}}$-null set $A$. Furthermore, $\left\vert \mathbf{1}_{B_{m^{\prime
\prime }}}(z)h_{m^{\prime \prime }}(z)\right\vert \leq h^{\ast }(z)$ holds
for every $m^{\prime \prime }$ and for every $z\notin A$. But then the
Dominated Convergence Theorem gives 
\begin{equation*}
\Pr (\mathbf{z}_{m^{\prime \prime }}\in B_{m^{\prime \prime }})=\int_{%
\mathbb{R}^{n}}\mathbf{1}_{B_{m^{\prime \prime }}}(z)h_{m^{\prime \prime
}}(z)d\lambda _{\mathbb{R}^{n}}(z)\rightarrow 0.
\end{equation*}%
Since the subsequence $m^{\prime }$ was arbitrary, we conclude that $\Pr (%
\mathbf{z}_{m}\in B_{m})$ converges to zero. $\blacksquare $

As a consequence of this lemma, versions of all the size control results in
the paper (except for the exact size control results) can be given under the
maintained assumptions of this subsection if $H$ satisfies the assumption in
Lemma \ref{Lcont2}. We illustrate this exemplarily with the following
version of Theorem \ref{HAR_F_all}, which is an immediate consequence of the
just mentioned theorem combined with Lemma \ref{Lcont2}. Observe that --
other than with the extensions discussed in the preceding subsection -- the
critical values $C^{\prime }(\alpha )$ in the subsequent theorem may now
differ from the critical values one obtains in the Gaussian case.

\begin{theorem}
\label{HAR_F_all_general} Suppose the maintained assumptions of this
subsection are satisfied with $H$ having a $\lambda _{\mathbb{R}^{n}}$%
-integrable envelope $h^{\ast }$. Suppose Assumptions \ref{AW} and \ref%
{R_and_X} are satisfied and $T_{w}$ is defined by (\ref{tlr}). Then for
every $0<\alpha <1$ there exists a real number $C^{\prime }(\alpha )$ such
that%
\begin{equation}
\sup_{(h,U)\in \mathcal{D}}\sup_{\mu _{0}\in \mathfrak{M}_{0}}\sup_{0<\sigma
^{2}<\infty }\sup_{f\in \mathfrak{F}_{\mathrm{all}}}Q_{\mu _{0},\sigma
^{2}\Sigma (f),h}(T_{w}\geq C^{\prime }(\alpha ))\leq \alpha  \label{supQ1.5}
\end{equation}%
holds, provided that 
\begin{equation}
\limfunc{span}(E_{n,\rho (\gamma ,\mathfrak{M}_{0}^{lin})}(\gamma
))\not\subseteq \mathsf{B}\text{ \ \ for every \ }\gamma \in \lbrack 0,\pi ].
\label{non-incl**}
\end{equation}%
In case the set $\mathsf{B}$ coincides with $\func{span}(X)$, condition (\ref%
{non-incl**}) can equivalently be expressed as 
\begin{equation*}
\func{rank}(X,E_{n,\rho (\gamma ,\mathfrak{M}_{0}^{lin})}(\gamma ))>k\text{
\ \ for every \ }\gamma \in \lbrack 0,\pi ].
\end{equation*}
\end{theorem}

Extensions of the other size control results in the paper to the present
setting follow a similar pattern and will not be given.

\begin{remark}
For the computation of critical values the cumbersome optimization over $h$
in (\ref{supQ1.5}) can in principle be avoided by determining $C^{\prime
}(\alpha )$ from%
\begin{equation*}
\sup_{U\in \mathcal{U}}\sup_{f\in \mathfrak{F}_{\mathrm{all}}}\int 
\boldsymbol{1}(T_{w}(z)\geq C^{\prime }(\alpha ))h^{\ast }(U^{\prime }\Sigma
^{-1/2}(f)(z-\mu _{0}))\det (\Sigma ^{-1/2}(f))d\lambda _{\mathbb{R}%
^{n}}(z)\leq \alpha
\end{equation*}%
for some $\mu _{0}\in \mathfrak{M}_{0}$, since the l.h.s. in the preceding
display is easily seen to be an upper bound for the l.h.s. of (\ref{supQ1.5}%
) in view of $G(\mathfrak{M}_{0})$-invariance of $T_{w}$ (cf. Remark \ref%
{special cases}). However, this will often lead to a quite conservative
choice for $C^{\prime }(\alpha )$. A similar remark applies to the more
general result indicated by the lemma given below.
\end{remark}

We next show that the above reasoning based on Lemma \ref{Lcont2} can
actually be extended to even larger classes of distributions, including
cases where $\mathbf{z}$ need not have a density. Consider again the linear
model (\ref{lm}) with all the assumptions as in Section \ref{frame}, except
that $\mathbf{U}$ is now assumed to be distributed as $\sigma \Sigma ^{1/2}U%
\mathbf{z}$ where $U$ is an orthogonal matrix and $\mathbf{z}$ has a
distribution $\zeta $ with the following properties: the pair $(\zeta ,U)$
belongs to a set $\mathcal{E}\subseteq Z\times \mathcal{U}$, where $Z$ is a
set of distributions with the property that (i) no $\zeta \in Z$ has an atom
at the origin and (ii) that the distribution of the random vector $\mathbf{z}%
/\left\Vert \mathbf{z}\right\Vert $ under each $\zeta \in Z$ has a density $%
\bar{h}_{\zeta }$ w.r.t. the uniform distribution $\upsilon _{S^{n-1}}$ on
the Borel-sets of the unit sphere $S^{n-1}$.\footnote{%
Without a further assumption such as that $\mathbf{z}$ has mean zero and
unit covariance matrix under each $\zeta \in Z$ the interpretation of the
parameters in the model is somewhat obscure and they are not guaranteed to
be identified.} Let $Q_{\mu ,\sigma ^{2}\Sigma ,\zeta ,U}$ denote the
distribution of $\mathbf{Y}$ resulting from the preceding assumptions. Then
we have the following lemma.

\begin{lemma}
\label{Lcont3} Suppose the assumptions in the preceding paragraph are
satisfied. Furthermore, assume that there is a $\upsilon _{S^{n-1}}$%
-integrable envelope $\bar{h}^{\ast }$ for $\left\{ \bar{h}_{\zeta }:\zeta
\in Z\right\} $ (i.e., $\bar{h}_{\zeta }(s)\leq \bar{h}^{\ast }(s)$ holds $%
\upsilon _{S^{n-1}}$-a.e. for every $\zeta \in Z$ and $\int_{S^{n-1}}\bar{h}%
^{\ast }(s)d\upsilon _{S^{n-1}}(s)<\infty $). Let $W_{m}$ be a sequence of $%
G(\mathfrak{M}_{0})$-invariant rejection regions such that%
\begin{equation}
\sup_{\mu _{0}\in \mathfrak{M}_{0}}\sup_{0<\sigma ^{2}<\infty }\sup_{\Sigma
\in \mathfrak{C}}P_{\mu _{0},\sigma ^{2}\Sigma }(W_{m})\rightarrow 0\text{ \
\ for \ }m\rightarrow \infty .  \label{supP2}
\end{equation}%
Then%
\begin{equation}
\sup_{(\zeta ,U)\in \mathcal{E}}\sup_{\mu _{0}\in \mathfrak{M}%
_{0}}\sup_{0<\sigma ^{2}<\infty }\sup_{\Sigma \in \mathfrak{C}}Q_{\mu
_{0},\sigma ^{2}\Sigma ,\zeta ,U}(W_{m})\rightarrow 0\text{ \ \ for \ }%
m\rightarrow \infty .  \label{supQ2}
\end{equation}
\end{lemma}

\textbf{Proof:} Let $Q_{\mu _{0,m},\sigma _{m}^{2}\Sigma _{m},\zeta
_{m},U_{m}}$ with $\mu _{0,m}\in \mathfrak{M}_{0}$, $0<\sigma
_{m}^{2}<\infty $, $\Sigma _{m}\in \mathfrak{C}$, and $(\zeta _{m},U_{m})\in 
\mathcal{E}$ be a sequence such that $Q_{\mu _{0,m},\sigma _{m}^{2}\Sigma
_{m},\zeta _{m},U_{m}}(W_{m})$ differs from the multiple supremum in (\ref%
{supQ2}) only by a null sequence. Then we have $P_{\mu _{0,m},\sigma
_{m}^{2}\Sigma _{m}}(W_{m})\rightarrow 0$ as $m\rightarrow \infty $ as a
consequence of (\ref{supP2}). Furthermore,%
\begin{eqnarray}
P_{\mu _{0,m},\sigma _{m}^{2}\Sigma _{m}}(W_{m}) &=&\Pr (\mu _{0,m}+\sigma
_{m}\Sigma _{m}^{1/2}\mathbf{G}\in W_{m})=\Pr (\mu _{0,m}+\sigma _{m}\Sigma
_{m}^{1/2}\mathbf{G/}\left\Vert \mathbf{G}\right\Vert \in W_{m})
\label{PPPP} \\
&=&\Pr (U_{m}^{\prime }\mathbf{G/}\left\Vert \mathbf{G}\right\Vert \in
B_{m})=\Pr (\mathbf{G/}\left\Vert \mathbf{G}\right\Vert \in
B_{m})=\int_{S^{n-1}}\mathbf{1}_{B_{m}}(s)d\upsilon _{S^{n-1}}(s)  \notag
\end{eqnarray}%
where $B_{m}=U_{m}^{\prime }\sigma _{m}^{-1}\Sigma _{m}^{-1/2}(W_{m}-\mu
_{0,m})$, where we have used $G(\mathfrak{M}_{0})$-invariance in the second
step, and where $\mathbf{G}$ is a standard Gaussian $n\times 1$ random
vector. Similarly, 
\begin{eqnarray*}
Q_{\mu _{0,m},\sigma _{m}^{2}\Sigma _{m},\zeta _{m},U_{m}}(W_{m}) &=&\Pr
(\mu _{0,m}+\sigma _{m}\Sigma _{m}^{1/2}U_{m}\mathbf{z}_{m}\in W_{m}) \\
&=&\Pr (\mu _{0,m}+\sigma _{m}\Sigma _{m}^{1/2}U_{m}\mathbf{z}_{m}\mathbf{/}%
\left\Vert \mathbf{z}_{m}\right\Vert \in W_{m}) \\
&=&\Pr (\mathbf{z}_{m}\mathbf{/}\left\Vert \mathbf{z}_{m}\right\Vert \in
B_{m})=\int_{S^{n-1}}\mathbf{1}_{B_{m}}(s)\bar{h}_{\zeta _{m}}(s)d\upsilon
_{S^{n-1}}(s),
\end{eqnarray*}%
where $\mathbf{z}_{m}$ is a random vector with distribution $\zeta _{m}$.
Let $m^{\prime }$ be an arbitrary subsequence of $m$. Since a fortiori the
integral in (\ref{PPPP}) converges to zero along the subsequence $m^{\prime
} $, there exists a subsequence $m^{\prime \prime }$ of $m^{\prime }$ and a $%
\upsilon _{S^{n-1}}$-null set $A_{1}\subseteq S^{n-1}$ such that $\mathbf{1}%
_{B_{m^{\prime \prime }}}(s)\rightarrow 0$ for every $s\in S^{n-1}\backslash
A_{1}$ (e.g., Theorem 3.12 in \cite{rudin}). Since the envelope $\bar{h}%
^{\ast }$ is $\upsilon _{S^{n-1}}$-integrable, it is finite everywhere,
except possibly on a $\upsilon _{S^{n-1}}$-null set $A_{2}\subseteq S^{n-1}$%
. By assumption, $0\leq \bar{h}_{\zeta _{m^{\prime \prime }}}(s)\leq \bar{h}%
^{\ast }(s)$ holds for all $m^{\prime \prime }$ outside of a $\upsilon
_{S^{n-1}}$-null set $A_{3}\subseteq S^{n-1}$. Set $A=A_{1}\cup A_{2}\cup
A_{3}$. But then $\mathbf{1}_{B_{m^{\prime \prime }}}(s)\bar{h}_{\zeta
_{m^{\prime \prime }}}(s)\rightarrow 0$ holds for every $s$ outside the $%
\upsilon _{S^{n-1}}$-null set $A$. Furthermore, $\left\vert \mathbf{1}%
_{B_{m^{\prime \prime }}}(s)\bar{h}_{\zeta _{m^{\prime \prime
}}}(s)\right\vert \leq \bar{h}^{\ast }(s)$ holds for every $m^{\prime \prime
}$ and for every $s\notin A$. But then the Dominated Convergence Theorem
gives 
\begin{equation*}
\Pr (\mathbf{z}_{m^{\prime \prime }}\mathbf{/}\left\Vert \mathbf{z}%
_{m^{\prime \prime }}\right\Vert \in B_{m^{\prime \prime }})=\int_{S^{n-1}}%
\mathbf{1}_{B_{m^{\prime \prime }}}(s)\bar{h}_{\zeta _{m^{\prime \prime
}}}(s)d\upsilon _{S^{n-1}}(s)\rightarrow 0.
\end{equation*}%
Since the subsequence $m^{\prime }$ was arbitrary, we conclude that $Q_{\mu
_{0,m},\sigma _{m}^{2}\Sigma _{m},\zeta _{m},U_{m}}(W_{m})$, and hence (\ref%
{supQ2}), converges to zero. $\blacksquare $

Since all rejection regions in the size control results in this paper are $G(%
\mathfrak{M}_{0})$-invariant, it is now obvious how this lemma can be used
to transfer the size control results for the Gaussian case to the setup
considered here. We abstain from spelling out the details.

\begin{remark}
(i) Restricted to $G(\mathfrak{M}_{0})$-invariant rejection regions, Lemma %
\ref{Lcont2} is indeed a special case of Lemma \ref{Lcont3}. This follows
since Lemma D.1 in \cite{PP17} applied to $h\in H$ shows that $\mathbf{z}%
/\left\Vert \mathbf{z}\right\Vert $ has a density $\bar{h}$; furthermore,
applying that lemma once again, but now to $h^{\ast }$ (which is -- up to a
normalization -- a probability density), produces an $\upsilon _{S^{n-1}}$%
-integrable envelope $\bar{h}^{\ast }$ for the collection of densities $\bar{%
h}$.

(ii) Similar to Remark \ref{cross} one can also allow for cross-restrictions
between $(\beta ,\sigma ^{2},\Sigma )$ and $(h,U)$ ($(\zeta ,U)$,
respectively) here.

(iii) The exact size control results can also be generalized beyond
Gaussianity and elliptical symmetry under appropriate assumptions, but we do
not discuss this here.
\end{remark}

\section{Appendix\label{App F}: Description of the algorithms}

Subsequently, for a symmetric and positive definite matrix $\Sigma $ we
denote by $\mathrm{chol}(\Sigma )$ the (unique) lower triangular matrix that
satisfies $\Sigma =\mathrm{chol}(\Sigma )\mathrm{chol}(\Sigma )^{\prime }$.

\begin{center}
\begin{algorithm}
\caption{Numerical approximation of $C_{\Diamond}(\alpha, p) =
\sup_{\rho \in (-1, 1)^p} F_{\rho}^{-1}(1-\alpha)$.} \label{alg:AR} 
\begin{algorithmic}[1]
\State \textbf{Input} Positive integers $M_0 \geq M_1 \geq M_2$, $N_0 \leq  N_1 \leq N_2$.
\State \textbf{Stage 0: Initial value search}
\State Generate a pseudorandom sample $Z_{1}, \hdots, Z_{N_0}$ from $P_{0, I_n}$. 
\For{$j = 1$ \texttt{to} $j = M_0$} 
\State Obtain a candidate $\rho_j \in (-1, 1)^p$. 
\State Compute $\tilde{F}_{j}^{-1}(1-\alpha)$ where $\tilde{F}_j(x) = N_0^{-1} \sum_{i = 1}^{N_0} \mathbf{1}_{(-\infty, x]}(T_w(\mu_0 + \mathrm{chol}(\Sigma(f_{\rho_j}))Z_{i}))$ for $x \in \mathbb{R}$.
\EndFor
\State Rank the candidates $\rho_j$ according to the value (from largest to smallest) of the corresponding quantities $\tilde{F}_{j}^{-1}(1-\alpha)$ to obtain $\rho_{1:M_0}, \hdots, \rho_{M_1:M_0}$, the initial values for the next stage.
\State \textbf{Stage 1: Coarse localized optimizations}
\For{$j = 1$ \texttt{to} $j = M_1$} 
\State Generate a pseudorandom sample $Z_{1}, \hdots, Z_{N_1}$ from $P_{0, I_n}$. 
\State Let $\bar{F}_{j, \rho}(x) = N_1^{-1} \sum_{i = 1}^{N_1} \mathbf{1}_{(-\infty, x]}(T_w(\mu_0 + \mathrm{chol}(\Sigma(f_{\rho}))Z_{i}))$ for $x \in \mathbb{R}$ and $\rho \in (-1, 1)^p$.
\State Obtain $\rho^*_{j}$ by running a numerical optimization algorithm for the problem $\sup_{\rho \in (-1, 1)^p} \bar{F}^{-1}_{j, \rho}(1-\alpha)$ initialized at $\rho_{j:M_0}$.
\EndFor
\State Rank the obtained numbers $\rho^*_{j}$ according to the value (from largest to smallest) of the corresponding $\bar{F}^{-1}_{j, \rho^*_{j}}(1-\alpha)$ to obtain $\rho^*_{1:M_1}, \hdots, \rho^*_{M_2:M_1}$, the initial values for the next stage.
\State \textbf{Stage 2: Refined localized optimization}
\For{$j = 1$ \texttt{to} $j = M_2$}
\State Generate a pseudorandom sample $Z_{1}, \hdots, Z_{N_2}$ from $P_{0, I_n}$. 
\State Let $\bar{\bar{F}}_{j, \rho}(x) = N_2^{-1} \sum_{i = 1}^{N_2} \mathbf{1}_{(-\infty, x]}(T_w(\mu_0 + \mathrm{chol}(\Sigma(f_{\rho}))Z_{i}))$ for $x \in \mathbb{R}$ and $\rho \in (-1, 1)^p$.
\State Obtain $\rho^{**}_{j}$ by running a numerical optimization algorithm for the problem $\sup_{\rho \in (-1, 1)^p} \bar{\bar{F}}^{-1}_{j, \rho}(1-\alpha)$ initialized at $\rho^*_{j:M_1}$. 
\EndFor
\State \textbf{Return} $\max_{j = 1, \hdots, M_2} \bar{\bar{F}}^{-1}_{j, \rho^{**}_{j}}(1-\alpha)$. 
\end{algorithmic}
\end{algorithm}\pagebreak
\end{center}

\begin{remark}
\emph{(Other test statistics)} Algorithm \ref{alg:AR} has been formulated
for the test statistic $T_{w}$, but clearly works as it stands also for any
other test statistic $T$ (upon replacing $T_{w}$ by $T$) provided $T$
satisfies the conditions given in Remark \ref{other} in Section \ref%
{numerical}.
\end{remark}

\begin{remark}
\label{rem:gencand} \emph{(Generation of candidates in Stage 0 of Algorithm %
\ref{alg:AR})} In Stage 0 of Algorithm \ref{alg:AR} candidates $\rho _{j}\in
(-1,1)^{p}$ need to be obtained, from which the best ones (in the sense of
giving the highest values of the (numerically approximated) objective
function) are then used in Stage 1 as starting values for a numerical
optimization procedure. The main purpose of Stage 0 is to decrease the risk
of running into a local, but not global, optimum. Different approaches can
be used to generate these candidates: If the autoregressive order $p$ is
relatively small, one could use an equally spaced grid of starting values.
Since the dimension of the feasible set is growing linearly in $p$, this is
not feasible for moderate to large autoregressive orders, and in particular
not feasible for the case where $\mathfrak{F}=\mathfrak{F}_{\mathrm{AR(}n-1%
\mathrm{)}}$ (which is equivalent to using $\mathfrak{F}=\mathfrak{F}_{%
\mathrm{all}}$) for typical sample sizes $n$. In such cases, one can
generate the candidates by drawing partial autocorrelation coefficients from
a distribution that induces a uniform distribution on the stationarity
region of $\mathrm{AR(}p\mathrm{)}$ coefficients as described in \cite{Jones}%
. One could also think of many variants of this approach that are designed
to more thoroughly exhaust subsets of the feasible set corresponding to
lower-dimensional autoregressive coefficient vectors. One variant is as
follows: One generates starting values on $(-1,1)^{p_{1}},\ldots
,(-1,1)^{p_{l}}$, respectively, for $1\leq p_{1}<p_{2}<\ldots <p_{l}=p$
following the method described in \cite{Jones}, and then converts the vector
of partial autocorrelation coefficients in $(-1,1)^{p_{i}}$ for $i=1,\ldots
,l-1$ into a vector of partial autocorrelation coefficients in $(-1,1)^{p}$
by setting the remaining $p-p_{i}$ coefficients to $0$.
\end{remark}

\begin{remark}
\emph{(Numerical optimization in Stages 1 and 2 of Algorithm \ref{alg:AR}) }%
Since numerical computation of derivatives would be computationally
intensive, we use derivative-free optimization methods as, e.g., variants of
the \cite{nelder} algorithm in Stages 1 and 2.
\end{remark}

\begin{remark}
\label{rem:restrpart} \emph{(Critical values if }$\mathfrak{F}\ $\emph{is a
subset of }$\mathfrak{F}_{\mathrm{AR(}p\mathrm{)}}$\emph{\ described by
restrictions on the partial autocorrelation coefficients)} Suppose it is
desired to solve the problem (\ref{eqn:def}) with $\mathfrak{F}$ a subset of 
$\mathfrak{F}_{\mathrm{AR(}p\mathrm{)}}$ described by restrictions on the
partial autocorrelation coefficients such as $\rho \in (-1+\varepsilon
,1-\varepsilon )^{p}$ for some $0<\varepsilon <1$. Algorithm \ref{alg:AR}
can easily be modified to accommodate such a situation, by choosing
candidates in Stage 0 from the set $(-1+\varepsilon ,1-\varepsilon )^{p}$
(e.g., by suitably modified versions of the procedures discussed in Remark %
\ref{rem:gencand}), and by solving the optimization problems in Stages 1 and
2 over the set $(-1+\varepsilon ,1-\varepsilon )^{p}$ instead of $(-1,1)^{p}$%
. The so-obtained critical value then, of course, additionally also depends
on $\varepsilon $.
\end{remark}

\begin{remark}
\emph{(Computation of Cholesky factorization)} As is well known, $\mathrm{%
chol}(\Sigma (f_{\rho }))$ can be efficiently obtained from the partial
autocorrelation coefficients through a variant of the Durbin-Levinson
recursion (e.g., by combining \cite{porat} Table 6.2 on p. 159 and Theorem
2.13 in the same reference).\pagebreak
\end{remark}

\begin{center}
\begin{algorithm}
\caption{Numerical approximation of $\sup_{\mu_0 \in
\mathfrak{M}_0} \sup_{0 < \sigma^2 < \infty} \sup_{f \in
\mathfrak{F}_{\mathrm{AR(p)}}} P_{\mu_0, \sigma^2 \Sigma(f)} \left( T_w \geq
C \right)$.} \label{alg:size} 
\begin{algorithmic}[1]
	\State \textbf{Input} A real number $C$ and positive integers $M_0 \geq M_1 \geq M_2$, $N_0 \leq  N_1 \leq N_2$.
		\State \textbf{Stage 0: Initial value search}
                                \State Generate a pseudorandom sample $Z_{1}, \hdots, Z_{N_0}$ from $P_{0, I_n}$.
		\For{$j = 1$ \texttt{to} $j = M_0$} 
		\State Obtain a candidate $\rho_j \in (-1, 1)^p$. 
		\State Compute $\tilde{p}_j = N_0^{-1} \sum_{i = 1}^{N_0} \mathbf{1}_{[C, \infty)}(T_w(\mu_0 + \mathrm{chol}(\Sigma(f_{\rho_j}))Z_{i}))$.
		\EndFor
		\State Rank the candidates $\rho_j$ according to the value (from largest to smallest) of the corresponding quantities $\tilde{p}_j $ to obtain $\rho_{1:M_0}, \hdots, \rho_{M_1:M_0}$, the initial values for the next stage.
		\State \textbf{Stage 1: Coarse localized optimizations}
		\For{$j = 1$ \texttt{to} $j = M_1$} 
		\State Generate a pseudorandom sample $Z_{1}, \hdots, Z_{N_1}$ from $P_{0, I_n}$. 
		\State Let $\bar{p}_{j, \rho} = N_1^{-1} \sum_{i = 1}^{N_1} \mathbf{1}_{[C, \infty)}(T_w(\mu_0 + \mathrm{chol}(\Sigma(f_{\rho}))Z_{i}))$ for $\rho \in (-1, 1)^p$.
		\State Obtain $\rho^*_{j}$ by running a numerical optimization algorithm for the problem $\sup_{\rho \in (-1, 1)^p} \bar{p}_{j, \rho}$ initialized at $\rho_{j:M_0}$.
		\EndFor
		\State Rank the obtained numbers $\rho^*_{j}$ according to the value (from largest to smallest) of the corresponding $\bar{p}_{j, \rho^*_{j}}$ to obtain $\rho^*_{1:M_1}, \hdots, \rho^*_{M_2:M_1}$, the initial values for the next stage.
		\State \textbf{Stage 2: Refined localized optimization}
		\For{$j = 1$ \texttt{to} $j = M_2$}
		\State Generate a pseudorandom sample $Z_{1}, \hdots, Z_{N_2}$ from $P_{0, I_n}$. 
		\State Let $\bar{\bar{p}}_{j, \rho} = N_2^{-1} \sum_{i = 1}^{N_2} \mathbf{1}_{[C, \infty)}(T_w(\mu_0 + \mathrm{chol}(\Sigma(f_{\rho}))Z_{i}))$ for $\rho \in (-1, 1)^p$.
		\State Obtain $\rho^{**}_{j}$ by running a numerical optimization algorithm for the problem $\sup_{\rho \in (-1, 1)^p} \bar{\bar{p}}_{j, \rho}$ initialized at $\rho^*_{j:M_1}$. 
		\EndFor
		\State \textbf{Return} $\max_{j = 1, \hdots, M_2} \bar{\bar{p}}_{j, \rho^{**}_{j}} $.
	\end{algorithmic}
\end{algorithm}
\end{center}

\begin{remark}
Remarks analogous to the ones given after Algorithm \ref{alg:AR} also apply
to Algorithm \ref{alg:size}.
\end{remark}

\section{Appendix\label{App G}: Numerically checking Condition (\protect\ref%
{non-incl}) for the test problems in Section \protect\ref{sec:macro}}

First observe that for the hypothesis considered we have that $\mathfrak{M}%
_{0}^{lin}$ is the span of $e_{+}$ and $(1,2,3,\ldots ,100)^{\prime }$, and
hence coincides with the span of $(E_{100,0}(0),E_{100,1}(0))$. It follows
that $\rho (0,\mathfrak{M}_{0}^{lin})=2$, since $\func{span}(E_{100,2}(0))=%
\func{span}((1,2^{2},3^{2},\ldots ,100^{2})^{\prime })$ is not contained in $%
\mathfrak{M}_{0}^{lin}$. We now show that $\rho (\gamma ,\mathfrak{M}%
_{0}^{lin})=0$ for every $\gamma \in (0,\pi ]$: Note that $\func{rank}%
(E_{100,0}(0),E_{100,1}(0),E_{100,0}(\gamma ))=\min (100,4)=4$ if $0<\gamma
<\pi $ and equals $\min (100,3)=3$ if $\gamma =\pi $, because of Lemma \ref%
{fullrank} in Appendix \ref{App C}. But this implies that the span of $%
E_{100,0}(\gamma )$ is not contained in $\mathfrak{M}_{0}^{lin}$ whenever $%
\gamma \in (0,\pi ]$, establishing the claim.

Second, we now verify that condition (\ref{non-incl}) holds for each of the
128 design matrices considered. That, is we show that $\func{span}%
(E_{100,\rho (\gamma ,\mathfrak{M}_{0}^{lin})}(\gamma ))\not\subseteq 
\mathsf{B}$ for every $\gamma \in \lbrack 0,\pi ]$ where $\rho (\gamma ,%
\mathfrak{M}_{0}^{lin})=2$ for $\gamma =0$ and $\rho (\gamma ,\mathfrak{M}%
_{0}^{lin})=0$ for $\gamma \neq 0$. Consider first the case where $\gamma =0$%
: Since $\func{span}(E_{100,2}(0))=\func{span}((1,2^{2},3^{2},\ldots
,100^{2})^{\prime })$ it suffices to show that $B((1,2^{2},3^{2},\ldots
,100^{2})^{\prime })$ is nonzero where $B$ is given in (\ref{Def_B}). Since $%
\func{rank}(R(X^{\prime }X)^{-1}X^{\prime })=q=1$ holds, it is in turn
sufficient to show that each coordinate of the residual vector obtained from
regressing $(1,2^{2},3^{2},\ldots ,100^{2})^{\prime }$ onto the design
matrix $X$ is nonzero. For each of the 128 design matrices considered this
has been numerically confirmed. Next, to check condition (\ref{non-incl})
for $\gamma \in (0,\pi ]$ it suffices to verify that for each of the 128
cases $c(\gamma )\notin \mathsf{B}$ for every $\gamma \in (0,\pi ]$ where $%
c(\gamma )=(\cos (\gamma ),\cos (2\gamma ),\ldots ,\cos (100\gamma
))^{\prime }$ is the first column of $E_{100,0}(\gamma )$. Since $c(\gamma
)\in \mathsf{B}$ is equivalent to $B(c(\gamma ))=0$ (as $q=1$ holds), we
compute for each of the $128$ cases the function 
\begin{equation}
\gamma \mapsto \Vert B(c(\gamma ))\Vert _{\infty }\quad \text{ where }\quad
c(\gamma )=(\cos (\gamma ),\cos (2\gamma ),\ldots ,\cos (100\gamma
))^{\prime }  \label{eqn:deteval}
\end{equation}%
(a) on a grid of $100~000$ equally spaced points in $[0,\pi ]$ and (b) also
on a grid of $100~000$ equally spaced points in $[0,10^{-6}]$ to get a more
refined resolution in this region (note that the function in the preceding
display has a root at $0$ in each of the 128 cases since $\rho (0,\mathfrak{M%
}_{0}^{lin})>0$). Then we plot for each value $\gamma $ in these grids the
smallest of the 128 norms in (\ref{eqn:deteval}) corresponding to the 128
cases considered. These plots are shown in Figures 3(a) and 3(b),
respectively. The figures suggest that $\Vert B(c(\gamma ))\Vert \neq 0$
holds for $\gamma \in (0,\pi ]$ for each of the 128 cases, implying that for
each of these cases we have $\func{span}(E_{100,0}(\gamma ))\notin \mathsf{B}
$ for $\gamma \in (0,\pi ]$.

\begin{figure}[tbp]
\centering
\includegraphics[scale = .4]{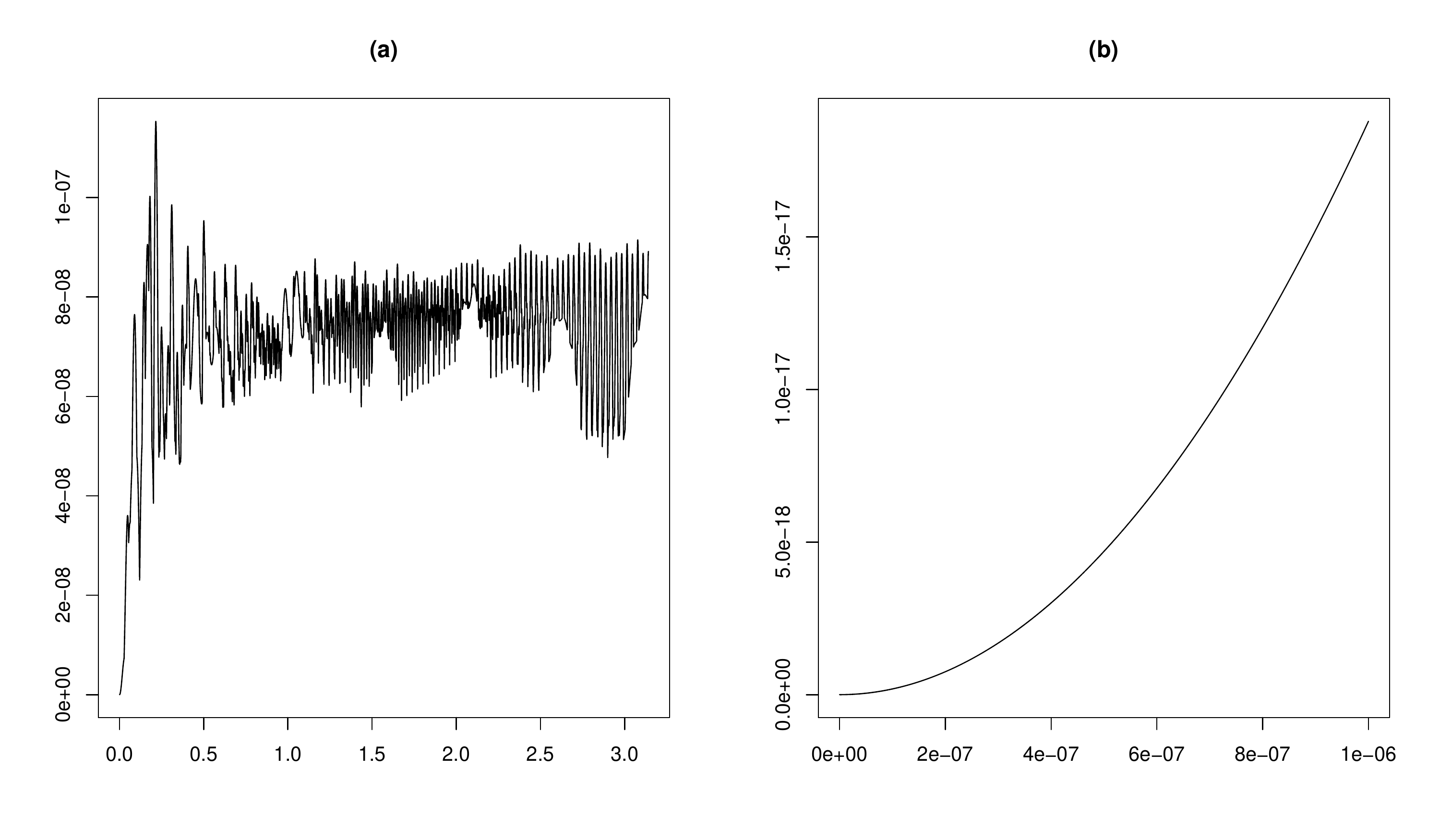}
\caption{(a) The point-wise minimum over all 128 regression models of the
function defined in (\protect\ref{eqn:deteval}) on the interval $[0,\protect%
\pi ]$. (b) The point-wise minimum over all 128 regression models of the
function defined in (\protect\ref{eqn:deteval}) on the interval $[0,10^{-6}]$%
.}
\end{figure}

%

\section{Appendix\label{App H}: Settings for Algorithms \protect\ref{alg:AR}
and \protect\ref{alg:size} used in Section \protect\ref{sec:macro}}

Here we describe the settings for Algorithms \ref{alg:AR} and \ref{alg:size}
used in the computations for Section \ref{sec:macro}. The actual
computations were performed using the implementations of Algorithms \ref%
{alg:AR} and \ref{alg:size} provided in the \textsf{R}-package \textbf{acrt}
(\cite{csart}) through the functions \textsf{critical.value} and \textsf{size%
}, respectively. Initial values in Stage 0 are generated as follows, cf.
Remark \ref{rem:gencand} in Appendix \ref{App F}:

\begin{enumerate}
\item For $p\in \{1,2\}$ the initial values are a pseudorandom sample of
size $5000$ in $(-1,1)^{p}$ drawn according to the distribution generating a
uniform distribution on the stationarity region of $\mathrm{AR(}p\mathrm{)}$
coefficients following Jones (1987).

\item For $p\in \{5,10,25,50,99\}$ we proceed as follows: For every $%
p_{i}\in \left\{ 2,5,10,25,50,99\right\} $ that does not exceed $p$, we
generate pseudorandom samples in $(-1,1)^{p_{i}}$ of size $5000$, each
according to the distribution that generates a uniform distribution on the
stationarity region of $\mathrm{AR(}p_{i}\mathrm{)}$ parameters. Then $%
\mathrm{AR(}p\mathrm{)}$ partial autocorrelation coefficients are obtained
from these pseudorandom samples in $(-1,1)^{p_{i}}$ by appending $p-p_{i}$
zeros.
\end{enumerate}

Furthermore we use $M_{1}=10$, $M_{2}=2$ and $N_{0}=1000$, $N_{1}=10~000$, $%
N_{2}=50~000$. Note that the number of replications in the Monte-Carlo
algorithms of $10~000$ in the first stage, and of $50~000$ in the second
stage are (at least) of the same order of magnitude as the number of
replications used in contemporary simulation studies concerning rejection
probabilities of autocorrelation robust tests (e.g., \cite{SPJ11} use 10 000
replications). [Of course, in a particular application, where one needs to
determine the critical value for \emph{one} model and for \emph{one}
parameter $p$ only, one could choose parameters $M_{i}$ and $N_{i}$ that
provide an even higher level of accuracy.] The optimization algorithm we
employ in Stages 1 and 2 is a Nelder-Mead algorithm with default settings
concerning the reflection, contraction, and expansion coefficients as
implemented in the \textsf{optim} function in \textsf{R}. The relative
convergence tolerance was set equal to $N_{1}^{-1/2}$ and to $N_{2}^{-1/2}$
in Stages 1 and 2, respectively. Furthermore, the maximal number of
iterations in Stages 1 and 2 were set equal to $20n$ and $30n$, respectively
(recall $n=100$ here). Since the \textsf{optim} function supplies an
implementation of the Nelder-Mead algorithm that optimizes functions from a
Euclidean space to the real numbers, we rephrased our optimization problems
as unrestricted optimization problems over $\mathbb{R}^{p}$ using the
function $(2/\pi )\arctan $. For the i.i.d. case (i.e., $p=0$) the problem
considerably simplifies as noted earlier, since the distribution of the test
statistic under the null does then not depend on any nuisance parameter; in
this case the maximal rejection probabilities and the $1-\alpha $ quantiles
were in each scenario obtained from a Monte Carlo sample of size $50~000$.

\section{Appendix\label{App I}: A random-walk-based critical value}

\begin{figure}[tbp]
\centering
\includegraphics[scale = .4]{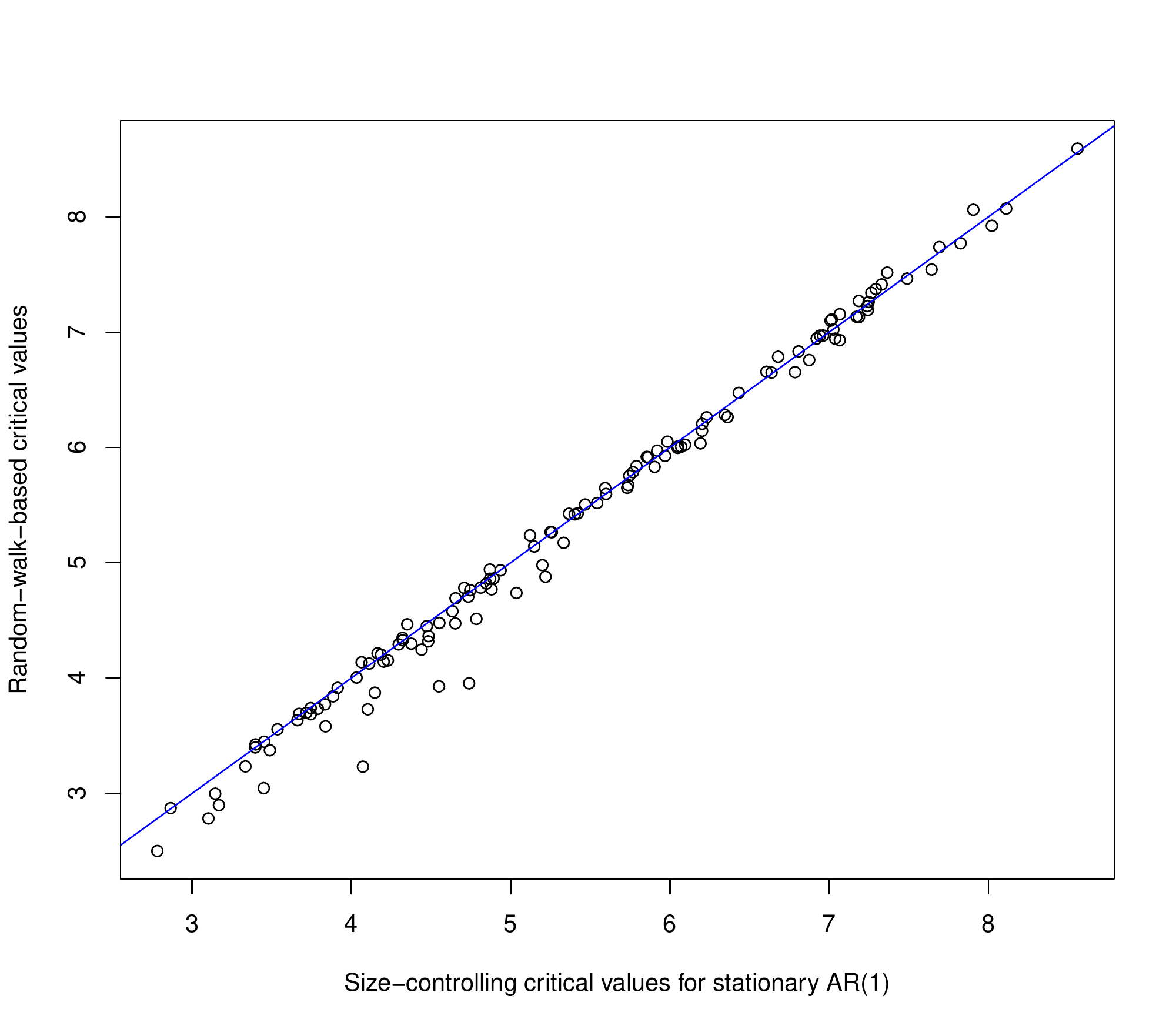}
\caption{Scatterplot of random-walk-based critical values versus critical
values controlling size over $\mathfrak{F}_{\mathrm{AR(}1\mathrm{)}}$, i.e.,
over the class of stationary $\mathrm{AR(}1\mathrm{)}$ processes. The line
represents the graph of the identity function.}
\label{fig:random_walk_vs_AR(1)}
\end{figure}

\begin{figure}[tbp]
\centering
\includegraphics[scale = .4]{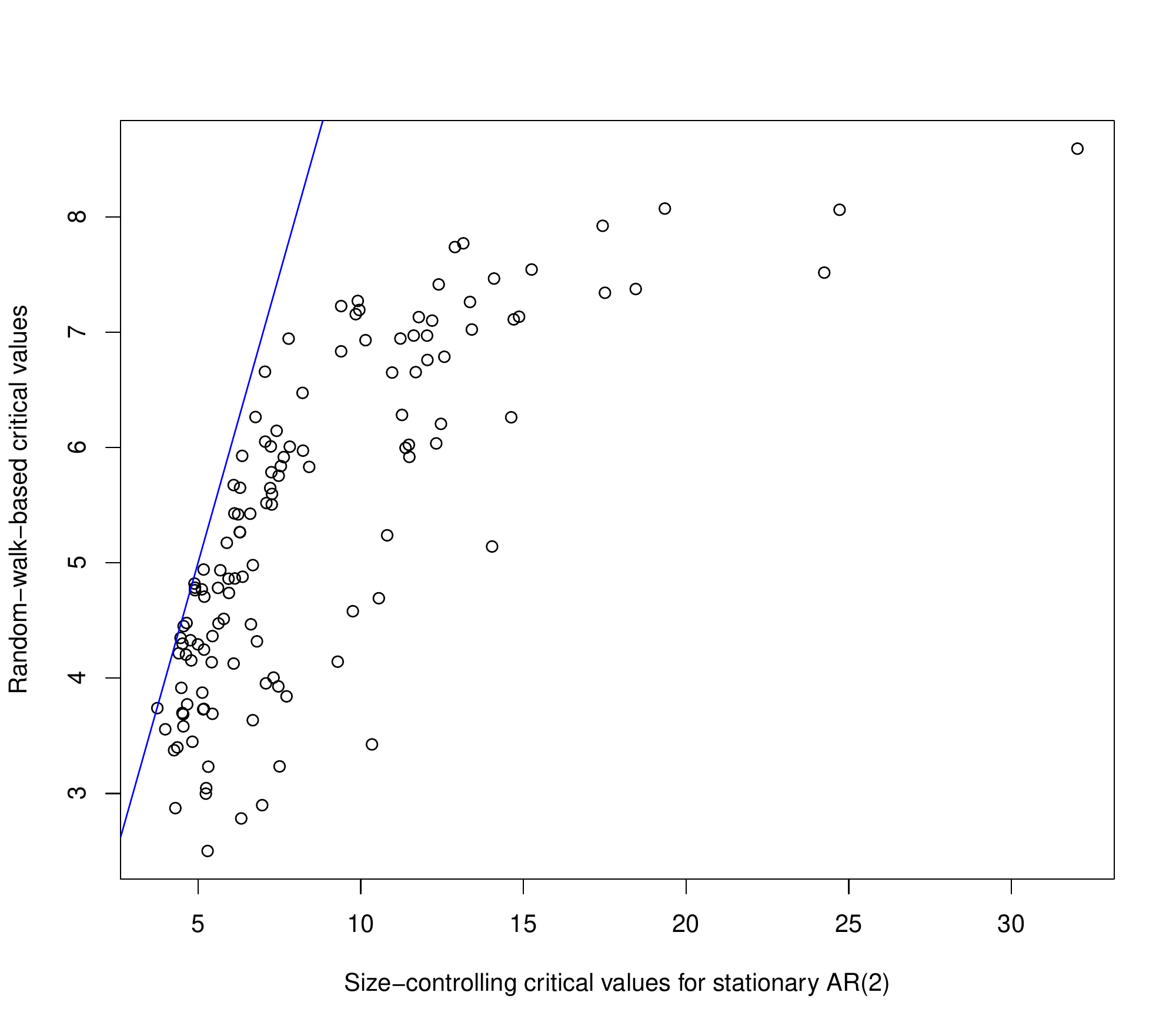}
\caption{Scatterplot of random-walk-based critical values versus critical
values controlling size over $\mathfrak{F}_{\mathrm{AR(}2\mathrm{)}} $,
i.e., over the class of stationary $\mathrm{AR(}2\mathrm{)}$ processes. The
line represents the graph of the identity function.}
\label{fig:random_walk_vs_AR(2)}
\end{figure}

As suggested by a referee, we compute -- in the context of the null
hypotheses and the $128$ models considered in Section \ref{sec:macro} -- the
critical values for $\left\vert t_{w}\right\vert $ that result from the
presumption that the errors follow a Gaussian random walk and where we use
the same weights as in Section \ref{sec:macro}. Since we may set $\beta
_{1}=\beta _{2}=\beta _{3}=0$ and the innovation variance $\sigma ^{2}=1$ by
invariance, the computation of such a critical value reduces to determining
the $(1-\alpha )$-quantile of the distribution of $\left\vert
t_{w}\right\vert $ under a fixed Gaussian distribution, which can easily be
achieved by Monte-Carlo. For $\alpha =0.05$ and for each of the $128$ models
considered in Section \ref{sec:macro} we have computed these critical values
from $10~000$ Monte-Carlo samples. We show a scatter-plot of these
random-walked-based critical values versus the critical values that control
size over $\mathfrak{F}_{\mathrm{AR(}1\mathrm{)}}$ (which have been obtained
in Section \ref{sec:macro}) in Figure \ref{fig:random_walk_vs_AR(1)},
whereas in Figure \ref{fig:random_walk_vs_AR(2)} we plot these versus the
critical values that control size over $\mathfrak{F}_{\mathrm{AR(}2\mathrm{)}%
}$ (which also have been obtained in Section \ref{sec:macro}). We can draw
the following conclusions:

1. The random-walk-based critical values are -- for the majority of the $128$
models considered -- roughly of the same magnitude as the critical values
that guarantee size control over $\mathfrak{F}_{\mathrm{AR(}1\mathrm{)}}$,
although for some models they are too small.

2. Compared to the critical values guaranteeing size control over $\mathfrak{%
F}_{\mathrm{AR(}2\mathrm{)}}$ the random-walk-based critical values are way
too small, and hence will not control size over $\mathfrak{F}_{\mathrm{AR(}2%
\mathrm{)}}$; and a fortiori not over $\mathfrak{F}_{\mathrm{AR(}p\mathrm{)}%
} $ with $p\geq 2$ or $\mathfrak{F}_{\mathrm{all}}$.

\begin{figure}[tbp]
\centering
\includegraphics[scale = .4]{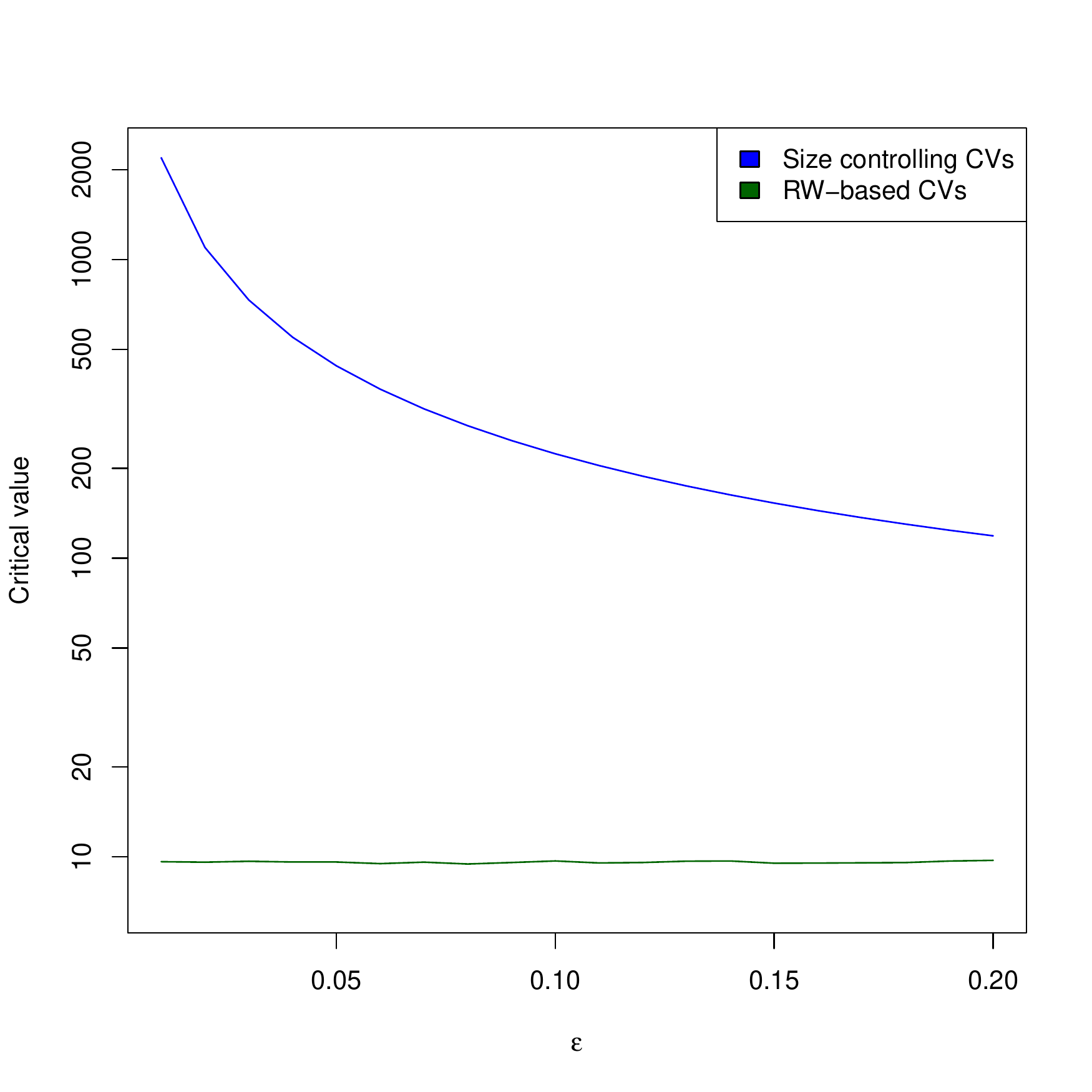}
\caption{Lower curve shows random-walk-based critical values as a function
of $\protect\varepsilon $. Upper curve shows critical values controlling
size over $\mathfrak{F}_{\mathrm{AR(}1\mathrm{)}}$, i.e., over the class of
stationary $\mathrm{AR(}1\mathrm{)}$ processes, as a function of $\protect%
\varepsilon $.}
\label{fig:perturbedloc}
\end{figure}

As a consequence, the random-walk-based critical values are certainly no
substitute for the size-controlling critical values whenever one wants to
allow for correlation structures richer than stationary $\mathrm{AR(}1%
\mathrm{)}$. If one is willing to only maintain stationary $\mathrm{AR(}1%
\mathrm{)}$ correlations, the first conclusion above may lead one to believe
that the random-walk-based critical values roughly deliver size control over 
$\mathfrak{F}_{\mathrm{AR(}1\mathrm{)}}$. However, this is not true either
in general. Note that Conclusion 1 above is based only on computations
involving the $128$ models (regressors) considered. It does \emph{not}
generalize to other models (regressors) as is easily seen by the following
example where size control over $\mathfrak{F}_{\mathrm{AR(}1\mathrm{)}}$ is
possible, but where the random-walk-based critical value is much too small.
To this end consider a model with only one regressor given by 
\begin{equation}
X=(1+\varepsilon ,1+\varepsilon ,1-\varepsilon ,1-\varepsilon ,1+\varepsilon
,1+\varepsilon ,\ldots ,1+\varepsilon ,1+\varepsilon )^{\prime },  \label{X}
\end{equation}%
where $n=100$ and where we vary $\varepsilon $ from $0.01$ to $0.20$ in
steps of size $0.01$. The null hypothesis is that the coefficient of the
regressor is equal to zero and again the corresponding test statistic $%
\left\vert t_{w}\right\vert $ with the same weights as before is used. In
Figure \ref{fig:perturbedloc} we present the random-walk-based critical
values (computed from $10~000$ Monte-Carlo samples) as well as the critical
values that control size over $\mathfrak{F}_{\mathrm{AR(}1\mathrm{)}}$
(computed via Algorithm \ref{alg:AR}) as a function of $\varepsilon $.%
\footnote{%
The settings used here for Algorithm \ref{alg:AR} are similar as in Appendix %
\ref{App H} with $p=1$, but with $M_{0}=100$, $M_{1}=M_{2}=1$, and $%
N_{0}=N_{1}=1~000$, $N_{2}=10~000$.} It is apparent from that figure that
the random-walk-based critical values are way to small. As an additional
observation we note that for the location model (i.e., $\varepsilon =0$ in (%
\ref{X})) the random-walk-based critical value can be computed to be $9.6$.
However, as discussed earlier, for the location model with \emph{stationary} 
$\mathrm{AR(}1\mathrm{)}$ errors no critical value exists that leads to size
control for $\left\vert t_{w}\right\vert $ (since size is equal to $1$ for
every choice of critical value). Hence the random-walk-based critical value
is completely misleading in this model.

\section{Appendix\label{App J}: Covariance models corresponding to
starting-value solutions of autoregressive models}

A referee has asked what happens if instead of the covariance model $%
\mathfrak{C}(\mathfrak{F}_{\mathrm{AR(}1\mathrm{)}})$ derived from
stationary autoregressive processes of order $1$ one considers the
covariance model $\mathfrak{C}_{\mathrm{AR(}1\mathrm{)}}^{0}$ generated by $(%
\mathbf{u}_{1},\ldots ,\mathbf{u}_{n})$ where $\mathbf{u}_{t}=\rho \mathbf{u}%
_{t-1}+\varepsilon _{t}$ for $t=1,\ldots ,n$ with starting-value $\mathbf{u}%
_{0}=0$, the innovations $\varepsilon _{t}$ are distributed independently as 
$N(0,\sigma _{\varepsilon }^{2})$, $0<\sigma _{\varepsilon }^{2}<\infty $,
and $\rho \in (-1,1)$.\footnote{%
This model has a problematic aspect to it in that it assigns a special
meaning to the time point $t=0$ which is hardly justifiable; cf. the
discussion in Section 3.2.2 of \cite{PP2016}.} More precisely, $\mathfrak{C}%
_{\mathrm{AR(}1\mathrm{)}}^{0}$ consists of the $n\times n$ matrices with $%
(i,j)$-th entry given by%
\begin{equation*}
\rho ^{\max (i,j)-\min (i,j)}(1-\rho ^{2\min (i,j)})/(1-\rho ^{2}).
\end{equation*}%
It is easy to see that the covariance model $\mathfrak{C}_{\mathrm{AR(}1%
\mathrm{)}}^{0}$ is norm bounded and has no singular limit points (and the
same is true for higher-order analoga) and hence one would, given the
results in \cite{PP2016}, intuitively expect that size control for test
statistics like $T_{w}$ (or equivalently $\left\vert t_{w}\right\vert $) is 
\emph{always} possible. This is indeed the case and follows from Remark \ref%
{bounded away}(ii), cf. also Remark \ref{star}. However, this does not mean
that standard critical values suggested in the literature based on
asymptotic considerations will come close to providing size control over the
model $\mathfrak{C}_{\mathrm{AR(}1\mathrm{)}}^{0}$ (or over higher-order
analoga). In fact, as we shall show below, such critical values will often
be much too small and will lead to considerable size distortions.

\begin{figure}[tbp]
\centering
\includegraphics[scale = .4]{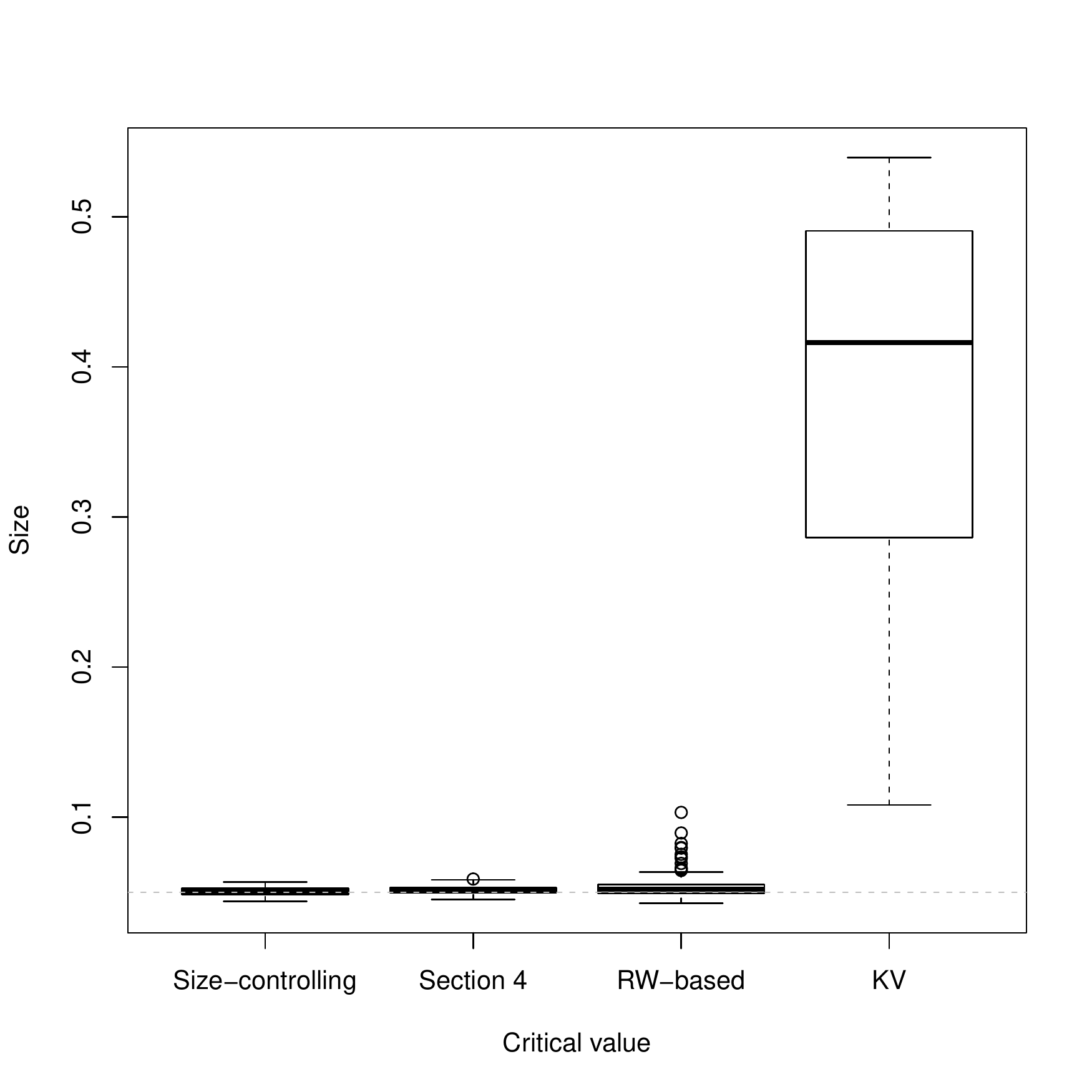}
\caption{The boxplots labeled \textquotedblleft
Size-controlling\textquotedblright , \textquotedblleft Section
4\textquotedblright , \textquotedblleft RW-based\textquotedblright , and
\textquotedblleft KV\textquotedblright , respectively, are each the boxplots
of $\ $the $128$ sizes corresponding to the critical values $c_{1i}$, $c_{2i}
$, $c_{i3}$, and $2.260568$, respectively.}
\label{fig:sizeboxplot}
\end{figure}

We first consider again the $128$ models from Section \ref{sec:macro} with
the only difference that we now assume that the errors follow the
starting-value solution $\mathbf{u}_{t}$ described in the preceding
paragraph. For the test statistic $\left\vert t_{w}\right\vert $ (with the
same weights as in Section \ref{sec:macro}) and for each of the $128$ models
we then numerically computed the critical value $c_{1i}$, $i=1,\ldots ,128$,
that guarantees size control (at $\alpha =0.05$) over the covariance model $%
\mathfrak{C}_{\mathrm{AR(}1\mathrm{)}}^{0}$ by a suitable variant of
Algorithm \ref{alg:AR}. We then consider the following four tests for each
of the 128 models: (i) Reject if $\left\vert t_{w}\right\vert \geq c_{1i}$,
(ii) reject if $\left\vert t_{w}\right\vert \geq c_{2i}$, where $c_{2i}$ is
the critical value computed in Section \ref{sec:macro} (i.e., the critical
value that would control size if the errors were stationary $\mathrm{AR(}1%
\mathrm{)}$-processes), (iii) $\left\vert t_{w}\right\vert \geq c_{3i}$
where $c_{3i}$ is the random-walk-based critical value computed in Appendix %
\ref{App I}, (iv) $\left\vert t_{w}\right\vert \geq 2.260568$
(Kiefer-Vogelsang critical value). For each of the four tests we then
computed the size of the test (using $\mathfrak{C}_{\mathrm{AR(}1\mathrm{)}%
}^{0}$ as the underlying covariance model of course). Note that by
construction of $c_{1i}$ the size of the test in (i) should be $0.05$. Since 
$c_{1i}$ has been determined only numerically, we recomputed the size also
in this case. Also observe that, by construction, the random-walk-based
critical value $c_{3i}$ -- ignoring numerical error -- can not be larger
than the size-controlling critical value $c_{1i}$. The results are shown in
Figure \ref{fig:sizeboxplot} in the form of boxplots, each boxplot
representing the size of one of the tests over the $128$ models.

%

It transpires that $c_{1i}$ indeed controls size (as it should by
construction) and that the critical values $c_{2i}$ (computed to deliver
size control over the related covariance model $\mathfrak{C}(\mathfrak{F}_{%
\mathrm{AR(}1\mathrm{)}})$) pretty much work also under the covariance model
considered here. The random-walk-based critical values do a reasonable job
in most, but not all of the $128$ cases. Finally, the Kiefer-Vogelsang
critical value is seen to be way too small and leads to quite dramatic size
distortions.

\begin{figure}[tbp]
\centering
\includegraphics[scale = .4]{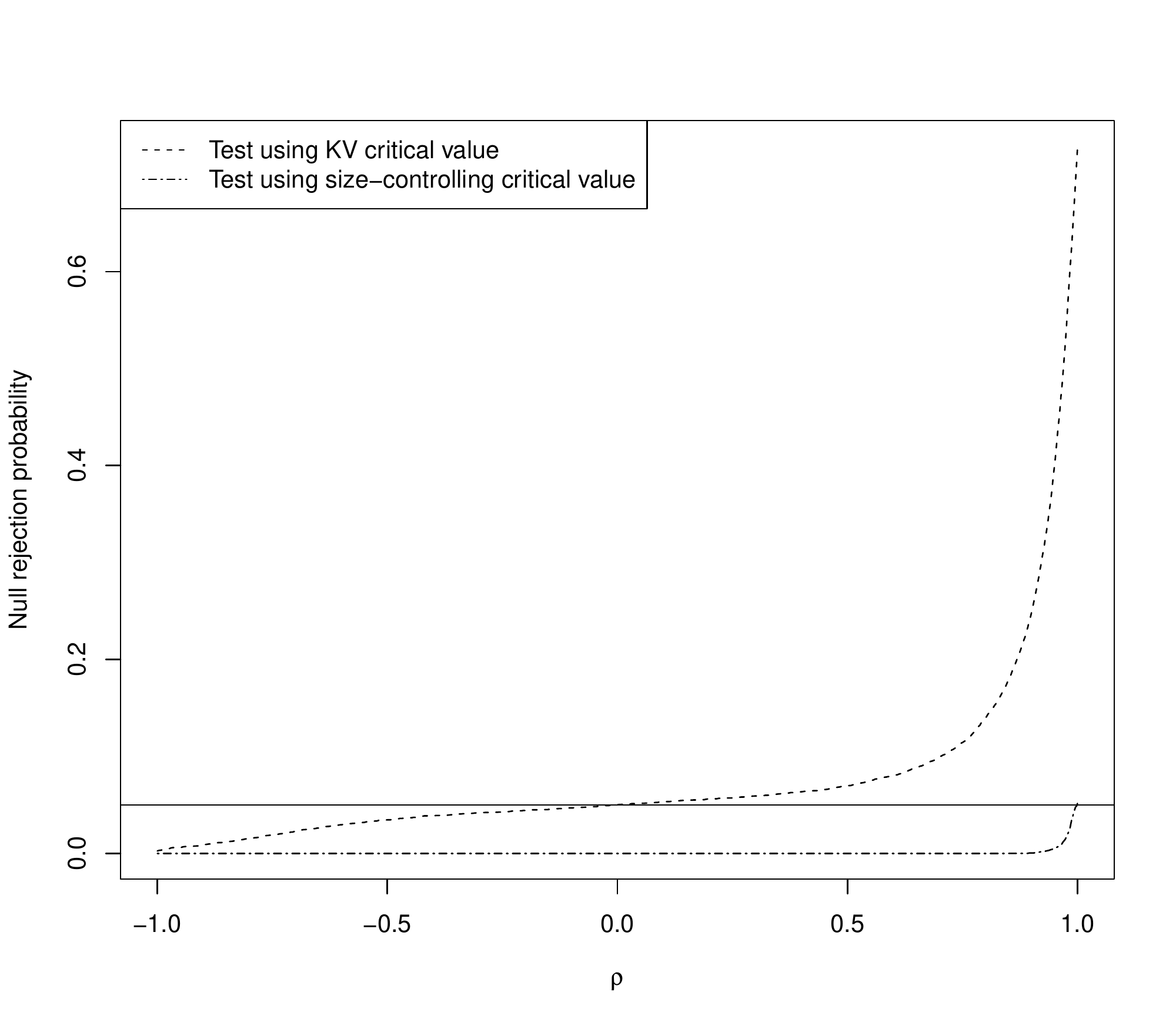}
\caption{Null-rejection probabilities in location model and for $\mathfrak{C}%
_{\mathrm{AR(}1\mathrm{)}}^{0}$. The horizontal line corresponds to $\protect%
\alpha =0.05$.}
\label{fig:ar1locmod}
\end{figure}

We also considered the location model with the errors following the process
as described in the first paragraph. The null hypothesis is that the
location parameter is zero and the test statistic considered is the
corresponding test statistic $\left\vert t_{w}\right\vert $ (with the same
weights as before). While no size-controlling critical value exists for this
problem if the underlying covariance model is $\mathfrak{C}(\mathfrak{F}_{%
\mathrm{AR(}1\mathrm{)}})$ as already noted earlier, this is different if
the covariance model $\mathfrak{C}_{\mathrm{AR(}1\mathrm{)}}^{0}$ is
maintained as is done here (since under $\mathfrak{C}_{\mathrm{AR(}1\mathrm{)%
}}^{0}$ size control is always possible as noted above). We computed this
size-controlling value $c_{1,loc}$ (corresponding to $\alpha =0.05$) by a
suitable variant of Algorithm \ref{alg:AR}. We then computed the
null-rejection probabilities of the tests (i) $\left\vert t_{w}\right\vert
\geq c_{1,loc}$ and (ii) $\left\vert t_{w}\right\vert \geq 2.260568$
(Kiefer-Vogelsang critical value) as a function of $\rho $ (note that these
probabilities do not depend on $\sigma _{\varepsilon }^{2}$). The graphs of
these two functions are given in Figure \ref{fig:ar1locmod}. By construction
of $c_{1,loc}$, the corresponding null-rejection probabilities never exceed $%
\alpha =0.05$, and reach that value at the right endpoint of the parameter
interval for $\rho $. The Kiefer-Vogelsang critical values again lead to
substantial size distortions.

%

Finally, critical values such as $c_{1i}$, $c_{2i}$, or $c_{3i}$ defined
above (and a fortiori the Kiefer-Vogelsang critical value) will typically be
too small to guarantee size control once we pass from $\mathrm{AR(}1\mathrm{)%
}$-models to higher order models.

\section{Appendix\label{App K}: Comments on stochastic regressors}

The assumption of nonstochastic regressors made in the present paper may be
considered restrictive for some applications and has been criticized by
referees. It should be noted, however, that the results obtained under this
assumption are quite strong, in that we are able to obtain exact
finite-sample size results and, in particular, size guarantee results, which
is in contrast to merely (pointwise) asymptotic results in the literature
that carry no size guarantee (in fact, are often plagued by considerable
size distortions). Furthermore, as already mentioned in Section \ref{frame},
the case where $X$ is random and independent of $\mathbf{U}$ can easily be
accommodated by our theory through conditioning on $X$. Of course, the
resulting conditional size control results then a fortiori imply
unconditional size control.\footnote{%
The implied critical values will -- by construction -- typically depend on
the observed values of $X$, and thus will be \textquotedblleft
random\textquotedblright\ critical values when viewed from an unconditional
viewpoint.} Conditioning on $X$ in such a scenario makes perfect sense, as
one can argue that values of $X$ other than the observed ones should be
irrelevant. For a more detailed discussion of arguments supporting
conditional inference see, e.g., \cite{RO1979}. Again, this form of strict
exogeneity (i.e., independence of $X$ and $\mathbf{U}$) may be considered
restrictive for some applications.

We now turn to the case where $X$ is random, with $X$ and $\mathbf{U}$ being
dependent. While our current theory is then not applicable as it stands, we
show by means of numerical examples that the critical values obtained from a
naive application of the algorithm proposed in Section \ref{numerical}
(acting as if $X$ and $\mathbf{U}$ were independent) -- while not guaranteed
to deliver size control -- lead to tests that have better size properties
than tests that are based on standard critical values suggested in the
literature. And this is so despite the fact that the asymptotic theory used
in the literature to justify the latter critical values \emph{is} applicable
to cases where $X$ and $\mathbf{U}$ are dependent. Of course, the subsequent
discussion is based only on a limited Monte-Carlo study and the case of
stochastic regressors needs further study which is in progress.

In a first step we shall now formulate a suitable framework for
investigating size properties of OLS-based tests allowing for dependence
between $X$ and $\mathbf{U}$. To this end consider the linear regression
model $y_{t}=a+bx_{t}+u_{t}$ where $x_{t}$ follows a stationary Gaussian
process and $u_{t}$ is given by the stationary solution of the model $%
u_{t}=\rho u_{t-1}+\varepsilon _{t}$ with $\rho $ varying in $(-1,1)$ and
Gaussian white noise $\varepsilon _{t}$ with positive variance $\sigma
_{\varepsilon }^{2}$. We furthermore assume that $x_{t}$ and $\varepsilon
_{t}$ (and thus $x_{t}$ and $u_{t}$) are jointly Gaussian, ergodic, and
stationary and that $x_{t}$ has positive variance. The probability law of
the process $(x_{t},\varepsilon _{t}/\sigma _{\varepsilon })$ may also
depend on some (possibly infinite dimensional) parameter $\eta $, say, where
we assume that $\eta $ and $\rho $ are variation-free. [We work with $%
(x_{t},\varepsilon _{t}/\sigma _{\varepsilon })$ instead of $%
(x_{t},\varepsilon _{t})$ here, because this conveniently eliminates the
dependence on $\sigma _{\varepsilon }^{2}$.] Autocorrelation robust test
statistics such as $T_{w}$ are based on the ordinary least squares estimator 
$(\hat{a},\hat{b})$. In order for such inference to be meaningful, we need
to ensure at a minimum that this estimator is consistent for $(a,b)$. By the
very notion of consistency this means that $(\hat{a},\hat{b})$ has to
converge to $(a,b)$ in probability regardless of what the true values of $a$%
, $b$, and of the nuisance parameters $\rho $, $\sigma _{\varepsilon }^{2}$,
and $\eta $ are. But as is easy to see, this means that $x_{t}$ and $u_{t}$
have to be uncorrelated whatever the value of $\rho \in (-1,1)$ (and of the
other nuisance parameters) is. This uncorrelatedness condition on $x_{t}$
and $u_{t}$ imposes a restriction on the structure of $x_{t}$ as we show in
the following proposition. Recall that $u_{t}=\sum_{j=0}^{\infty }\rho
^{j}\varepsilon _{t-j}$ holds, and that the uncorrelatedness condition just
mentioned can equivalently be expressed as uncorrelatedness of $x_{t}$ and $%
\sum_{j=0}^{\infty }\rho ^{j}(\varepsilon _{t-j}/\sigma _{\varepsilon })$
(which has the advantage of making the problem free of $\sigma _{\varepsilon
}^{2}$). Let $P_{\eta }$ denote the probability law of the process $%
(x_{t},\varepsilon _{t}/\sigma _{\varepsilon })$ and let $E_{\eta }$ denote
the corresponding expectation operator.

\begin{proposition}
Suppose $E_{\eta }(x_{t}\sum_{j=0}^{\infty }\rho ^{j}(\varepsilon
_{t-j}/\sigma _{\varepsilon }))=0$ holds for every $\rho \in (-1,1)$ and
some (every, respectively) $\eta $. Then $x_{t}$ is independent of $%
(\varepsilon _{t}/\sigma _{\varepsilon }$,$\varepsilon _{t-1}/\sigma
_{\varepsilon },\varepsilon _{t-2}/\sigma _{\varepsilon },\ldots )$ under $%
P_{\eta }$ for this (every, respectively) $\eta $.
\end{proposition}

\textbf{Proof:} By the Gaussianity assumption it suffices to show that $%
E_{\eta }(x_{t}(\varepsilon _{t-j}/\sigma _{\varepsilon }))=0$ for every $%
j\geq 0$ and some (every, respectively) $\eta $. The assumption can be
written as $\sum_{j=0}^{\infty }\rho ^{j}E_{\eta }(x_{t}\varepsilon
_{t-j}/\sigma _{\varepsilon })=0$ for some (every, respectively) $\eta $,
the sum being guaranteed to converge absolutely for $\rho \in (-1,1)$, since
the coefficients $E_{\eta }(x_{t}\varepsilon _{t-j}/\sigma _{\varepsilon })$
are a bounded sequence for given $\eta $ in view of Cauchy-Schwartz. Also
note that the coefficients do not depend on $\rho $ by the
variation-freeness assumption. Hence by the identity theorem for analytic
functions the claim follows. \ $\blacksquare $

Before we proceed with the numerical examples we note the following upshot
of the preceding result: In the context of the model considered here, we
either must be willing to assume that $x_{t}$ is dependent on some future
innovations $\varepsilon _{t+i}$, $i>0$ (in which case $X$ and $\mathbf{U}$
will typically be dependent), or otherwise arrive at independence of $X$ and 
$\mathbf{U}$ (in which case our theory can be applied via conditioning).

The numerical examples are now constructed in such a way that the processes $%
(x_{t})$ and $(u_{t})$ (and in particular $X$ and $\mathbf{U}$) are indeed
dependent in order to generate a scenario that is \emph{unfavorable} to our
procedure in the sense that applying our \emph{theory} via conditioning is
not feasible. In light of the proposition this means that we need to let $%
x_{t}$ depend on future innovations (since $x_{t}$ must be independent of
current and past innovations to justify ordinary least squares estimation in
light of the preceding proposition). It should be noted, however, that the
so-constructed examples are favorable to the Kiefer-Vogelsang approach in
the sense that the asymptotic theory developed in those papers applies.

Let now $x_{t}$ be given as $x_{t}=\gamma (\varepsilon _{t+1}/\sigma
_{\varepsilon })+(1-\gamma )z_{t}$ with $0\leq \gamma \leq 1$ and with $%
z_{t} $ the stationary solution of $z_{t}=\delta z_{t-1}+\tilde{\varepsilon}%
_{t}$, where $(\tilde{\varepsilon}_{t})$ is Gaussian white noise that is
independent of the process $(\varepsilon _{t})$ and where $\delta \in (-1,1)$%
. The variance of $\tilde{\varepsilon}_{t}$ is denoted by $\sigma _{\tilde{%
\varepsilon}}^{2}$ and is assumed to be positive. Note that we have $\eta
=(\gamma ,\delta ,\sigma _{\tilde{\varepsilon}}^{2})$. We are interested in
testing the null hypothesis $b=0$ versus the alternative $b\neq 0$ in the
above model at the $0.05$ significance level. [In fact, by invariance
considerations it follows that the results are the same as for testing the
null hypothesis $b=b_{\ast }$.] We use the test-statistic $T_{w}$ (or
equivalently the corresponding root $t_{w}$). As in Section \ref{numerical},
we use sample size $n=100$ and the Bartlett weights $w(j,n)=(1-\left\vert
j\right\vert /M_{n})\mathbf{1}_{(-1,1)}(j/M_{n})$ with $M_{n}=n/10$ (i.e.,
bandwidth parameter $M_{n}$ equal to $10$). For a range of choices for the
parameters $\gamma $, $\delta $, and $\rho $ (and setting $\sigma _{\tilde{%
\varepsilon}}^{2}=1-\delta ^{2}$, $\sigma _{\varepsilon }^{2}=1-\rho ^{2}$)
we obtain a Monte-Carlo sample for $(x_{1},\ldots ,x_{n})^{\prime }$ and $%
\mathbf{U=}(u_{1},\ldots ,u_{n})^{\prime }$ from which we then build $X$ in
an obvious way. We proceed to computing the dependent variable $\mathbf{Y}$
as $X(0,0)^{\prime }+\mathbf{U}$ (observe that also setting $a=0$ has no
effect on the null-rejection probabilities in view of invariance
considerations). We then do two things: First, we compute the test statistic 
$t_{w}$ from the data $\mathbf{Y}$ and $X$ and compare its absolute value to
the \cite{KV2005} critical value $2.260568$ (corresponding to a nominal
significance level of $0.05$). By repeating this for $1~000$ Monte-Carlo
samples we obtain a Monte-Carlo estimate for the null-rejection probability
when using this Kiefer-Vogelsang critical value. We report these in the
right panel of the tables given below. Second, for a given Monte-Carlo
sample, we use the matrix $X$ and employ our Algorithm \ref{alg:AR} (with $%
p=1$ and tuning parameters $M_{0}=20$, $M_{1}=M_{2}=1$, $N_{0}=N_{1}=1~000$, 
$N_{2}=10~000$, cf. Appendix \ref{App F}) to compute critical values for $%
\left\vert t_{w}\right\vert $ that would result in size control over $%
\mathfrak{F}_{\mathrm{AR(}1\mathrm{)}}$ if the regressors matrix were
nonrandom and equal to the given $X$. This results in a critical value $c(X)$%
. We then compare the test statistic $\left\vert t_{w}\right\vert $ computed
from the given Monte-Carlo sample with $c(X)$ and record whether we reject
or not. Repeating over the Monte-Carlo samples this gives us a Monte-Carlo
estimate of the null-rejection probability of the so defined procedure. We
report these in the left panel of the tables given below.

\begin{center}
{\small 
\begin{tabular}{l|llllll||llllll}
\hline\hline
$\gamma =0$ &  &  &  &  &  &  &  &  &  &  &  &  \\ \hline\hline
Method &  &  & Alg &  &  &  &  &  &  & KV &  &  \\ \hline\hline
$\rho ~\big \backslash~\delta $ & 0 & 0.5 & 0.75 & 0.9 & 0.95 & 0.99 & 0 & 
0.5 & 0.75 & 0.9 & 0.95 & 0.99 \\ \hline\hline
0 & 0.04 & 0.04 & 0.02 & 0.01 & 0.00 & 0.00 & 0.06 & 0.07 & 0.07 & 0.07 & 
0.07 & 0.07 \\ 
0.5 & 0.04 & 0.04 & 0.02 & 0.01 & 0.01 & 0.00 & 0.05 & 0.08 & 0.07 & 0.09 & 
0.11 & 0.11 \\ 
0.75 & 0.04 & 0.04 & 0.03 & 0.02 & 0.01 & 0.01 & 0.05 & 0.09 & 0.11 & 0.15 & 
0.16 & 0.17 \\ 
0.9 & 0.03 & 0.04 & 0.05 & 0.03 & 0.03 & 0.02 & 0.05 & 0.10 & 0.15 & 0.20 & 
0.26 & 0.24 \\ 
0.95 & 0.03 & 0.04 & 0.04 & 0.03 & 0.03 & 0.03 & 0.05 & 0.09 & 0.15 & 0.24 & 
0.28 & 0.33 \\ 
0.99 & 0.02 & 0.03 & 0.03 & 0.04 & 0.05 & 0.04 & 0.03 & 0.09 & 0.15 & 0.28 & 
0.37 & 0.41 \\ \hline
\end{tabular}%
}

{\small 
\begin{tabular}{l|llllll||llllll}
\hline\hline
$\gamma =\frac{1}{6}$ &  &  &  &  &  &  &  &  &  &  &  &  \\ \hline\hline
Method &  &  & Alg &  &  &  &  &  &  & KV &  &  \\ \hline\hline
$\rho ~\big \backslash~\delta $ & 0 & 0.5 & 0.75 & 0.9 & 0.95 & 0.99 & 0 & 
0.5 & 0.75 & 0.9 & 0.95 & 0.99 \\ \hline\hline
0 & 0.04 & 0.03 & 0.02 & 0.01 & 0.00 & 0.00 & 0.06 & 0.04 & 0.07 & 0.08 & 
0.09 & 0.07 \\ 
0.5 & 0.04 & 0.03 & 0.03 & 0.01 & 0.01 & 0.01 & 0.06 & 0.06 & 0.09 & 0.10 & 
0.11 & 0.10 \\ 
0.75 & 0.04 & 0.04 & 0.04 & 0.02 & 0.01 & 0.01 & 0.06 & 0.07 & 0.11 & 0.15 & 
0.16 & 0.17 \\ 
0.9 & 0.04 & 0.05 & 0.06 & 0.03 & 0.03 & 0.02 & 0.06 & 0.09 & 0.15 & 0.18 & 
0.25 & 0.27 \\ 
0.95 & 0.04 & 0.03 & 0.04 & 0.04 & 0.04 & 0.03 & 0.06 & 0.09 & 0.14 & 0.25 & 
0.31 & 0.34 \\ 
0.99 & 0.03 & 0.04 & 0.04 & 0.04 & 0.06 & 0.06 & 0.05 & 0.08 & 0.14 & 0.26 & 
0.35 & 0.44 \\ \hline
\end{tabular}%
}

{\small 
\begin{tabular}{l|llllll||llllll}
\hline\hline
$\gamma =\frac{1}{3}$ &  &  &  &  &  &  &  &  &  &  &  &  \\ \hline\hline
Method &  &  & Alg &  &  &  &  &  &  & KV &  &  \\ \hline\hline
$\rho ~\big \backslash~\delta $ & 0 & 0.5 & 0.75 & 0.9 & 0.95 & 0.99 & 0 & 
0.5 & 0.75 & 0.9 & 0.95 & 0.99 \\ \hline\hline
0 & 0.05 & 0.04 & 0.02 & 0.01 & 0.01 & 0.01 & 0.06 & 0.07 & 0.06 & 0.06 & 
0.07 & 0.06 \\ 
0.5 & 0.05 & 0.06 & 0.02 & 0.01 & 0.02 & 0.02 & 0.08 & 0.07 & 0.08 & 0.09 & 
0.10 & 0.11 \\ 
0.75 & 0.06 & 0.05 & 0.04 & 0.02 & 0.01 & 0.02 & 0.07 & 0.09 & 0.12 & 0.14 & 
0.12 & 0.14 \\ 
0.9 & 0.06 & 0.06 & 0.06 & 0.05 & 0.04 & 0.05 & 0.08 & 0.11 & 0.13 & 0.22 & 
0.23 & 0.24 \\ 
0.95 & 0.06 & 0.04 & 0.06 & 0.05 & 0.05 & 0.08 & 0.09 & 0.08 & 0.15 & 0.26 & 
0.28 & 0.30 \\ 
0.99 & 0.07 & 0.05 & 0.05 & 0.06 & 0.07 & 0.09 & 0.10 & 0.10 & 0.17 & 0.25 & 
0.37 & 0.39 \\ \hline
\end{tabular}%
}

{\small 
\begin{tabular}{l|llllll||llllll}
\hline\hline
$\gamma =\frac{1}{2}$ &  &  &  &  &  &  &  &  &  &  &  &  \\ \hline\hline
Method &  &  & Alg &  &  &  &  &  &  & KV &  &  \\ \hline\hline
$\rho ~\big \backslash~\delta $ & 0 & 0.5 & 0.75 & 0.9 & 0.95 & 0.99 & 0 & 
0.5 & 0.75 & 0.9 & 0.95 & 0.99 \\ \hline\hline
0 & 0.04 & 0.04 & 0.03 & 0.02 & 0.02 & 0.02 & 0.06 & 0.07 & 0.07 & 0.06 & 
0.07 & 0.06 \\ 
0.5 & 0.06 & 0.05 & 0.05 & 0.04 & 0.03 & 0.04 & 0.07 & 0.08 & 0.09 & 0.11 & 
0.10 & 0.10 \\ 
0.75 & 0.07 & 0.08 & 0.07 & 0.06 & 0.06 & 0.08 & 0.09 & 0.12 & 0.11 & 0.12 & 
0.15 & 0.16 \\ 
0.9 & 0.09 & 0.07 & 0.08 & 0.07 & 0.07 & 0.11 & 0.11 & 0.09 & 0.14 & 0.21 & 
0.21 & 0.22 \\ 
0.95 & 0.11 & 0.11 & 0.08 & 0.09 & 0.07 & 0.13 & 0.14 & 0.14 & 0.15 & 0.22 & 
0.23 & 0.26 \\ 
0.99 & 0.15 & 0.12 & 0.11 & 0.10 & 0.11 & 0.18 & 0.18 & 0.16 & 0.19 & 0.28 & 
0.32 & 0.38 \\ \hline
\end{tabular}%
}

{\small 
\begin{tabular}{l|llllll||llllll}
\hline\hline
$\gamma =\frac{2}{3}$ &  &  &  &  &  &  &  &  &  &  &  &  \\ \hline\hline
Method &  &  & Alg &  &  &  &  &  &  & KV &  &  \\ \hline\hline
$\rho ~\big \backslash~\delta $ & 0 & 0.5 & 0.75 & 0.9 & 0.95 & 0.99 & 0 & 
0.5 & 0.75 & 0.9 & 0.95 & 0.99 \\ \hline\hline
0 & 0.04 & 0.04 & 0.04 & 0.02 & 0.04 & 0.03 & 0.05 & 0.06 & 0.05 & 0.06 & 
0.07 & 0.06 \\ 
0.5 & 0.07 & 0.06 & 0.06 & 0.05 & 0.05 & 0.06 & 0.08 & 0.08 & 0.08 & 0.08 & 
0.09 & 0.10 \\ 
0.75 & 0.09 & 0.09 & 0.10 & 0.10 & 0.09 & 0.11 & 0.11 & 0.12 & 0.12 & 0.14 & 
0.15 & 0.14 \\ 
0.9 & 0.14 & 0.13 & 0.14 & 0.15 & 0.14 & 0.17 & 0.17 & 0.15 & 0.17 & 0.20 & 
0.20 & 0.23 \\ 
0.95 & 0.19 & 0.15 & 0.14 & 0.16 & 0.17 & 0.22 & 0.21 & 0.17 & 0.19 & 0.24 & 
0.26 & 0.28 \\ 
0.99 & 0.26 & 0.19 & 0.18 & 0.18 & 0.21 & 0.27 & 0.29 & 0.22 & 0.25 & 0.27 & 
0.33 & 0.37 \\ \hline
\end{tabular}%
}

{\small 
\begin{tabular}{l|llllll||llllll}
\hline\hline
$\gamma =\frac{5}{6}$ &  &  &  &  &  &  &  &  &  &  &  &  \\ \hline\hline
Method &  &  & Alg &  &  &  &  &  &  & KV &  &  \\ \hline\hline
$\rho ~\big \backslash~\delta $ & 0 & 0.5 & 0.75 & 0.9 & 0.95 & 0.99 & 0 & 
0.5 & 0.75 & 0.9 & 0.95 & 0.99 \\ \hline\hline
0 & 0.04 & 0.04 & 0.05 & 0.04 & 0.04 & 0.04 & 0.06 & 0.06 & 0.07 & 0.06 & 
0.05 & 0.06 \\ 
0.5 & 0.07 & 0.07 & 0.08 & 0.06 & 0.07 & 0.07 & 0.09 & 0.09 & 0.09 & 0.08 & 
0.10 & 0.08 \\ 
0.75 & 0.12 & 0.12 & 0.12 & 0.12 & 0.11 & 0.11 & 0.14 & 0.13 & 0.13 & 0.14 & 
0.13 & 0.13 \\ 
0.9 & 0.17 & 0.18 & 0.18 & 0.20 & 0.17 & 0.18 & 0.19 & 0.20 & 0.19 & 0.23 & 
0.19 & 0.20 \\ 
0.95 & 0.22 & 0.21 & 0.18 & 0.21 & 0.23 & 0.21 & 0.23 & 0.23 & 0.19 & 0.24 & 
0.26 & 0.23 \\ 
0.99 & 0.29 & 0.29 & 0.29 & 0.32 & 0.30 & 0.30 & 0.32 & 0.33 & 0.34 & 0.36 & 
0.35 & 0.35 \\ \hline
\end{tabular}%
}

{\small 
\begin{tabular}{l|llllll||llllll}
\hline\hline
$\gamma =1$ &  &  &  &  &  &  &  &  &  &  &  &  \\ \hline\hline
Method &  &  & Alg &  &  &  &  &  &  & KV &  &  \\ \hline\hline
$\rho ~\big \backslash~\delta $ & 0 & 0.5 & 0.75 & 0.9 & 0.95 & 0.99 & 0 & 
0.5 & 0.75 & 0.9 & 0.95 & 0.99 \\ \hline\hline
0 & 0.07 & 0.05 & 0.06 & 0.05 & 0.05 & 0.05 & 0.09 & 0.07 & 0.07 & 0.06 & 
0.07 & 0.06 \\ 
0.5 & 0.07 & 0.07 & 0.08 & 0.07 & 0.09 & 0.07 & 0.09 & 0.08 & 0.10 & 0.08 & 
0.11 & 0.08 \\ 
0.75 & 0.12 & 0.11 & 0.12 & 0.12 & 0.10 & 0.11 & 0.14 & 0.13 & 0.14 & 0.13 & 
0.12 & 0.13 \\ 
0.9 & 0.20 & 0.19 & 0.17 & 0.17 & 0.17 & 0.19 & 0.21 & 0.21 & 0.18 & 0.19 & 
0.19 & 0.21 \\ 
0.95 & 0.25 & 0.21 & 0.24 & 0.24 & 0.23 & 0.22 & 0.27 & 0.23 & 0.26 & 0.25 & 
0.25 & 0.23 \\ 
0.99 & 0.33 & 0.32 & 0.32 & 0.31 & 0.30 & 0.32 & 0.36 & 0.35 & 0.35 & 0.34 & 
0.35 & 0.36 \\ \hline
\end{tabular}%
}
\end{center}

\bigskip

The results in the preceding tables clearly show that our method based on
Algorithm \ref{alg:AR}, which presently has no theoretical justification in
the example considered (except in the case where $\gamma =0$), typically
performs better in terms of null-rejection probabilities than the competitor
method based on the asymptotic theory developed in \cite{KV2005}, often by a
considerable margin. In particular, the range of parameters $\gamma, \rho,
\delta$, for which the null-rejection probabilities of our method do not
exceed the prescribed significance level of $0.05$, is considerably larger
than the corresponding range for the competitor method (and includes cases
with $\gamma \neq 0$); for the other cases, where the null-rejection
probabilities of our method exceed $0.05$, the null-rejection probabilities
of the competitor method are even larger, sometimes considerably. [Note that
in case $\gamma =1$ the results do not depend on the value of $\delta $, the
variation in the rows of these subtables hence only reflecting Monte-Carlo
uncertainty.] A fortiori the same conclusion applies if the Kiefer-Vogelsang
critical value is replaced by a smaller critical value (as, e.g., the
standard normal critical value $1.96$ suggested by standard
bandwidth-to-zero asymptotics).

\section{Appendix\label{app:tables}: Tables}

In the following we present tables for the numerical results underlying
Figure 1 in Section \ref{sec:macro}. For a detailed description of the
FRED-MD database, including description of the variables and the
transformations applied to each time series to achieve stationarity, we
refer the reader to \cite{mccracken}.

{\footnotesize 
\begin{longtable}{lrrrrrrrr}	
	\caption{Sizes of  the t-type tests corresponding to $T_w$ using the Kiefer-Vogelsang critical value and  obtained by an application of Algorithm \ref{alg:size}.} \label{tab:size} \\
	\hline \hline
	\multicolumn{1}{l}{\textbf{Regressor}} & 
	\multicolumn{1}{c}{\textbf{i.i.d.}} & 
	\multicolumn{1}{c}{\textbf{AR(1)}} &
	\multicolumn{1}{c}{\textbf{AR(2)}} &
	\multicolumn{1}{c}{\textbf{AR(5)}} &
	\multicolumn{1}{c}{\textbf{AR(10)}} &
	\multicolumn{1}{c}{\textbf{AR(25)}} &
	\multicolumn{1}{c}{\textbf{AR(50)}} &
	\multicolumn{1}{c}{\textbf{AR(99)}} \\
	\endfirsthead
	
	\multicolumn{9}{r}{{Table \ref{tab:size} continues on the next page.}} \\ 
	\endfoot
	
	\hline \hline
	\endlastfoot
	
	RPI & 0.11 & 0.40 & 0.97 & 0.97 & 0.96 & 0.94 & 0.94 & 0.92 \\ 
	W875RX1 & 0.11 & 0.41 & 0.96 & 0.99 & 0.97 & 0.95 & 0.92 & 0.91 \\ 
	DPCERA3M086SBEA & 0.11 & 0.45 & 0.97 & 0.94 & 0.95 & 0.93 & 0.94 & 0.94 \\ 
	CMRMTSPLx & 0.19 & 0.46 & 0.95 & 0.99 & 0.98 & 0.88 & 0.89 & 0.90 \\ 
	RETAILx & 0.15 & 0.24 & 0.52 & 0.53 & 0.53 & 0.53 & 0.53 & 0.53 \\ 
	INDPRO & 0.15 & 0.25 & 0.63 & 0.64 & 0.64 & 0.64 & 0.64 & 0.64 \\ 
	IPFPNSS & 0.18 & 0.32 & 0.58 & 0.59 & 0.60 & 0.59 & 0.60 & 0.60 \\ 
	IPFINAL & 0.16 & 0.28 & 0.60 & 0.60 & 0.60 & 0.60 & 0.60 & 0.60 \\ 
	IPCONGD & 0.13 & 0.37 & 0.96 & 0.90 & 0.90 & 0.95 & 0.96 & 0.94 \\ 
	IPDCONGD & 0.17 & 0.32 & 0.68 & 0.69 & 0.69 & 0.80 & 0.83 & 0.87 \\ 
	IPNCONGD & 0.08 & 0.38 & 0.97 & 0.99 & 0.94 & 0.92 & 0.92 & 0.98 \\ 
	IPBUSEQ & 0.12 & 0.23 & 0.57 & 0.57 & 0.57 & 0.57 & 0.58 & 0.57 \\ 
	IPMAT & 0.12 & 0.22 & 0.68 & 0.68 & 0.68 & 0.68 & 0.68 & 0.68 \\ 
	IPDMAT & 0.13 & 0.23 & 0.61 & 0.63 & 0.63 & 0.62 & 0.63 & 0.63 \\ 
	IPNMAT & 0.16 & 0.22 & 0.68 & 0.69 & 0.69 & 0.69 & 0.69 & 0.69 \\ 
	IPMANSICS & 0.18 & 0.29 & 0.61 & 0.62 & 0.62 & 0.62 & 0.62 & 0.62 \\ 
	IPB51222S & 0.07 & 0.14 & 0.65 & 0.69 & 0.69 & 0.69 & 0.69 & 0.69 \\ 
	IPFUELS & 0.07 & 0.11 & 0.78 & 0.78 & 0.78 & 0.78 & 0.78 & 0.78 \\ 
	CUMFNS & 0.13 & 0.23 & 0.59 & 0.60 & 0.60 & 0.60 & 0.60 & 0.60 \\ 
	HWI & 0.12 & 0.42 & 0.95 & 0.99 & 0.77 & 0.83 & 0.82 & 0.90 \\ 
	HWIURATIO & 0.12 & 0.50 & 0.98 & 0.98 & 0.97 & 0.95 & 0.96 & 0.95 \\ 
	CLF16OV & 0.09 & 0.49 & 0.97 & 0.93 & 0.94 & 0.94 & 0.96 & 0.94 \\ 
	CE16OV & 0.12 & 0.52 & 0.97 & 0.94 & 0.94 & 0.95 & 0.96 & 0.96 \\ 
	UNRATE & 0.14 & 0.52 & 0.97 & 0.95 & 0.96 & 0.96 & 0.95 & 0.97 \\ 
	UEMPMEAN & 0.09 & 0.53 & 0.98 & 0.95 & 0.96 & 0.97 & 0.97 & 0.96 \\ 
	UEMPLT5 & 0.08 & 0.37 & 0.96 & 0.98 & 0.86 & 0.91 & 0.90 & 0.93 \\ 
	UEMP5TO14 & 0.14 & 0.26 & 0.75 & 0.76 & 0.76 & 0.76 & 0.76 & 0.77 \\ 
	UEMP15OV & 0.12 & 0.52 & 0.97 & 0.96 & 0.95 & 0.95 & 0.96 & 0.96 \\ 
	UEMP15T26 & 0.16 & 0.32 & 0.64 & 0.65 & 0.65 & 0.65 & 0.65 & 0.65 \\ 
	UEMP27OV & 0.11 & 0.53 & 0.98 & 0.95 & 0.97 & 0.95 & 0.95 & 0.98 \\ 
	CLAIMSx & 0.17 & 0.28 & 0.71 & 0.71 & 0.71 & 0.71 & 0.71 & 0.71 \\ 
	PAYEMS & 0.13 & 0.52 & 0.97 & 0.95 & 0.95 & 0.95 & 0.95 & 0.98 \\ 
	USGOOD & 0.14 & 0.53 & 0.97 & 0.94 & 0.99 & 0.97 & 0.96 & 0.96 \\ 
	CES1021000001 & 0.09 & 0.43 & 0.94 & 0.87 & 0.88 & 0.91 & 0.87 & 0.91 \\ 
	USCONS & 0.11 & 0.52 & 0.97 & 0.95 & 0.96 & 0.95 & 0.98 & 0.96 \\ 
	MANEMP & 0.15 & 0.53 & 0.97 & 0.94 & 0.95 & 0.93 & 0.95 & 0.96 \\ 
	DMANEMP & 0.15 & 0.51 & 0.97 & 0.99 & 0.99 & 0.93 & 0.93 & 0.95 \\ 
	NDMANEMP & 0.13 & 0.50 & 0.97 & 0.95 & 0.97 & 0.97 & 0.96 & 0.98 \\ 
	SRVPRD & 0.12 & 0.51 & 0.97 & 0.94 & 0.98 & 0.94 & 0.95 & 0.96 \\ 
	USTPU & 0.13 & 0.52 & 0.97 & 0.96 & 0.92 & 0.94 & 0.95 & 0.96 \\ 
	USWTRADE & 0.14 & 0.53 & 0.97 & 0.95 & 0.95 & 0.96 & 0.95 & 0.95 \\ 
	USTRADE & 0.13 & 0.50 & 0.97 & 0.95 & 0.94 & 0.93 & 0.96 & 0.95 \\ 
	USFIRE & 0.11 & 0.52 & 0.98 & 0.95 & 0.95 & 0.96 & 0.95 & 0.96 \\ 
	USGOVT & 0.08 & 0.47 & 0.97 & 0.93 & 0.94 & 0.93 & 0.93 & 0.96 \\ 
	CES0600000007 & 0.09 & 0.41 & 0.96 & 0.89 & 0.89 & 0.89 & 0.89 & 0.88 \\ 
	AWOTMAN & 0.10 & 0.38 & 0.95 & 0.99 & 0.89 & 0.87 & 0.91 & 0.88 \\ 
	AWHMAN & 0.10 & 0.40 & 0.95 & 0.89 & 0.90 & 0.90 & 0.93 & 0.94 \\ 
	HOUST & 0.11 & 0.51 & 0.97 & 0.94 & 1.00 & 0.97 & 0.96 & 0.95 \\ 
	HOUSTNE & 0.11 & 0.44 & 0.97 & 0.96 & 0.97 & 0.96 & 0.95 & 0.96 \\ 
	HOUSTMW & 0.07 & 0.16 & 0.52 & 0.62 & 0.63 & 0.63 & 0.63 & 0.62 \\ 
	HOUSTS & 0.09 & 0.47 & 0.97 & 0.91 & 0.98 & 0.93 & 0.92 & 0.95 \\ 
	HOUSTW & 0.08 & 0.47 & 0.97 & 0.94 & 0.92 & 0.95 & 0.94 & 0.93 \\ 
	PERMIT & 0.13 & 0.46 & 0.96 & 0.99 & 0.93 & 0.92 & 0.94 & 0.93 \\ 
	PERMITNE & 0.15 & 0.29 & 0.68 & 0.80 & 0.81 & 0.79 & 0.81 & 0.80 \\ 
	PERMITMW & 0.08 & 0.36 & 0.94 & 0.85 & 0.88 & 0.86 & 0.89 & 0.87 \\ 
	PERMITS & 0.11 & 0.45 & 0.95 & 0.99 & 0.82 & 0.90 & 0.92 & 0.88 \\ 
	PERMITW & 0.10 & 0.47 & 0.96 & 0.91 & 0.90 & 0.92 & 0.92 & 0.94 \\ 
	ACOGNO & 0.08 & 0.33 & 0.61 & 0.62 & 0.62 & 0.62 & 0.62 & 0.64 \\ 
	AMDMNOx & 0.10 & 0.21 & 0.61 & 0.62 & 0.61 & 0.62 & 0.62 & 0.62 \\ 
	ANDENOx & 0.08 & 0.99 & 1.00 & 1.00 & 1.00 & 1.00 & 1.00 & 1.00 \\ 
	AMDMUOx & 0.11 & 0.46 & 0.93 & 0.86 & 0.92 & 0.83 & 0.91 & 0.95 \\ 
	BUSINVx & 0.14 & 0.39 & 0.63 & 0.64 & 0.64 & 0.63 & 0.64 & 0.64 \\ 
	ISRATIOx & 0.12 & 0.46 & 0.97 & 0.99 & 0.94 & 0.93 & 0.96 & 0.96 \\ 
	M1SL & 0.10 & 0.41 & 0.70 & 0.71 & 0.71 & 0.71 & 0.71 & 0.71 \\ 
	M2SL & 0.11 & 0.50 & 0.97 & 0.94 & 1.00 & 0.93 & 0.95 & 0.96 \\ 
	M2REAL & 0.11 & 0.48 & 0.97 & 0.93 & 0.93 & 0.94 & 0.95 & 0.97 \\ 
	AMBSL & 0.08 & 0.19 & 0.67 & 0.67 & 0.67 & 0.67 & 0.67 & 0.67 \\ 
	TOTRESNS & 0.08 & 0.21 & 0.66 & 0.68 & 0.68 & 0.68 & 0.68 & 0.68 \\ 
	NONBORRES & 0.09 & 0.36 & 0.96 & 0.98 & 0.91 & 0.95 & 0.95 & 0.96 \\ 
	BUSLOANS & 0.10 & 0.52 & 0.98 & 0.96 & 0.96 & 0.96 & 0.96 & 0.96 \\ 
	REALLN & 0.09 & 0.51 & 0.97 & 0.96 & 0.96 & 0.96 & 0.96 & 0.96 \\ 
	NONREVSL & 0.11 & 0.49 & 0.97 & 0.99 & 0.94 & 0.97 & 0.95 & 0.95 \\ 
	CONSPI & 0.09 & 0.43 & 0.96 & 0.91 & 0.90 & 0.92 & 0.93 & 0.93 \\ 
	S.P.500 & 0.11 & 0.32 & 0.94 & 0.99 & 0.71 & 0.84 & 0.83 & 0.68 \\ 
	S.P..indust & 0.10 & 0.31 & 0.61 & 0.62 & 0.62 & 0.62 & 0.65 & 0.63 \\ 
	S.P.div.yield & 0.14 & 0.39 & 0.96 & 0.99 & 0.91 & 0.87 & 0.88 & 0.93 \\ 
	S.P.PE.ratio & 0.21 & 0.38 & 0.71 & 0.71 & 0.71 & 0.72 & 0.72 & 0.72 \\ 
	FEDFUNDS & 0.22 & 0.46 & 0.96 & 0.91 & 0.92 & 0.88 & 0.93 & 0.96 \\ 
	CP3Mx & 0.17 & 0.47 & 0.97 & 0.93 & 0.94 & 0.92 & 0.93 & 0.95 \\ 
	TB3MS & 0.22 & 0.45 & 0.97 & 0.94 & 0.92 & 0.92 & 0.94 & 0.94 \\ 
	TB6MS & 0.18 & 0.47 & 0.97 & 0.99 & 0.92 & 0.94 & 0.95 & 0.94 \\ 
	GS1 & 0.16 & 0.46 & 0.97 & 0.93 & 0.93 & 0.94 & 0.97 & 0.96 \\ 
	GS5 & 0.08 & 0.47 & 0.98 & 0.95 & 0.95 & 0.94 & 0.96 & 0.95 \\ 
	GS10 & 0.07 & 0.23 & 0.52 & 0.53 & 0.53 & 0.53 & 0.53 & 0.52 \\ 
	AAA & 0.08 & 0.41 & 0.95 & 0.91 & 0.97 & 0.91 & 0.93 & 0.95 \\ 
	BAA & 0.10 & 0.50 & 0.98 & 0.95 & 0.96 & 0.97 & 0.96 & 0.96 \\ 
	COMPAPFFx & 0.19 & 0.42 & 0.96 & 0.91 & 0.91 & 0.89 & 0.92 & 0.96 \\ 
	TB3SMFFM & 0.22 & 0.39 & 0.95 & 0.86 & 0.86 & 0.87 & 0.89 & 0.90 \\ 
	TB6SMFFM & 0.12 & 0.17 & 0.69 & 0.76 & 0.76 & 0.76 & 0.76 & 0.76 \\ 
	T1YFFM & 0.09 & 0.50 & 0.98 & 0.95 & 0.95 & 0.95 & 0.95 & 0.96 \\ 
	T5YFFM & 0.08 & 0.30 & 0.69 & 0.71 & 0.71 & 0.71 & 0.71 & 0.70 \\ 
	T10YFFM & 0.07 & 0.17 & 0.68 & 0.69 & 0.69 & 0.69 & 0.69 & 0.69 \\ 
	AAAFFM & 0.08 & 0.14 & 0.60 & 0.61 & 0.62 & 0.61 & 0.61 & 0.61 \\ 
	BAAFFM & 0.12 & 0.18 & 0.66 & 0.66 & 0.66 & 0.66 & 0.66 & 0.66 \\ 
	TWEXMMTH & 0.08 & 0.51 & 0.98 & 0.95 & 0.95 & 0.96 & 0.95 & 0.97 \\ 
	EXSZUSx & 0.09 & 0.47 & 0.96 & 0.93 & 0.94 & 0.95 & 0.95 & 0.96 \\ 
	EXJPUSx & 0.09 & 0.53 & 0.98 & 0.96 & 0.97 & 0.95 & 0.97 & 0.97 \\ 
	EXUSUKx & 0.11 & 0.16 & 0.72 & 0.71 & 0.72 & 0.72 & 0.72 & 0.72 \\ 
	EXCAUSx & 0.09 & 0.52 & 0.98 & 0.99 & 0.97 & 0.96 & 0.96 & 0.96 \\ 
	WPSFD49207 & 0.08 & 0.46 & 0.97 & 0.95 & 0.99 & 0.95 & 0.93 & 0.95 \\ 
	WPSFD49502 & 0.08 & 0.47 & 0.97 & 0.94 & 0.97 & 0.94 & 0.96 & 0.96 \\ 
	WPSID61 & 0.08 & 0.44 & 0.96 & 0.91 & 0.95 & 1.00 & 0.93 & 0.96 \\ 
	WPSID62 & 0.08 & 0.34 & 0.85 & 0.62 & 0.63 & 0.63 & 0.63 & 0.64 \\ 
	OILPRICEx & 0.08 & 0.37 & 0.96 & 0.99 & 0.88 & 0.88 & 0.90 & 0.90 \\ 
	PPICMM & 0.11 & 0.22 & 0.57 & 0.58 & 0.58 & 0.58 & 0.58 & 0.58 \\ 
	CPIAUCSL & 0.08 & 0.25 & 0.57 & 0.58 & 0.58 & 0.58 & 0.58 & 0.58 \\ 
	CPIAPPSL & 0.08 & 0.41 & 0.95 & 0.99 & 0.80 & 0.83 & 0.83 & 0.88 \\ 
	CPITRNSL & 0.08 & 0.47 & 0.98 & 0.94 & 0.94 & 0.97 & 0.96 & 0.95 \\ 
	CPIMEDSL & 0.10 & 0.26 & 0.68 & 0.69 & 0.69 & 0.69 & 0.69 & 0.69 \\ 
	CUSR0000SAC & 0.08 & 0.47 & 0.97 & 0.94 & 0.96 & 0.96 & 0.94 & 0.95 \\ 
	CUUR0000SAD & 0.08 & 0.54 & 0.98 & 0.96 & 0.96 & 0.97 & 0.96 & 0.96 \\ 
	CUSR0000SAS & 0.09 & 0.52 & 0.98 & 0.95 & 0.96 & 0.95 & 0.96 & 0.96 \\ 
	CPIULFSL & 0.08 & 0.24 & 0.55 & 0.57 & 0.57 & 0.57 & 0.56 & 0.56 \\ 
	CUUR0000SA0L2 & 0.07 & 0.42 & 0.96 & 0.91 & 0.93 & 0.91 & 0.92 & 0.92 \\ 
	CUSR0000SA0L5 & 0.08 & 0.25 & 0.58 & 0.58 & 0.59 & 0.58 & 0.58 & 0.59 \\ 
	PCEPI & 0.08 & 0.26 & 0.58 & 0.59 & 0.59 & 0.59 & 0.59 & 0.59 \\ 
	DDURRG3M086SBEA & 0.08 & 0.52 & 0.97 & 0.95 & 0.96 & 0.96 & 0.96 & 0.96 \\ 
	DNDGRG3M086SBEA & 0.08 & 0.47 & 0.96 & 0.95 & 0.96 & 0.94 & 0.96 & 0.95 \\ 
	DSERRG3M086SBEA & 0.11 & 0.52 & 0.97 & 0.96 & 0.96 & 0.96 & 0.97 & 0.96 \\ 
	CES0600000008 & 0.09 & 0.28 & 0.68 & 0.70 & 0.70 & 0.70 & 0.70 & 0.70 \\ 
	CES2000000008 & 0.08 & 0.29 & 0.74 & 0.74 & 0.74 & 0.74 & 0.74 & 0.74 \\ 
	CES3000000008 & 0.08 & 0.31 & 0.72 & 0.72 & 0.72 & 0.72 & 0.72 & 0.72 \\ 
	UMCSENTx & 0.07 & 0.24 & 0.61 & 0.63 & 0.63 & 0.63 & 0.63 & 0.62 \\ 
	MZMSL & 0.10 & 0.48 & 0.97 & 0.95 & 0.95 & 0.97 & 0.95 & 0.95 \\ 
	DTCOLNVHFNM & 0.13 & 0.33 & 0.66 & 0.66 & 0.66 & 0.66 & 0.66 & 0.66 \\ 
	DTCTHFNM & 0.15 & 0.49 & 0.93 & 0.88 & 0.85 & 0.87 & 0.90 & 0.95 \\ 
	INVEST & 0.08 & 0.29 & 0.90 & 0.72 & 0.69 & 0.69 & 0.78 & 0.74 \\ 
	VXOCLSx & 0.12 & 0.27 & 0.94 & 0.99 & 0.87 & 0.85 & 0.91 & 0.89 \\
	\hline
\end{longtable}}

{\footnotesize 
\begin{longtable}{lrrrrrrrr}
	\caption{Critical Values guaranteeing size $\leq $ 0.05 for the t-type tests corresponding to $T_w$ obtained by an application of Algorithm \ref{alg:AR}.} 	\label{tab:cv} \\
	\hline \hline 
	\multicolumn{1}{l}{\textbf{Regressor}} & 
	\multicolumn{1}{c}{\textbf{i.i.d.}} & 
	\multicolumn{1}{c}{\textbf{AR(1)}} &
	\multicolumn{1}{c}{\textbf{AR(2)}} &
	\multicolumn{1}{c}{\textbf{AR(5)}} &
	\multicolumn{1}{c}{\textbf{AR(10)}} &
	\multicolumn{1}{c}{\textbf{AR(25)}} &
	\multicolumn{1}{c}{\textbf{AR(50)}} &
	\multicolumn{1}{c}{\textbf{AR(99)}} \\
	\endfirsthead
	\multicolumn{9}{r}{{Table \ref{tab:cv} continues on the next page.}} \\ 
	\endfoot
	
	\hline \hline
	\endlastfoot
	
	RPI & 2.81 & 5.26 & 6.29 & 6.68 & 8.13 & 8.27 & 8.29 & 8.37 \\ 
	W875RX1 & 2.83 & 5.37 & 6.60 & 6.95 & 6.97 & 7.69 & 7.91 & 7.94 \\ 
	DPCERA3M086SBEA & 2.84 & 6.34 & 11.27 & 11.21 & 11.17 & 12.56 & 12.53 & 12.59 \\ 
	CMRMTSPLx & 3.65 & 6.61 & 7.06 & 8.61 & 8.44 & 8.58 & 8.65 & 8.67 \\ 
	RETAILx & 3.21 & 4.32 & 4.46 & 4.89 & 4.90 & 5.97 & 5.95 & 6.02 \\ 
	INDPRO & 3.24 & 4.44 & 5.18 & 5.34 & 5.90 & 6.07 & 6.29 & 6.13 \\ 
	IPFPNSS & 3.61 & 5.33 & 5.88 & 6.18 & 6.51 & 8.03 & 7.92 & 8.02 \\ 
	IPFINAL & 3.34 & 4.79 & 5.79 & 5.83 & 6.63 & 7.09 & 7.29 & 7.26 \\ 
	IPCONGD & 3.04 & 4.81 & 4.91 & 6.09 & 6.61 & 6.56 & 6.69 & 6.71 \\ 
	IPDCONGD & 3.36 & 4.65 & 5.63 & 6.40 & 7.10 & 7.11 & 7.15 & 6.94 \\ 
	IPNCONGD & 2.60 & 4.87 & 5.94 & 7.08 & 7.21 & 7.38 & 7.49 & 7.48 \\ 
	IPBUSEQ & 2.96 & 4.23 & 4.79 & 4.81 & 4.95 & 5.43 & 5.92 & 5.82 \\ 
	IPMAT & 2.96 & 4.06 & 5.42 & 5.44 & 5.48 & 6.22 & 6.24 & 6.21 \\ 
	IPDMAT & 3.08 & 4.15 & 5.13 & 5.43 & 6.00 & 6.13 & 6.14 & 6.06 \\ 
	IPNMAT & 3.37 & 4.07 & 5.31 & 5.35 & 6.45 & 6.55 & 6.44 & 6.46 \\ 
	IPMANSICS & 3.56 & 5.04 & 5.95 & 6.27 & 6.78 & 7.79 & 7.61 & 7.76 \\ 
	IPB51222S & 2.46 & 3.17 & 6.97 & 7.11 & 7.84 & 7.73 & 7.87 & 8.34 \\ 
	IPFUELS & 2.53 & 2.78 & 5.29 & 7.18 & 8.50 & 7.52 & 8.23 & 8.31 \\ 
	CUMFNS & 3.02 & 4.10 & 5.16 & 5.13 & 6.61 & 6.80 & 6.69 & 6.68 \\ 
	HWI & 3.01 & 5.74 & 6.09 & 6.24 & 6.65 & 7.08 & 7.09 & 7.06 \\ 
	HWIURATIO & 2.92 & 7.01 & 12.20 & 12.47 & 12.55 & 12.58 & 12.51 & 12.55 \\ 
	CLF16OV & 2.66 & 7.33 & 12.40 & 12.45 & 12.56 & 17.39 & 18.71 & 16.62 \\ 
	CE16OV & 2.96 & 7.17 & 14.87 & 15.46 & 15.26 & 15.52 & 15.42 & 15.34 \\ 
	UNRATE & 3.10 & 6.96 & 11.63 & 11.84 & 11.71 & 11.77 & 12.10 & 11.94 \\ 
	UEMPMEAN & 2.65 & 8.56 & 32.03 & 33.13 & 32.98 & 32.75 & 33.15 & 33.13 \\ 
	UEMPLT5 & 2.50 & 3.75 & 3.75 & 5.62 & 6.27 & 6.26 & 5.87 & 6.08 \\ 
	UEMP5TO14 & 3.15 & 4.55 & 7.47 & 7.46 & 7.49 & 9.73 & 9.64 & 9.64 \\ 
	UEMP15OV & 3.00 & 7.25 & 13.36 & 13.84 & 13.90 & 14.01 & 13.82 & 14.11 \\ 
	UEMP15T26 & 3.33 & 5.20 & 6.68 & 6.75 & 8.25 & 8.60 & 8.34 & 8.61 \\ 
	UEMP27OV & 2.85 & 7.64 & 15.25 & 15.86 & 15.82 & 15.88 & 15.87 & 15.86 \\ 
	CLAIMSx & 3.51 & 4.74 & 7.09 & 7.15 & 7.16 & 8.48 & 8.53 & 8.56 \\ 
	PAYEMS & 2.99 & 7.03 & 13.41 & 13.85 & 13.92 & 14.84 & 14.04 & 13.87 \\ 
	USGOOD & 3.07 & 7.19 & 11.79 & 11.98 & 12.11 & 12.03 & 12.19 & 12.24 \\ 
	CES1021000001 & 2.67 & 6.05 & 11.38 & 11.65 & 11.63 & 11.65 & 12.50 & 12.50 \\ 
	USCONS & 2.90 & 7.26 & 17.51 & 18.48 & 18.39 & 18.37 & 18.36 & 18.47 \\ 
	MANEMP & 3.25 & 7.19 & 9.91 & 10.11 & 10.05 & 10.00 & 10.29 & 10.25 \\ 
	DMANEMP & 3.28 & 7.24 & 9.40 & 9.40 & 9.46 & 9.98 & 9.95 & 10.01 \\ 
	NDMANEMP & 3.03 & 6.88 & 12.05 & 12.38 & 12.35 & 12.32 & 12.58 & 12.58 \\ 
	SRVPRD & 2.92 & 7.02 & 14.70 & 15.24 & 15.09 & 15.21 & 15.40 & 15.75 \\ 
	USTPU & 3.03 & 6.94 & 12.04 & 12.50 & 12.48 & 12.41 & 12.52 & 12.53 \\ 
	USWTRADE & 3.13 & 7.69 & 12.89 & 13.30 & 13.33 & 13.15 & 13.36 & 13.45 \\ 
	USTRADE & 3.05 & 6.79 & 11.69 & 11.68 & 12.06 & 12.08 & 12.03 & 12.33 \\ 
	USFIRE & 2.88 & 7.29 & 18.45 & 19.43 & 18.99 & 19.46 & 19.61 & 19.80 \\ 
	USGOVT & 2.61 & 5.54 & 7.10 & 7.15 & 7.19 & 7.84 & 7.91 & 7.94 \\ 
	CES0600000007 & 2.65 & 4.87 & 5.17 & 5.28 & 5.84 & 6.88 & 6.85 & 6.81 \\ 
	AWOTMAN & 2.78 & 4.38 & 4.52 & 4.86 & 4.92 & 4.93 & 5.47 & 5.46 \\ 
	AWHMAN & 2.76 & 4.55 & 4.65 & 4.89 & 5.03 & 5.48 & 5.48 & 5.53 \\ 
	HOUST & 2.88 & 6.64 & 10.96 & 11.27 & 11.20 & 12.11 & 13.56 & 12.27 \\ 
	HOUSTNE & 2.92 & 4.75 & 4.91 & 7.34 & 7.32 & 7.99 & 8.06 & 7.76 \\ 
	HOUSTMW & 2.44 & 2.87 & 4.30 & 4.79 & 6.29 & 6.08 & 6.30 & 6.29 \\ 
	HOUSTS & 2.64 & 4.94 & 5.68 & 8.66 & 8.90 & 8.94 & 8.73 & 8.92 \\ 
	HOUSTW & 2.57 & 4.89 & 6.13 & 8.22 & 8.94 & 8.93 & 9.22 & 9.38 \\ 
	PERMIT & 3.07 & 6.20 & 12.46 & 12.68 & 12.56 & 12.50 & 12.58 & 12.76 \\ 
	PERMITNE & 3.44 & 4.03 & 7.32 & 7.35 & 8.21 & 7.76 & 8.01 & 8.07 \\ 
	PERMITMW & 2.61 & 4.32 & 4.77 & 5.91 & 6.69 & 6.51 & 6.13 & 6.58 \\ 
	PERMITS & 2.90 & 6.23 & 14.62 & 14.86 & 14.86 & 14.87 & 14.78 & 14.84 \\ 
	PERMITW & 2.78 & 6.19 & 12.32 & 12.40 & 12.46 & 14.85 & 14.42 & 14.82 \\ 
	ACOGNO & 2.56 & 4.30 & 5.00 & 5.00 & 4.99 & 5.88 & 5.79 & 5.80 \\ 
	AMDMNOx & 2.79 & 3.79 & 5.18 & 5.21 & 6.12 & 6.42 & 6.22 & 6.44 \\ 
	ANDENOx & 2.62 & 3.40 & 4.36 & 5.61 & 6.15 & 6.10 & 5.94 & 6.06 \\ 
	AMDMUOx & 2.87 & 6.68 & 12.57 & 12.98 & 12.90 & 12.98 & 13.02 & 13.47 \\ 
	BUSINVx & 3.07 & 6.05 & 7.24 & 7.31 & 7.39 & 7.40 & 7.75 & 7.65 \\ 
	ISRATIOx & 2.89 & 4.85 & 4.89 & 4.97 & 4.92 & 5.71 & 5.86 & 5.69 \\ 
	M1SL & 2.75 & 6.10 & 11.48 & 11.55 & 11.60 & 11.55 & 11.61 & 11.72 \\ 
	M2SL & 2.89 & 7.04 & 7.78 & 8.84 & 8.83 & 10.32 & 10.28 & 10.56 \\ 
	M2REAL & 2.84 & 5.97 & 6.36 & 6.31 & 7.23 & 7.05 & 7.05 & 7.08 \\ 
	AMBSL & 2.58 & 3.66 & 6.68 & 6.73 & 6.73 & 7.15 & 7.47 & 7.49 \\ 
	TOTRESNS & 2.60 & 3.88 & 7.72 & 7.75 & 7.75 & 8.40 & 8.34 & 8.07 \\ 
	NONBORRES & 2.66 & 5.12 & 10.81 & 10.99 & 11.03 & 11.02 & 12.67 & 11.42 \\ 
	BUSLOANS & 2.72 & 7.83 & 13.15 & 13.16 & 13.43 & 13.75 & 13.64 & 13.65 \\ 
	REALLN & 2.65 & 7.49 & 14.10 & 14.13 & 14.13 & 16.08 & 15.86 & 16.02 \\ 
	NONREVSL & 2.86 & 7.07 & 10.15 & 10.18 & 10.18 & 10.10 & 10.73 & 10.79 \\ 
	CONSPI & 2.72 & 5.42 & 6.12 & 6.67 & 7.94 & 8.25 & 8.17 & 8.10 \\ 
	S.P.500 & 2.85 & 4.48 & 6.81 & 6.83 & 6.84 & 7.22 & 7.11 & 7.03 \\ 
	S.P..indust & 2.79 & 4.35 & 6.62 & 6.80 & 6.84 & 6.84 & 6.93 & 6.95 \\ 
	S.P.div.yield & 3.17 & 4.88 & 5.12 & 5.23 & 6.63 & 6.70 & 6.62 & 6.63 \\ 
	S.P.PE.ratio & 3.75 & 5.22 & 6.37 & 6.39 & 6.70 & 7.24 & 7.18 & 7.20 \\ 
	FEDFUNDS & 3.84 & 5.73 & 6.29 & 6.95 & 7.04 & 7.06 & 7.02 & 7.14 \\ 
	CP3Mx & 3.41 & 5.99 & 7.06 & 7.09 & 7.06 & 7.48 & 7.55 & 7.53 \\ 
	TB3MS & 3.83 & 5.77 & 7.25 & 7.28 & 8.00 & 8.12 & 8.42 & 7.84 \\ 
	TB6MS & 3.52 & 6.07 & 7.82 & 7.76 & 7.94 & 8.08 & 7.97 & 8.13 \\ 
	GS1 & 3.37 & 5.92 & 8.22 & 8.36 & 8.18 & 8.28 & 8.34 & 8.35 \\ 
	GS5 & 2.57 & 5.90 & 8.42 & 9.12 & 8.87 & 10.15 & 10.24 & 10.18 \\ 
	GS10 & 2.48 & 3.49 & 4.26 & 4.32 & 5.54 & 5.59 & 5.57 & 5.54 \\ 
	AAA & 2.57 & 4.73 & 5.19 & 5.27 & 5.56 & 6.10 & 6.08 & 6.05 \\ 
	BAA & 2.80 & 6.43 & 8.21 & 8.29 & 9.22 & 9.36 & 9.44 & 9.69 \\ 
	COMPAPFFx & 3.56 & 4.47 & 4.56 & 5.32 & 5.67 & 5.68 & 5.86 & 5.49 \\ 
	TB3SMFFM & 3.82 & 4.49 & 5.44 & 6.80 & 7.27 & 7.11 & 6.81 & 6.86 \\ 
	TB6SMFFM & 2.97 & 3.10 & 6.33 & 6.71 & 7.41 & 7.39 & 7.27 & 8.07 \\ 
	T1YFFM & 2.65 & 5.60 & 7.27 & 10.46 & 11.02 & 11.24 & 11.96 & 10.96 \\ 
	T5YFFM & 2.52 & 4.20 & 9.29 & 9.36 & 9.38 & 10.03 & 9.55 & 9.68 \\ 
	T10YFFM & 2.49 & 3.40 & 10.35 & 10.46 & 10.34 & 10.46 & 10.45 & 10.35 \\ 
	AAAFFM & 2.55 & 3.14 & 5.24 & 5.35 & 5.28 & 5.61 & 6.01 & 6.16 \\ 
	BAAFFM & 2.94 & 3.45 & 5.25 & 5.24 & 5.33 & 5.90 & 5.87 & 5.93 \\ 
	TWEXMMTH & 2.62 & 6.81 & 9.40 & 9.37 & 10.61 & 11.56 & 11.52 & 11.56 \\ 
	EXSZUSx & 2.69 & 6.20 & 7.42 & 7.48 & 7.47 & 8.03 & 7.99 & 8.05 \\ 
	EXJPUSx & 2.66 & 8.02 & 17.44 & 17.67 & 17.65 & 17.70 & 18.56 & 17.69 \\ 
	EXUSUKx & 2.84 & 3.33 & 7.50 & 7.59 & 8.03 & 8.26 & 8.19 & 8.20 \\ 
	EXCAUSx & 2.67 & 7.07 & 9.85 & 9.83 & 11.44 & 11.86 & 11.36 & 11.11 \\ 
	WPSFD49207 & 2.58 & 5.79 & 7.54 & 7.55 & 7.57 & 8.62 & 8.62 & 8.57 \\ 
	WPSFD49502 & 2.57 & 5.86 & 7.63 & 8.50 & 8.41 & 8.74 & 8.66 & 8.66 \\ 
	WPSID61 & 2.58 & 5.25 & 6.28 & 6.28 & 6.34 & 6.61 & 7.34 & 7.36 \\ 
	WPSID62 & 2.60 & 4.19 & 4.63 & 4.61 & 5.25 & 5.42 & 5.67 & 5.68 \\ 
	OILPRICEx & 2.57 & 4.16 & 4.41 & 5.49 & 5.89 & 6.16 & 6.12 & 5.89 \\ 
	PPICMM & 2.89 & 3.84 & 4.55 & 4.65 & 5.25 & 5.81 & 5.77 & 5.73 \\ 
	CPIAUCSL & 2.58 & 3.72 & 4.52 & 4.81 & 5.10 & 5.57 & 5.48 & 5.58 \\ 
	CPIAPPSL & 2.59 & 5.85 & 11.49 & 11.61 & 11.63 & 11.61 & 11.64 & 11.62 \\ 
	CPITRNSL & 2.57 & 5.40 & 6.23 & 6.23 & 6.93 & 7.09 & 7.11 & 7.11 \\ 
	CPIMEDSL & 2.72 & 4.11 & 6.09 & 6.13 & 6.12 & 6.49 & 6.89 & 7.09 \\ 
	CUSR0000SAC & 2.52 & 5.75 & 7.48 & 8.10 & 8.03 & 8.38 & 8.25 & 8.24 \\ 
	CUUR0000SAD & 2.62 & 8.11 & 19.35 & 19.41 & 19.48 & 19.43 & 19.48 & 20.10 \\ 
	CUSR0000SAS & 2.68 & 7.91 & 24.72 & 25.02 & 24.85 & 25.76 & 25.95 & 26.74 \\ 
	CPIULFSL & 2.57 & 3.54 & 3.99 & 4.03 & 5.11 & 5.24 & 5.16 & 5.03 \\ 
	CUUR0000SA0L2 & 2.52 & 4.71 & 5.61 & 6.73 & 7.07 & 7.12 & 7.16 & 7.30 \\ 
	CUSR0000SA0L5 & 2.56 & 3.75 & 4.54 & 5.40 & 4.59 & 5.47 & 5.42 & 5.42 \\ 
	PCEPI & 2.56 & 3.83 & 4.67 & 5.26 & 5.33 & 5.52 & 5.37 & 5.42 \\ 
	DDURRG3M086SBEA & 2.52 & 6.92 & 11.22 & 11.39 & 11.37 & 11.30 & 11.39 & 11.54 \\ 
	DNDGRG3M086SBEA & 2.54 & 5.59 & 7.22 & 7.83 & 7.85 & 8.23 & 8.23 & 8.17 \\ 
	DSERRG3M086SBEA & 2.83 & 7.37 & 24.25 & 24.92 & 24.89 & 24.36 & 24.94 & 26.66 \\ 
	CES0600000008 & 2.64 & 4.65 & 10.56 & 11.00 & 11.08 & 11.10 & 11.07 & 11.38 \\ 
	CES2000000008 & 2.52 & 4.63 & 9.76 & 9.91 & 9.90 & 9.97 & 12.22 & 10.34 \\ 
	CES3000000008 & 2.56 & 5.15 & 14.04 & 14.28 & 14.22 & 14.19 & 14.16 & 14.75 \\ 
	UMCSENTx & 2.53 & 3.67 & 5.44 & 5.46 & 5.69 & 6.21 & 6.11 & 6.28 \\ 
	MZMSL & 2.81 & 6.36 & 6.77 & 6.96 & 8.63 & 8.81 & 8.62 & 8.57 \\ 
	DTCOLNVHFNM & 3.04 & 5.47 & 7.27 & 7.26 & 7.43 & 8.93 & 9.00 & 8.98 \\ 
	DTCTHFNM & 3.26 & 7.24 & 9.96 & 10.07 & 10.14 & 10.66 & 10.83 & 10.82 \\ 
	INVEST & 2.53 & 3.91 & 4.49 & 4.61 & 4.98 & 5.40 & 5.45 & 5.46 \\ 
	VXOCLSx & 2.90 & 3.45 & 4.83 & 4.86 & 4.87 & 5.34 & 5.30 & 5.38 \\ 
	\hline
\end{longtable}}

{\footnotesize 
\begin{longtable}{lrrrrrrrr}
	\caption{Critical Values guaranteeing size $\leq $ 0.05 for the t-type tests corresponding to $T_{E, W}$ obtained by an application of Algorithm \ref{alg:AR}.}	\label{tab:cve} \\
	\hline \hline
	\multicolumn{1}{l}{\textbf{Regressor}} & 
	\multicolumn{1}{c}{\textbf{i.i.d.}} & 
	\multicolumn{1}{c}{\textbf{AR(1)}} &
	\multicolumn{1}{c}{\textbf{AR(2)}} &
	\multicolumn{1}{c}{\textbf{AR(5)}} &
	\multicolumn{1}{c}{\textbf{AR(10)}} &
	\multicolumn{1}{c}{\textbf{AR(25)}} &
	\multicolumn{1}{c}{\textbf{AR(50)}} &
	\multicolumn{1}{c}{\textbf{AR(99)}} \\
	\endfirsthead
	\multicolumn{9}{r}{{Table \ref{tab:cve} continues on the next page.}} \\ 
	\endfoot
	
	\hline \hline
	\endlastfoot
	
	RPI & 2.24 & 4.23 & 5.20 & 5.27 & 7.18 & 7.36 & 7.29 & 7.14 \\ 
	W875RX1 & 2.28 & 4.08 & 4.95 & 5.26 & 6.07 & 6.10 & 6.09 & 6.04 \\ 
	DPCERA3M086SBEA & 2.28 & 5.21 & 7.61 & 7.54 & 7.53 & 10.61 & 10.51 & 10.71 \\ 
	CMRMTSPLx & 2.22 & 3.58 & 4.43 & 4.48 & 4.49 & 4.48 & 4.57 & 4.50 \\ 
	RETAILx & 2.19 & 3.23 & 3.81 & 3.87 & 3.85 & 3.88 & 3.93 & 3.99 \\ 
	INDPRO & 2.23 & 3.16 & 3.81 & 3.89 & 3.88 & 3.97 & 3.99 & 3.98 \\ 
	IPFPNSS & 2.24 & 3.19 & 3.60 & 3.64 & 3.61 & 3.63 & 3.66 & 3.67 \\ 
	IPFINAL & 2.20 & 3.08 & 3.42 & 3.45 & 3.49 & 3.73 & 3.79 & 3.77 \\ 
	IPCONGD & 2.20 & 3.74 & 4.58 & 4.58 & 5.08 & 5.12 & 5.21 & 5.21 \\ 
	IPDCONGD & 2.16 & 3.15 & 3.77 & 3.77 & 3.84 & 3.95 & 3.99 & 3.97 \\ 
	IPNCONGD & 2.21 & 4.27 & 5.36 & 5.30 & 6.49 & 6.80 & 6.77 & 6.62 \\ 
	IPBUSEQ & 2.24 & 3.46 & 4.18 & 4.20 & 4.21 & 4.62 & 4.59 & 4.60 \\ 
	IPMAT & 2.24 & 3.27 & 4.22 & 4.33 & 4.33 & 4.60 & 4.67 & 4.62 \\ 
	IPDMAT & 2.27 & 3.17 & 3.55 & 3.66 & 3.67 & 4.15 & 4.19 & 4.16 \\ 
	IPNMAT & 2.15 & 2.54 & 4.07 & 4.06 & 4.25 & 4.29 & 4.32 & 4.34 \\ 
	IPMANSICS & 2.23 & 2.94 & 3.27 & 3.28 & 3.32 & 3.33 & 3.38 & 3.40 \\ 
	IPB51222S & 2.08 & 2.44 & 3.59 & 3.61 & 4.64 & 4.66 & 4.67 & 4.69 \\ 
	IPFUELS & 2.10 & 2.50 & 5.55 & 5.57 & 5.58 & 5.61 & 5.73 & 5.64 \\ 
	CUMFNS & 2.25 & 3.06 & 3.70 & 3.78 & 3.80 & 4.22 & 4.25 & 4.26 \\ 
	HWI & 2.29 & 4.28 & 5.72 & 5.83 & 5.70 & 5.82 & 5.84 & 5.98 \\ 
	HWIURATIO & 2.34 & 5.60 & 10.47 & 10.66 & 10.66 & 10.64 & 10.81 & 10.78 \\ 
	CLF16OV & 2.29 & 5.69 & 9.47 & 9.26 & 9.28 & 11.80 & 12.00 & 12.29 \\ 
	CE16OV & 2.33 & 5.77 & 11.70 & 11.69 & 11.78 & 11.85 & 11.84 & 12.23 \\ 
	UNRATE & 2.30 & 5.18 & 9.48 & 9.78 & 9.75 & 9.77 & 9.64 & 9.78 \\ 
	UEMPMEAN & 2.41 & 7.79 & 29.75 & 31.05 & 31.16 & 31.12 & 31.09 & 31.05 \\ 
	UEMPLT5 & 2.08 & 2.83 & 3.00 & 3.66 & 3.66 & 3.70 & 3.77 & 3.75 \\ 
	UEMP5TO14 & 2.15 & 2.56 & 3.24 & 3.28 & 3.37 & 3.59 & 3.66 & 3.66 \\ 
	UEMP15OV & 2.35 & 5.53 & 11.44 & 11.97 & 11.89 & 11.67 & 11.82 & 11.88 \\ 
	UEMP15T26 & 2.25 & 3.23 & 3.53 & 3.54 & 3.62 & 3.79 & 3.81 & 3.83 \\ 
	UEMP27OV & 2.38 & 6.11 & 14.43 & 15.28 & 15.24 & 15.17 & 14.98 & 15.24 \\ 
	CLAIMSx & 2.18 & 2.61 & 3.06 & 3.04 & 3.04 & 3.33 & 3.36 & 3.35 \\ 
	PAYEMS & 2.34 & 5.71 & 11.68 & 12.07 & 11.88 & 11.84 & 11.97 & 12.11 \\ 
	USGOOD & 2.32 & 5.29 & 9.52 & 9.88 & 9.78 & 9.93 & 9.99 & 9.92 \\ 
	CES1021000001 & 2.38 & 5.22 & 10.61 & 10.69 & 10.84 & 10.83 & 10.82 & 10.92 \\ 
	USCONS & 2.36 & 6.20 & 16.16 & 17.26 & 17.25 & 17.27 & 16.92 & 17.07 \\ 
	MANEMP & 2.31 & 4.75 & 7.32 & 7.54 & 7.55 & 7.50 & 7.53 & 7.59 \\ 
	DMANEMP & 2.31 & 4.54 & 6.66 & 6.86 & 6.84 & 6.86 & 6.88 & 6.94 \\ 
	NDMANEMP & 2.32 & 5.55 & 10.60 & 10.78 & 10.65 & 10.71 & 10.65 & 11.03 \\ 
	SRVPRD & 2.34 & 5.97 & 13.18 & 13.26 & 13.50 & 13.50 & 13.69 & 13.69 \\ 
	USTPU & 2.32 & 5.43 & 9.96 & 10.21 & 10.04 & 10.22 & 10.24 & 10.38 \\ 
	USWTRADE & 2.33 & 5.21 & 8.66 & 8.92 & 8.88 & 8.92 & 9.12 & 9.14 \\ 
	USTRADE & 2.32 & 5.46 & 9.90 & 10.04 & 9.97 & 10.01 & 10.20 & 10.28 \\ 
	USFIRE & 2.34 & 6.43 & 16.43 & 16.32 & 16.61 & 16.85 & 17.20 & 17.68 \\ 
	USGOVT & 2.33 & 6.25 & 11.95 & 12.00 & 12.01 & 11.99 & 11.98 & 12.48 \\ 
	CES0600000007 & 2.27 & 4.20 & 4.56 & 4.56 & 5.15 & 5.33 & 5.39 & 5.37 \\ 
	AWOTMAN & 2.31 & 3.94 & 4.12 & 4.21 & 4.62 & 5.02 & 5.06 & 5.05 \\ 
	AWHMAN & 2.29 & 3.87 & 4.01 & 4.00 & 4.00 & 4.64 & 4.64 & 4.62 \\ 
	HOUST & 2.17 & 4.38 & 6.78 & 6.91 & 6.94 & 6.96 & 6.94 & 6.91 \\ 
	HOUSTNE & 2.07 & 3.38 & 4.44 & 4.42 & 4.50 & 4.51 & 5.16 & 4.75 \\ 
	HOUSTMW & 2.06 & 3.03 & 3.48 & 3.54 & 3.92 & 3.88 & 3.98 & 4.01 \\ 
	HOUSTS & 2.14 & 4.09 & 5.72 & 5.83 & 5.72 & 5.82 & 7.04 & 5.96 \\ 
	HOUSTW & 2.14 & 4.25 & 6.41 & 6.53 & 6.41 & 6.94 & 7.78 & 6.86 \\ 
	PERMIT & 2.19 & 4.29 & 7.17 & 7.21 & 7.20 & 7.18 & 7.31 & 7.22 \\ 
	PERMITNE & 2.06 & 3.04 & 3.97 & 4.01 & 4.01 & 4.08 & 4.06 & 5.00 \\ 
	PERMITMW & 2.16 & 3.75 & 4.61 & 4.60 & 4.67 & 5.89 & 4.82 & 4.92 \\ 
	PERMITS & 2.20 & 4.29 & 6.61 & 6.69 & 6.80 & 6.78 & 6.89 & 6.81 \\ 
	PERMITW & 2.16 & 4.32 & 6.78 & 6.96 & 7.00 & 7.20 & 7.45 & 7.47 \\ 
	ACOGNO & 2.32 & 4.50 & 6.15 & 6.14 & 6.15 & 6.47 & 6.66 & 6.61 \\ 
	AMDMNOx & 2.18 & 3.14 & 3.40 & 3.48 & 3.51 & 4.09 & 4.09 & 4.09 \\ 
	ANDENOx & 2.15 & 3.75 & 4.57 & 4.65 & 4.60 & 4.87 & 4.89 & 4.86 \\ 
	AMDMUOx & 2.36 & 4.93 & 9.68 & 9.83 & 9.83 & 9.81 & 9.84 & 9.84 \\ 
	BUSINVx & 2.34 & 4.19 & 5.45 & 5.49 & 5.62 & 5.60 & 5.65 & 5.62 \\ 
	ISRATIOx & 2.32 & 4.34 & 4.74 & 4.74 & 5.10 & 5.16 & 5.28 & 5.23 \\ 
	M1SL & 2.31 & 4.13 & 5.12 & 5.10 & 5.12 & 5.62 & 5.71 & 5.69 \\ 
	M2SL & 2.34 & 5.03 & 5.89 & 5.98 & 6.75 & 7.15 & 7.22 & 7.16 \\ 
	M2REAL & 2.30 & 4.49 & 4.83 & 4.85 & 4.87 & 5.65 & 5.69 & 5.69 \\ 
	AMBSL & 2.25 & 3.30 & 4.39 & 4.40 & 5.20 & 5.28 & 5.30 & 5.31 \\ 
	TOTRESNS & 2.25 & 3.40 & 4.51 & 4.57 & 5.16 & 5.29 & 5.26 & 5.29 \\ 
	NONBORRES & 2.31 & 4.19 & 5.07 & 5.38 & 5.42 & 6.81 & 6.87 & 7.00 \\ 
	BUSLOANS & 2.42 & 6.73 & 14.66 & 14.81 & 14.85 & 14.86 & 14.77 & 15.11 \\ 
	REALLN & 2.40 & 6.50 & 11.81 & 11.89 & 11.86 & 12.55 & 12.63 & 12.73 \\ 
	NONREVSL & 2.33 & 5.75 & 10.07 & 9.95 & 10.23 & 10.26 & 10.26 & 10.34 \\ 
	CONSPI & 2.25 & 4.48 & 5.27 & 5.98 & 5.28 & 6.42 & 6.40 & 6.32 \\ 
	S.P.500 & 2.18 & 3.53 & 4.36 & 4.42 & 4.42 & 4.48 & 4.91 & 4.88 \\ 
	S.P..indust & 2.22 & 3.46 & 4.30 & 4.30 & 4.34 & 4.78 & 4.84 & 4.77 \\ 
	S.P.div.yield & 2.30 & 3.11 & 3.81 & 3.87 & 3.87 & 3.88 & 3.98 & 4.00 \\ 
	S.P.PE.ratio & 2.28 & 3.06 & 3.59 & 3.63 & 3.65 & 4.16 & 4.19 & 4.16 \\ 
	FEDFUNDS & 2.20 & 3.77 & 5.06 & 5.04 & 5.08 & 5.33 & 5.33 & 5.38 \\ 
	CP3Mx & 2.27 & 4.32 & 6.17 & 6.66 & 6.19 & 6.59 & 6.62 & 6.63 \\ 
	TB3MS & 2.18 & 3.73 & 5.08 & 5.13 & 5.13 & 5.38 & 5.44 & 5.44 \\ 
	TB6MS & 2.23 & 4.22 & 6.25 & 6.18 & 6.11 & 6.53 & 6.62 & 6.59 \\ 
	GS1 & 2.24 & 4.64 & 7.32 & 7.29 & 7.31 & 7.56 & 7.46 & 7.66 \\ 
	GS5 & 2.25 & 4.64 & 5.41 & 5.47 & 5.62 & 5.70 & 5.72 & 5.86 \\ 
	GS10 & 2.21 & 3.33 & 3.82 & 3.83 & 4.25 & 4.34 & 4.43 & 4.36 \\ 
	AAA & 2.25 & 3.71 & 3.86 & 3.84 & 3.97 & 5.11 & 5.17 & 5.12 \\ 
	BAA & 2.28 & 4.44 & 5.08 & 5.85 & 6.66 & 6.61 & 6.64 & 6.72 \\ 
	COMPAPFFx & 2.18 & 3.95 & 5.19 & 5.21 & 5.16 & 5.25 & 5.35 & 5.36 \\ 
	TB3SMFFM & 2.12 & 3.09 & 3.72 & 3.72 & 3.73 & 3.80 & 3.77 & 3.79 \\ 
	TB6SMFFM & 2.09 & 2.98 & 3.02 & 3.13 & 3.64 & 3.75 & 3.81 & 3.78 \\ 
	T1YFFM & 2.22 & 5.43 & 7.25 & 7.52 & 7.53 & 7.65 & 7.63 & 7.60 \\ 
	T5YFFM & 2.27 & 3.74 & 8.02 & 8.04 & 8.05 & 8.05 & 8.04 & 8.09 \\ 
	T10YFFM & 2.24 & 3.23 & 6.37 & 6.42 & 6.41 & 7.18 & 7.09 & 7.07 \\ 
	AAAFFM & 2.21 & 2.90 & 4.16 & 4.17 & 4.18 & 5.99 & 5.96 & 5.94 \\ 
	BAAFFM & 2.18 & 2.75 & 3.54 & 3.66 & 4.37 & 4.42 & 4.42 & 4.35 \\ 
	TWEXMMTH & 2.36 & 5.72 & 7.62 & 7.67 & 7.68 & 8.27 & 8.30 & 8.17 \\ 
	EXSZUSx & 2.31 & 4.82 & 5.41 & 5.99 & 6.06 & 6.18 & 6.16 & 6.18 \\ 
	EXJPUSx & 2.39 & 6.16 & 11.80 & 11.88 & 11.91 & 11.88 & 11.89 & 13.36 \\ 
	EXUSUKx & 2.14 & 2.85 & 3.81 & 3.85 & 3.92 & 4.55 & 4.46 & 4.55 \\ 
	EXCAUSx & 2.38 & 6.25 & 9.42 & 9.50 & 9.50 & 9.93 & 9.93 & 9.91 \\ 
	WPSFD49207 & 2.35 & 5.43 & 8.39 & 8.43 & 8.44 & 8.75 & 8.63 & 8.75 \\ 
	WPSFD49502 & 2.34 & 5.52 & 8.33 & 8.43 & 8.42 & 8.83 & 8.89 & 8.87 \\ 
	WPSID61 & 2.35 & 5.17 & 6.96 & 6.98 & 6.98 & 7.77 & 7.77 & 7.83 \\ 
	WPSID62 & 2.31 & 4.37 & 5.40 & 5.40 & 6.28 & 6.39 & 6.30 & 6.29 \\ 
	OILPRICEx & 2.28 & 4.31 & 5.38 & 5.44 & 5.44 & 5.91 & 6.08 & 6.09 \\ 
	PPICMM & 2.23 & 3.02 & 3.44 & 3.53 & 3.90 & 4.04 & 4.03 & 4.03 \\ 
	CPIAUCSL & 2.27 & 4.12 & 5.96 & 6.01 & 6.01 & 6.39 & 6.15 & 6.13 \\ 
	CPIAPPSL & 2.35 & 5.35 & 8.58 & 8.66 & 8.67 & 8.66 & 8.85 & 8.78 \\ 
	CPITRNSL & 2.34 & 5.48 & 8.04 & 8.21 & 8.21 & 8.35 & 8.38 & 8.32 \\ 
	CPIMEDSL & 2.21 & 3.51 & 4.68 & 4.67 & 4.70 & 5.09 & 5.14 & 5.13 \\ 
	CUSR0000SAC & 2.36 & 5.69 & 9.44 & 9.55 & 9.49 & 9.57 & 9.54 & 9.55 \\ 
	CUUR0000SAD & 2.35 & 6.43 & 10.58 & 10.73 & 10.72 & 10.90 & 10.93 & 11.38 \\ 
	CUSR0000SAS & 2.38 & 7.33 & 24.37 & 24.85 & 24.62 & 25.00 & 25.95 & 27.12 \\ 
	CPIULFSL & 2.24 & 3.98 & 5.62 & 5.68 & 5.68 & 6.08 & 6.06 & 5.86 \\ 
	CUUR0000SA0L2 & 2.31 & 4.93 & 7.21 & 7.30 & 7.29 & 7.57 & 7.29 & 7.53 \\ 
	CUSR0000SA0L5 & 2.27 & 4.16 & 5.96 & 6.00 & 6.00 & 6.34 & 6.36 & 6.11 \\ 
	PCEPI & 2.29 & 4.26 & 6.19 & 6.26 & 6.27 & 6.44 & 6.45 & 6.42 \\ 
	DDURRG3M086SBEA & 2.31 & 6.11 & 9.59 & 9.66 & 9.66 & 9.66 & 10.04 & 9.66 \\ 
	DNDGRG3M086SBEA & 2.35 & 5.66 & 9.13 & 9.26 & 9.27 & 9.26 & 9.28 & 9.30 \\ 
	DSERRG3M086SBEA & 2.35 & 6.46 & 18.82 & 19.26 & 19.12 & 19.04 & 19.66 & 20.80 \\ 
	CES0600000008 & 2.23 & 4.24 & 8.65 & 8.78 & 8.65 & 8.77 & 8.71 & 9.32 \\ 
	CES2000000008 & 2.19 & 4.31 & 8.97 & 9.23 & 9.19 & 9.20 & 10.28 & 10.31 \\ 
	CES3000000008 & 2.29 & 4.78 & 12.70 & 12.96 & 13.03 & 13.01 & 13.02 & 13.63 \\ 
	UMCSENTx & 2.10 & 2.65 & 3.43 & 3.58 & 4.09 & 4.14 & 4.11 & 4.13 \\ 
	MZMSL & 2.39 & 5.28 & 6.18 & 6.27 & 7.22 & 6.95 & 7.11 & 7.19 \\ 
	DTCOLNVHFNM & 2.23 & 3.47 & 3.95 & 3.98 & 4.02 & 3.99 & 4.08 & 4.12 \\ 
	DTCTHFNM & 2.29 & 4.32 & 5.17 & 5.27 & 5.38 & 5.52 & 5.48 & 5.54 \\ 
	INVEST & 2.22 & 3.34 & 3.90 & 4.00 & 4.43 & 4.43 & 4.43 & 4.39 \\ 
	VXOCLSx & 2.16 & 2.82 & 2.91 & 3.59 & 3.07 & 3.76 & 3.80 & 3.91 \\ 
	\hline
\end{longtable}}

\bibliographystyle{ims}
\bibliography{refs}

\end{document}